\documentclass{article}
\usepackage[T1]{fontenc}
\usepackage{amsthm, fullpage}
\usepackage{amsmath}
\usepackage{amssymb}
\usepackage{amsfonts}
\usepackage{eucal}
\usepackage[all]{xy}
\usepackage[frenchb, english]{babel}
\usepackage{pslatex}
\usepackage{hyperref}

\newtheorem{prop}{Proposition}[subsection]
\newtheorem{theo}[prop]{Théor\`eme}
\newtheorem*{theo**}{Théorème}
\newtheorem{coro}[prop]{Corollaire}
\newtheorem*{conj*}{Conjecture}
\newtheorem{lemm}[prop]{Lemme}
\newtheorem{lemm*}{Lemme}[prop]

\theoremstyle{definition}

\newtheorem{vide}[prop]{}
\newtheorem{defi}[prop]{Définition}
\newtheorem*{defi*}{Définition}
\newtheorem{nota}[prop]{Notations}

\theoremstyle{remark}
\newtheorem{rema}[prop]{Remarques}

\numberwithin{equation}{prop}

\newcommand{\riso}{ \overset{\sim}{\longrightarrow}\, }
\newcommand{\liso}{ \overset{\sim}{\longleftarrow}\, }

\newcommand{\Spf}{\mathrm{Spf}\,}

\renewcommand{\sp}{\mathrm{sp}}

\newcommand{\FF}{{\mathcal{F}}}

\newcommand{\E}{{\mathcal{E}}}
\newcommand{\G}{{\mathcal{G}}}

\newcommand{\D}{{\mathcal{D}}}

\newcommand{\PP}{{\mathcal{P}}}

\renewcommand{\O}{{\mathcal{O}}}

\newcommand{\V}{\mathcal{V}}

\newcommand{\Y}{\mathcal{Y}}
\newcommand{\ZZ}{\mathcal{Z}}
\newcommand{\X}{\mathfrak{X}}

\newcommand{\U}{\mathfrak{U}}

\renewcommand{\P}{\mathbb{P}}

\newcommand{\DD}{\mathbb{D}}

\newcommand{\R}{\mathbb{R}}
\newcommand{\Q}{\mathbb{Q}}
\newcommand{\Z}{\mathbb{Z}}
\newcommand{\N}{\mathbb{N}}

\newcommand{\hdag}{  \phantom{}{^{\dag} }    }

\begin{document}
\selectlanguage{frenchb}

\title{Pleine fidélité sans structure de Frobenius et isocristaux partiellement surconvergents\\
(Full faithfulness without Frobenius structure and partially overcoherent isocrystals)}
\author{Daniel Caro \footnote{L'auteur a bénéficié du soutien du réseau européen TMR \textit{Arithmetic Algebraic Geometry}
(contrat numéro UE MRTN-CT-2003-504917).}}


\maketitle

\selectlanguage{english}
\begin{abstract}
Let $K$ be a mixed characteristic complete discrete valuation field with perfect residue field $k$.
Let $X$ be a variety over $k$, 
$Y$ be an open of $X$, 
$Y'$ be an open of $Y$ dense in $X$.
We extend Kedlaya's full faithfulness theorem as follows (we do not suppose $Y$ to be smooth): 
the canonical functor 
$F\text{-}\mathrm{Isoc} ^{\dag} (Y,X/K) \to F\text{-}\mathrm{Isoc} ^{\dag} (Y,Y/K) $ 
is fully faithfull.

Suppose now $Y$ smooth. We construct the category of 
partially overcoherent isocrystals over $(Y,X)$ denoted by $\mathrm{Isoc} ^{\dag\dag} (Y,X/K) $
whose objects are some particular arithmetic $\D$-modules. 
Furthermore, we check the equivalence of categories 
$\sp _{(Y,X),+}\,:\,
\mathrm{Isoc} ^{\dag} (Y,X/K) \cong 
\mathrm{Isoc} ^{\dag\dag} (Y,X/K)$.

\end{abstract}

\selectlanguage{frenchb}
\begin{abstract}

Soit $\V$ un anneau de valuation discrète complet d'inégales caractéristiques $(0,p)$, de corps résiduel parfait $k$, de corps des fractions $K$.
Soient $X$ une variété sur $k$, 
$Y$ un ouvert de $X$, 
$Y'$ un ouvert de $Y$ dense dans $X$.
Nous prolongeons le théorème de pleine fidélité de Kedlaya de la manière suivante (en effet, nous ne supposons pas $Y$ lisse): 
le foncteur canonique
$F\text{-}\mathrm{Isoc} ^{\dag} (Y,X/K) \to F\text{-}\mathrm{Isoc} ^{\dag} (Y,Y/K) $ 
est pleinement fidèle.

Supposons à présent $Y$ lisse. Nous construisons la catégorie $\mathrm{Isoc} ^{\dag\dag} (Y,X/K) $ des isocristaux partiellement surcohérents sur $(Y,X)$  dont les objets sont certains $\D$-modules arithmétiques.
De plus, nous vérifions l'équivalence de catégories
$\sp _{(Y,X),+}\,:\,
\mathrm{Isoc} ^{\dag} (Y,X/K) \cong 
\mathrm{Isoc} ^{\dag\dag} (Y,X/K)$.

\end{abstract}

\date
\tableofcontents

\section*{Introduction}
Soit $\V$ un anneau de valuation discrète complet d'inégales caractéristiques $(0,p)$, de corps résiduel parfait $k$, de corps des fractions $K$.
Soient $\PP $
un $\V$-schéma formel  séparé et lisse,
$P$ sa fibre spéciale, 
$T$ un diviseur de $P$, $X $ un sous-schéma fermé de $P$
et $Y := X \setminus T$.
On note alors $(F\text{-})\mathrm{Isoc} ^{\dag} (Y,X/K) $ la catégorie des ($F$-)isocristaux sur $Y$ surconvergents le long de $X \setminus Y$ (voir le livre \cite{LeStum-livreRigCoh} ou \cite{Berig}). Lorsque $X$ est propre, cette catégorie se note $(F\text{-})\mathrm{Isoc} ^{\dag} (Y/K) $ car elle 
ne dépend pas, à isomorphisme canonique près, du choix de $X$.

Lorsque $X$ est lisse, nous avons construit dans \cite{caro-construction}
un foncteur pleinement fidèle noté $\sp _{X \hookrightarrow \PP,T,+}$ (lorsque le diviseur $T$ est vide, on ne l'indique pas)
de la catégorie $(F\text{-})\mathrm{Isoc} ^{\dag} (Y,X/K) $ dans celle 
des $(F\text{-})\D ^\dag _{\PP} (\hdag T) _\Q$-modules cohérents $\E$ à support dans $X$.
On note $(F\text{-})\mathrm{Isoc} ^{\dag \dag} (\PP, T, X/K)$ cette image essentielle. 
La catégorie $(F\text{-})\mathrm{Isoc} ^{\dag \dag} (\PP, T, X/K)$ est stable par le foncteur dual $\D ^\dag _{\PP} (\hdag T) _\Q$-linéaire noté $\DD _{T}$ et, d'après \cite[6.1.4]{caro_devissge_surcoh}, 
un objet $\E$ de $(F\text{-})\mathrm{Isoc} ^{\dag \dag} (\PP, T, X/K)$ est 
$\D ^\dag _{\PP} (\hdag T) _\Q$-surcohérent (en gros, un module cohérent est surcohérent si sa cohérence est préservée par image inverse extraordinaire : voir la définition de \cite{caro_surcoherent} ou sa caractérisation \cite{caro_caract-surcoh}).

Lorsque $Y$ est lisse (mais $X$ quelconque), 
la construction de la catégorie $(F\text{-})\mathrm{Isoc} ^{\dag \dag} (\PP, T, X/K)$ des $(F\text{-})$isocristaux partiellement surcohérents sur $(Y,X)$
( ou sur $(\PP, T, X/K)$ si on veut préciser les choix faits) 
se généralise 
(voir la définition \cite[6.2.1]{caro_devissge_surcoh}) de la manière suivante : 
les objets de $(F\text{-})\mathrm{Isoc} ^{\dag \dag} (\PP, T, X/K)$ sont 
les $(F\text{-})\D ^\dag _{\PP} (\hdag T) _\Q$-modules cohérents $\E$ à support dans $X$
 tels que 
  $\E$ et $\DD _{T} (\E)$ soient $\D ^\dag _{\PP} (\hdag T) _\Q$-surcohérents,
 $\E |\U$ soit dans l'image essentielle de $\sp _{Y\hookrightarrow \U, +}$ (cela a un sens car $Y$ est lisse).
Lorsque $\PP$ est propre, nous avons vérifié en nous ramenant par recollement au cas où $Y$ est affine et lisse (voir \cite[1.3.6 et 2.3.1]{caro-2006-surcoh-surcv}) l'équivalence de catégories 
$\sp _{Y,+}\,:\, F\text{-}\mathrm{Isoc} ^{\dag} (Y/K) \cong F\text{-}\mathrm{Isoc} ^{\dag \dag} (\PP, T, X/K)$.
Dans ce papier, nous étendons cette équivalence au cas général sans structure de Frobenius. 
Comme d'habitude, l'idée de la preuve est de se ramener grâce au théorème de désingularisation de de Jong (voir \cite{Ber-alterationdejong}) au cas où $X$ est lisse via un théorème de descente. 
Néanmoins, comme le théorème de pleine fidélité de Kedlaya du foncteur restriction
$F\text{-}\mathrm{Isoc} ^{\dag} (Y,X/K) \to F\text{-}\mathrm{Isoc} ^{\dag} (Y,Y/K)$ (voir \cite[4.2.1]{kedlaya-semistableII} ou \cite{kedlaya_full_faithfull})
est, si on enlève les structures de Frobenius, une conjecture et comme 
ce théorème est à la base des résultats de finitude de \cite{caro_devissge_surcoh} ou \cite{caro-2006-surcoh-surcv} 
(plus précisément, voir la preuve de \cite[6.3.1]{caro_devissge_surcoh}),
l'obtention de cette extension nécessite un nouvel angle d'approche. 
Cette approche est expliquée plus précisément dans la description des chapitres ci-dessous.

Voici plus précisément le contenu de ce travail.

Nous étudions dans une premième partie la notion d'holonomie sans structure de Frobenius. 
Supposons dans ce paragraphe $P$ intègre. 
Soit $\E$ un $\D ^{\dag } _{\PP , \Q }$-module cohérent. 
Nous définissons la dimension et la codimension de $\E$.
Ces définitions coïncident à celles de Berthelot lorsque $\E$ est muni d'une structure de Frobenius.
L'inégalité de Bernstein reste vraie sans structure de Frobenius : si $\E$ est non nul alors
$ \dim (\E )  \geq \dim P$ (ce qui équivaut à $\mathrm{codim}  (\E)\leq  \dim P$). 
Comme d'habitude, $\E$ est par définition holonome si $\mathrm{codim}  (\E)\geq  \dim P$.
On en déduit alors,
pour tout $i > \dim X$, 
$\mathcal{E} xt ^{i} _{\D ^{\dag } _{\PP , \Q }} (\E, \D ^{\dag } _{\PP , \Q })=0$
(voir \ref{dim-hom-m-dag}).
Le {\og critère homologique d'holonomie\fg} de Virrion (voir \cite{virrion}) reste valable :
$\E$ est holonome si et seulement si, pour tout $i \not = \dim P$, 
$\mathcal{E} xt ^{i} _{\D ^{\dag } _{\PP , \Q }} (\E, \D ^{\dag } _{\PP , \Q })=0$ (\ref{III.3.5}).
Puis nous établissons la version sans structure de Frobenius de \cite[2.2.12]{caro_courbe-nouveau} (c'est d'ailleurs la première étape pour
dévisser en isocristaux surconvergents) :
si pour tout point fermé $x$ de
$P$, $ i ^* _x ( \E ) $ est un $K$-espace vectoriel de dimension finie, 
il existe alors un diviseur $T$ de $X$ tel que $(\hdag T) (\E)$ soit un isocristal surconvergent le long de $T$
(voir \ref{borel9.3}).
On en déduit le {\og critère surcohérent d'holonomie\fg} 
(qui généralise \cite[2]{caro_surholonome} à la situation sans structure de Frobenius) : 
si les espaces de cohomologie de $\DD (\E)$ sont
  à fibres extraordinaires finies (e.g., si $\DD  (\E)$ est $\D ^\dag _{\PP, \Q}$-surcohérent),
  alors $\E$ est holonome (voir \ref{surcoh=>hol}).

Soient $a\,:\, X^{(0)} \to X$ 
un morphisme propre, surjectif, génériquement étale de $k$-variétés intègres, $Y$ un ouvert dense de $X$, 
$\widetilde{Y}$ un ouvert dense de $Y$, 
$\widetilde{j}\,:\,\widetilde{Y} \hookrightarrow X$ l'immersion ouverte correspondante, 
$Y^{(0)}:= a ^{-1} (Y)$, 
$\widetilde{Y}^{(0)}:= a ^{-1} (\widetilde{Y})$. On suppose $Y$ et $Y ^{(0)}$ lisses.
On dispose alors du foncteur canonique pleinement fidèle : 
\begin{equation}
\notag
(*)\hspace{1cm}
(a ^{*}, \widetilde{j} ^{\dag})\,:\,
\mathrm{Isoc} ^{\dag} (Y,X/K)
\to 
\mathrm{Isoc} ^{\dag} (Y^{(0)},X^{(0)}/K) 
\times _{\mathrm{Isoc} ^{\dag} (\widetilde{Y}^{(0)},X^{(0)}/K)} \mathrm{Isoc} ^{\dag} (\widetilde{Y},X/K).
\end{equation}
En effet, 
puisque $Y$ est lisse, cela résulte du théorème de  
pleine fidélité du foncteur extension canonique 
$\widetilde{j} ^{\dag}\,:\,
\mathrm{Isoc} ^{\dag} (Y,X/K)
\to 
\mathrm{Isoc} ^{\dag} (\widetilde{Y},X/K)$
 (voir \cite[4.1.1]{tsumono}, \cite[5.2.1]{kedlaya-semistableI} ou encore \cite{Tsuzuki-NormVar09} pour la version la plus générale). 
Dans la seconde partie, on étend de plus le théorème \cite[4.2.1]{kedlaya-semistableII} ou \cite{kedlaya_full_faithfull} de Kedlaya au cas où $Y$ est non nécessairement lisse, i.e.,  
le foncteur canonique $F\text{-}\mathrm{Isoc} ^{\dag} (Y,X/K) \to F\text{-}\mathrm{Isoc} ^{\dag} (Y,Y/K)$ est pleinement fidèle (voir \ref{EqCat-rig-coro-Frob}). L'idée est de se ramener au cas lisse 
via le théorème de descente propre de Shiho (voir \cite[7.3]{shiho-logRC-RCII}).

Abordons à présent le troisième chapitre. 
Afin d'établir l'équivalence entre isocristaux surconvergents sur $(Y,X)$ et
 isocristaux surcohérents sur $(Y,X)$, l'idée est d'obtenir l'analogue du foncteur $(*)$. 
 On vérifiera cet analogue dans cette partie dans le contexte suivant :
 supposons en outre $X ^{(0)}$ lisse et qu'il existe
$\PP $ un $\V$-schéma formel  séparé et lisse,
$P$ sa fibre spéciale, 
$T$ et $\widetilde{T}$ deux diviseurs de $P$
tels que $X $ soit un sous-schéma fermé de $P$,
$Y = X \setminus T$,
$\widetilde{Y} = X \setminus \widetilde{T}$ ; de même en ajoutant des exposants $(0)$.
De plus, on suppose qu'il existe un morphisme propre et lisse $f\,:\, \PP ^{(0)} \to \PP$ prolongeant $a$
tel que $T ^{(0)}=f ^{-1}(T)$, $\widetilde{T} ^{(0)}=f ^{-1}(\widetilde{T})$, 
le morphisme 
$\widetilde{Y} ^{(0)}\to \widetilde{Y}$ induit par $a$ soit fini et étale
(on obtient ainsi un diagramme de la forme \ref{formel-II-diag}).
D'après \ref{Th2}, 
le foncteur canonique 
\begin{equation}
\notag
(**)\hspace{1cm}
(a ^{+},\, (\hdag \widetilde{T}) )\,:\,
\mathrm{Isoc} ^{*} (\PP, T, X/K) 
\to 
\mathrm{Isoc} ^{*} (\PP ^{(0)}, T ^{(0)}, X^{(0)}/K) 
\times _{\mathrm{Isoc} ^{*} (\PP ^{(0)}, \widetilde{T} ^{(0)}, X^{(0)}/K)}
\mathrm{Isoc} ^{*} (\PP, \widetilde{T}, X/K),
\end{equation}
où $\mathrm{Isoc} ^{*}$ est une catégorie ad hoc 
disposant d'un foncteur image inverse $a ^{+}$ (voir \ref{stabIsoc*inv-ii})
et qui s'avérera égale à $\mathrm{Isoc} ^{\dag \dag}$ (voir \ref{Isoc*=dagdag}),
est un foncteur pleinement fidèle. 
On donne de plus une description de l'image essentielle de ce foncteur
$(**)$ grâce au critère surcohérent d'holonomie et 
au résultat sur la dimension cohomologique des $\D ^{\dag} _{\PP,\Q}$-modules cohérents 
établis dans le premier chapitre.

Dans la quatrième partie, on vérifie l'équivalence de catégories
notée
$\sp _{X \hookrightarrow \PP, T,+}\,:\,\mathrm{Isoc} ^{\dag} (Y,X/K)\cong
\mathrm{Isoc} ^{\dag \dag} (\PP, T, X/K) $ (voir le théorème \ref{diag-eqcat-dagdag-ii}).
La preuve est grosso modo la suivante : grâce au cas déjà connu où $X$ est lisse (voir \cite{caro-construction}),
grâce aux foncteurs pleinement fidèles $(*)$ et $(**)$ ci-dessus et au théorème de désingularisation de de Jong,
on se ramène au cas où il existe une désingularisation finie et étale en dehors du diviseur $T$ (i.e., avec les notations précédentes, on peut supposer
$T= \widetilde{T}$).
Ce cas est traité de la manière suivante (pour plus de détails, voir la preuve de \ref{Isoc*=dagdagFét}). 
D'un côté, on dispose du théorème de descente de Shiho pour les isocristaux partiellement surconvergent.
De l'autre côté, on établit l'analogue pour les isocristaux partiellement surcohérents dans le contexte de la désingularisation finie et étale en dehors du diviseur $T$ (techniquement, celui-ci permet de définir des images directes d'isocristaux surcohérents) :  la donnée d'un isocristal surcohérent sur ($Y ^{(0)}, X ^{(0)}$) muni d'une donnée de descente est équivalente à celle d'un isocristal surcohérent sur $(Y,X)$ (voir \ref{desc-fini-ét}). 
Ces deux théorèmes de descente nous permettent alors d'obtenir l'équivalence de catégories par $\sp _{X \hookrightarrow \PP, \widetilde{T},+}$.

Dans une cinquième partie, nous vérifions que la catégorie $\mathrm{Isoc} ^{\dag \dag} (\PP, T, X/K) $
ne dépend canoniquement que de $(Y,X/K)$ (les équivalences de catégories sont explicites). On la note alors
$\mathrm{Isoc} ^{\dag \dag} (Y, X/K) $. Lorsque $X$ est propre, elle ne dépend pas d'un tel choix de $X$ et se note
$\mathrm{Isoc} ^{\dag \dag} (Y/K) $.

Enfin, dans une dernière partie, on établit que le foncteur explicite $\sp _{Y,+}$ construit à la main dans \cite{caro_devissge_surcoh} lorsque $Y$ est une $k$-variété affine 
induit une équivalence de catégories
$\sp _{Y,+}\,:\, \mathrm{Isoc} ^{ \dag} (Y/K) \cong \mathrm{Isoc} ^{\dag \dag} (Y/K)$, 
ce qui étend le cas déjà connu où l'on disposait d'une structure de Frobenius (voir \cite{caro-2006-surcoh-surcv}).
Cette équivalence de catégories sans structure de Frobenius
permet d'éviter les problèmes de compatibilité à Frobenius.
Par exemple, cette équivalence $\sp _{Y,+}$ sera concrètement utilisée dans la vérification de l'égalité entre holonomie et surholonomie pour les
variétés quasi-projectives dans \cite{caro-stab-holo}. 
Enfin, on obtient la notion de dévissage en isocristaux partiellement surconvergents.

\subsection*{Remerciement}
Je remercie Nobuo Tsuzuki pour m'avoir communiqué son papier \cite{Tsuzuki-NormVar09}, pour une erreur décelée dans une précédente version, son invitation au Japon et notamment à une conférence à Sendai très fructueuse. 
Je remercie aussi Atsushi Shiho pour m'avoir indiqué lors de ce passage au Japon son théorème de descente propre (un ingrédient fondamental de ce travail) et pour une autre erreur trouvée dans une version précédente de ce papier. 
Je remercie Pierre Berthelot et Tomoyuki Abe pour leur contre-exemple sur les inclusions de variétés caractéristiques lorsque l'on fait varier le niveau (plus précisément, voir la remarque \ref{rema-m-geq-m'}).

\section*{Notations et rappels}

\indent 

$\bullet$ On désigne par $\V$ un anneau de valuation discrète complet d'inégale caractéristique $(0,p)$ de corps résiduel parfait $k$, de corps des fractions $K$. Soit $\pi $ une uniformisante de $\V$. On fixe $s\geq 1$ un entier naturel et $F$ désigne la puissance
$s$-ième de l'endomorphisme de Frobenius.

$\bullet$ 
Les $\V$-schémas formels (resp. les $\V$-schémas formels faibles, resp. les $k$-variétés) seront notés via des lettres calligraphique ou gothiques
(resp. des lettres droites surmontés du symbole $\dag$, resp. des lettres droites).
Les fibres spéciales des $\V$-schémas formels lisses seront désignées par les lettres droites correspondantes.
Les $\V$-schémas formels déduits par complétion $p$-adique d'un $\V$-schéma formel faible
seront notés par la lettre calligraphique ou gothique correspondante, e.g., 
si $P ^\dag$ un $\V$-schémas formels faible, on pose
$\PP := \widehat{P ^{\dag}}$ 
et $P:= \PP \otimes _{\V} k$.
De plus, si $u $ est un morphisme de $\V$-schémas formels lisses, le morphisme induit
au niveau des fibres spéciales sera noté $u _0$ ou par abus de notations $u$.

$\bullet$ Les modules sont par défaut à gauche.

$\bullet$  En général, dans les notations qui font intervenir un diviseur (noté par défaut $T$), lorsque celui-ci est vide, on omet alors de l'indiquer.

$\bullet$ Nous reprenons les notations usuelles concernant les opérations cohomologiques de la théorie des $\D$-modules arithmétiques (voir \cite{Beintro2} et \cite[1]{caro_surholonome}) :
soit $f$ : $\PP' \rightarrow \PP$ un morphisme de $\V$-schémas formels lisses.
Le foncteur image inverse extraordinaire par $f$ est noté $f ^{!}$.
Le foncteur image directe par $f$ est noté $f _{+}$.

Soit $X$ un sous-schéma fermé de $P$. On notera par 
$\R \underline{\Gamma} ^\dag _{X}$ le foncteur cohomologique local défini dans \cite[2]{caro_surcoherent} (conjecturalement égal à celui défini par Berthelot dans \cite{Beintro2}). On note $(\hdag X)$ le foncteur localisation en dehors de $X$ (voir \cite[2]{caro_surcoherent}). 
Rappelons que ces deux foncteurs commutent canoniquement aux images directes et images inverses extraordinaires. 

Soit $T$ un diviseur de $P$. 
 Le foncteur  $\DD _{T}$ désigne le foncteur dual $\D ^\dag _{\PP} (\hdag T) _\Q$-linéaire au sens de Virrion (voir \cite{virrion}).
Rappelons que le foncteur $\DD _{T}$ préserve la $\D ^\dag _{\PP} (\hdag T) _\Q$-cohérence (et la perfection) car le faisceau $\D ^\dag _{\PP} (\hdag T) _\Q$ est cohérent et de dimension cohomologique finie (voir \cite{huyghe_finitude_coho}).

\section{Holonomie sans structure de Frobenius et critères d'holonomie}

Nous vérifions dans cette section que la notion d'holonomie s'étend de façon naturelle sans structure de Frobenius.
Cette section sera utilisée via \ref{dim-hom-m-dag} et \ref{borel9.3} pour valider le théorème \ref{Th2} (ou plutôt le lemme \ref{lemm2Th2-cas lisse}).
Dans ce chapitre, $\X$ sera un $\V$-schéma formel lisse.

\subsection{Préliminaire sur les variétés caractéristiques de niveau $m$}

On suppose $\X$ intègre.

\begin{nota}
\label{var-car-nivm}

Soit $\E ^{(m)}$ un $\smash{\widehat{\D}} ^{(m)} _{\X,\,\Q}$-module cohérent.

$\bullet$
D'après \cite[5.2.5]{Beintro2}, on lui associe de manière canonique une variété caractéristique notée $\mathrm{Car} ^{(m)}  (\E  ^{(m)})$. 
Notons $F ^{m}$ la puissance $m$-ième de l'endomorphisme de Frobenius $X \to X$.
Posons $\E ^{(0,m)} := (F ^{m} ) ^{\flat}\smash{\widehat{\D}} ^{(0)} _{\X,\,\Q} \otimes _{\smash{\widehat{\D}} ^{(m)} _{\X,\,\Q}} \E  ^{(m)}$.
D'après \cite[2.1.4.1]{Beintro2}, 
$\E ^{(0,m)}$ est un $\smash{\widehat{\D}} ^{(0)} _{\X,\,\Q}$-module cohérent induisant l'isomorphisme canonique 
$F ^{m*} (\E ^{(0,m)}) \riso \E ^{(m)}$, 
où  et $F ^{m*}$ le foncteur image inverse par $F ^{m}$.
Via \cite[5.2.2.1]{Beintro2}, on identifiera $\mathrm{Car} ^{(0)} (\E ^{(0,m)})$ avec $ \mathrm{Car} ^{(m)} (\E ^{(m)})$ (en effet, si on choisit un modèle sans $p$-torsion de $\E ^{(0,m)}$ alors son image par $F ^{m*} $ donne un modèle sans $p$-torsion de $\E ^{(m)}$).

$\bullet$
La {\og dimension de niveau $m$ de $\E ^{(m)}$\fg} est définie en posant 
$\mathrm{dim} ^{(m)} (\E^{(m)})  := \dim \mathrm{Car} ^{(m)} (\E^{(m)})$.
La {\og codimension de niveau $m$ de $\E ^{(m)}$\fg} est définie en posant 
$\mathrm{codim} ^{(m)} (\E^{(m)}) := 2 \dim X - \mathrm{dim} ^{(m)} (\E^{(m)}) $.

S'il n'y a pas d'ambiguïté sur le niveau (e.g. il se peut a priori qu'un $\smash{\widehat{\D}} ^{(m+1)} _{\X,\,\Q}$-module cohérent
soit aussi un $\smash{\widehat{\D}} ^{(m)} _{\X,\,\Q}$-module cohérent), nous noterons plus simplement
$\mathrm{dim} (\E^{(m)})$ et $\mathrm{codim} (\E^{(m)})$.
De plus, on dispose des encadrements $0\leq \mathrm{dim} (\E^{(m)}) \leq 2 \dim X$
et $0\leq \mathrm{codim} (\E^{(m)})\leq 2 \dim X$.

$\bullet$
Avec l'équivalence de catégories canonique entre $\smash{\widehat{\D}} ^{(m)} _{\X,\,\Q}$-modules à gauche cohérents 
et $\smash{\widehat{\D}} ^{(m)} _{\X,\,\Q}$-modules à droite cohérents induit par le foncteur $-\otimes _{\O _{\X}} \omega _{\X}$ (on rappelle que 
$\omega _{\X}:= \wedge ^{d _X} \Omega _{\X/\V} $), on obtient des constructions et définitions analogues pour les modules à droite.
\end{nota}

\begin{vide}
\label{3.1.3.(ii)Beintro}
D'après \cite[3.1.3.(ii)]{Beintro2}, si $\X$ est affine alors la dimension homologique de 
$\Gamma (\X, \smash{\widehat{\D}} ^{(m)} _{\X,\,\Q})$ est comprise entre $\dim X$ et $2 \dim X$.
\end{vide}

\begin{prop}
\label{III.2.4Virrion}
Soit $\E ^{(m)}$ un $\smash{\widehat{\D}} ^{(m)} _{\X,\,\Q}$-module cohérent. Alors : 
\begin{enumerate}
\item  \label{m-codimgeqi}
$\mathrm{codim} ^{(m)} (\mathcal{E} xt ^{i} _{\smash{\widehat{\D}} ^{(m)} _{\X,\,\Q}} (\E ^{(m)}, \smash{\widehat{\D}} ^{(m)} _{\X,\,\Q}))
\geq i$, pour tout $i \geq 0$.

\item
\label{m-codimgeqii}
$\mathcal{E} xt ^{i} _{\smash{\widehat{\D}} ^{(m)} _{\X,\,\Q}} (\E ^{(m)}, \smash{\widehat{\D}} ^{(m)} _{\X,\,\Q})=0$, pour tout $i < \mathrm{codim} ^{(m)} (\E ^{(m)})$. 
 
\end{enumerate}

\end{prop}

\begin{proof}
On procède de manière analogue à \cite[III.2.4]{virrion} : soit $\E ^{(0,m)}:= (F ^{m} ) ^{\flat}\smash{\widehat{\D}} ^{(0)} _{\X,\,\Q} \otimes _{\smash{\widehat{\D}} ^{(m)} _{\X,\,\Q}} \E  ^{(m)}$.
D'après Virrion (voir \cite[II.3]{virrion}), on dispose de l'isomorphisme canonique
\begin{equation}
\label{FflatF*}
(F ^{m})^{\flat} 
\mathcal{E} xt ^{i} _{\smash{\widehat{\D}} ^{(0)} _{\X,\,\Q}} (\E ^{(0,m)}, \smash{\widehat{\D}} ^{(0)} _{\X,\,\Q})
\riso 
\mathcal{E} xt ^{i} _{\smash{\widehat{\D}} ^{(m)} _{\X,\,\Q}} (F ^{m*} \E ^{(0,m)}, \smash{\widehat{\D}} ^{(m)} _{\X,\,\Q})
\riso 
\mathcal{E} xt ^{i} _{\smash{\widehat{\D}} ^{(m)} _{\X,\,\Q}} (\E ^{(m)}, \smash{\widehat{\D}} ^{(m)} _{\X,\,\Q}) .
\end{equation}
Modulo l'identification \cite[5.2.2.1]{Beintro2}, il en résulte que
$\mathrm{Car} ^{(0)} (\mathcal{E} xt ^{i} _{\smash{\widehat{\D}} ^{(0)} _{\X,\,\Q}} (\E ^{(0,m)}, \smash{\widehat{\D}} ^{(0)} _{\X,\,\Q}))$
et
$\mathrm{Car} ^{(m)} (\mathcal{E} xt ^{i} _{\smash{\widehat{\D}} ^{(m)} _{\X,\,\Q}} (\E ^{(m)}, \smash{\widehat{\D}} ^{(m)} _{\X,\,\Q}) )$
sont isomorphes.
Grâce à \cite[III.2.3.(i)]{virrion}, il en résulte \ref{m-codimgeqi}.

De plus, comme $\mathrm{Car} ^{(0)} (\E ^{(0,m)})= \mathrm{Car} ^{(m)} (\E ^{(m)})$,
alors $\mathrm{codim} (\E ^{(0,m)} )=\mathrm{codim}  (\E ^{(m)})$. 
D'après \cite[III.2.3.(ii)]{virrion}, on en déduit que
$\mathcal{E} xt ^{i} _{\smash{\widehat{\D}} ^{(m)} _{\X,\,\Q}} (\E ^{(m)}, \smash{\widehat{\D}} ^{(m)} _{\X,\,\Q})=0$
pour tout $i < \mathrm{codim} ^{(m)} (\E)$.
\end{proof}

\begin{prop}
\label{Vir-III.3.5}
Soit $\E ^{(m)}$ un $\smash{\widehat{\D}} ^{(m)} _{\X,\,\Q}$-module cohérent non nul. Alors: 
\begin{equation}
\mathrm{codim} ^{(m)} (\E ^{(m)}) = \min \left \{i \,|\, \mathcal{E} xt ^{i}  _{\smash{\widehat{\D}} ^{(m)} _{\X,\,\Q}} ( \E ^{(m)},\smash{\widehat{\D}} ^{(m)} _{\X,\,\Q}) \not = 0\right \}.
\end{equation}
\end{prop}

\begin{proof}
Il s'agit de reprendre \cite[III.3.4--5]{virrion} en remplaçant $\D ^{\dag}$ par $\widehat{\D} ^{(m)}$ et (grâce à \ref{3.1.3.(ii)Beintro}) $\dim X$ par $2\dim X$. 
Pour la commodité du lecteur, voici les grandes lignes de la preuve : 
posons $k:=\mathrm{codim} ^{(m)} (\E ^{(m)}) $.
Grâce à l'isomorphisme de bidualité, on dispose de la suite spectrale 
$$E _{2} ^{i,j}= \mathcal{E} xt ^{i}  _{\smash{\widehat{\D}} ^{(m)} _{\X,\,\Q}} ( \mathcal{E} xt ^{-j}  _{\smash{\widehat{\D}} ^{(m)} _{\X,\,\Q}} ( \E ^{(m)},\smash{\widehat{\D}} ^{(m)} _{\X,\,\Q}),\smash{\widehat{\D}} ^{(m)} _{\X,\,\Q})
\Rightarrow 
E ^{N}= \mathcal{H} ^{i+j} (\E ^{(m)}).
$$
D'après \ref{III.2.4Virrion}.\ref{m-codimgeqii}, $E _{2} ^{i,-i}=0$ pour $i< k$. 
Par l'absurde, supposons $E _{2} ^{k,-k}=0$.
En outre, via \ref{III.2.4Virrion}.\ref{m-codimgeqi}, $\mathrm{codim}(E _{2} ^{i,-i})\geq k +1$ pour $i \geq k+1$.
Comme la filtration de $\E ^{(m)}$ induite par la suite spectrale est finie, on obtient alors :
$\mathrm{codim} (\E ^{(m)} ) \geq k +1$, ce qui est absurde. 
\end{proof}

\begin{coro}
\label{m-geq-m'}
Soit $\E ^{(m)}$ un $\smash{\widehat{\D}} ^{(m)} _{\X,\,\Q}$-module cohérent. 
Pour $m' \geq m$, posons $\E ^{(m')} :=\smash{\widehat{\D}} ^{(m')} _{\X,\,\Q} \otimes _{\smash{\widehat{\D}} ^{(m)} _{\X,\,\Q}} \E ^{(m)}$.
On dispose de l'inégalité $\mathrm{codim} ^{(m')} (\E ^{(m')}) \geq \mathrm{codim} ^{(m)} (\E ^{(m)})$.
\end{coro}

\begin{proof}
Si $\E ^{(m')}=0$ alors l'assertion est immédiate. Supposons donc $\E ^{(m')} \not =0$.
Comme l'extension $\smash{\widehat{\D}} ^{(m)} _{\X,\,\Q} \to \smash{\widehat{\D}} ^{(m')} _{\X,\,\Q}$ est plate, pour tout $i < \mathrm{codim} ^{(m)} (\E ^{(m)})$, on obtient :
$$0=\mathcal{E} xt ^{i}  _{\smash{\widehat{\D}} ^{(m)} _{\X,\,\Q}} ( \E ^{(m)},\smash{\widehat{\D}} ^{(m)} _{\X,\,\Q}) \otimes _{\smash{\widehat{\D}} ^{(m)} _{\X,\,\Q}} \smash{\widehat{\D}} ^{(m')} _{\X,\,\Q}
\riso 
\mathcal{E} xt ^{i}  _{\smash{\widehat{\D}} ^{(m')} _{\X,\,\Q}} ( \E ^{(m')},\smash{\widehat{\D}} ^{(m')} _{\X,\,\Q})
$$
D'où le résultat via \ref{Vir-III.3.5}.
\end{proof}

\begin{rema}
\label{rema-m-geq-m'}
Avec les notations de \ref{m-geq-m'},
d'après des contre-exemples de Berthelot et de Abe, il n'existe pas d'inclusion
entre les variétés caractéristiques de respectivement $\E ^{(m)}$ et $\E ^{(m')}$.
\end{rema}

\subsection{Inégalité de Bernstein, holonomie et son critère homologique}

Pour simplifier l'exposé, supposons $X $ intègre (voir aussi la remarque \ref{rema-Xgen}).

\begin{nota}
\label{var-car-dag}
$\bullet$ 
Soit  $\E$ un $\D ^{\dag } _{\X , \Q }$-module cohérent. 
Par commodité notons $\mathrm{niv} (\E)$ le plus petit entier $m\geq 0$ 
tel qu'il existe un $\smash{\widehat{\D}} ^{(m)} _{\X,\,\Q}$-module cohérent 
$\E ^{(m)}$ et un isomorphisme $\D ^{\dag } _{\X , \Q }$-linéaire de la forme 
$ \D ^{\dag } _{\X , \Q } \otimes _{\smash{\widehat{\D}} ^{(m )} _{\X,\,\Q}} \E ^{(m)} \riso \E$.
D'après \cite[3.6.2]{Be1}, cet entier est bien défini. De plus, 
pour tout entier $m\geq \mathrm{niv} (\E)$, 
un tel $\smash{\widehat{\D}} ^{(m)} _{\X,\,\Q}$-module cohérent 
ne dépend, à isomorphisme non canonique près, que du niveau $m$ choisi. On le notera $\E ^{(m)}$.
On définit alors, à isomorphisme non canonique près, la variété caractéristique de niveau $m$ associée à $\E$ en posant 
$\mathrm{Car} ^{(m)}  (\E):= \mathrm{Car} ^{(m)}  (\E^{(m)})$.

$\bullet$ 
La {\og dimension de niveau $m$ de $\E$\fg} est définie en posant 
$\mathrm{dim} ^{(m)} (\E)  := \dim \mathrm{Car} ^{(m)} (\E)$.
Il résulte de \ref{m-geq-m'} que, pour tous $m' \geq m \geq \mathrm{niv} (\E)$, on a $\mathrm{dim} ^{(m)} (\E)  \geq \mathrm{dim} ^{(m')} (\E) $.

$\bullet$ On note $\mathrm{dim}  (\E) $ le minimum des $\mathrm{dim} ^{(m)} (\E) $, i.e., la limite de la suite stationnaire 
$(\mathrm{dim} ^{(m)} (\E) ) _{m\geq \mathrm{niv}(\E)}$.
On pose aussi $\mathrm{codim}  (\E) := 2 \dim X - \mathrm{dim} (\E)$.
\end{nota}

\begin{rema}
Avec les notations de \ref{var-car-dag}, 
si $\E$ est un $\D ^{\dag } _{\X , \Q }$-module cohérent muni d'une structure de Frobenius, pour tout niveau $m$ (assez grand si $p=2$), 
les variétés caractéristiques $\mathrm{Car} ^{(m)}  (\E)$ ne dépendent pas du niveau $m$ et sont égales à la variété caractéristique 
définie par Berthelot (voir \cite[5.2.7]{Beintro2}).
\end{rema}

Étendons d'abord via les deux propositions suivantes les résultats de la section précédente:
\begin{prop}
\label{III.2.4Virriondag}
Soit $\E$ un $\D ^{\dag } _{\X , \Q }$-module cohérent.  Alors, pour tout $m \geq \mathrm{niv}(\E)$, on a  : 
\begin{enumerate}
\item  \label{codimgeqi}
$\mathrm{codim} ^{(m)} (\mathcal{E} xt ^{i} _{\D ^{\dag } _{\X , \Q }} (\E, \D ^{\dag } _{\X , \Q }))
\geq i$, pour tout $i \geq 0$.

\item \label{codimgeqii} $\mathcal{E} xt ^{i} _{\D ^{\dag } _{\X , \Q }} (\E, \D ^{\dag } _{\X , \Q })=0$, pour tout $i < \mathrm{codim} ^{(m)} (\E)$. 
 
\end{enumerate}

\end{prop}

\begin{proof}
Comme l'extension $\smash{\widehat{\D}} ^{(m)} _{\X,\,\Q}\to \D ^{\dag } _{\X , \Q }$ est plate (voir \cite{Be1}), comme le foncteur dual commute aux extensions (voir \cite[I.4]{virrion}), on obtient le second isomorphisme :
\begin{equation}
\label{Exti-m->dag}
\mathcal{E} xt ^{i} _{\D ^{\dag } _{\X , \Q }} (\E, \D ^{\dag } _{\X , \Q })
\riso
\mathcal{E} xt ^{i} _{\D ^{\dag } _{\X , \Q }} (\D ^{\dag } _{\X , \Q }
\otimes _{\smash{\widehat{\D}} ^{(m)} _{\X,\,\Q}}
\E ^{(m)}, \D ^{\dag } _{\X , \Q })
\riso 
\mathcal{E} xt ^{i} _{\smash{\widehat{\D}} ^{(m)} _{\X,\,\Q}} (\E ^{(m)}, \smash{\widehat{\D}} ^{(m)} _{\X,\,\Q}) 
\otimes _{\smash{\widehat{\D}} ^{(m)} _{\X,\,\Q}}
\D ^{\dag } _{\X , \Q }.
\end{equation}

D'après \ref{III.2.4Virrion}, il en résulte que $\mathcal{E} xt ^{i} _{\D ^{\dag } _{\X , \Q }} (\E, \D ^{\dag } _{\X , \Q })= 0$
pour tout $i < \mathrm{codim} ^{(m)} (\E)$.
Par définition des variétés caractéristiques de niveau $m$,
l'isomorphisme \ref{Exti-m->dag} se traduit par 
$\mathrm{Car} ^{(m)} (\mathcal{E} xt ^{i} _{\D ^{\dag } _{\X , \Q }} (\E, \D ^{\dag } _{\X , \Q }))
=
\mathrm{Car} ^{(m)}(\mathcal{E} xt ^{i} _{\smash{\widehat{\D}} ^{(m)} _{\X,\,\Q}} (\E ^{(m)}, \smash{\widehat{\D}} ^{(m)} _{\X,\,\Q}))$.
D'après \ref{III.2.4Virrion}, il en résulte alors \ref{codimgeqi}.
\end{proof}

\begin{prop}
\label{Vir-III.3.5dag}
Soit $\E$ un $\D ^{\dag } _{\X , \Q }$-module cohérent non nul. Alors: 
\begin{equation}
\mathrm{codim} (\E) = \min \left \{i \,|\, 
\mathcal{E} xt ^{i} _{\D ^{\dag } _{\X , \Q }} (\E, \D ^{\dag } _{\X , \Q }) \not = 0\right \}.
\end{equation}
\end{prop}

\begin{proof}
Il suffit d'utiliser \ref{Vir-III.3.5} et les isomorphismes \ref{Exti-m->dag}.
\end{proof}

\begin{theo}
[Inégalité de Bernstein]
Pour tout $\D ^{\dag } _{\X , \Q }$-module cohérent $\E$  non nul, 
on a 
\begin{gather}
\label{Inég-Bernstein}
\dim \mathrm{Car}  ^{(m)} (\E ^{(m)}) 
\geq 
\dim (\E )  \geq \dim X.
\end{gather}
\end{theo}

\begin{proof}
La première inégalité est déjà connue. Il reste à vérifier 
$\dim (\E) \geq \dim X$. 
La preuve est identique à celle de \cite[5.3.4]{Beintro2}. Pour la commodité du lecteur, nous donnons néanmoins ici ses grandes lignes : 
comme $X$ est somme de ses composantes irréductibles, on peut supposer $X$ intègre. 
On procède alors à une récurrence sur la dimension de $X$. 
Supposons que le support de $\E$ soit $X$.  Soit $\smash{\overset{_{\circ}}{\E}}^{(m )}$ un $\smash{\widehat{\D}} _{\X} ^{(m )}$-module cohérent  sans $p$-torsion tel que 
$\smash{\widehat{\D}} ^{(m)} _{\X , \Q } \otimes _{\smash{\widehat{\D}} ^{(m)} _{\X}} \smash{\overset{_{\circ}}{\E}}^{(m ) } \riso \E ^{(m )} $.
Alors le support de $\mathrm{gr} (\smash{\overset{_{\circ}}{\E}}^{(m)}/\pi  \smash{\overset{_{\circ}}{\E}}^{(m)})$ 
contient $X$ et l'inégalité est vérifiée.
Sinon, le support $Z$ de $\E$ est un fermé de dimension strictement plus petite à celle de $X$. 
Quitte à changer $X$ par un ouvert affine $U$ tel que $Z \cap U$ soit lisse et dense dans $Z$, on peut supposer $Z$ lisse
et que $Z \subset X$ se relève en une immersion fermée
$u \,:\, \ZZ \hookrightarrow \X$ de $\V$-schémas formels lisses.
Comme $\E ^{(m )}$ est à support dans $Z$, d'après \cite[5.3.2]{Beintro2}, quitte à augmenter $m$, il existe 
un $\smash{\widehat{\D}} ^{(m)} _{\ZZ,\,\Q}$-module cohérent 
$\FF ^{(m)}$ tel que $ u _{+} ^{(m)} (\FF ^{(m)}) \riso \E ^{(m )}$.
On pose $\E ^{(0,m)}:= (F ^{m} ) ^{\flat}\smash{\widehat{\D}} ^{(0)} _{\X,\,\Q} \otimes _{\smash{\widehat{\D}} ^{(m)} _{\X,\,\Q}} \E  ^{(m)}$
et 
$\FF ^{(0,m)}:= (F ^{m} ) ^{\flat}\smash{\widehat{\D}} ^{(0)} _{\X,\,\Q} \otimes _{\smash{\widehat{\D}} ^{(m)} _{\X,\,\Q}} \FF  ^{(m)}$
Soit  
 $\smash{\overset{_{\circ}}{\FF}}^{(0,m )}$ un $\smash{\widehat{\D}} _{\X} ^{(0)}$-module cohérent  sans $p$-torsion tel que 
$\smash{\widehat{\D}} ^{(0)} _{\X , \Q } \otimes _{\smash{\widehat{\D}} ^{(0)} _{\X}} \smash{\overset{_{\circ}}{\FF}}^{(0,m ) } \riso \FF ^{(0,m )} $.
On pose 
$\smash{\overline{\FF}}^{(0,m )}:=\smash{\D} ^{(0)} _{X } \otimes _{\smash{\widehat{\D}} ^{(0)} _{\X}} \smash{\overset{_{\circ}}{\FF}}^{(0,m ) } $.
Comme dans le cas classique (voir ce qui suit \cite[5.3.3]{Beintro2}),
on vérifie la formule 
$\mathrm{Car}  ^{(0)} (u ^{(0)} _{+}(\smash{\overline{\FF}}^{(0,m )})) = v ( q ^{-1} \mathrm{Car}  ^{(0)} (\smash{\overline{\FF}}^{(0,m )}) )$, où 
$ T ^{*} Z \overset{q}{\longleftarrow} Z \times _{X} T ^{*} X \overset{v}{\longrightarrow}  T ^{*} X $ sont les applications canoniques.
Comme $u ^{(0)} _{+}(\FF^{(0,m )})$ donne un modèle sans $p$-torsion de $\E ^{(0,m)}$, 
on conclut alors par hypothèse de récurrence en remarquant que $q$ est plat et $v$ est une immersion fermée.
\end{proof}

\begin{theo}
\label{dim-hom-m-dag}
Soit 
$\E$ un $\D ^{\dag } _{\X , \Q }$-module cohérent. 
Alors, pour tout $i > \dim X$, 
$$
\mathcal{E} xt ^{i} _{\D ^{\dag } _{\X , \Q }} (\E, \D ^{\dag } _{\X , \Q })=0.$$
\end{theo}

\begin{proof}
 Par l'absurde, supposons qu'il existe $i > \dim X$ tel que $\mathcal{E} xt ^{i} _{\D ^{\dag } _{\X , \Q }} (\E, \D ^{\dag } _{\X , \Q })\not =0$.
L'inégalité de Bernstein se traduit par
$\mathrm{codim} ^{(m)} (\mathcal{E} xt ^{i} _{\D ^{\dag } _{\X , \Q }} (\E, \D ^{\dag } _{\X , \Q }))\leq \dim X$, ce qui contredit 
\ref{III.2.4Virriondag}.\ref{codimgeqi}.
\end{proof}

\begin{defi}
\label{defi-hol}
Soit $\E$ un $\D ^{\dag } _{\X , \Q }$-module cohérent. 
On dit que $\E$ est holonome si $\mathrm{codim}  (\E)\geq  \dim X$. 
\end{defi}

\begin{prop}
\label{III.3.5}
Soit $\E$ un $\D ^{\dag } _{\X , \Q }$-module cohérent. Alors $\E$ est holonome si et seulement si, pour tout $i \not = 0$, 
$\mathcal{E} xt ^{i} _{\D ^{\dag } _{\X , \Q }} (\E, \D ^{\dag } _{\X , \Q }) [\dim X]=0$.
\end{prop}

\begin{proof}
Cela se vérifie comme \cite[III.4.2]{virrion} en utilisant 
\ref{Vir-III.3.5dag} et de \ref{dim-hom-m-dag}.
\end{proof}

\begin{rema}
\label{rema-Xgen}
Lorsque $\X$ n'est plus intègre, la définition \ref{defi-hol} et le critère \ref{III.3.5} sont encore vrais de la manière suivante. 
Comme $\X$ est lisse, $\X$ est somme de ses composantes connexes notées $(\X _{r}) _r$.
Soit $\E$ un $\D ^{\dag } _{\X , \Q }$-module cohérent. 
On désigne alors par $\dim (\E) $ (resp. $\mathrm{codim} (\E)$, resp. $\dim X$) le faisceau sur $X$ localement constant 
égal à $\dim (\E |\X _{r}) $ (resp. $\mathrm{codim} (\E|\X _{r})$, resp. $\dim X _{r}$) sur $X _{r}$.
\end{rema}

La proposition \ref{III.3.5} nous permet d'étendre de manière naturelle la notion d'holonomie de la manière suivante :

\begin{defi}
\label{defi-hol-div}
Soient $T$ un diviseur de $X$ et $\E$ un $\D ^{\dag } _{\X } (\hdag T) _{\Q }$-module cohérent. 
On dit que $\E$ est {\og $\D ^{\dag } _{\X } (\hdag T) _{\Q }$-holonome\fg} ou {\og holonome en tant que 
$\D ^{\dag } _{\X } (\hdag T) _{\Q }$-module\fg} 
si, pour tout $i \not = 0$, 
$\mathcal{H}^{i} \DD _{T} (\E) =0$. 
\end{defi}

\begin{prop}
\label{hol-suiteexacte}
Soient $T$ un diviseur de $X$ 
et 
$0\to \E' \to \E \to \E'' \to 0$ une suite exacte de $\D ^{\dag } _{\X } (\hdag T) _{\Q }$-modules cohérents. 
Alors $\E$ est $\D ^{\dag } _{\X } (\hdag T) _{\Q }$-holonome si et seulement si 
$\E'$ et $\E''$ sont $\D ^{\dag } _{\X } (\hdag T) _{\Q }$-holonomes.
\end{prop}

\begin{proof}
Si $\E'$ et $\E''$ sont $\D ^{\dag } _{\X } (\hdag T) _{\Q }$-holonomes, il est immédiat que $\E$ l'est aussi.
Traitons la réciproque. On suppose ainsi que $\E$ est $\D ^{\dag } _{\X } (\hdag T) _{\Q }$-holonome.
Par \cite[4.3.12]{Be1}, on se ramène au cas où le diviseur $T$ est vide.
Pour $m$ assez grand, via \cite[5.3.2]{Be1} (et \cite[3.4]{Be1} pour l'existence d'un modèle de $\E''$), 
on se ramène à supposer que la suite exacte $0\to \E' \to \E \to \E'' \to 0$
provient par extension d'une suite exacte de $\widehat{\D} ^{(m)} _{\X}$-modules cohérents sans $p$-torsion.
En considérant la suite exacte déduite par réduction modulo $\pi$,
la proposition résulte alors de \cite[5.2.4]{Beintro2} (et de la proposition \ref{III.3.5}).
\end{proof}

\begin{prop}
\label{Stab-hol-dual}
Soit $T$ un diviseur de $X$. 
Le foncteur $\DD _{T}$ induit une auto-équivalence de la catégories des 
$\D ^{\dag } _{\X } (\hdag T) _{\Q }$-modules holonomes.
\end{prop}

\begin{proof}
Cela résulte aussitôt de l'isomorphisme de bidualité et de notre définition \ref{defi-hol-div}.
\end{proof}

\subsection{Surconvergence générique d'un $\D$-module cohérent à fibres finies, second critère d'holonomie}

Le résultat principal de cette section est le théorème 
\ref{borel9.3}. 
Il correspond à une version sans structure de Frobenius du théorème \cite[2.2.17]{caro_courbe-nouveau}.
On établira à partir de ce théorème que les complexes partiellements surcohérents se dévissent en isocristaux surconvergents (voir \ref{caradev}).

Nous avons incorporé ce théorème \ref{borel9.3} dans ce chapitre  
car il nous permet d'établir le critère d'holonomie \ref{surcoh=>hol} que l'on utilisera  
pour vérifier le théorème \ref{Th2} (très précisément au cours de la preuve du lemme \ref{lemm2Th2-cas lisse}).

Les ingrédients qui nous avaient permis d'obtenir \cite[2.2.17]{caro_courbe-nouveau} seront réutilisés. Comme nous nous contenterons de les citer, le lecteur pourra commencer par la section \cite[2.2]{caro_courbe-nouveau} avant de lire cette section.

\begin{nota}
Pour tout point fermé $x$ de $X$, on notera par
$i _x$ : $ \Spf \V (x)\hookrightarrow \X$, un relèvement de l'immersion fermée canonique induite par $x$.
On notera $i _{x} ^{*}$ le foncteur image inverse (isomorphe à $\mathcal{H} ^{d _{X}} \circ i _{x} ^{!}$).
On remarque que $\V (x)$ est un anneau de valuation discrète complet d'inégales caractéristiques $(0,p)$.
Son corps des fractions sera désigné par $K(x)$.
\end{nota}

\begin{lemm}
\label{isoc-dense-libre}
Soit $\E$ un $\D ^{\dag } _{\X , \Q }$-module cohérent, $\O _{\X,\,\Q}$-cohérent.
Il existe alors un ouvert dense $\Y$ de $\X$ tel que $\E |\Y$ soit un $\O _{\Y,\,\Q}$-module libre de type fini. 
\end{lemm}

\begin{proof}
On peut supposer $\X$ affine et muni de coordonnées locales. 
Soit $m _{0}$ un entier. D'après \cite[3.1]{Be0}, il existe un $\smash{\widehat{\D}} _{\X} ^{(m _{0} )}$-module cohérent $\smash{\overset{_{\circ}}{\E}}^{(m _{0})}$ sans $p$-torsion, $\O _{\X}$-cohérent et un isomorphisme $\D ^{\dag } _{\X , \Q }$-linéaire  : 
$ \D ^{\dag } _{\X , \Q } \otimes _{\smash{\widehat{\D}} ^{(m _0)} _{\X}} \smash{\overset{_{\circ}}{\E}}^{(m _{0} ) } \riso \E $.
De manière analogue à \cite[VII.9.3]{borel}, on vérifie qu'il existe alors un ouvert affine et dense $\Y$ de $\X$ tel que $\O _X \otimes _{\O_{\X}} \smash{\overset{_{\circ}}{\E}}^{(m _{0})} |Y $ soit un $\O _{Y}$-module libre (forcément de type fini).
On termine la preuve via par exemple \cite[2.2.16]{caro_courbe-nouveau}.

\end{proof}

Soit $\E$ un $\D ^{\dag } _{\X , \Q }$-module cohérent tel que, pour tout point fermé $x$ de $X$, le $K$-espace vectoriel $ i ^* _x ( \E ) $ est de dimension finie. Un problème supplémentaire technique qui apparaît contrairement à la situation avec structure de Frobenius est que lorsque le niveau $m$ s'élève, l'ouvert au-dessus duquel le $\smash{\widehat{\D}} ^{(m)} _{\X,\,\Q}$-module cohérent associé à $\E$ devient $\O _{\X,\Q}$-cohérent 
rétrécit a priori. 
Le phénomène de contagiosité du lemme \ref{m->m'} ci-dessous nous permettra de résoudre cet obstacle.

\begin{lemm}
\label{m->m'}
On suppose $\X$ affine et muni de coordonnées locales. 
Soient $m'\geq m$ deux entiers, $\U$ un ouvert affine dense dans $\X$, 
$\E$ un $\smash{\widehat{\D}} ^{(m)} _{\X,\,\Q}$-module cohérent 
tel que $\E$ soit un $\O _{\X,\Q}$-module projectif de type fini. 

On suppose que $\Gamma (\U,\E )$ est muni d'une structure de $\Gamma (\U,\smash{\widehat{\D}} ^{(m')} _{\X,\,\Q})$-module prolongeant sa structure de $\Gamma (\U,\smash{\widehat{\D}} ^{(m)} _{\X,\,\Q})$-module. 
Alors $\Gamma (\X,\E )$ est muni d'une et d'une seule structure de $\Gamma (\X,\smash{\widehat{\D}} ^{(m')} _{\X,\,\Q})$-module prolongeant sa structure de $\Gamma (\X,\smash{\widehat{\D}} ^{(m)} _{\X,\,\Q})$-module.
\end{lemm}

\begin{proof}
L'unicité provient du fait que $\Gamma (\X,\smash{\widehat{\D}} ^{(m)} _{\X,\,\Q})$ est dense dans $\Gamma (\X,\smash{\widehat{\D}} ^{(m')} _{\X,\,\Q})$ (e.g. voir \cite[2.2.9]{caro_courbe-nouveau}).
Soient $e \in \Gamma (\X,\E )$ et $1 \otimes e$ l'élément induit de $\Gamma (\U,\E )$.
Comme les éléments $\underline{\partial} ^{<\underline{k} >_{(m')}} $ sont de norme $1$ dans $\Gamma (\X,\smash{\widehat{\D}} ^{(m')} _{\X,\,\Q})$ (pour la norme de Gauss), 
la famille $\underline{\partial} ^{<\underline{k} >_{(m')}} \cdot (1\otimes e)$ avec $\underline{k}$ parcourant $\N ^{d}$ est bornée pour la topologie de $\Gamma (\U,\E )$ induite par sa structure de $\Gamma (\U,\smash{\widehat{\D}} ^{(m')} _{\X,\,\Q})$-module de type fini.
Or, d'après \cite[4.1.2]{Be1}, cette topologie est identique à celle induite par la structure de $\Gamma (\U , \O _{\X,\Q})$-module de type fini sur $\Gamma (\U,\E )$. 

Soient $N$ un entier et $\FF$ un $\O _{\X,\Q}$-module cohérent tel que $\E \oplus \FF \riso \O _{\X,\Q} ^{N}$.
On obtient le diagramme commutatif : 
$$\xymatrix @R=0,3cm {
{\Gamma (\X , \O _{\X,\Q}) ^{N}} 
\ar[r] ^-{\sim}
\ar@{^{(}->}[d] ^-{}
&
{\Gamma (\X , \E) \oplus \Gamma (\X , \FF)} 
\ar@{->>}[r] ^-{}
& 
{\Gamma (\X , \E) } 
\ar@{^{(}->}[d] ^-{}
\\ 
{\Gamma (\U , \O _{\X,\Q})^{N}} 
\ar[r] ^-{\sim}
&
{\Gamma (\U , \E) \oplus \Gamma (\U , \FF)} 
\ar@{->>}[r] ^-{}
&
{\Gamma (\U , \E) }
}$$
dont la flèche injective de droite envoie $\underline{\partial} ^{<\underline{k} >_{(m')}}.e$ sur $\underline{\partial} ^{<\underline{k} >_{(m')}}.(1\otimes e)$
(cela a bien un sens car $\underline{\partial} ^{<\underline{k} >_{(m')}} \in \Gamma (\X,\smash{\widehat{\D}} ^{(m)} _{\X,\,\Q})$).
Comme $\Gamma (\U , \O _{\X})\cap \Gamma (\X , \O _{\X,\Q}) = \Gamma (\X , \O _{\X})$ (en effet, comme $\X$ est lisse et $U$ est dense dans $X$, 
le morphisme $\Gamma (X , \O _{X}) \to \Gamma (U , \O _{X})$ est injectif),
on vérifie alors que la famille $\{\underline{\partial} ^{<\underline{k} >_{(m')}} \cdot e \,|\, \underline{k} \in \N ^{d}\}$ est un sous-ensemble borné de $\Gamma (\X,\E )$ pour la topologie induite par sa structure de $\Gamma (\X , \O _{\X,\Q})$-module de type fini.
Soit $P = \sum _{\underline{k}\in \N ^{d}} a _{\underline{k}} \underline{\partial} ^{<\underline{k} >_{(m')}}$ un élément de $\Gamma (\X,\smash{\widehat{\D}} ^{(m')} _{\X,\,\Q})$. 
Il en résulte la convergence de la somme  
$\sum _{\underline{k}\in \N ^{d}} a _{\underline{k}} (\underline{\partial} ^{<\underline{k} >_{(m')}}.e)$.
On pose alors $P.e := \sum _{\underline{k}\in \N ^{d}} a _{\underline{k}} (\underline{\partial} ^{<\underline{k} >_{(m')}}.e)$. 

Il reste à vérifier que cela munit bien $\Gamma (\X,\E)$ d'une structure de 
$\Gamma (\X,\smash{\widehat{\D}} ^{(m')} _{\X,\,\Q})$-module. 
En multipliant par $a _{\underline{k}}$ les égalités $\underline{\partial} ^{<\underline{k} >_{(m')}}.e = \underline{\partial} ^{<\underline{k} >_{(m')}}.(1\otimes e)$, puis en les sommant, 
on obtient
$P.e = P.(1\otimes e)$.
Ainsi, $\Gamma (\X,\E)$ est un sous-$\Gamma (\X,\smash{\widehat{\D}} ^{(m')} _{\X,\,\Q})$-module de $\Gamma (\U,\E)$.
\end{proof}

\begin{theo}
\label{preborel9.3} 
Soit $\E$ un $\D ^{\dag } _{\X , \Q }$-module cohérent satisfaisant la propriété suivante:  
pour tout 
point fermé $x$ de $X$, le $K$-espace vectoriel $ i ^* _x ( \E ) $ est de dimension finie.

Il existe alors un ouvert affine dense $\Y$ de $\X$ tel que $\E |\Y$ soit un $\O _{\Y,\Q}$-module libre de rang fini. 
\end{theo}

\begin{proof}
I) Supposons dans un premier temps que le cardinal de $k$ n'est pas dénombrable. 

Le théorème étant local, supposons $\X$ affine, intègre et muni de coordonnées locales.
Soit  $m _0\geq 0$ un entier tel qu'il existe un $\smash{\widehat{\D}} ^{(m _0)} _{\X,\,\Q}$-module cohérent 
$\E ^{(m _0)}$ et un isomorphisme $\D ^{\dag } _{\X , \Q }$-linéaire de la forme 
$ \D ^{\dag } _{\X , \Q } \otimes _{\smash{\widehat{\D}} ^{(m _0)} _{\X,\,\Q}} \E ^{(m _0)} \riso \E$
(voir \cite[3.6.2]{Be1}).
Pour tout entier $m\geq m _0$, on pose $\E ^{(m)}:= \smash{\widehat{\D}} _{\X,\,\Q} ^{(m)} \otimes _{\smash{\widehat{\D}} _{\X,\,\Q} ^{(m _{0})}} \E ^{(m _0)} $,
  $E ^{(m)} := \Gamma (\X,\,\E ^{(m)})$ et $E := \Gamma (\X,\,\E)$.
  Soit $\smash{\overset{_{\circ}}{\E}}^{(m _{0})}$ un $\smash{\widehat{\D}} _{\X} ^{(m _{0} )}$-module cohérent sans $p$-torsion tel que l'on ait un isomorphisme
$\smash{\widehat{\D}} _{\X,\,\Q} ^{(m _{0} )}$-linéaire : $\smash{\overset{_{\circ}}{\E}}^{(m _{0} )}_\Q \riso \E ^{(m_0 )}$ (voir \cite[3.4.5]{Be1}).
On note $\smash{\overset{_{\circ}}{\E}}^{(m )}$ le quotient de 
$\smash{\widehat{\D}} _{\X} ^{(m )} 
\otimes _{\smash{\widehat{\D}} _{\X} ^{(m _{0} )}}
\smash{\overset{_{\circ}}{\E}}^{(m _{0})}$ 
par sa partie de $p$-torsion. D'après \cite[3.4.4]{Be1}, $\smash{\overset{_{\circ}}{\E}}^{(m )}$ est un $\smash{\widehat{\D}} _{\X} ^{(m )} $-module cohérent. On pose $\smash{\overset{_{\circ}}{E}} ^{(m)} : = \Gamma (\X,\,\smash{\overset{_{\circ}}{\E}}^{(m)})$. La preuve du théorème se décompose en six étapes: 

\begin{enumerate}
\item[1)] Pour tout entier $m\geq m _{0}$, il existe un ouvert affine dense $\Y _{(m)}$ de $\X$ 
tel que $T _{m}:= X \setminus Y _{(m)}$ soit le support d'un diviseur de $X$ et 
tel que $\smash{\overset{_{\circ}}{\E}}^{(m)}| \Y _{(m)}$
soit  isomorphe au complété $p$-adique d'un $\O _{\Y _{(m)}}$-module libre.
\end{enumerate}

Le faisceau $ \O _X \otimes _{\O_{\X}} \smash{\overset{_{\circ}}{\E}}^{(m)}$
est un $\D ^{(m)} _{X}$-module cohérent.
De manière analogue à \cite[VII.9.3]{borel},
il existe alors un ouvert affine et dense $\Y _{(m)} $ de $\X$ tel que $\O _X \otimes _{\O_{\X}} \smash{\overset{_{\circ}}{\E}}^{(m)} |Y _{(m)}$ soit un $\O _{Y _{(m)}}$-module libre.
Il résulte alors de \cite[2.2.16]{caro_courbe-nouveau} 
que
$\Gamma (\Y _{(m)},\, \smash{\overset{_{\circ}}{\E}}^{(m)})$ est isomorphe au complété $p$-adique d'un
$\Gamma (\Y_{(m)},\,\O _{\Y _{(m)}})$-module libre.
Or, par passage à la limite projective de \cite[2.3.5.2]{Be1}, on obtient l'isomorphisme
$\O _{\Y _{(m)}} \widehat{\otimes} _{\Gamma (\Y _{(m)},\,\O _{\X})} \Gamma (\Y _{(m)},\, \smash{\overset{_{\circ}}{\E}}^{(m)})
\riso
\smash{\widehat{\D}} ^{(m )} _{\Y_{(m)}}
\widehat{\otimes} _{\Gamma (\Y_{(m)},\, \smash{\widehat{\D}} ^{(m )} _{\X})}
\Gamma (\Y_{(m)},\, \smash{\overset{_{\circ}}{\E}}^{(m)})$.
Comme $\Gamma (\Y_{(m)},\, \smash{\overset{_{\circ}}{\E}}^{(m)})$ est un $\Gamma (\Y_{(m)},\, \smash{\widehat{\D}} ^{(m )} _{\X})$-module de type fini et comme 
$\smash{\overset{_{\circ}}{\E}}^{(m)} |\Y _{(m)}$ est un $\smash{\widehat{\D}} ^{(m )} _{\Y _{(m)}}$-module cohérent,  via le théorème de type $A$ de Berthelot \cite[3.3.9]{Be1}, 
on obtient alors les isomorphismes :
$$
\smash{\widehat{\D}} ^{(m )} _{\Y_{(m)}}
\widehat{\otimes} _{\Gamma (\Y_{(m)},\, \smash{\widehat{\D}} ^{(m )} _{\X})}
\Gamma (\Y_{(m)},\, \smash{\overset{_{\circ}}{\E}}^{(m)})
\tilde{\leftarrow}
\smash{\widehat{\D}} ^{(m )} _{\Y_{(m)}}
\otimes _{\Gamma (\Y_{(m)},\, \smash{\widehat{\D}} ^{(m )} _{\X})}
\Gamma (\Y_{(m)},\, \smash{\overset{_{\circ}}{\E}}^{(m)})
\riso \smash{\overset{_{\circ}}{\E}}^{(m)} |\Y _{(m)}.$$
Ainsi, $\smash{\overset{_{\circ}}{\E}}^{(m)} |\Y _{(m)}$
est isomorphe au complété $p$-adique d'un $\O _{\Y _{(m)}}$-module libre.

\begin{enumerate}
\item[2)]Comme $k$ n'est pas de cardinal dénombrable, 
il existe alors un point fermé $x$ de $X$ n'appartenant pas à $\cup _{m\in\N} T _{m}$.
Quitte à rétrécir $\Y _{(m)}$, on peut supposer $\Y _{(m)} \supset \Y _{(m+1)} $.
Notons $\mathcal{I} \subset \O_\X$ l'idéal définissant $i _x$
 et $I := \Gamma(\X,\, \mathcal{I})$.
\end{enumerate}

  \begin{enumerate}
\item[3)] 
Il existe $m _{1} \geq m _{0}$ tel que, pour tout $m \geq m _{1}$, le faisceau $\smash{\overset{_{\circ}}{\E}} ^{(m)} | \Y _{(m)}$ soit un $\O _{\Y _{(m)}}$-module libre de type fini. 
\end{enumerate}
  Par l'absurde, supposons que, pour tout entier $m\geq m _{0}$,
  $E ^{(m)}/I E ^{(m)}$ soient des
  ${K (x)}$-espaces vectoriels de dimension infinie.
D'après la proposition \cite[2.2.6.(i)]{caro_courbe-nouveau},
  $\smash{\overset{_{\circ}}{E}} ^{(m)} /I\smash{\overset{_{\circ}}{E}} ^{(m)}
  \riso  
  i ^{*} _{x} (\smash{\overset{_{\circ}}{\E}} ^{(m)}) $.
Notons $i _{x} ^{(m)}$ l'immersion fermée $\Spf \V (x) \hookrightarrow \Y _{(m)}$ induite par $i _{x}$. Comme $i ^{*} _{x} (\smash{\overset{_{\circ}}{\E}} ^{(m)}) \riso 
  i ^{(m)*} _{x} (\smash{\overset{_{\circ}}{\E}} ^{(m)} | \Y _{(m)})$, 
   $\smash{\overset{_{\circ}}{E}} ^{(m)} /I\smash{\overset{_{\circ}}{E}} ^{(m)}$ est donc, d'après l'étape $1$, isomorphe au complété
  $p$-adique d'un ${\V (x)}$-module libre (voir \cite[2.2.6.(ii)]{caro_courbe-nouveau}).
On peut alors reprendre la preuve de \cite[2.2.10]{caro_courbe-nouveau} à partir de {\og Par l'absurde, supposons que...\fg} pour aboutir à une contradiction.
On a ainsi vérifié qu'il existe un entier $m _{1} \geq m _{0}$ tel que $E ^{(m _{1})}/I E ^{(m _{1})}$ soit un ${K (x)}$-espace vectoriel de dimension finie. 
Avec \cite[2.2.8]{caro_courbe-nouveau}, il en résulte alors qu'il en ait de même
de $E ^{(m)}/I E ^{(m )}$ pour tout $m \geq m _{1}$. D'où le résultat.

 \begin{enumerate}
\item[4)] 
Soit $m $ un entier tel que $m\geq m _{1}$.
Le morphisme canonique de $\O _{ \Y _{(m+1)} ,\Q}$-modules libres de type fini : 
$\E ^{(m)} | \Y _{(m+1)} \to 
\E ^{(m+1)} | \Y _{(m+1)} $
est surjectif.
\end{enumerate}

En effet, comme l'assertion est locale, il suffit de l'établir pour un ouvert affine $\U$ de $\Y _{(m+1)}$. 
  D'après \cite[3.7.3.1]{bosch},
  en munissant $\Gamma (\U,\, \E ^{(m)}  )$ et $\Gamma (\U,\, \E ^{(m+1)}  )$ de la topologie induite par leur structure de
  $\Gamma (\U,\, \O _{\X,\Q}  )$-module de type fini,
  l'image de $\Gamma (\U,\, \E ^{(m)} ) $ dans $\Gamma (\U,\, \E ^{(m+1)}  )$ est fermée.
  Or, d'après \cite[2.2.9]{caro_courbe-nouveau},
  l'image
  de $\Gamma (\U,\, \E ^{(m)}) $ dans $\Gamma (\U,\, \E ^{(m+1)}  )$ est dense,
   lorsque $\Gamma (\U,\, \E ^{(m+1)}  )$ est muni de la topologie induite par sa structure
  de $\Gamma (\U,\,\smash{\widehat{\D}} _{\X,\Q} ^{(m +1)})$-module de type fini.
  Mais, grâce à \cite[4.1.2]{Be1}, ces deux topologies sur $\Gamma (\U,\, \E ^{(m+1)}  )$ coïncident.
  On en déduit que le morphisme canonique
  $\Gamma (\U,\, \E ^{(m)}) \rightarrow \Gamma (\U,\, \E ^{(m+1)} ) $ est surjectif.

\begin{enumerate}
\item[5)]

Si $E /I E$ est un ${K (x) }$-espace vectoriel de dimension finie égale à $r$, 
 alors il existe un entier $m _{2} \geq m _{1}$ tel que pour tout $m \geq m _{2}$, 
 le morphisme canonique $\E ^{(m)} | 
 \Y _{(m+1)}  \to  \E ^{(m+1)} | \Y _{(m+1)} $
 soit un isomorphisme de $\O _{\Y _{(m+1)},\Q}$-modules libres de rang $r$. 
\end{enumerate}

Il résulte alors de l'étape $4)$ que le morphisme canonique $E ^{(m)}/I E ^{(m )} \to E ^{(m+1)}/I E ^{(m +1)}$ est surjectif.
Comme $E /IE \riso \underset{\longrightarrow}{\lim}\, _{ _m} \, E ^{(m)} / I E ^{(m)}$ (e.g., voir la fin de la preuve de \cite[2.2.10]{caro_courbe-nouveau}),  il existe un entier $m _{2} \geq m _{1}$ tel que pour tout $m \geq m _{2}$, la dimension
du ${K (x) }$-espace vectoriel $E ^{(m)}/I E ^{(m )} $ est $r$.
Or, on remarque que le rang de $\E ^{(m)} | \Y _{(m)}$ est la dimension du ${K (x) }$-espace vectoriel
$E ^{(m)}/I E ^{(m )}$. 
Pour tout $m \geq m _{2}$,
le rang de $\E ^{(m)} | \Y _{(m)}  $ vaut donc $r$. 
Pour $m \geq m _{2}$, le morphisme canonique $\E ^{(m)} | \Y _{(m+1)} \to \E ^{(m+1)} | \Y _{(m+1)} $ est donc un morphisme de $\O _{\Y _{(m+1)},\Q}$-modules libres de rang $r$. 
D'après l'étape $4)$, ce morphisme est surjectif. Comme pour tout ouvert affine $\U$ de $\Y _{(m+1)}$ l'anneau $\Gamma (\U, \O _{\X,\Q})$ est intègre, ce morphisme devient injectif après extension par l'homomorphisme de $\Gamma (\U, \O _{\X,\Q})$ dans son corps des fractions.
On obtient donc son injectivité.

\begin{enumerate}
\item[6)] Le morphisme canonique $\E ^{(m _{2})} |\Y _{(m _{2})} \to \E  |\Y _{(m _{2})}$
est un isomorphisme. En particulier, $\E  |\Y _{(m _{2})}$ est un $\O _{\Y _{(m _{2})},\Q}$-module libre de rang $r$.
\end{enumerate}

Soit $m \geq m _{2}$ en entier. 
D'après l'étape $5$, le morphisme canonique $\E ^{(m _{2})} | 
 \Y _{(m)}  \to \E ^{(m)} | \Y _{(m)} $
 est un isomorphisme
de $\O _{\Y _{(m)},\Q}$-modules libres de rang $r$. 
Il en dérive que, pour pour tout ouvert affine non vide $\U$ de $\Y _{(m)}$, 
$\Gamma (\U  , \E ^{(m _{2})})$ est 
muni d'une structure de $\Gamma (\U,\smash{\widehat{\D}} ^{(m)} _{\X,\,\Q})$-module.
Avec \ref{m->m'}, il en résulte que
$\Gamma (\Y _{(m _{2})}  , \E ^{(m _{2})})$ est 
muni d'une (unique) structure de $\Gamma (\Y _{(m _{2})},\smash{\widehat{\D}} ^{(m)} _{\X,\,\Q})$-module prolongeant sa structure de $\Gamma (\Y _{(m _{2})},\smash{\widehat{\D}} ^{(m _{2})} _{\X,\,\Q})$.
Cela implique que le morphisme $\Gamma (\Y _{(m _{2})}  , \E ^{(m _{2})}) \to \Gamma (\Y _{(m _{2})}  , \E ^{(m)})$ admet une rétractation canonique. Celui-ci est donc injectif. 
Vérifions à présent sa surjectivité : le morphisme canonique 
$\Gamma (\Y _{(m _{2})}  , \E ^{(m _{2})}) \to \Gamma (\Y _{(m _{2})}  , \E ^{(m)})$, où, pour $m'=m$ et $m'= m _{2}$,
$\Gamma (\Y _{(m _{2})}  , \E ^{(m')})$ est muni de la topologie induite par sa structure de 
$\Gamma (\Y _{(m _{2})},\smash{\widehat{\D}} ^{(m')} _{\X,\,\Q})$-module de type fini, est continue.
Or, il résulte de \cite[4.1.2]{Be1} que la topologie de $\Gamma (\Y _{(m _{2})}  , \E ^{(m _{2})})$
induite par sa structure de $\Gamma (\Y _{(m _{2})},\smash{\widehat{\D}} ^{(m _{2})} _{\X,\,\Q})$-module de type fini et celle induite par sa structure de $\Gamma (\Y _{(m _{2})},\smash{\widehat{\D}} ^{(m)} _{\X,\,\Q})$-module de type fini sont identiques.
Comme $\Gamma (\Y _{(m _{2})},\smash{\widehat{\D}} ^{(m _{2})} _{\X,\,\Q})$ est dense dans
$\Gamma (\Y _{(m _{2})},\smash{\widehat{\D}} ^{(m)} _{\X,\,\Q})$, il en résulte que le morphisme canonique
$\Gamma (\Y _{(m _{2})}  , \E ^{(m _{2})}) \to \Gamma (\Y _{(m _{2})}  , \E ^{(m)})$
est 
$\Gamma (\Y _{(m _{2})},\smash{\widehat{\D}} ^{(m)} _{\X,\,\Q})$-linéaire. 
Son image est donc fermée. Comme elle est aussi dense, ce morphisme est surjectif. 
L'injectivité étant aisée (analogue à la fin de l'étape 5), le morphisme canonique 
$\Gamma (\Y _{(m _{2})}  , \E ^{(m _{2})}) \to \Gamma (\Y _{(m _{2})}  , \E ^{(m)})$
est donc un isomorphisme.
En passant à la limite inductive sur $m$, il en découle l'isomorphisme canonique 
$\Gamma (\Y _{(m _{2})}  , \E ^{(m _{2})}) \riso \Gamma (\Y _{(m _{2})}  , \E)$. 
D'après \cite[2.2.13]{caro_courbe-nouveau}, on en déduit que 
$\E | \Y _{(m _{2})} $ est $\O _{\Y _{(m _{2})},\Q}$-cohérent, i.e., est un isocristal convergent sur 
$\Y _{(m _{2})}$. Dans ce cas, on sait alors que le morphisme canonique
$\E ^{(m _{2})} |\Y _{(m _{2})} \to \E  |\Y _{(m _{2})}$ est un isomorphisme (voir la preuve de \ref{isocniveaum}
ou \cite[3.1.4.1 et 3.1.4.2]{Be0}).

\begin{enumerate}
\item[II)] On se ramène au cas où le cardinal de $k$ est quelconque grâce au lemme \ref{chgt-base-denombrable} 
ci-dessous. 
\end{enumerate}
 
\end{proof}

\begin{vide}
[Changement de base]
\label{nota-chgt-base}
Soient $\V \to \V'$ un morphisme d'anneaux de valuation discrète complets d'inégales caractéristiques $(0,p)$,
$\X ' := \X \times _{\Spf (\V)} \Spf (\V ')$, 
$f\,:\,\X' \to \X$ la projection canonique.
 
Par passage à la limite projective, on dispose d'après \cite[2.2.2]{Be1} et \cite[3.2.1]{Be2} des homomorphismes canoniques
d'anneaux
\begin{equation}
\label{iso-nivm}
f ^{-1} \widehat{\D} ^{(m)} _{\X/\Spf (\V)}
\to 
f ^{*} \widehat{\D} ^{(m)} _{\X/\Spf (\V)}
\riso 
\widehat{\D} ^{(m)} _{\X'/\Spf (\V')}
\end{equation}
dont le dernier est un isomorphisme. 
Il en résulte que le foncteur image inverse $f ^{*}$ se factorise en un foncteur de la catégorie des 
$\widehat{\D} ^{(m)} _{\X/\Spf (\V)}$-modules cohérents dans celle des 
$\widehat{\D} ^{(m)} _{\X'/\Spf (\V')}$-modules cohérents.
Par tensorisation par $\Q$, 
on obtient alors un foncteur de la catégorie des
$\widehat{\D} ^{(m)} _{\X/\Spf (\V),\Q}$-modules cohérents dans celle des 
$\widehat{\D} ^{(m)} _{\X'/\Spf (\V'),\Q}$-modules cohérents.
Pour préciser le niveau $m$, on le notera $f ^{(m) *}$. 
Ce foncteur $f ^{(m) *}$ est le foncteur changement de base (via $\V \to \V'$) de niveau $m$. 

Comme les homomorphismes \ref{iso-nivm} commutent aux changements de niveaux, on obtient 
par passage à la limite inductive sur le niveau les homomorphismes d'anneaux 
$f ^{-1} \D ^{\dag } _{\X/\Spf (\V),\Q}
\to 
f ^{*} \D ^{\dag } _{\X/\Spf (\V),\Q}
\riso 
\D ^{\dag } _{\X'/\Spf (\V'),\Q}$.
Le foncteur $f ^{*}$ se factorise donc aussi en un foncteur dit de changement de base (via $\V \to \V'$) de la catégorie 
des $\D ^{\dag } _{\X /\Spf (\V), \Q }$-modules cohérents dans celle des $\D ^{\dag } _{\X' /\Spf (\V'), \Q }$-modules cohérents.
De plus, via cette compatibilité au changement de niveaux, on vérifie que
si $\E ^{(m)}$ est un $\widehat{\D} ^{(m)} _{\X/\Spf (\V),\Q}$-module cohérent, on dispose alors de l'isomorphisme canonique 
\begin{equation}
\label{iso-nivm-dag}
\D ^{\dag } _{\X' /\Spf (\V'), \Q }\otimes _{\widehat{\D} ^{(m)} _{\X'/\Spf (\V'),\Q}} f ^{(m)*} (\E ^{(m)})
\riso
f ^{*} (\D ^{\dag } _{\X /\Spf (\V), \Q }\otimes _{\widehat{\D} ^{(m)} _{\X/\Spf (\V),\Q}} \E ^{(m)}).
\end{equation}

\end{vide}

\begin{lemm}
\label{chgt-base-denombrable}
Avec les notations de \ref{nota-chgt-base},
soient $\E$ un $\D ^{\dag } _{\X /\Spf (\V), \Q }$-module cohérent
et $\E':= f ^{*} (\E)$ le $\D ^{\dag } _{\X' /\Spf (\V'), \Q }$-module cohérent déduit par changement de base. 
Les deux assertions sont alors équivalentes: 
\begin{enumerate}
\item \label{chgt-base-denombrable-1} Il existe un ouvert dense $\Y$ de $\X$ tel que $\E |\Y$ soit un $\O _{\Y,\,\Q}$-module libre de type fini. 
\item \label{chgt-base-denombrable-2} Il existe un ouvert dense $\Y'$ de $\X'$ tel que $\E '|\Y'$ soit un $\O _{\Y',\,\Q}$-module libre de type fini. 
\end{enumerate}

\end{lemm}

\begin{proof}
L'implication  $\ref{chgt-base-denombrable-1} \Rightarrow \ref{chgt-base-denombrable-2}$ est immédiate.
Réciproquement, supposons qu'il existe un ouvert dense $\Y'$ de $\X'$ tel que $\E '|\Y'$ soit un $\O _{\Y',\,\Q}$-module libre de rang fini $r$.

D'après l'étape $1$ de la preuve de \ref{preborel9.3} et avec ses notations, 
pour tout $m \geq m _{0}$, 
il existe un ouvert affine dense $\Y _{(m )}$ de $\X$ 
tel que $\smash{\overset{_{\circ}}{\E}}^{(m)}| \Y _{(m)}$
soit  isomorphe au complété $p$-adique d'un $\O _{\Y _{(m)}}$-module libre.
Quitte à rétrécir $\Y _{(m)}$, on peut d'ailleurs supposer $\Y _{(m)} \supset \Y _{(m+1)} $.
Le rang (un entier s'il est fini, sinon on le définit égal à $+ \infty$) de $\smash{\overset{_{\circ}}{\E}}^{(m)}| \Y _{(m)}$ en tant que $\O _{\Y _{(m)}}$-module est alors le même que celui de $\E ^{(m)}| \Y _{(m)}$ en tant que $\O _{\Y _{(m)},\Q}$-module.
Notons $ g\,:\, f ^{-1} ( \Y _{(m)}) \to  \Y _{(m)}$ le morphisme induit par $f$, $g ^{(m) *}$ le foncteur de changement de base de niveaux $m$, $\Y ' _{(m)}:=\Y' \cap f ^{-1} ( \Y _{(m)})$ et $h ^{(m) *}:=|\Y ' _{(m)} \circ g ^{(m) *}$.

Comme le foncteur $h ^{(m) *}$ préserve la cohérence et avec \ref{iso-nivm-dag}, 
on obtient l'isomorphisme de $\D ^{\dag } _{\Y ' _{(m)} , \Q }$-modules cohérents: 
$ \D ^{\dag } _{\Y ' _{(m)}, \Q } \otimes _{\smash{\widehat{\D}} ^{(m )} _{\Y ' _{(m)},\Q}} (h ^{(m) *} (\E ^{(m)}| \Y _{(m)})) \riso \E '| \Y ' _{(m)} $.
Par \cite[3.1]{Be0}, comme $\E '| \Y ' _{(m)} $ est en outre $\O _{\Y ' _{(m)},\Q}$-cohérent, 
il en résulte 
$h ^{(m) *} (\E ^{(m)}| \Y _{(m)}) \riso \E '| \Y ' _{(m)} $.
Comme le rang de $h ^{(m) *} (\E ^{(m)}| \Y _{(m)})$ en tant que $\O _{\Y ' _{(m)},\Q}$-module
est le même que celui de $\E ^{(m)}| \Y _{(m)}$ en tant que $\O _{\Y _{(m)},\Q}$-module, il en résulte que $\smash{\overset{_{\circ}}{\E}}^{(m)}| \Y _{(m)}$ est un $\O _{\Y _{(m)}}$-module 
de rang $r$ (i.e., on a établit l'étape $3$ de \ref{preborel9.3} pour $m _{1}=m _{0}$).
Or, en reprenant la preuve de l'étape $4$ de \ref{preborel9.3}, on vérifie alors que
le morphisme canonique  
$\E ^{(m)} | \Y _{(m+1)} \to 
\E ^{(m+1)} | \Y _{(m+1)} $
est surjectif.
Comme c'est un morphisme de $\O _{\Y _{(m+1)} ,\Q}$-modules libres de rang $r$, celui-ci est un isomorphisme (i.e., on a établit l'étape $5$ de \ref{preborel9.3} pour $m _{2}=m _{0}$). En procédant comme pour l'étape $6$ de la preuve de \ref{preborel9.3}, on vérifie alors que $\E  | \Y _{(m _{0})} $ est un $\O _{ \Y _{(m _{0})},\Q}$-module libre de rang $r$.

\end{proof}

\begin{theo}
\label{borel9.3}  
Soit $\E$ un $\D ^{\dag } _{\X , \Q }$-module cohérent tel que, pour tout point fermé $x$ de
$X$, $ i ^* _x ( \E ) $ soit un $K$-espace vectoriel de dimension finie.
Il existe alors un diviseur $T$ de $X$ tel que $(\hdag T) (\E)$ soit un isocristal surconvergent le long de $T$.
\end{theo}

\begin{proof}
Cela résulte aussitôt de \ref{preborel9.3} et du théorème de Berthelot énoncé dans \cite[2.2.12]{caro_courbe-nouveau}. 

\end{proof}

Grâce au théorème \ref{borel9.3} qui étend la version précédente au cas sans structure de Frobenius, la section \cite[2]{caro_surholonome} se généralise. En particulier: 
\begin{theo}
[Critère d'holonomie]
\label{surcoh=>hol}
  Soient $T$ un diviseur de $X$ et $\E$ un $\D ^\dag _{\X} (\hdag T) _{\Q}$-module cohérent. Si les espaces de cohomologie de $\DD _{T}(\E)$ sont
  à fibres extraordinaires finies (e.g., si $\DD _{T} (\E)$ est $\D ^\dag _{\X} (\hdag T) _{\Q}$-surcohérent),
  alors $\E$ est {\og $\D ^\dag _{\X} (\hdag T) _{\Q}$-holonome\fg}, i.e., pour tout $l\not = 0$, $\mathcal{H} ^{l} \DD _{T}(\E)=0$.
\end{theo}

\begin{proof}
Grâce à \cite[4.3.12]{Be1}, on se ramène au cas où le diviseur $T$ est vide. 
Il suffit alors de reprendre les preuves de la section \cite[2]{caro_surholonome}, en particulier celle de \cite[2.5]{caro_surholonome}, en remplaçant \cite[2.2.17]{caro_courbe-nouveau} par sa version sans structure de Frobenius \ref{borel9.3}. 
\end{proof}

\begin{rema}
Avec ses notations, il résulte du critère d'holonomie \ref{surcoh=>hol}, qu'un $\D ^{\dag } _{\X , \Q }$-module surholonome est holonome.
En particulier, le foncteur $\DD$ induit une autoéquivalence de la catégorie des $\D ^{\dag } _{\X , \Q }$-modules surholonomes.
De même, on vérifie que le foncteur $\DD _{T} $
induit une autoéquivalence de la catégorie des $\D ^{\dag } _{\X , \Q }$-modules surholonomes qui sont en outre 
munis d'une structure de $\D ^{\dag } _{\X } (\hdag T) _{\Q }$-modules cohérents.
\end{rema}

Dans notre lancée, terminons par la proposition suivante qui donne une version sans Frobenius de \cite[2.2.14]{caro_courbe-nouveau}.
\begin{prop}
\label{isocniveaum}
Soit $\E$ un $\D ^{\dag } _{\X , \Q }$-module cohérent.
Soit  $m _0\geq 0$ un entier tel qu'il existe un $\smash{\widehat{\D}} ^{(m _0)} _{\X,\,\Q}$-module cohérent 
$\E ^{(m _0)}$ induisant $\E$ par extension (\cite[3.6.2]{Be1}). Pour tout entier $m\geq m _0$, notons $\E ^{(m)}:= \smash{\widehat{\D}} _{\X,\,\Q} ^{(m)} \otimes _{\smash{\widehat{\D}} _{\X,\,\Q} ^{(m _{0})}} \E ^{(m _0)} $.

Le faisceau $\E $ est $\O _{\X,\,\Q}$-cohérent si et seulement s'il existe une suite strictement croissante $(m _{n}) _{n\in \N}$ telle que 
  $\E ^{(m _{n})} $ soit $\O _{\X,\,\Q}$-cohérent.
\end{prop}

\begin{proof} L'assertion est locale et l'on peut supposer que $\X$ est affine.
  Si $\E $ est $\O _{\X,\,\Q}$-cohérent alors $\E$ est associé à un isocristal surconvergent
  (\cite[4.1.4]{Be1}). D'après \cite[3.1.4.1 et 3.1.4.2]{Be0}, on obtient alors que
  $\E ^{(m_0)}$ est isomorphe à $\E$.

  Réciproquement, supposons qu'il existe une suite strictement croissante $(m _{n}) _{n\in \N}$ telle que 
  $\E ^{(m _{n})} $ soit $\O _{\X,\,\Q}$-cohérent.
  Le morphisme canonique
  $\Gamma (\X,\, \E ^{(m _{n})}) \rightarrow \Gamma (\X,\, \E ^{(m _{n+1})} ) $ est alors surjectif (comme pour l'étape $4$ de la preuve de \ref{preborel9.3}).
  En passant à la limite sur $n$, il en résulte la surjectivité de
$\Gamma (\X,\, \E ^{(m _{n})} ) \rightarrow \Gamma (\X,\, \E ) $. 
D'après \cite[2.2.13]{caro_courbe-nouveau},
$\E$ est donc $\O _{\X,\,\Q}$-cohérent.
\end{proof}

\begin{vide}
Avec les notations du théorème \ref{isocniveaum}, le fait que $\E ^{(m _{0})}$ soit $\O _{\X,\,\Q}$-cohérent ne devrait pas a priori impliquer que $\E$ soit $\O _{\X,\,\Q}$-cohérent. 
\end{vide}

\section{Théorèmes de pleine fidélité pour les isocristaux surconvergents}
\subsection{Foncteur image inverse-restriction}

Le théorème de descente propre de Shiho 
 (voir \cite[7.3]{shiho-logRC-RCII}) implique que le foncteur canonique du théorème \ref{plfid-dag} ci-dessous est pleinement fidèle. 
Lorsque l'on dispose d'une structure de Frobenius, nous établirons que ce foncteur induit une équivalence de catégories (voir \ref{plfid-dag-Frob}).

 \begin{theo}
\label{plfid-dag}

Soient $a\,:\, X^{(0)} \to X$ un morphisme propre surjectif de $k$-variétés intègres, $Y$ un ouvert dense de $X$, $j\,:\, Y \hookrightarrow X$ l'immersion ouverte, $Y^{(0)}:= a ^{-1} (Y)$.
Le foncteur canonique 
\begin{equation}
\label{diag-plfid-dag}
(a ^{*}, j ^{*}) \,:\, \mathrm{Isoc} ^{\dag} (Y,X/K)
\to 
\mathrm{Isoc} ^{\dag} (Y^{(0)},X^{(0)}/K) \times _{\mathrm{Isoc} ^{\dag} (Y^{(0)},Y^{(0)}/K)} \mathrm{Isoc} ^{\dag} (Y,Y/K)
\end{equation}
est pleinement fidèle.
\end{theo}

\begin{proof}
Comme le foncteur canonique 
$\mathrm{Isoc} ^{\dag} (Y,X/K) \to \mathrm{Isoc} ^{\dag} (Y,Y/K)$ 
est fidèle, il en résulte qu'il en est de même du foncteur canonique de 
\ref{diag-plfid-dag}.

Soient $E  _{1}, E  _{2}$ deux objets de $\mathrm{Isoc} ^{\dag} (Y,X/K)$. Notons 
$E  ^{(0)}_{1}, E ^{(0)} _{2}$ (resp. $\widehat{E}  _{1}$, $\widehat{E}  _{2}$) les objets de $\mathrm{Isoc} ^{\dag} (Y^{(0)},X^{(0)}/K)$ (resp. $\mathrm{Isoc} ^{\dag} (Y,Y/K)$)
correspondants. 
Soient $\alpha^{(0)}\,:\,E ^{(0)} _{1} \to E ^{(0)} _{2}$ un morphisme de $\mathrm{Isoc} ^{\dag} (Y^{(0)},X^{(0)}/K)$, 
$\beta\,:\,\widehat{E}  _{1} \to \widehat{E}  _{2}$ un morphisme de $\mathrm{Isoc} ^{\dag} (Y,Y/K)$ tels que $\alpha ^{(0)}$ et $\beta$ induisent canoniquement le même morphisme de $\mathrm{Isoc} ^{\dag} (Y^{(0)},Y^{(0)}/K)$ de la forme $ b ^{*} (\widehat{E}  _{2} ) \to b ^{*} (\widehat{E}  _{1} )$.
Il s'agit de construire un morphisme $\alpha\,:\,E  _{1} \to E  _{2}$ de $\mathrm{Isoc} ^{\dag} (Y,X/K)$
induisant $\alpha ^{(0)}$ et $\beta$.

Notons $b\,:\, Y^{(0)} \to Y$ le morphisme induit par $a$ et, pour $i=1,2$, $a _{i}\,:\, X^{(0)} \times _{X} X ^{(0)} \to X^{(0)}$ 
et
$b _{i}\,:\, Y^{(0)} \times _{Y} Y ^{(0)} \to Y^{(0)}$
les projections canoniques respectives à gauche et à droite.
Considérons le diagramme de $\mathrm{Isoc} ^{\dag} (Y^{(0)}\times _{Y} Y^{(0)},X^{(0)}\times _{X} X^{(0)}/K)$ de gauche ci-dessous 
\begin{equation}
\label{diag-plfid-preuve-dag}
\xymatrix{
{a _{1} ^{*} a ^{*} (E  _{1})} 
\ar@{=}[r] 
\ar[d] ^-{\sim}
& 
{a _{1} ^{*} (E ^{(0)} _{1})} 
\ar[r] ^-{a ^{*} _{1} (\alpha^{(0)})}
\ar[d] ^-{\sim}
& 
{a _{1} ^{*} (E ^{(0)} _{2})} 
\ar[d] ^-{\sim}
&
{a _{1} ^{*} a ^{*} (E  _{2})} 
\ar@{=}[l] 
\ar[d] ^-{\sim}
\\ 
{a _{2} ^{*} a ^{*} (E  _{1})} 
\ar@{=}[r]
&
{a _{2} ^{*} (E ^{(0)} _{1})}  
\ar[r] ^-{a ^{*} _{2} (\alpha^{(0)})}
& 
{a _{2} ^{*} (E ^{(0)} _{2})} 
&
{a _{2} ^{*} a ^{*} (E  _{2}),} 
\ar@{=}[l] 
}
\xymatrix{
{b _{1} ^{*} b ^{*} (\widehat{E}  _{1})} 
\ar[r] ^-{b _{1} ^{*} b ^{*} (\beta)}
\ar[d] ^-{\sim}
&
{b _{1} ^{*} b ^{*} (\widehat{E}  _{2})} 
\ar[d] ^-{\sim}
\\ 
{b _{2} ^{*} b ^{*} (\widehat{E}  _{1})} 
\ar[r] ^-{b _{2} ^{*} b ^{*} (\beta)}
&
{b _{2} ^{*} b ^{*} (\widehat{E}  _{2})} 
}
\end{equation}
dont les isomorphismes verticaux résultent de l'égalité $a \circ a _{1}=a \circ a _{2}$.
En appliquant le foncteur canonique de restriction
$\mathrm{Isoc} ^{\dag} (Y^{(0)}\times _{Y} Y^{(0)},X^{(0)}\times _{X} X^{(0)}/K) \to \mathrm{Isoc} ^{\dag} (Y^{(0)}\times _{Y} Y^{(0)},Y^{(0)}\times _{Y} Y^{(0)}/K)$ 
au diagramme de gauche de \ref{diag-plfid-preuve-dag}, on obtient le diagramme à sa droite qui est commutatif par fonctorialité. 
Comme ce foncteur est fidèle, il en résulte alors la commutativité du diagramme de gauche de \ref{diag-plfid-preuve-dag}.
Grâce au théorème de descente propre de Shiho (voir \cite[7.3]{shiho-logRC-RCII}), il en résulte l'existence d'un morphisme $\alpha\,:\,E  _{1} \to E  _{2}$ de $\mathrm{Isoc} ^{\dag} (Y,X/K)$
induisant $\alpha ^{(0)}$ et $\beta$.
 \end{proof}

Avec des structures de Frobenius, on obtient grâce au théorème de pleine fidélité de Kedlaya 
le théorème ci-dessous qui étend \cite[Théorème 2]{Etesse-descente-etale} :
\begin{theo}
\label{plfid-dag-Frob}
Soient $a\,:\, X^{(0)} \to X$ 
un morphisme propre surjectif de $k$-variétés intègres, $Y$ un ouvert dense de $X$, $Y^{(0)}:= a ^{-1} (Y)$.
Le foncteur canonique 
\begin{equation}
\label{diag-plfid-dagbid}
(a ^{*}, j ^{*}) \,:\, 
F \text{-}\mathrm{Isoc} ^{\dag} (Y,X/K)
\to 
F \text{-}\mathrm{Isoc} ^{\dag} (Y^{(0)},X^{(0)}/K) 
\times _{F \text{-}\mathrm{Isoc} ^{\dag} (Y^{(0)},Y^{(0)}/K)} 
F \text{-}\mathrm{Isoc} ^{\dag} (Y,Y/K)
\end{equation}
est une équivalence de catégories.
\end{theo}

\begin{proof}
(I) D'après \ref{plfid-dag}, le foncteur est pleinement fidèle. 
Notons $X^{(1)}:=X^{(0)}\times _{X} X^{(0)}$, $Y^{(1)}:=Y^{(0)}\times _{Y} Y^{(0)}$. 
Pour $i=1,2$, notons $a _{i}\,:\, X^{(1)}\to X^{(0)}$ 
les projections canoniques respectives à gauche et à droite.

(II) Établissons à présent que ce foncteur est essentiellement surjectif. 

Soit $( E ^{(0)}, \widehat{E},\rho)$ un objet de $F\text{-}\mathrm{Isoc} ^{\dag} (Y^{(0)},X^{(0)}/K) \times _{F\text{-}\mathrm{Isoc} ^{\dag} (Y^{(0)},Y^{(0)}/K)} F\text{-}\mathrm{Isoc} ^{\dag} (Y,Y/K)$. 
D'après le théorème de désingularisation de de Jong, il existe $a '\,:\,X''\to X ^{(1)}$ un morphisme projectif, surjectif, génériquement fini et étale tel que $X''$ soit lisse. Notons $Y'':= X''\setminus a ^{\prime -1} (X ^{(1)} \cap T ^{(1)})$, $b' \,:\, Y '' \to Y ^{(1)}$ le morphisme induit par $a'$, $j ''\,:\, Y'' \subset X''$ 
et, pour $i =0,1,2$, $j ^{(i)}\,:\, Y^{(i)} \subset X^{(i)}$ les immersions ouvertes canoniques. 
On désigne par $j ^{\prime \prime *}\,:\, F\text{-}\mathrm{Isoc} ^{\dag} (Y'',X''/K) \to F\text{-}\mathrm{Isoc} ^{\dag} (Y'',Y''/K) $
et $j ^{(i) *}\,:\, F\text{-}\mathrm{Isoc} ^{\dag} (Y^{(i)},X^{(i)}/K) \to F\text{-}\mathrm{Isoc} ^{\dag} (Y^{(i)},Y^{(i)}/K) $
les foncteurs restrictions.

i) Construction de l'isomorphisme de recollement 
$\theta\,:\, a _{2} ^{*} ( E ^{(0)}) \riso a _{1} ^{*} ( E ^{(0)})$.  Nous allons pour cela utiliser la pleine fidélité (prouvée dans la partie (I) de la preuve) du foncteur canonique suivant noté
\begin{equation}
\label{pfFidX''X(1)}
\phi:=( a ^{\prime *},  j ^{(1)*})\,:\, F\text{-}\mathrm{Isoc} ^{\dag} (Y^{(1)},X^{(1)}/K)
\to 
F\text{-}\mathrm{Isoc} ^{\dag} (Y'',X''/K) \times _{F\text{-}\mathrm{Isoc} ^{\dag} (Y'',Y''/K)} F\text{-}\mathrm{Isoc} ^{\dag} (Y^{(1)},Y^{(1)}/K).
\end{equation}
Pour $i = 1,2$, notons $\rho _{i}$ l'isomorphisme canonique 
$j ^{\prime \prime *} \circ a ^{\prime *} [a _{i} ^{*} ( E ^{(0)})] \riso 
b ^{\prime * }\circ  j ^{(1)*} [ a _{i} ^{*} ( E ^{(0)})]$.
L'image de $a _{i} ^{*} ( E ^{(0)})$ par le foncteur \ref{pfFidX''X(1)} est donc
$\phi (a _{i} ^{*} ( E ^{(0)}))= (  a ^{\prime *} [a _{i} ^{*} ( E ^{(0)})],   j ^{(1)*} [ a _{i} ^{*} ( E ^{(0)})], \rho _{i})$.

On définit canoniquement les isomorphismes 
$\theta _{Y ^{(1)}}\,:\, j ^{(1)*} [ a _{2} ^{*} ( E ^{(0)})]\riso j ^{(1)*} [ a _{1} ^{*} ( E ^{(0)})]$ 
et 
$\widehat{\theta}\,:\, b _{2} ^{*} ( j ^{(0) *} ( E ^{(0)})) \riso b _{1} ^{*} ( j ^{(0) *} ( E ^{(0)}))$
à partir de $\rho$ comme ceux rendant commutatif le diagramme canonique suivant : 
\begin{equation}
\label{diag-defi-thetaU1}
\xymatrix{
{j ^{(1)*} [ a _{2} ^{*} ( E ^{(0)})]} 
\ar@{.>}[d] ^-{\theta _{Y ^{(1)}}}
\ar[r] ^-{\sim}
& 
{b _{2} ^{*} ( j ^{(0) *} ( E ^{(0)}))}  
\ar[r] ^-{\sim} _-{b _{2} ^{*}(\rho)}
\ar@{.>}[d] ^-{\widehat{\theta}}
& 
{b _{2} ^{*} ( b ^{*} ( \widehat{E}   ))}  
\ar[d] ^-{\sim}
\\ 
{j ^{(1)*} [ a _{1} ^{*} ( E ^{(0)})]} 
\ar[r] ^-{\sim}
& 
{b _{1} ^{*} ( j ^{(0) *} ( E ^{(0)}))}  
\ar[r] ^-{\sim} _-{b _{1} ^{*}(\rho)}
& 
{b _{1} ^{*} ( b ^{*} ( \widehat{E}   )).}  
}
\end{equation}
On définit canoniquement l'isomorphisme 
$\theta _{Y''}  \,:\, j ^{\prime \prime*} a ^{\prime *} [ a _{2} ^{*} ( E ^{(0)})]\riso j ^{\prime \prime*} a ^{\prime *}  [ a _{1} ^{*} ( E ^{(0)})]$ à partir de $\rho$ via le diagramme commutatif: 
\begin{equation}
\label{diag-defi-theta''}
\xymatrix{
{j ^{\prime \prime *} a ^{\prime *} [ a _{2} ^{*} ( E ^{(0)})]}
\ar@{.>}[d] ^-{\theta _{Y''} }
\ar[r] ^-{\sim} _-{\rho _{2}}
& 
{b ^{\prime *} j ^{(1)*} [ a _{2} ^{*} ( E ^{(0)})]} 
\ar[d] ^-{\sim} _-{b ^{\prime * }(\theta _{Y ^{(1)}})}
\\ 
{j ^{\prime \prime *} a ^{\prime *} [ a _{1} ^{*} ( E ^{(0)})]}
\ar[r] ^-{\sim} _-{\rho _{1}}
& 
{b ^{\prime * } j ^{(1)*} [ a _{1} ^{*} ( E ^{(0)})].}  
}
\end{equation}
Comme $X''$ est lisse, d'après le théorème \cite[4.2.1]{kedlaya-semistableII}, le foncteur $j ^{\prime \prime *}$ est pleinement fidèle.
Il existe donc un et un seul isomorphisme $\theta ''\,:\, a ^{\prime *} a _{2} ^{*} ( E ^{(0)} )\riso a ^{\prime *} a _{1} ^{*} ( E ^{(0)})$
tel que $j ^{\prime \prime *} (\theta '') =\theta _{Y''} $.
Il découle de la commutativité de \ref{diag-defi-theta''}, que $(\theta '', \theta _{Y ^{(1)}})$ induit un isomorphisme de la forme
$\phi (a _{2} ^{*} ( E ^{(0)})) \riso \phi (a _{1} ^{*} ( E ^{(0)}))$.
Par pleine fidélité de $\phi$, on obtient un isomorphisme 
$\theta \,:\, a _{2} ^{*} ( E ^{(0)}) \riso a _{1} ^{*} ( E ^{(0)})$
tel que $a ^{\prime *} (\theta ) = \theta ''$ et $j ^{(1)*} (\theta) = \theta _{Y ^{(1)}}$.

ii) Vérifions à présent que $\theta$ vérifie la condition de cocycle. 
Or, il suffit de le vérifier après application du foncteur $j ^{(2)*} $ car celui-ci est fidèle. 
 Comme $j ^{(1)*} (\theta) = \theta _{Y ^{(1)}}$, cela résulte aussitôt de la construction de 
$\theta _{Y ^{(1)}}$ (voir \ref{diag-defi-thetaU1}).

iii) Grâce au théorème de descente propre de Shiho (voir \cite[7.3]{shiho-logRC-RCII}), il en résulte l'existence d'un (et d'un seul à isomorphisme canonique près) objet $E \in F\text{-}\mathrm{Isoc} ^{\dag} (Y,X/K)$ 
 et d'un isomorphisme $ \rho _{(Y ^{(0)}, X ^{(0)})} \,:\, a ^{*} (E) \riso E ^{(0)}$ s'inscrivant dans le diagramme commutatif
\begin{equation}
\label{rhoY0X0}
\xymatrix @C1,8cm {
{a ^{*} _{2} a ^{*} (E)} 
\ar[r] ^-{\sim} _-{a ^{*} _{2}\rho _{(Y ^{(0)}, X ^{(0)})}}
\ar[d] ^-{\sim}
& 
{a ^{*} _{2}(E ^{(0)} )}
\ar[d] ^-{\sim} _-{\theta}
\\ 
{a ^{*} _{1}a ^{*} (E)} 
\ar[r] ^-{\sim} _-{a ^{*} _{1}\rho _{(Y ^{(0)}, X ^{(0)})}}
& 
{a ^{*} _{1}( E ^{(0)} ).} 
}
\end{equation}
On définit l'isomorphisme $\rho _{(Y ^{(0)}, Y ^{(0)})}\,:\, b ^{*} j ^{*} (E) \riso b ^{*}\widehat{E}$ via le diagramme commutatif : 
\begin{equation}
\label{rhoY0Y0}
\xymatrix @C1,8cm{
{j ^{(0)*} a ^{*} (E)} 
\ar[r] ^-{\sim} _-{j ^{(0)*}\rho _{(Y ^{(0)}, X ^{(0)})}}
\ar[d] ^-{\sim} _{\rho _{\mathrm{can}}}
& 
{j ^{(0)*}(E ^{(0)} )}
\ar[d] ^-{\sim} _-{\rho}
\\ 
{b ^{*} j ^{*} (E)} 
\ar@{.>}[r] ^-{\sim} _-{\rho _{(Y ^{(0)}, Y ^{(0)})}}
&
{b ^{*}\widehat{E},} 
}
\end{equation}
où $\rho _{\mathrm{can}}$ désigne l'isomorphisme canonique.
Considérons à présent le {\og cube $A$\fg} dont la face de devant est l'image par le foncteur $ b _{2} ^{*}$ de \ref{rhoY0Y0},
la face de derrière est l'image par le foncteur $ b _{1} ^{*}$ de \ref{rhoY0Y0}, la face de droite correspond
au carré de droite de \ref{diag-defi-thetaU1} 
et dont les isomorphismes de la face de gauche sont canoniques. 
La face de gauche du cube $A$ est canoniquement commutative. Celle de droite (resp. de devant, resp. de derrière) l'est par définition. 
Pour s'assurer de la commutativité du cube $A$, il suffit alors de valider celle de la face du haut. 
Pour cela considérons le {\og cube $B$\fg} dont la face du bas est la face du haut du cube $A$, 
dont la face du haut est l'image par $ j ^{(1)*}$ du carré commutatif \ref{rhoY0X0}, 
dont la face de devant (resp. de derrière) est le carré commutatif par fonctorialité en l'isomorphisme 
$ j ^{(1)*} \circ a _{2} ^{*} \riso b _{2} ^{*} \circ j ^{(0)*} $
(resp. $ j ^{(1)*} \circ a _{1} ^{*} \riso b _{1} ^{*} \circ j ^{(0)*} $).
On remarque que la face de droite du cube $B$ est le carré de gauche de \ref{diag-defi-thetaU1} (cela a un sens car $j ^{(1)*} (\theta) = \theta _{Y ^{(1)}}$). Il est donc commutatif. Celle de la face de gauche du cube $B$ se vérifie canoniquement. Le cube $B$ est donc commutatif. D'où celle du cube $A$.
Ainsi, la face du bas du cube $A$ est commutatif, i.e., l'isomorphisme $\rho _{(Y ^{(0)}, Y ^{(0)})}\,:\, b ^{*} j ^{*} (E) \riso b ^{*}\widehat{E}$ commute aux données de recollement. Il en résulte qu'il existe un isomorphisme $\rho _{(Y , Y)}\,:\,  j ^{*} (E) \riso\widehat{E}$ tel que
$\rho _{(Y , Y)}= b ^{*}(\rho _{(Y ^{(0)}, Y ^{(0)})})$.
La commutativité du diagramme \ref{rhoY0Y0} signifie alors que l'on dispose 
d'un isomorphisme $(\rho _{(Y ^{(0)}, X ^{(0)})}, \rho _{(Y , Y)})\,:\,( a ^{*} (E) ,  j ^{*} (E), \rho _{\mathrm{can}}) \riso 
( E ^{(0)}, \widehat{E},\rho)$. D'où le résultat.

\end{proof}

\subsection{Foncteur restriction}

\begin{theo}
\label{EqCat-rig-coro-Frob}
Soient $X$ une $k$-variété intègre, $Y$ un ouvert dense de $X$ et
$j\,:\,Y\hookrightarrow X$ l'immersion ouverte correspondante.
Le foncteur $j ^{*}\,:\,
F \text{-}\mathrm{Isoc} ^{\dag} (Y,X/K)
\to 
F \text{-}\mathrm{Isoc} ^{\dag} (Y,Y/K)$
est pleinement fidèle. 
\end{theo}

\begin{proof}
La fidélité étant connue, il reste à établir que celle-ci est pleine. 
D'après le théorème de désingularisation de de Jong (voir \cite{Ber-alterationdejong}), 
il existe un morphisme projectif, surjectif, génériquement fini et étale de $k$-variétés intègres de la forme $a\,:\, X^{(0)} \to X$ tel que $X ^{(0)}$ soit lisse.
Notons
$Y^{(0)}:= a ^{-1} (Y)$.
Soient $E _{1},E _{2} \in F \text{-}\mathrm{Isoc} ^{\dag} (Y,X/K)$ et $\widehat{\phi}\,:\, j ^{*} E _{1} \to j ^{*} E _{2}$.
Or, comme $X ^{(0)}$ est lisse, 
le foncteur canonique 
$F \text{-}\mathrm{Isoc} ^{\dag} (Y^{(0)},X^{(0)}/K)
\to 
F \text{-}\mathrm{Isoc} ^{\dag} (Y ^{(0)},Y^{(0)}/K)$
est pleinement fidèle (voir \cite[4.2.1]{kedlaya-semistableII}). 
Le morphisme $\widehat{\phi}$ induit donc canoniquement un morphisme 
$\phi ^{(0)}\,:\,a ^{*}(E _{1}) \to a ^{*} (E _{2})$.
Il résulte du théorème \ref{plfid-dag-Frob} qu'il existe un morphisme 
$\phi \,:\,E _{1} \to E _{2}$ induisant $\phi ^{(0)}$ et $\widehat{\phi}$.

\end{proof}

\begin{rema}
N. Tsuzuki conjecture (voir la conjecture \cite[1.2.1]{tsumono}) que le foncteur $j ^{*}$ de \ref{EqCat-rig-coro-Frob}
reste pleinement fidèle sans structure de Frobenius. 
\end{rema}

\subsection{Foncteur extension et image inverse-extension}

\begin{vide}
[Rappels sur le foncteur extension]
\label{jdag-pl-fid}
Soient $X$ une variété sur $k$, $Y$ un ouvert de $X$ avec $Y$ normal, 
$\widetilde{Y}$ un ouvert dense de $Y$, 
$\widetilde{j}\,:\,\widetilde{Y} \hookrightarrow X$ l'immersion ouverte correspondante.
 Tsuzuki a démontré récemment dans \cite{Tsuzuki-NormVar09} que 
 le foncteur canonique
$\widetilde{j} ^{\dag}\,:\,\mathrm{Isoc} ^{\dag} (Y,X/K) \to \mathrm{Isoc} ^{\dag} (\widetilde{Y},X/K) $ est pleinement fidèle. 
Ce théorème a d'abord été établi par Tsuzuki lorsque $X$ est lisse (voir \cite[4.1.1]{tsumono}), puis étendu par 
Kedlaya lorsque $Y$ est lisse dans \cite[5.2.1]{kedlaya-semistableI} (nous n'aurons besoin que du cas où la variété est lisse).
On note $\mathrm{Isoc} ^{\dag} (\widetilde{Y}\supset Y,X/K) $ l'image essentielle du foncteur 
$\widetilde{j} ^{\dag}\,:\,\mathrm{Isoc} ^{\dag} (Y,X/K) \to \mathrm{Isoc} ^{\dag} (\widetilde{Y},X/K) $.
\end{vide}

Le lemme suivant est évident compte tenu du paragraphe \ref{jdag-pl-fid} et de ses notations.

\begin{lemm}
\label{EqCat-rig}
Soient $a\,:\, X^{(0)} \to X$ 
un morphisme propre, surjectif, génériquement fini et étale de $k$-variétés intègres, $Y$ un ouvert de $X$, 
$\widetilde{Y}$ un ouvert dense de $Y$, 
$\widetilde{j}\,:\,\widetilde{Y} \hookrightarrow X$ l'immersion ouverte correspondante, 
$Y^{(0)}:= a ^{-1} (Y)$, 
$\widetilde{Y}^{(0)}:= a ^{-1} (\widetilde{Y})$.
On suppose de plus $Y$ est $Y ^{(0)}$ lisses.
Le foncteur canonique 
\begin{equation}
\label{diag-plfid-dag2}
(a ^{*}, \widetilde{j} ^{\dag})\,:\,
\mathrm{Isoc} ^{\dag} (Y,X/K)
\to 
\mathrm{Isoc} ^{\dag} (Y^{(0)},X^{(0)}/K) 
\times _{\mathrm{Isoc} ^{\dag} (\widetilde{Y}^{(0)},X^{(0)}/K)} \mathrm{Isoc} ^{\dag} (\widetilde{Y}\supset Y,X/K) 
\end{equation}
est une équivalence de catégories.

\end{lemm}

\section{Isocristaux partiellement surcohérents}

\subsection{Définitions et premières propriétés}
Soient $\PP$ un $\V$-schéma formel séparé et lisse, $T$ un diviseur de $P$, $X$ un sous-schéma fermé de $P$, $\U$ l'ouvert de $\PP$ complémentaire de $T$. On suppose que $Y:= X \setminus T$ est un $k$-schéma lisse. 
\begin{vide}
[Rappels]
\label{nota-6.2.1dev-pre}

D'après \cite{caro-construction},
comme $Y$ est lisse, on dispose du foncteur canonique pleinement fidèle $\sp _{Y\hookrightarrow \U, +}$
de la catégorie $(F\text{-})\mathrm{Isoc} (Y/K)$ des $(F\text{-})$isocristaux convergents sur $Y$ dans celle des $(F\text{-})\D ^\dag _{\U,\Q} $-modules surcohérents (et même surholonomes d'après \cite{caro_surholonome}) 
à support dans $Y$. Ce foncteur commute aux foncteurs duaux ainsi qu'aux {\og images inverses (extraordinaires)\fg} respectives (voir \cite{caro-construction}). 
 
Nous avons défini dans \cite[6.2.1]{caro_devissge_surcoh} la
catégorie $(F\text{-})\mathrm{Isoc} ^{\dag \dag} (\PP, T, X/K)$
des $(F\text{-})\D ^\dag _{\PP} (\hdag T) _\Q$-modules cohérents $\E$ à support dans $X$
 tels que $\E$ et $\DD _{T} (\E)$ soient $\D ^\dag _{\PP} (\hdag T) _\Q$-surcohérents
 et $\E |\U$ soit dans l'image essentielle de $\sp _{Y\hookrightarrow \U, +}$.

\end{vide}

Sans structure de Frobenius, la technique de descente utilisée lors de la preuve de \cite[6.3.1]{caro_devissge_surcoh} ne fonctionne plus 
car on ne peut plus utiliser le théorème de pleine fidélité de Kedlaya (voir \cite[4.2.1]{kedlaya-semistableII} ou \cite{kedlaya_full_faithfull}).
Il en résulte que la stabilité par image inverse des catégories de la forme $\mathrm{Isoc} ^{\dag \dag} (\PP, T, X/K)$ n'est plus aussi directe. 
Pour contourner cette difficulté technique, 
nous allons introduire les catégories de \ref{nota-6.2.1dev} (qui s'avéreront d'après \ref{Isoc*=dagdag} être égales à $\mathrm{Isoc} ^{\dag \dag} (\PP, T, X/K)$).
Leur stabilité par images inverses (extraordinaires) est quasiment tautologique (voir la proposition \ref{stabIsoc*inv}). 
Par contre, il faudra jusqu'à la fin de ce chapitre (plus précisément en \ref{Isoc*=dagdag}) se passer de la stabilité par foncteur dual (voir aussi le deuxième point de la remarque \ref{DDIsoc*=mod}).
\begin{nota}
\label{nota-6.2.1dev}

 \begin{itemize}

\item On note $\mathrm{Coh} (\PP, T,X/K)$ la catégorie des $\D ^{\dag} _{\PP}(\hdag T) _{\Q}$-modules cohérents à support dans $X$. 

\item On note $\mathrm{Isoc} ^{*} (\PP, T, X/K)$ la sous-catégorie pleine de $\mathrm{Coh}  (\PP, T, X/K)$
des $\D ^{\dag} _{\PP}(\hdag T) _{\Q}$-modules cohérents $\E$ tels que $\DD _{T} (\E)$ soit $\D ^{\dag} _{\PP}(\hdag T) _{\Q}$-surcohérent et $\E |\U$ soit dans l'image essentielle du foncteur $\sp _{Y\hookrightarrow \U, +}$.
Les objets de $\mathrm{Isoc} ^{*} (\PP, T, X/K)$ sont {\og les isocristaux partiellement surcohérents sur $(\PP, T, X/K)$\fg}.
Lorsque $\PP$ est propre, on omet le qualificatif {\og partiellement\fg}. 

\item On note $\mathrm{Isoc} ^{**} (\PP, T, X/K)$ la sous-catégorie pleine de $\mathrm{Coh}  (\PP, T, X/K)$
des $\D ^{\dag} _{\PP}(\hdag T) _{\Q}$-modules surcohérents $\E$ tels que 
$\E |\U$ soit dans l'image essentielle du foncteur $\sp _{Y\hookrightarrow \U, +}$.

\end{itemize}

Lorsque le diviseur $T$ est vide, on ne l'indique pas dans les notations \ref{nota-6.2.1dev} ci-dessus. 
\end{nota}

\begin{vide}
[Cas de la compactification partielle lisse]
\label{Xlisse-Eqsp+}
 Lorsque $X$ est lisse, d'après \cite[6.1.4]{caro_devissge_surcoh}, on dispose 
 des égalités 
 $\mathrm{Isoc} ^{*} (\PP, T, X/K) =\mathrm{Isoc} ^{**} (\PP, T, X/K) \mathrm{Isoc} ^{\dag\dag} (\PP, T, X/K)$. 
 De plus, avec le théorème \cite[2.5.10]{caro-construction} (voir aussi ci-après \ref{prop-donnederecol-dag}), 
 on bénéficie de l'équivalence de catégories 
$\sp _{X \hookrightarrow \PP, T,+}\,:\,\mathrm{Isoc} ^{\dag} (Y,X/K) 
\cong \mathrm{Isoc} ^{*} (\PP, T, X/K) $.
Soient $T \subset \widetilde{T}$ un second diviseur de $P$, 
$\widetilde{Y}:= X \setminus \widetilde{T}$, $\widetilde{j}\,:\,\widetilde{Y} \hookrightarrow X$ l'immersion ouverte correspondante.
D'après \cite[6.2.1]{caro_devissge_surcoh}, \cite[6.2.2]{caro_devissge_surcoh} et \cite{caro-construction}, 
on dispose des diagrammes commutatifs suivants:
\begin{equation}
\label{Xlisse-Eqsp+-diag}
\xymatrix{
{\mathrm{Isoc} ^{\dag} (Y,X/K)} 
\ar[d] ^-{\widetilde{j} ^{\dag}}
\ar[r] ^-{\sp _{X \hookrightarrow \PP, T,+}} _-{\cong}
& 
{\mathrm{Isoc} ^{\dag \dag} (\PP, T, X/K)  } 
\ar[d] ^-{(\hdag \widetilde{T})}
\\ 
{ \mathrm{Isoc} ^{\dag} (\widetilde{Y},X/K)} 
\ar[r] ^-{\sp _{X \hookrightarrow \PP, \widetilde{T},+}} _-{\cong}
& 
{ \mathrm{Isoc} ^{\dag \dag} (\PP, \widetilde{T}, X/K)}, 
}
\xymatrix{
{\mathrm{Isoc} ^{\dag} (Y,X/K)} 
\ar[d] ^-{\overset{\vee}{}}
\ar[r] ^-{\sp _{X \hookrightarrow \PP, T,+}} _-{\cong}
& 
{\mathrm{Isoc} ^{\dag \dag} (\PP, T, X/K)  } 
\ar[d] ^-{\DD _{T}}
\\ 
{ \mathrm{Isoc} ^{\dag} (Y,X/K)} 
\ar[r] ^-{\sp _{X \hookrightarrow \PP, T,+}} _-{\cong}
& 
{ \mathrm{Isoc} ^{\dag \dag} (\PP, T, X/K)}. 
}
\end{equation}
D'après \cite{caro-construction}, on dispose aussi d'un isomorphisme de commutation des foncteurs de la forme $\sp _{X \hookrightarrow \PP, T,+}$ 
aux images inverses (extraordinaires).
\end{vide}

\begin{rema}
\label{DDIsoc*=mod}
Avec les notations \ref{nota-6.2.1dev}, on dispose des deux propriétés immédiates ci-dessous. 
\begin{enumerate}
\item Soit $\E \in \mathrm{Isoc} ^{*} (\PP, T, X/K)$.
Comme le foncteur $\sp _{Y\hookrightarrow \U, +}$ commute aux foncteurs duaux respectifs (voir \cite{caro_comparaison}), 
 pour tout entier $j\not = 0$, $\mathcal{H} ^{j} (\DD  (\E |\U) )=0$.
Par \cite[4.3.12]{Be1}, il en résulte que, pour tout entier $j\not = 0$, $\mathcal{H} ^{j} (\DD _{T} (\E) )=0$.
Enfin, d'après \ref{Xlisse-Eqsp+}, $\mathrm{Isoc} ^{*} (\U, Y/K) = 
\mathrm{Isoc} ^{**} (\U, Y/K)  =\mathrm{Isoc} ^{\dag \dag} (\U, Y/K)$.

\item Comme le foncteur $\sp _{Y\hookrightarrow \U, +}$ commute aux foncteurs duaux respectifs (voir \cite{caro_comparaison}) et via le théorème de bidualité (voir \cite{virrion}), 
on vérifie que le foncteur $\DD _{T}$ induit une équivalence entre les catégories 
$\mathrm{Isoc} ^{*} (\PP, T, X/K)$ et $\mathrm{Isoc} ^{**} (\PP, T, X/K)$.
En fait, nous établirons que ces deux catégories sont égales (voir \ref{Isoc*=dagdag}) mais avant nous prendrons garde
de les distinguer.
\end{enumerate}
\end{rema}

\begin{rema}
Grâce à \cite{caro-Tsuzuki-2008}, un élément de
$F\text{-}\mathrm{Isoc} ^{\dag \dag} (\PP, T, X/K)$
est aussi un $F\text{-}\D ^\dag _{\PP,\Q}$-module surholonome.
Notons que la structure de Frobenius est indispensable puisque les isocristaux surconvergents ne sont pas en général $\D ^\dag _{\PP,\Q}$-cohérents ni a fortiori $\D ^\dag _{\PP,\Q}$-surholonomes (e.g., voir le contre-exemple de Berthelot donné à la fin de \cite{Be1}).
\end{rema}

\begin{defi}
\label{defi-triplet}
Un {\og triplet $(\PP, T,X)$ lisse en dehors du diviseur\fg} est la donnée d'un $\V$-schéma formel séparé et lisse $\PP$, 
d'un diviseur $T$ de $P$, d'un sous-schéma fermé $X$ de $P$ tels que $Y:= X \setminus T$ soit lisse. 
Un tel quadruplet $(\PP, T,X,Y)$ est {\og un quadruplet lisse en dehors du diviseur\fg}.
Le quadruplet canonique induit par $(\PP, T,X)$ est $(\PP, T,X,Y)$.

Un morphisme de triplets lisses en dehors du diviseur $\theta \,:\, (\PP', T',X')\to (\PP, T,X)$ ou 
un morphisme de quadruplets lisses en dehors du diviseur
$\theta \,:\, (\PP', T',X',Y')\to (\PP, T,X,Y)$
est la donnée d'un morphisme $f\,:\, \PP' \to \PP$ tel que $f$ induise un morphisme $X'\to X$ tel que $f (Y') \subset Y$.
Si $a\,:\,X' \to X$, $b \,:\, Y' \to Y$ sont les morphismes induits par $f$, le morphisme $\theta$ se note aussi $(f,a,b)$.
\end{defi}

\begin{prop}
\label{stabIsoc*inv}
Soit $\theta = (f,a,b)\,:\, (\PP', T',X',Y')\to (\PP, T,X,Y)$ un morphisme de quadruplets lisses en dehors du diviseur.
Soient $\E \in (F\text{-})\mathrm{Isoc} ^{*} (\PP, T, X/K)$,
$\FF \in (F\text{-})\mathrm{Isoc} ^{**} (\PP, T, X/K)$.
Pour tout entier $j \in \Z\setminus \{0\}$, les égalités suivantes sont alors vérifiées :
$\mathcal{H} ^j (\DD _{T'} (\hdag T') \R \underline{\Gamma} ^\dag _{X'} f ^! \DD _{T} (\E)   [-d _{X'/X}] )=0,
\
\mathcal{H} ^j ((\hdag T')  \R \underline{\Gamma} ^\dag _{X'} f ^! (\FF)   [-d _{X'/X}] )=0.$

Plus précisément, 
on dispose des factorisations 
\begin{align}
\label{stabIsoc*inv-i} 
\theta ^{!}\,:=\,(\hdag T') \circ \R \underline{\Gamma} ^\dag _{X'} \circ f ^! [-d _{X'/X}] \,:\,&(F\text{-})\mathrm{Isoc} ^{**} (\PP, T, X/K) \to (F\text{-})\mathrm{Isoc} ^{**} (\PP', T', X'/K)\\
\label{stabIsoc*inv-ii}
\theta ^{+}\,:=\,\DD _{T'}\circ \R \underline{\Gamma} ^\dag _{X'}\circ (\hdag T') \circ f ^! [-d _{X'/X}] \circ \DD _{T} \,:\,&
(F\text{-})\mathrm{Isoc} ^{*} (\PP, T, X/K)\to (F\text{-})\mathrm{Isoc} ^{*} (\PP', T', X'/K)).
\end{align}
De plus ces foncteurs sont transitives pour la composition : 
si $\theta ^{\prime}= (f',a',b')\,:\, (\PP'', T'',X'',Y'')\to (\PP', T',X',Y')$ est un morphisme de quadruplets lisses en dehors du diviseur, 
on dispose des isomorphismes canoniques 
$\theta ^{\prime !} \circ \theta ^{!} \riso (\theta \circ \theta ^{\prime}) ^{!}$
et
$\theta ^{\prime +} \circ \theta ^{+} \riso (\theta \circ \theta ^{\prime}) ^{+}$. 
\end{prop}

\begin{proof}
Pour établir les deux annulations, par \cite[4.3.12]{Be1}, il suffit de le voir en dehors des diviseurs, ce qui nous ramène au 
cas déjà étudié de \cite{caro-construction} où la compactification partielle est lisse.
Les factorisations s'en déduisent par stabilité de la surcohérence par image inverse extraordinaire et foncteur cohomologique local
à support propre. La transitivité pour la composition de morphismes résulte de celle des foncteurs images inverses extraordinaires
et des isomorphismes de commutation de l'image inverse extraordinaire aux foncteurs extensions et restrictions (voir \cite{caro_surcoherent}).
\end{proof}

\begin{nota}
\label{stabIsoc*inv-rema}
Soit $\theta = (f,a,b)\,:\, (\PP', T',X',Y')\to (\PP, T,X,Y)$ un morphisme de quadruplets lisses en dehors du diviseur.

$\bullet$ On dispose du foncteur image inverse 
\begin{equation}
\label{theta*} 
\theta ^{*}\,:\, \mathrm{Isoc} ^{\dag} (Y,X/K) \to \mathrm{Isoc} ^{\dag} (Y',X'/K).
\end{equation}
qui d'ailleurs ne dépend canoniquement que de $(a,b)$ (voir \cite{LeStum-livreRigCoh}). 
Lorsque $a ^{-1} (Y) = Y'$, ce foncteur ne dépend que de $a$ et se notera
$a ^{*}$.

$\bullet$ Avec les notations de la proposition \ref{stabIsoc*inv},
lorsque $T'= f ^{-1} (T)$, on remarque que les foncteurs $(\hdag T') $ dans les expressions \ref{stabIsoc*inv-i} et \ref{stabIsoc*inv-ii} sont inutiles. 
Dans ce cas, s'il n'y a aucune ambiguïté sur $b$ et $f$, on notera alors respectivement $a ^{!}$ et $a ^{+}$ à la place 
$\theta ^{!}$ et $\theta ^{+}$. 
\end{nota}

\begin{lemm}
\label{rema-isocDagDagDense}
Avec les notations de \ref{nota-6.2.1dev}
notons {\og$?$\fg} l'un des symboles {\og$*$\fg}, {\og$**$\fg} ou {\og$\dag \dag$\fg}. 
Soit $\overline{Y}$ l'adhérence de $Y$ dans $X$.
Soient $Y_1,\dots ,Y _{N}$ les composantes irréductibles de $Y$ et $\overline{Y}_1,\dots, \overline{Y} _{N}$ leur adhérence dans $X$. 

\begin{enumerate}
\item On dispose de l'égalité
$(F\text{-})\mathrm{Isoc} ^{?} (\PP, T, \overline{Y}/K)=(F\text{-})\mathrm{Isoc} ^{?} (\PP, T, X/K)$.

\item On dispose d'une équivalence canonique de catégories : 
$(F\text{-})\mathrm{Isoc} ^{?} (\PP, T, X/K)\cong \oplus _{r=1} ^{N} (F\text{-})\mathrm{Isoc} ^{?} (\PP, T, \overline{Y} _{r}/K)$.
\end{enumerate} 
\end{lemm}

\begin{proof}
Traitons d'abord la première égalité. 
Commençons par le cas $?=**$. 
L'inclusion de la première égalité est évidente. Réciproquement, soit $\E \in (F\text{-})\mathrm{Isoc} ^{**} (\PP, T, X/K)$.
Par stabilité de la surcohérence par foncteur cohomologique local (voir \cite{caro_surcoherent}), on dispose alors du morphisme canonique $\R \underline{\Gamma} ^{\dag} _{\overline{Y}}(\E) \to \E$ de 
complexes à cohomologie bornée et $\D ^\dag _{\PP} (\hdag T) _\Q$-surcohérente. 
Comme c'est le cas en dehors de $T$ (car $\overline{Y}\setminus T =Y$), avec \cite[4.3.12]{Be1},
ce morphisme est alors un isomorphisme.
Le faisceau $\E$ est donc à support dans $\overline{Y}$. D'où l'inclusion inverse.
Avec la deuxième remarque de \ref{DDIsoc*=mod}, on en déduit par dualité la première égalité lorsque $?=*$.
Il résulte des deux premiers cas celui où $?=\dag \dag$.

En utilisant de même le théorème \cite[4.3.12]{Be1}, on vérifie que le foncteur $\oplus _{r=1} ^{N} \R \underline{\Gamma} ^{\dag} _{\overline{Y} _{r}}$ induit 
l'équivalence de catégories :
$(F\text{-})\mathrm{Isoc} ^{?} (\PP, T, X/K)\cong \oplus _{r=1} ^{N} (F\text{-})\mathrm{Isoc} ^{?} (\PP, T, \overline{Y} _{r}/K)$
pour $?=**$ ou pour $?=\dag\dag$.
Le dernier cas s'en déduit par dualité (plus précisément, l'équivalence de catégories est donnée par 
$\DD _{T} \circ (\oplus _{r=1} ^{N} \R \underline{\Gamma} ^{\dag} _{\overline{Y} _{r}})\circ \DD _{T}$).
\end{proof}

\begin{rema}
\label{rema-isocDagDense}
Le version rigide du lemme \ref{rema-isocDagDagDense} est déjà connue (voir \cite{LeStum-livreRigCoh} ou \cite{Berig}). 
Plus précisément, avec les notations de \ref{rema-isocDagDagDense}, 
on dispose de l'égalité
 $(F\text{-})\mathrm{Isoc} ^{\dag} (Y, \overline{Y}/K)=(F\text{-})\mathrm{Isoc} ^{\dag} (Y, X/K)$
 ainsi que de l'équivalence canonique de catégories : 
$(F\text{-})\mathrm{Isoc} ^{\dag} (Y, X/K)\cong \oplus _{r=1} ^{N} (F\text{-})\mathrm{Isoc} ^{\dag} (Y _{r}, \overline{Y} _{r}/K)$
\end{rema}

\begin{lemm}
\label{*dagdag-coh-PXTindtP}
  Soient $f$ : $\PP' \rightarrow \PP$ un morphisme propre de $\V$-schémas formels séparés et lisses,
  $X$ un sous-schéma fermé de $P'$ tel que le morphisme induit $X \rightarrow P$ soit une immersion fermée,
  $Y$ un ouvert  de $X$, 
  $T$ un diviseur de $P$ tel que $T' := f ^{-1} (T)$  soit un diviseur de $P'$ et tel que $Y = X \setminus T$
  (et donc $Y = X \setminus T'$).
  On suppose $Y$ lisse. 

\begin{enumerate} 
\item 
\label{*dagdag-coh-PXTindtP-i}
Soient $\E \in (F\text{-})\mathrm{Isoc} ^{*} (\PP, T, X/K)$,
$\FF \in (F\text{-})\mathrm{Isoc} ^{**} (\PP, T, X/K)$.
Soit $\E '\in (F\text{-})\mathrm{Isoc} ^{?} (\PP', T', X'/K)$), avec $? =*$ ou $? = **$.
Pour tout entier $j \in \Z\setminus \{0\}$, les égalités suivantes sont alors vérifiées :
$$\mathcal{H} ^j (\DD _{T'} \R \underline{\Gamma} ^\dag _X f ^! \DD _{T} (\E) ) =0,
\hspace{2cm}
\mathcal{H} ^j (\R \underline{\Gamma} ^\dag _X f ^! (\FF) ) =0,
\hspace{2cm}
\mathcal{H} ^j (f_+(\E')) =0.$$

\item \label{*dagdag-coh-PXTindtP-ii} Les foncteurs $\R \underline{\Gamma} ^\dag _X  f ^! $ (resp. $\DD _{T'} \R \underline{\Gamma} ^\dag _X f ^! \DD _{T} $) et
$f _+$ induisent des équivalences quasi-inverses 
entre les catégories $(F\text{-})\mathrm{Isoc} ^{**} (\PP, T, X/K)$ et $(F\text{-})\mathrm{Isoc} ^{**} (\PP', T', X/K)$
(resp. $(F\text{-})\mathrm{Isoc} ^{*} (\PP, T, X/K)$ et $(F\text{-})\mathrm{Isoc} ^{*} (\PP', T', X/K)$).

\end{enumerate}
\end{lemm}

\begin{proof}
Les autres annulations de (\ref{*dagdag-coh-PXTindtP-i}) étant déjà connues (voir \ref{stabIsoc*inv}), traitons 
la dernière. Par \cite[4.3.12]{Be1}, il suffit de le vérifier en dehors de $T$, ce qui nous ramène au cas $X$ est lisse et $T$ est vide. 
Comme cela est local en $\PP$, supposons $P$ affine. Dans ce cas, on dispose de relèvements $u '\,:\,\X \to \PP'$ 
et $u\,:\, \X \to \PP$ de $X \to P$ et $X' \to P'$. Par hypothèse, il existe $\G$ un $\D ^{\dag} _{\X,\Q}$-module cohérent, $\O _{\X,\Q}$-cohérent
tel que $\E' \riso u ^{\prime} _{+} (\G)$. D'où $f _{+} (\E') \riso u  _{+} (\G)$. Comme $u$ est une immersion fermée, 
le foncteur $ u _{+}$ est exact. 
Cela entraîne donc que, pour tout $j \in \Z\setminus \{0\}$,
$\mathcal{H} ^j (f_+(\E')) =0.$

Traitons à présent (\ref{*dagdag-coh-PXTindtP-ii}). Supposons $?= **$. Comme $f$ est propre, on dispose des morphismes d'adjonction
$f _{+} \circ \R \underline{\Gamma} ^\dag _X  f ^! (\FF) \to \FF$ et 
$\E' \to \R \underline{\Gamma} ^\dag _X  f ^! (\E') \circ f _{+}$.
Ces morphismes sont des isomorphismes en dehors respectivement du diviseur $T$ et $T'$ (en effet, comme cela est local, 
on se ramène alors au cas où l'on dispose de relèvements de $X \hookrightarrow P$), 
par \cite[4.3.12]{Be1}, ce sont donc des isomorphismes. D'où le résultat. 
\end{proof}

\begin{rema}
Le lemme \ref{*dagdag-coh-PXTindtP} sera généralisé (le morphisme $f$ n'y étant plus forcément propre, le diviseur $T'$ étant plus indépendant de $T$) en
\ref{gen-coh-PXTindtP}.
Néanmoins, avant d'établir ce cas plus général, cette première étape sera utilisée 
pour construire les morphismes d'adjonction au cours de la preuve du lemme \ref{theo-adj-fini-étale}
(plus précisément en fait, au cours de la preuve du lemme \ref{lemm-iso-adj-cube} qui permet d'établir \ref{theo-adj-fini-étale}).
On en déduira 
que la catégorie $(F\text{-})\mathrm{Isoc} ^{*} (\PP, T, X/K)$ 
ne dépend que du couple $(Y,X)$.
\end{rema}

\begin{lemm}
\label{Isoc*stblKer}
La catégorie $\mathrm{Isoc} ^{*} (\PP, T, X/K)$ est stable par noyau, image, conoyau. 
\end{lemm}

\begin{proof}
Soit $\phi$ un morphisme de $\mathrm{Isoc} ^{*} (\PP, T, X/K)$.
Comme la catégorie $\mathrm{Isoc} ^{\dag} (Y, Y/K)$ est stable par noyau, image et conoyau, via l'équivalence de catégories 
$\sp _{Y\hookrightarrow \U,+}\,:\, \mathrm{Isoc} ^{\dag} (Y, Y/K) \cong \mathrm{Isoc} ^{\dag \dag} (\U, Y/K)$ de \cite[6.2.2]{caro_devissge_surcoh}, il en résulte qu'il en est de même pour $\mathrm{Isoc} ^{\dag \dag} (\U, Y/K)$. 
Si $\E$ est le noyau ou l'image ou le conoyau de $\phi$, il en résulte que, pour tout entier $j\not = 0$, le faisceau $\mathcal{H} ^{j}(\DD  (\E| \U)) =0$ (voir \ref{DDIsoc*=mod}).
Par \cite[4.3.12]{Be1}, il en résulte que, pour tout entier $j\not = 0$, le faisceau $\mathcal{H} ^{j} \DD _{T} (\E) =0$. On en déduit que le noyau (resp. l'image, resp. le conoyau) de $\phi$ est le conoyau (resp. l'image, resp. le noyau) de $\DD _{T} (\phi)$. Or, comme $\DD _{T} (\phi)$ est un morphisme de $\D ^{\dag} _{\PP}(\hdag T) _{\Q}$-modules surcohérents, son noyau, son image, son conoyau sont $\D ^{\dag} _{\PP}(\hdag T) _{\Q}$-surcohérents.
\end{proof}

\begin{nota}
[Morphisme fini et étale en dehors des singularités]
\label{nota-diag-pre}
Soit le diagramme commutatif
\begin{equation}
\label{nota-diag-deux}
\xymatrix{
{Y^{(0)}} \ar@{^{(}->}[r]  ^-{j ^{(0)}} \ar[d] ^-{b}
& {X^{(0)} } \ar@{^{(}->}[r] ^-{u^{(0)}} \ar[d] ^-{a}
& {\PP ^{(0)} }  \ar[d] ^-{f}
\\ 
{Y} \ar@{^{(}->}[r] ^-{j }& {X } \ar@{^{(}->}[r] ^-{u} & {\PP  ,}
}
\end{equation}
où le carré de gauche est cartésien, 
$f$ est un morphisme propre, lisse de $\V$-schémas formels séparés et lisses,
$a$ est un morphisme propre surjectif de $k$-variétés, 
$b$ est un morphisme fini et étale de $k$-variétés lisses, 
$j$ et $j^{(0)}$ sont des immersions ouvertes, $u$ et $u^{(0)}$ sont des immersions fermées.
On suppose en outre
qu'il existe un diviseur $T$ de $P$ tel que 
$Y = X \setminus T$. On note $\U := \PP \setminus T$, $T ^{(0)}:=f ^{-1}(T)$, $\U ^{(0)}:= \PP ^{(0)}\setminus T^{(0)}$ et $g \,:\, \U ^{(0)} \to \U$ le morphisme induit par $f$.
\end{nota}

\begin{vide}
[Foncteur $b _{+}$ adjoint à droite et à gauche de $b ^{+}$]
\label{b+b+}
On garde les notations de \ref{nota-diag-pre}.

$\bullet$ D'après \cite[6.2.2]{caro_devissge_surcoh}, on dispose de l'équivalence de catégories 
$\sp _{Y\hookrightarrow \U,+}\,:\, \mathrm{Isoc} ^{\dag} (Y, Y/K) \cong \mathrm{Isoc} ^{\dag \dag} (\U, Y/K)$. 
D'après \cite{caro-construction}, les foncteurs de la forme $\sp _{Y\hookrightarrow \U,+}$ commutent aux foncteurs duaux respectifs et aux images inverses (extraordinaires). 
Les deux foncteurs $b ^{!}:=\R \underline{\Gamma}  ^{\dag } _{Y ^{(0)} } \circ g ^{!}$
et
 $b ^{+}:=\DD \circ  \R \underline{\Gamma}  ^{\dag } _{Y ^{(0)} } \circ g ^{!} \circ \DD$
de
$\mathrm{Isoc} ^{\dag \dag} (\U, Y/K) \to \mathrm{Isoc} ^{\dag \dag} (\U  ^{(0)}, Y  ^{(0)}/K)$
sont donc isomorphes.
On dispose ainsi du diagramme commutatif (à isomorphisme canonique près) : 
\begin{equation}
\label{fini-étale-présdag-diag}
\xymatrix@C=1,8cm {
{\mathrm{Isoc} ^{\dag} (Y  ^{(0)}, Y  ^{(0)}/K)}
\ar[r] ^-{\sp _{Y  ^{(0)}\hookrightarrow \U  ^{(0)},+}} _-{\cong}
& 
{\mathrm{Isoc} ^{\dag \dag} (\U  ^{(0)}, Y/K) } 
\\ 
{\mathrm{Isoc} ^{\dag} (Y, Y/K)} 
\ar[r] ^-{\sp _{Y\hookrightarrow \U,+}} _-{\cong}
\ar[u] ^-{b ^{*}}
& 
{\mathrm{Isoc} ^{\dag \dag} (\U, Y/K) .} 
\ar[u] ^-{b ^{+}}
}
\end{equation}

$\bullet$ En outre, comme $b$ est fini et étale, le foncteur $g _{+}$ se factorise de la manière suivante
$g _{+}\,:\, \mathrm{Isoc} ^{\dag \dag} (\U  ^{(0)}, Y  ^{(0)}/K) \to \mathrm{Isoc} ^{\dag \dag} (\U, Y/K)$.
En effet, comme cela est local en $\U$, on peut supposer $\U$ affine. Via le théorème de Berthelot-Kashiwara, on se ramène au cas où $Y=U$ et $Y^{(0)}=U^{(0)}$. Dans ce cas, comme $g$ est fini et étale, le foncteur $g _{+}$ préserve la $\D ^{\dag}$-cohérence ainsi que la $\O$-cohérence. 
On notera $b _{+}$ cette factorisation. 

$\bullet$ Comme $g ^{!}$ est adjoint à droite de $g _{+}$, on vérifie que le foncteur $b ^{!}$ est adjoint à droite de $b _{+}$.
Comme $g$ est propre, grâce au théorème de dualité relative, on en déduit que le foncteur $b ^{+}$ est adjoint à gauche de $b _{+}$.

$\bullet$ Comme $b ^{!}\riso b ^{+}$, il en résulte que $b ^{!}$ est adjoint à gauche de $b _{+}$ et 
que $b ^{+}$ est adjoint à droite de $b _{+}$.
\end{vide}

\begin{rema}
Avec les notations de \ref{b+b+}, comme $b$ est fini et étale, on dispose du foncteur image directe
$b _{*}\,:\,\mathrm{Isoc} ^{\dag} (Y  ^{(0)}, Y  ^{(0)}/K) \to \mathrm{Isoc} ^{\dag} (Y, Y/K)$ (voir \cite[5]{tsumono}). 
Le foncteur $b ^{*}$ est adjoint à gauche du foncteur $b _{*}$.
Par unicité des foncteurs d'adjonction, on obtient 
le diagramme commutatif (à isomorphisme canonique près) :
\begin{equation}
\label{fini-étale-présdag-diag2}
\xymatrix @C=1,8 cm{
{\mathrm{Isoc} ^{\dag} (Y  ^{(0)}, Y  ^{(0)}/K)}
\ar[d] ^-{b _{*}}
\ar[r] ^-{\sp _{Y  ^{(0)}\hookrightarrow \U  ^{(0)},+}} _-{\cong}
& 
{\mathrm{Isoc} ^{\dag \dag} (\U  ^{(0)}, Y/K) }
\ar[d] ^-{b _{+}} 
\\ 
{\mathrm{Isoc} ^{\dag} (Y, Y/K)} 
\ar[r] ^-{\sp _{Y\hookrightarrow \U,+}} _-{\cong}
& 
{\mathrm{Isoc} ^{\dag \dag} (\U, Y/K) ,} 
}
\end{equation}
de telle manière que les morphismes d'adjonction soient compatibles.
Mais nous n'aurons pas besoin de ce fait. 
\end{rema}

Afin de définir \ref{defi-Isocbullet*}, nous aurons besoin de la notion d'image directe par un morphisme fini et étale en dehors des singularités 
surconvergentes définie dans la proposition \ref{Isoc*a+} ci-dessous.
\begin{lemm}
\label{Isoc*a+}
On reprend les notations de \ref{nota-diag-pre}.
\begin{enumerate}
\item \label{Isoc*a+1} On dispose des factorisations
$f _{+}\,:\, \mathrm{Isoc} ^{*} (\PP^{(0)}, T^{(0)}, X^{(0)}/K)
\to \mathrm{Isoc} ^{*} (\PP, T, X/K)$ 
et
$a ^{+}:=\DD _{T^{(0)}}\circ  \R \underline{\Gamma}  ^{\dag } _{X^{(0)}} \circ f ^{!} \circ \DD _{T}\,:\, \mathrm{Isoc} ^{*} (\PP, T, X/K)
\to \mathrm{Isoc} ^{*} (\PP^{(0)}, T^{(0)}, X^{(0)}/K)$. 
En posant $\theta = (f,a,b)$, on note $\theta _{+}:= f _{+}$.
S'il n'y a pas d'ambiguïté sur le prolongement $f$ de $a$, on pose $ a _{+}:= f _{+}$.

\item \label{Isoc*a+5} Le foncteur $ a _{+} $ est adjoint à droite de $a  ^{+}$.
On notera $\mathrm{adj} _{a}$ les morphismes  d'adjonctions $\mathrm{Id} \to a _{+} \circ a ^{+}$ et $a ^{+} \circ a _{+}\to \mathrm{Id} $.
Ces morphismes d'adjonctions sont transitifs pour la composition de diagramme de la forme \ref{nota-diag-deux} et
satisfaisant les conditions requises.
\end{enumerate}
\end{lemm}

\begin{proof}
Traitons d'abord \ref{Isoc*a+1}. La seconde factorisation est déjà connue (voir la factorisation \ref{stabIsoc*inv-ii} et les notations de \ref{stabIsoc*inv-rema}), 
vérifions à présent la première.
Soit $\E ^{(0)}\in \mathrm{Isoc} ^{*} (\PP^{(0)}, T^{(0)}, X^{(0)}/K)$.
Comme $f$ est propre, $f _{+}$ commute au dual et la (sur)cohérence est préservée par $f _{+}$.
On en déduit que
$f _{+} (\E ^{(0)}) $ est cohérent et $\DD _{T} \circ f _{+} (\E ^{(0)}) $ est surcohérent. 
Or, $f _{+} (\E ^{(0)}) |\U \riso g _{+} (\E ^{(0)} |\U ^{(0)} )$.
D'après le deuxième point de \ref{b+b+}, il en résulte que 
$f _{+} (\E ^{(0)}) |\U \in \mathrm{Isoc} ^{\dag \dag} (\U, Y/K)$, i.e., 
$f _{+} (\E ^{(0)}) |\U$ est dans l'image essentielle de $\sp _{Y\hookrightarrow \U,+}$.
Par \cite[4.3.12]{Be1}, 
il en résulte aussi que $f _{+} (\E ^{(0)})$ est un module.
On obtient ainsi la factorisation de $f _{+}$ demandée.

On obtient canoniquement les morphismes d'adjonction entre $a _{+}$ et $a ^{+}$ 
de la manière suivante : 
comme le foncteur $f _{+}$ est adjoint à gauche de $f ^{!}$ (pour les complexes surcohérents) on obtient le morphisme 
$f _{+} \circ f ^{!} \circ \DD _{T} \to \DD _{T} $. Via le morphisme canonique $ \R \underline{\Gamma}  ^{\dag } _{X^{(0)}} \to Id$, 
on obtient par fonctorialité $ f _{+} \R \underline{\Gamma}  ^{\dag } _{X^{(0)}} \circ f ^{!} \circ \DD _{T} \to \DD _{T}$.
Avec le théorème de bidualité et de dualité relative, il en résulte par application du foncteur $\DD _{T}$ le morphisme
$Id \riso \DD _{T} \circ \DD _{T} \to \DD _{T} \circ  f _{+} \R \underline{\Gamma}  ^{\dag } _{X^{(0)}} \circ f ^{!} \circ \DD _{T}
\riso 
 f _{+} \circ \DD _{T^{(0)}}  \circ\R \underline{\Gamma}  ^{\dag } _{X^{(0)}} \circ f ^{!} \circ \DD _{T}$. 
On obtient ainsi le morphisme $Id \to a _{+} \circ a ^{+}$.
En utilisant les mêmes propriétés, on construit le morphisme d'adjonction 
$a ^{+} \circ a _{+}\to \mathrm{Id} $. 
\end{proof}

\begin{rema}
\label{Isoc*a+-rem-hdag}
Nous conservons les notations de \ref{Isoc*a+}.
Soit $\widetilde{T}$ un diviseur de $\PP$ contenant $T$. Posons $\widetilde{T} ^{(0)}:= f ^{-1} (\widetilde{T})$,
$\widetilde{b}\,:\, (X ^{(0)} \setminus \widetilde{T} ^{(0)}) \to (X \setminus \widetilde{T})$ le morphisme induit par $a$,
$\theta = (f,a,b)$ et $\widetilde{\theta} = (f,a,\widetilde{b})$ les morphismes de triplets lisses en dehors du diviseur correspondants.
 
Comme les foncteurs duaux commutent aux extensions des scalaires (voir \cite{virrion}), on dispose des foncteurs $(\hdag \widetilde{T})\,:\, \mathrm{Isoc} ^{*} (\PP, T, X/K)\to \mathrm{Isoc} ^{*} (\PP, \widetilde{T}, X/K)$
et $(\hdag \widetilde{T} ^{(0)})\,:\, \mathrm{Isoc} ^{*} (\PP^{(0)}, T^{(0)}, X^{(0)}/K)\to \mathrm{Isoc} ^{*} (\PP^{(0)}, \widetilde{T}^{(0)}, X^{(0)}/K)$.

Or, les foncteurs images inverses extraordinaires et images directes commutent aux foncteurs extension (voir \cite[2.2.18]{caro-construction}).
On en déduit que les foncteurs $ \theta ^{+}$ et $\theta _{+}$ de \ref{stabIsoc*inv} et \ref{Isoc*a+} commutent aux extensions des scalaires, i.e., on dispose des
isomorphismes canoniques $\widetilde{\theta} ^{+} \circ (\hdag \widetilde{T}) \riso (\hdag \widetilde{T} ^{(0)})\circ \theta ^{+}$
et $\widetilde{\theta} _{+} \circ (\hdag \widetilde{T} ^{(0)}) \riso (\hdag \widetilde{T} )\circ \theta _{+}$. 

\end{rema}

\subsection{Descente finie et étale en dehors des singularités surconvergentes\label{DescFiniEta}}

\begin{nota}
\label{nota-diag}
Dans toute la suite de cette section \ref{DescFiniEta}, nous reprenons les notations et hypothèses de 
\ref{nota-diag-pre}.
De plus, notons $\PP ^{(1)}:=\PP ^{(0)} \times _{\PP} \PP ^{(0)}$, 
$\PP ^{(2)}:=\PP^{(0)} \times _{\PP} \PP^{(0)} \times _{\PP} \PP^{(0)}$, 
$ f _{1},f _{2}\,:\, \PP ^{(1)} \to \PP^{(0)}$ les projections respectives à gauche et à droite, 
$f _{ij }\,:\, \PP^{(2)} \to \PP^{(1)} $ les projections sur les facteurs d'indices $i, j$ pour $i <j$. 
De même, on note $X^{(1)}:=X^{(0)} \times _{X} X ^{(0)}$, $X^{(2)} :=X^{(0)} \times _{X} X^{(0)} \times _{X} X ^{(0)}$, $a _{i}\,:\, X^{(1)} \to X^{(0)}$, $a _{ij}\,:\, X^{(2)} \to X^{(1)}$ les projections canoniques ; de même en remplaçant respectivement $a$ par $b$ et $X$ par $Y$.
On note $T ^{(1)}:=f _{1} ^{-1}(T ^{(0)})= f _{2} ^{-1}(T ^{(0)})$,
$T ^{(2)}:=f _{12} ^{-1}(T ^{(1)})=f _{23} ^{-1}(T ^{(1)})=f _{13} ^{-1}(T ^{(1)})$. 
\end{nota}

\begin{defi}
\label{defi-Isocbullet*}
On définit la catégorie $\mathrm{Isoc} ^{*} (\PP^{(\bullet)}, T^{(\bullet)}, X^{(\bullet)}/K)$ de la manière suivante: 

-ses objets sont les faisceaux $\E ^{(0)}\in \mathrm{Isoc} ^{*} (\PP^{(0)}, T^{(0)}, X^{(0)}/K)$ munis d'une donnée de recollement, i.e., d'un isomorphisme dans $\mathrm{Isoc} ^{*} (\PP^{(1)}, T^{(1)}, X^{(1)}/K)$ de la forme 
$\theta \,:\,a _{2} ^{+} (\E ^{(0)}) \riso a _{1} ^{+} (\E ^{(0)}) $ satisfaisant la condition de cocycle : 
$a _{13} ^+ (\theta) = a _{12} ^+ (\theta) \circ a _{23} ^+ (\theta)$ (cela a bien un sens grâce à la transitivité de la proposition \ref{Isoc*a+}).

-ses morphismes $(\E^{(0)}, \theta )\to (\FF^{(0)}, \tau)$ sont les morphismes $\D ^{\dag} _{\PP}(\hdag T) _{\Q}$-linéaires $\phi\,:\,\E^{(0)}\to \FF^{(0)}$ 
commutant aux données de recollement, i.e. tels que 
$\tau \circ a _{2} ^{+} (\phi) =a _{1} ^{+} (\phi)\circ \theta$.

Lorsque les diviseurs sont vides, comme d'habitude on ne les indique pas.
\end{defi}

\begin{vide}
\label{DefLoc}
On dispose du foncteur 
$\mathcal{L}oc\,:\, \mathrm{Isoc} ^{*} (\PP, T, X/K)\to \mathrm{Isoc} ^{*} (\PP^{(\bullet)}, T^{(\bullet)}, X^{(\bullet)}/K)$
défini par $\E \mapsto (a ^{+} (\E), \theta)$, 
où $\theta$ est l'isomorphisme induit par la transitivité pour la composition des foncteurs de la forme $a ^{+}$ (voir \ref{Isoc*a+}).
\end{vide}

\begin{vide}
\label{def-recol}
On construit le foncteur 
$\mathcal{R}ecol\,:\,
\mathrm{Isoc} ^{*} (\PP^{(\bullet)}, T^{(\bullet)}, X^{(\bullet)}/K)\to \mathrm{Isoc} ^{*} (\PP, T, X/K)$ de la manière suivante: 
soit $(\E^{(0)},\theta) \in \mathrm{Isoc} ^{*} (\PP^{(\bullet)}, T^{(\bullet)}, X^{(\bullet)}/K)$. 
Posons $\widetilde{a}:=  a   \circ a _{1}= a   \circ a _{2}$.
On bénéficie du morphisme composé 
$\phi _{1}\,:\,a _{+} ( \E^{(0)} ) \overset{\mathrm{adj} _{a _{1}}}{\longrightarrow}  a _{+}  \circ a _{1+} \circ a _{1} ^{+} (\E^{(0)}) \riso \tilde{a} _{+}\circ a _{1} ^{+} (\E^{(0)}) $, où le premier morphisme est induit par adjonction (voir \ref{Isoc*a+}.\ref{Isoc*a+5}). On dispose d'un second morphisme composé :
$\phi _{2}\,:\,a _{+} ( \E^{(0)} ) \overset{\mathrm{adj} _{a _{2}}}{\longrightarrow}  a _{+}  \circ a _{2+} \circ a _{2} ^{+} (\E^{(0)}) \riso \tilde{a} _{+} \circ a _{2} ^{+} (\E^{(0)}) \underset{\mathrm{\theta}}{\riso} \tilde{a} _{+} \circ a _{1} ^{+} (\E^{(0)}) $.
On pose alors
$ \mathcal{R}ecol (\E^{(0)},\theta) =\ker \left ( \xymatrix { { a _{+} ( \E^{(0)} )} \ar@<1ex>[r]  ^-{\phi _{1}}\ar@<-1ex>[r] _-{\phi _{2}}
  &
{\tilde{a} _{+} \circ a _{1} ^{+} (\E^{(0)}) }} \right )$.
Avec \ref{Isoc*stblKer} et \ref{Isoc*a+}, on vérifie que $\mathcal{R}ecol (\E^{(0)},\theta) \in \mathrm{Isoc} ^{*} (\PP, T, X/K)$.
\end{vide}

\begin{rema}
\label{Loc-Recol-hdag}
Soient $\widetilde{T}$ un diviseur de $\PP$ contenant $T$ et $\widetilde{T} ^{(0)}:= f ^{-1} (\widetilde{T})$.
Via la remarque de \ref{Isoc*a+-rem-hdag}, on vérifie que les foncteurs 
$\mathcal{L}oc$ de \ref{DefLoc} et 
$\mathcal{R}ecol$ de \ref{def-recol}
commutent aux extensions des scalaires, i.e., on dispose des
isomorphismes canoniques $\mathcal{R}ecol \circ (\hdag \widetilde{T}^{(0)}) \riso (\hdag \widetilde{T} )\circ \mathcal{R}ecol$
et $\mathcal{L}oc \circ (\hdag \widetilde{T} ) \riso (\hdag \widetilde{T} ^{(0)})\circ \mathcal{L}oc$ (les foncteurs $\mathcal{L}oc$ et $\mathcal{R}ecol$ 
indiquent soit les foncteurs de \ref{DefLoc} et \ref{def-recol}, soit les mêmes en remplaçant $T$ par $\widetilde{T}$).

\end{rema}

Avant de considérer le cas général de cette sous-section, traitons d'abord la descente finie et étale des isocristaux convergents sur les $k$-schémas lisses, ce qui correspondra à ce qui se passe en dehors des singularités surconvergentes :
\begin{lemm}
[Situation en dehors des singularités surconvergentes]
\label{fini-étale-présdag}
Posons $\widetilde{b}=  b   \circ b _{1}= b   \circ b _{2}$.
Les foncteurs 
$\mathcal{L}oc\,:\, \mathrm{Isoc} ^{*} (\U, Y/K)\to \mathrm{Isoc} ^{*} (\U^{(\bullet)}, Y^{(\bullet)}/K)$
et $\mathcal{R}ecol\,:\,
\mathrm{Isoc} ^{*} (\U^{(\bullet)}, Y^{(\bullet)}/K)\to \mathrm{Isoc} ^{*} (\U, Y/K)$
définis respectivement en \ref{DefLoc} et \ref{def-recol} 
sont quasi-inverses. 
\end{lemm}

\begin{proof}
Soient $\E \in \mathrm{Isoc} ^{*} (\U, Y/K)$ et $(\E ^{(0)},\theta) \in  \mathrm{Isoc} ^{*} (\U^{(\bullet)}, Y^{(\bullet)}/K)$.
D'après les étapes I.1) et II.1) de la preuve de \ref{desc-fini-ét} (ces deux étapes n'utilisent pas \ref{fini-étale-présdag}), 
on dispose des morphismes 
$\E \to \mathcal{R}ecol \circ \mathcal{L}oc(\E)$
et $b ^{+} \circ \mathcal{R}ecol (\E ^{(0)},\theta)\to \E ^{(0)}$.
La vérification que le premier (resp. second) morphisme est un isomorphisme (resp. commute aux données de recollement et est un isomorphisme)
est locale. Via le théorème de Berthelot-Kashiwara, on se ramène alors à supposer $Y=U$ et $Y^{(0)}=U^{(0)}$. 
Dans ce cas $g=b$ est un morphisme fini et étale et $b _{+}$ (resp. $b ^{+}$) est canoniquement isomorphe à $b _{*}$ (resp. $b ^{*}$). De même en rajoutant des indices $i$ ou $ij$.

Dans le cas où $g =b$, via la théorie classique de la descente fidèlement plate, on vérifie les assertions suivantes : 
\begin{enumerate}
\item La catégorie $\mathrm{Isoc} ^{\dag\dag} (Y, Y/K)$ est équivalente à celle des objets $\E^{(0)} \in \mathrm{Isoc} ^{\dag\dag } (Y  ^{(0)}, Y  ^{(0)}/K)$ munis d'une donnée de recollement, i.e., d'un isomorphisme 
$\theta\,:\,b ^{*} _{2} (\E^{(0)}) \riso b ^{*} _{1} (\E^{(0)})$ satisfaisant la condition usuelle de cocycle.

\item Le foncteur quasi-inverse de cette équivalence de catégories est donné par le foncteur recollement défini en posant $ \mathcal{R}ecol (\E  ^{(0)},\theta) =\ker \left ( \xymatrix { { b _{*} ( \E  ^{(0)} )} \ar@<1ex>[r]  ^-{\phi _{1}}\ar@<-1ex>[r] _-{\phi _{2}}
  &
  {\widetilde{b} _{*}  \circ b _{1} ^{*} (\E  ^{(0)}) }} \right )$, où $\phi _{1}\,:\,b _{*} ( \E  ^{(0)}) \overset{\mathrm{adj} _{b _{1}}}{\longrightarrow}  b _{*}  \circ b _{1*} \circ b _{1} ^{*} (\E  ^{(0)}) \riso \tilde{b} _{*}\circ b _{1} ^{*} (\E  ^{(0)}) $
et où
$\phi _{2}\,:\,b _{*} ( \E  ^{(0)}) \overset{\mathrm{adj} _{b _{2}}}{\longrightarrow}  b _{*}  \circ b _{2*} \circ b _{2} ^{*} (\E  ^{(0)}) \riso \tilde{b} _{*} \circ b _{2} ^{*} (\E  ^{(0)}) \underset{\mathrm{\theta}}{\riso} \tilde{b} _{*} \circ b _{1} ^{*} (\E  ^{(0)}) $.

\item Le morphisme induit par adjonction 
$b ^{*}\mathcal{R}ecol (\E  ^{(0)},\theta) \to (\E  ^{(0)},\theta)$ est un isomorphisme qui commute aux données de recollement (voir par exemple, 
\cite[Début de la page 19]{milne}).
\end{enumerate}

\end{proof}

\begin{prop}
\label{desc-fini-ét}
Les foncteurs 
$\mathcal{L}oc$ et $\mathcal{R}ecol$ 
induisent des équivalences quasi-inverses entre les catégories 
$\mathrm{Isoc} ^{*} (\PP, T, X/K)$ et $\mathrm{Isoc} ^{*} (\PP^{(\bullet)}, T^{(\bullet)}, X^{(\bullet)}/K)$.
\end{prop}

\begin{proof}
I)  Soit $\E \in \mathrm{Isoc} ^{*} (\PP, T, X/K)$. 
Établissons l'isomorphisme canonique 
$\E \riso \mathcal{R}ecol \circ \mathcal{L}oc(\E)$.

1) Construction de ce morphisme : 
on pose $(\E^{(0)}, \theta):= \mathcal{L}oc(\E)$. 
Considérons le diagramme ci-dessous :
\begin{equation}
\label{recol-loc=id}
\xymatrix{
{\E}
\ar[r] ^-{\mathrm{adj _{a}}} 
\ar@{=}[dd]
& 
{a _{+} ( \E^{(0)} )} 
\ar[r] ^-{\mathrm{adj} _{a _{2}}}
\ar[d] ^-{\mathrm{adj} _{a _{1}}}
&
{a _{+}  \circ a _{2+} \circ a _{2} ^{+} (\E^{(0)})}
\ar[r] ^-{\sim}
&
{\widetilde{a} _{+}  \circ a _{2} ^{+} (\E^{(0)})}
\ar[dl] ^-{\sim} _-{\theta}
\ar[dd] ^-{\sim}
\\ 
{} 
&
{a _{+}  \circ a _{1+} \circ a _{1} ^{+} (\E^{(0)})}
\ar[r] ^-{\sim}
\ar[d] ^-{\sim}
&
{\widetilde{a} _{+}  \circ a _{1} ^{+} (\E^{(0)})}
\ar[dl] ^-{\sim}
& 
{ } 
\\ 
{\E} 
\ar[r] ^-{\mathrm{adj _{\widetilde{a}}}}
&
{\widetilde{a} _{+}  \circ \widetilde{a} ^{+} (\E)}
\ar@{=}[rr]
&&
{\widetilde{a} _{+}  \circ \widetilde{a} ^{+} (\E).}
}
\end{equation}
Le triangle du milieu du bas est commutatif par définition de $a _{+}  \circ a _{1+} \circ a _{1} ^{+} (\E^{(0)})\riso \widetilde{a} _{+}  \circ \widetilde{a} ^{+} (\E)$. Le triangle de droite du bas est commutatif par définition de l'isomorphisme $\theta$.
Par transitivité du morphisme d'adjonction, le rectangle gauche du diagramme est commutative. 
Pour la même raison, le contour de \ref{recol-loc=id} est commutatif. 
Dans le trapèze du haut à droite, on remarque que le composé
$a _{+} ( \E^{(0)} ) \to \widetilde{a} _{+}  \circ a _{1} ^{+} (\E^{(0)}) $ passant le bas (resp. le haut) est $\phi _{1}$ (resp. $\phi _{2}$). On en déduit la factorisation: 
$\E \to \mathcal{R}ecol \circ \mathcal{L}oc(\E)$.

2) Vérifions à présent que celui-ci est un isomorphisme.
D'après \cite[4.3.12]{Be1}, il suffit de l'établir en dehors de $T$. 
On se ramène ainsi à la situation où $T$ est vide, i.e. au cas où $\U =\PP$, i.e., à la situation sans singularités surconvergentes de \ref{fini-étale-présdag}. 
D'où le résultat.

II) Réciproquement, soit $(\E^{(0)}, \theta)\in \mathrm{Isoc} ^{*} (\PP^{(\bullet)}, T^{(\bullet)}, X^{(\bullet)}/K)$
et vérifions l'isomorphisme $\mathcal{L}oc \circ \mathcal{R}ecol(\E^{(0)}, \theta) \riso (\E^{(0)}, \theta)$.

1) Construisons canoniquement ce morphisme. Par définition, on dispose de l'inclusion
$\mathcal{R}ecol(\E^{(0)}, \theta) \subset  a _{+} (\E^{(0)})$. 
Par adjonction, il en résulte le morphisme
$\phi\,:\,a ^{+} \circ \mathcal{R}ecol(\E^{(0)}, \theta) \to \E^{(0)} $.

2) Notons $\tau$ la donnée de recollement de $a ^{+} \circ \mathcal{R}ecol(\E^{(0)}, \theta)$ telle que 
$(a ^{+} \circ \mathcal{R}ecol(\E^{(0)}, \theta) ,\tau)= \mathcal{L}oc \circ \mathcal{R}ecol(\E^{(0)}, \theta)$.
Il s'agit de vérifier que $\phi $ est un isomorphisme commutant aux données de recollement respectives, i.e., tel que $\tau \circ a _{2} ^{+} (\phi) =a _{1} ^{+} (\phi)\circ \theta$.
Par \cite[4.3.12]{Be1}, il suffit de l'établir en dehors de $T$, ce qui nous ramène 
à la situation de \ref{fini-étale-présdag}. 
\end{proof}

\begin{coro}
\label{coro*=dag-desc-fini-ét}
Avec les notations de \ref{nota-diag}, on suppose en outre $X^{(0)}$ lisse sur $k$. 
On bénéficie alors des égalités : 
$\mathrm{Isoc} ^{*} (\PP, T, X/K)=\mathrm{Isoc} ^{**} (\PP, T, X/K)=\mathrm{Isoc} ^{\dag \dag} (\PP, T, X/K)$.
\end{coro}

\begin{proof}
Avec la seconde remarque de \ref{DDIsoc*=mod}, 
il suffit de prouver l'égalité $\mathrm{Isoc} ^{*} (\PP, T, X/K)=\mathrm{Isoc} ^{\dag \dag} (\PP, T, X/K)$.
Comme le foncteur $ \sp _{Y \hookrightarrow \U, + }$ commute aux foncteurs duaux respectifs, on obtient l'inclusion (pleinement fidèle)
$\mathrm{Isoc} ^{*} (\PP, T, X/K)\supset\mathrm{Isoc} ^{\dag \dag} (\PP, T, X/K)$.
Réciproquement, soit $\E \in \mathrm{Isoc} ^{*} (\PP, T, X/K)$.
Posons $a ^{!}:=\R \underline{\Gamma}  ^{\dag } _{X^{(0)}} \circ f ^{!} $.
Comme $\DD _{T} (\E)$ est $\D ^{\dag} _{\PP} (\hdag T) _{\Q}$-surcohérent, 
par stabilité de la surcohérence par image inverse extraordinaire et foncteur cohomologique local, 
$a ^{!} (\DD _{T} (\E))$ est $\D ^{\dag} _{\PP^{(0)}} (\hdag T^{(0)}) _{\Q}$-(sur)cohérent.
En outre, le faisceau $a ^{!} (\DD _{T} (\E))| \U ^{(0)}$ est dans l'image essentielle de 
$ \sp _{Y ^{(0)}\hookrightarrow \U^{(0)}, + }$.
Or, comme $X^{(0)}$ est lisse, on dispose de la caractérisation de \cite[2.5.10]{caro-construction} de l'image essentielle
du foncteur $ \sp _{X ^{(0)}\hookrightarrow \PP^{(0)}, T ^{(0)}, + }$.
Il en résulte aussitôt que $a ^{!} (\DD _{T} (\E))$ est dans l'image essentielle de $ \sp _{X ^{(0)}\hookrightarrow \PP^{(0)}, T ^{(0)}, + }$.
D'après \cite[6.1.4]{caro_devissge_surcoh}, il en découle que 
$a ^{!} (\DD _{T} (\E))$ vérifie la propriété $P _{\PP^{(0)}, T ^{(0)}}$ définie dans \cite[6.1.1]{caro_devissge_surcoh}. Cela implique en particulier que 
$a ^{+} (\E) $ est $\D ^{\dag} _{\PP ^{(0)}} (\hdag T^{(0)}) _{\Q}$-surcohérent et 
(grâce de plus à l'isomorphisme de bidualité) que $a ^{+} _{1}( a ^{+} (\E) )$ est $\D ^{\dag} _{\PP ^{(1)}} (\hdag (T^{(1)})) _{\Q}$-surcohérent. 
Par préservation de la surcohérence par noyau et par image directe par un morphisme propre, il en résulte que $\mathcal{R}ecol \circ \mathcal{L}oc (\E)$ est $\D ^{\dag} _{\PP} (\hdag T) _{\Q}$-surcohérent. 
Or, d'après \ref{desc-fini-ét}, $\E\riso \mathcal{R}ecol \circ \mathcal{L}oc (\E)$.
D'où le résultat.
\end{proof}

\subsection{Pleine fidélité du foncteur extension}

Avec les notations \ref{nota-6.2.1dev}, le théorème suivant est l'analogue du théorème de Tsuzuki \cite[4.1.1]{tsumono} ou de la version étendue de Kedlaya de \cite[5.2.1]{kedlaya-semistableI}:

\begin{theo}
\label{dagT'pl-fid}
Soient $\PP$ un $\V$-schéma formel séparé et lisse, $T \subset T'$ deux diviseurs de $P$, $X$ un sous-schéma fermé de $P$. En posant $Y:= X \setminus T$, $Y':= X \setminus T'$, on suppose de plus $Y$ lisse et $Y'$ dense dans $Y$.
Le foncteur $(\hdag T')$ induit les foncteurs pleinement fidèle :
\begin{gather}
\label{dagT'pl-fid-fonct}
(\hdag T')\,:\, (F\text{-})\mathrm{Isoc} ^{*} (\PP, T, X/K) \to (F\text{-})\mathrm{Isoc} ^{*} (\PP, T', X/K),
\\
\label{dagT'pl-fid-fonctbis}
(\hdag T')\,:\, (F\text{-})\mathrm{Isoc} ^{**} (\PP, T, X/K) \to (F\text{-})\mathrm{Isoc} ^{**} (\PP, T', X/K),
\\
(\hdag T')\,:\, (F\text{-})\mathrm{Isoc} ^{\dag \dag} (\PP, T, X/K) \to (F\text{-})\mathrm{Isoc} ^{\dag\dag} (\PP, T', X/K).
\end{gather}
\end{theo}

\begin{proof}
Pour alléger les notations, on omettra d'indiquer {\og $(F\text{-})$\fg}.
Comme le foncteur $\DD _{T}$ (resp. $\DD _{T'}$) induit une équivalence entre $\mathrm{Isoc} ^{*} (\PP, T, X/K)$ et $\mathrm{Isoc} ^{**} (\PP, T, X/K)$
(resp. entre $\mathrm{Isoc} ^{*} (\PP, T', X/K)$ et $\mathrm{Isoc} ^{**} (\PP, T', X/K)$), comme on dispose de l'isomorphisme canonique 
$\DD _{T'} \circ (\hdag T') \riso (\hdag T')  \circ \DD _{T}$, il suffit de vérifier \ref{dagT'pl-fid-fonctbis}.

Soient $\E _{1}$, $\E _{2}$ deux objets de $(F\text{-})\mathrm{Isoc} ^{* *} (\PP, T, X/K)$.
Posons $\U:= \PP \setminus T$, $\U':= \PP \setminus T'$.
Via le lemme \ref{rema-isocDagDagDense} et la remarque \ref{rema-isocDagDense},
on se ramène au cas où $Y$ est intègre et dense dans $X$.

1) 
{\it Pour tout $i=1,2$, le morphisme canonique 
$\E _{i} \to \E _{i} (\hdag T')$
est injectif.} 

Preuve : soit $\E' _{i}$ le noyau de $\E _{i} \to \E _{i} (\hdag T')$.
Par \cite[4.3.12]{Be1}, comme $\E' _{i}$ est un $\D ^{\dag} _{\PP} (\hdag T) _{\Q}$-module (sur)cohérent, 
il suffit d'établir que $\E ' _{i}|\U =0$. On se ramène ainsi au cas où $T$ est vide. Dans ce cas comme $X=Y$ est lisse et que l'assertion est locale en $\PP$, on peut supposer qu'il existe un morphisme $\X \hookrightarrow \PP$ de $\V$-schémas formels affines et lisses relevant $X \hookrightarrow P$. 
Comme $Y'$ est dense dans $Y=X$, $T'$ ne contient pas $X$ et donc $T' \cap X$ est un diviseur de $X$ (car $X$ est intègre et lisse). 
En utilisant le théorème de Berthelot-Kashiwara, il en dérive que les foncteurs $u _{T',+}$ et $u ^{!} _{T'}$ (resp. $u _{+}$ et $u ^{!}$) induisent des équivalences de catégories quasi-inverses entre $\mathrm{Isoc} ^{* *} (\X, T'\cap X, X/K)$ et $\mathrm{Isoc} ^{* *} (\PP, T', X/K)$
(resp. entre $\mathrm{Isoc} ^{* *} (\X, X/K)$ et $\mathrm{Isoc} ^{* *} (\PP,X/K)$).
On se ramène ainsi à traiter le cas où $X =P$ ($P$ est toujours affine et $T$ vide). 
Or, la catégorie $\mathrm{Isoc} ^{* *} (\PP, P/K)$ est égale à la catégorie des $\D ^{\dag} _{\PP,\Q}$-modules cohérents, $\O _{\PP,\Q}$-cohérents.
Via l'équivalence canonique entre cette catégorie avec celle des isocristaux convergents sur $P$ (via l'image inverse par le morphisme de spécialisation : voir \cite[4.1.4]{Be1}), grâce aux théorèmes de type $A$ pour les $\O _{\PP,\Q}$-modules cohérents, on vérifie  
que les objets de $\mathrm{Isoc} ^{* *} (\PP, P/K)$ sont plus précisément des $\O _{\PP,\Q}$-modules projectifs de type fini.
En particulier, $\E _{1}$ et $\E _{2}$ sont des $\O _{\PP,\Q}$-modules projectifs. Comme le morphisme canonique 
$\O _{\PP,\Q} \to \O _{\PP} (\hdag T') _{\Q}$ est injectif, il en est alors de même de $\E _{i} \to \E _{i} (\hdag T')$.

2) Il résulte de l'étape $1)$ que le foncteur \ref{dagT'pl-fid-fonctbis} est fidèle. 

3) Soient $\phi\,:\, \E _{1} (\hdag T')\to \E _{2}(\hdag T')$ un morphisme de $\mathrm{Isoc} ^{* *} (\PP, T', X/K)$. Il reste à vérifier que $\phi$ provient par extension d'un morphisme de la forme $\E _{1}\to \E _{2}$. D'après l'étape $2)$, cela est local. On peut supposer $\PP$ affine et muni de coordonnées locales.

4) Notons $\G$ l'image du composé : $\E  _{1} \hookrightarrow   \E  _{1} (\hdag T') 
\overset{\phi}{\longrightarrow} \E  _{2} (\hdag T')$.
Comme $\G$, $\E _{2}$ sont des sous-$\D ^{\dag} _{\PP} (\hdag T) _{\Q}$-modules (sur)cohérents
du $\D ^{\dag} _{\PP} (\hdag T) _{\Q}$-module (sur)cohérent $ \E  _{2} (\hdag T')$, le faisceau
$\G \cap \E  _{2}$ (égal par définition au noyau de $\G \to \E  _{2} (\hdag T')/\E  _{2}$) est aussi un sous-$\D ^{\dag} _{\PP} (\hdag T) _{\Q}$-module (sur)cohérent
de $\E  _{2} (\hdag T')$.

5) {\it L'inclusion $\G \cap \E  _{2} \subset \G$ est un isomorphisme.}

Preuve: via \cite[4.3.12]{Be1}, on se ramène 
au cas où $T$ est vide.
De plus, via le théorème de Berthelot-Kashiwara, il ne coûte pas cher de supposer $X =P$.
En utilisant le théorème de type $A$ pour les $\D ^{\dag} _{\PP,\Q}$-modules cohérents, comme $\G$ est l'image du morphisme 
de $\D ^{\dag} _{\PP,\Q}$-modules cohérents de la forme $\E _{1}\to \E _{2} (\hdag T')$, 
on obtient une surjection $\Gamma (\PP,\D ^{\dag} _{\PP,\Q})$-linéaire
$\Gamma (\PP, \E _{1})\to \Gamma (\PP, \G)$. 
Or, via le théorème de type $A$ pour les $\O _{\PP,\Q}$-modules cohérents, 
$\Gamma (\PP, \E _{1})$ est $\Gamma (\PP, \O _{\PP,\Q})$-cohérent. 
Par nothérianité de $\Gamma (\PP, \O _{\PP,\Q})$, il en résulte que
$\Gamma (\PP, \G)$ et 
$\Gamma (\PP, \G \cap \E _{2})$
sont $\Gamma (\PP, \O _{\PP,\Q})$-cohérents. 
Comme $\G$ et $\G \cap \E _{2}$ sont en outre $\D ^{\dag} _{\PP,\Q}$-cohérents, on déduit de \cite[2.2.13]{caro_courbe-nouveau}
que $\G$ et $\G \cap \E _{2}$ sont $\O _{\PP,\Q}$-cohérents.
Le morphisme $\G \cap \E  _{2} \subset \G$ est ainsi un morphisme
de $\D ^{\dag} _{\PP,\Q}$-modules cohérents, $\O _{\PP,\Q}$-cohérents.
Donc, $\G / \G \cap \E  _{2}$ est aussi un $\D ^{\dag} _{\PP,\Q}$-module cohérent, $\O _{\PP,\Q}$-cohérent. 
Il en donc associé à un isocristal convergent $E$ sur $P$ (voir \cite[4.1.4]{Be1}). 
Comme $\G / \G \cap \E  _{2}| \U' =0$, l'isocristal convergent sur $U'$ déduit de $E$ est nul. Comme $\U'$ est dense dans $\PP$, il en résulte que $E=0$, i.e. $\G / \G \cap \E  _{2}=0$.

6) On obtient le morphisme $\theta\,:\,\E _{1} \twoheadrightarrow \G \liso \G \cap \E _{2} \hookrightarrow \E _{2}$ dont le composé avec $\E _{2} \to \E _{2} (\hdag T')$ est égal au morphisme composé
$\E _{1} \hookrightarrow \E _{1} (\hdag T') \overset{\phi}{\longrightarrow} \E  _{2} (\hdag T')$.
Comme  $(\hdag T') (\theta)|\U'= \phi |\U'$, il en résulte par fidélité de $|\U'$ l'égalité
$(\hdag T') (\theta)= \phi$.
\end{proof}

\begin{rema}
Le théorème \ref{dagT'pl-fid} est faux si $Y'$ n'est pas dense dans $Y$. 
Par exemple, si on prend $\PP$ une courbe, $Y$ un point et $Y'$ l'ensemble vide. 
\end{rema}

\subsection{Pleine fidélité du foncteur {\og restriction-image inverse\fg} (par un morphisme génériquement fini et étale)}
\label{formel-II}

Nous conserverons dans toute la suite de cette sous-section 
\ref{formel-II} 
les notations suivantes : soit le diagramme commutatif
\begin{equation}
\label{formel-II-diag}
\xymatrix{
{\widetilde{Y}^{(0)}} \ar@{^{(}->}[r]  ^-{l ^{(0)}} \ar[d] ^-{c}
& 
{Y^{(0)}} \ar@{^{(}->}[r]  ^-{j ^{(0)}} \ar[d] ^-{b}
& 
{X^{(0)} } \ar@{^{(}->}[r] ^-{u^{(0)}} \ar[d] ^-{a}
& 
{\PP ^{(0)} }  \ar[d] ^-{f}
\\ 
{\widetilde{Y}} \ar@{^{(}->}[r] ^-{l}
&
{Y} \ar@{^{(}->}[r] ^-{j }
&
 {X } \ar@{^{(}->}[r] ^-{u} & {\PP  ,}
}
\end{equation}
où les deux carrés de gauche sont cartésiens, 
$f$ est un morphisme propre et lisse de $\V$-schémas formels séparés et lisses,
$a$ est un morphisme propre, surjectif de $k$-variétés, 
$b$ est un morphisme de $k$-variétés lisses, 
$c$ est un morphisme fini et étale, 
$l$, $l^{(0)}$, $j$ et $j^{(0)}$ sont des immersions ouvertes, $u$ et $u^{(0)}$ sont des immersions fermées,
$\widetilde{Y}$ est dense dans $Y$
et 
$\widetilde{Y} ^{(0)}$ est dense dans $Y ^{(0)}$.
On note $\widetilde{j}\,:\, \widetilde{Y}\hookrightarrow X$ et $\widetilde{j}^{(0)}\,:\, \widetilde{Y} ^{(0)}\hookrightarrow X ^{(0)}$ les immersions ouvertes induites. 
On suppose en outre
qu'il existe un diviseur $T$ (resp. $\widetilde{T}$) de $P$ tel que 
$Y = X \setminus T$ (resp. $\widetilde{Y} = X \setminus  \widetilde{T}$). On note $\U := \PP \setminus T$, 
$T ^{(0)}:=f ^{-1}(T)$, $\U ^{(0)}:= \PP ^{(0)}\setminus T^{(0)}$ et $g \,:\, \U ^{(0)} \to \U$ le morphisme induit par $f$.
De même en rajoutant des tildes (e.g. on pose $\widetilde{\U} := \PP \setminus \widetilde{T}$ etc.).

Notons $\PP ^{(1)}:=\PP ^{(0)} \times _{\PP} \PP ^{(0)}$, 
$\PP ^{(2)}:=\PP^{(0)} \times _{\PP} \PP^{(0)} \times _{\PP} \PP^{(0)}$, 
$ f _{1},f _{2}\,:\, \PP ^{(1)} \to \PP^{(0)}$ les projections respectives à gauche et à droite, 
$f _{ij }\,:\, \PP^{(2)} \to \PP^{(1)} $ les projections sur les facteurs d'indices $i, j$ pour $i <j$. 
De même, on note $X^{(1)}:=X^{(0)} \times _{X} X ^{(0)}$, $X^{(2)} :=X^{(0)} \times _{X} X^{(0)} \times _{X} X ^{(0)}$, $a _{i}\,:\, X^{(1)} \to X^{(0)}$, $a _{ij}\,:\, X^{(2)} \to X^{(1)}$ les projections canoniques ; de même en remplaçant respectivement $a$ par $b$ 
ou $c$ et $X$ par $Y$ ou $\widetilde{Y}$.
On note $T ^{(1)}:=f _{1} ^{-1}(T ^{(0)})\cup f _{2} ^{-1}(T ^{(0)})$,
$T ^{(2)}:=f _{12} ^{-1}(T ^{(1)})\cup f _{23} ^{-1}(T ^{(1)})\cup f _{13} ^{-1}(T ^{(1)})$ ; de même avec des tildes.

\begin{lemm}
\label{a+|U-plfid-lemm}

Le foncteur canonique
\begin{equation}
(a ^{+},\, |\widetilde{\U})\,:\,
\mathrm{Isoc} ^{*} (\PP, \widetilde{T}, X/K) \to 
\mathrm{Isoc} ^{*} (\PP ^{(0)}, \widetilde{T} ^{(0)}, X^{(0)}K) 
\times _{\mathrm{Isoc} ^{*} (\widetilde{\U} ^{(0)}, \widetilde{Y} ^{(0)}/K)}
\mathrm{Isoc} ^{*} (\widetilde{\U}, \widetilde{Y}/K)
\end{equation}
est pleinement fidèle.
\end{lemm}

\begin{proof}
Comme le foncteur $|\widetilde{\U}$ est fidèle, le foncteur $(a ^{+},\, |\widetilde{\U})$ est fidèle. Il reste à établir que cette fidélité est pleine.
Soient $\E _{1}, \E _{2} \in \mathrm{Isoc} ^{*} (\PP, \widetilde{T}, X/K)$.
En posant $\E _{1} ^{(0)}:= a ^{+} (\E _{1})$, $\E _{2} ^{(0)}:= a ^{+} (\E _{2})$
 soient $\phi ^{(0)} \,:\, \E _{1} ^{(0)}  \to \E _{2}  ^{(0)}$ et
$\psi \,:\, \E _{1} |\widetilde{\U} \to \E _{2} |\widetilde{\U}$ deux morphismes induisant canoniquement le même morphisme dans 
$\mathrm{Isoc} ^{*} (\widetilde{\U} ^{(0)}, \widetilde{Y} ^{(0)}/K)$. Considérons le carré de 
$\mathrm{Isoc} ^{*} (\PP ^{(1)}, \widetilde{T} ^{(1)}, X^{(1)}/K)$ : 
\begin{equation}
\label{diag-lemm(a+,|U)}
\xymatrix{
{a _{1} ^{+} (\E _{1} ^{(0)})} 
\ar[r] ^-{\sim} _-{\theta _{1}} \ar[d] ^-{a _{1} ^{+} (\phi ^{(0)})}
& 
{a _{2} ^{+} (\E _{1} ^{(0)})} 
\ar[d] ^-{a _{2} ^{+} (\phi ^{(0)})}
\\ 
{a _{1} ^{+} (\E _{2} ^{(0)})} 
\ar[r] ^-{\sim} _-{\theta _{2}} 
& 
{a _{2} ^{+} (\E _{2} ^{(0)}),} 
} 
\end{equation}
où $\theta _{1}$ et $\theta _{2}$ sont les isomorphismes canoniques induits par transitivité (voir \ref{stabIsoc*inv}). 
Comme les morphismes $\phi ^{(0)}$ et $\psi$ sont compatibles, 
le diagramme \ref{diag-lemm(a+,|U)} devient commutatif après application du foncteur 
$| \U ^{(1)}$. 
Grâce à la proposition \cite[4.3.12]{Be1},
on en déduit la commutativité de \ref{diag-lemm(a+,|U)}.
On dispose donc du morphisme 
$\phi ^{(0)}\,:\,(\E _{1} ^{(0)}, \theta _{1}) \to \E _{2} ^{(0)}, \theta _{2})$ de 
$\mathrm{Isoc} ^{*} (\PP^{(\bullet)}, \widetilde{T}^{(\bullet)}, X^{(\bullet)}/K)$.
Or, pour $i=1,2$, on a par définition $\mathcal{L}oc (\E _{i} ) =(\E _{i} ^{(0)}, \theta _{i}) $.
Comme le foncteur $\mathcal{L}oc$ est pleinement fidèle (voir \ref{desc-fini-ét}), il existe 
un morphisme 
$\phi \,:\, \E _{1}  \to \E _{2}$ tel que $\mathcal{L}oc (\phi )= \phi ^{(0)}$, i.e., $a ^{+} (\phi) =\phi ^{(0)}$.  
Comme $\mathcal{L}oc ( \phi | \widetilde{\U}) = \phi ^{(0)} | \widetilde{\U} ^{(0)}=\mathcal{L}oc (\psi)$, 
par fidélité de $\mathcal{L}oc$,
il en résulte que $\phi | \widetilde{\U} = \psi$.
On a donc vérifié $(a ^{+},\, |\widetilde{\U}) (\phi)=(\phi ^{(0)}, \psi)$.
\end{proof}

\begin{prop}
\label{a+|U-plfid-prop}
Le foncteur canonique  
\begin{equation}
(a ^{+},\, |\U)\,:\,
\mathrm{Isoc} ^{*} (\PP, T, X/K) 
\to 
\mathrm{Isoc} ^{*} (\PP ^{(0)}, T ^{(0)}, X^{(0)}K) 
\times _{\mathrm{Isoc} ^{*} (\U ^{(0)}, Y ^{(0)}/K)}
\mathrm{Isoc} ^{*} (\U, Y/K)
\end{equation}
est pleinement fidèle.
\end{prop}

\begin{proof}
Afin de faire la distinction entre les deux foncteurs $a ^{+}$, notons $\theta = (f,a,b)$ et $\widetilde{\theta} = (f,a,c)$ les morphisme de triplets lisses en dehors du diviseur.
Considérons le diagramme suivant
\begin{equation}
\label{a+|U-plfid-prop-diag1}
\xymatrix{
{\mathrm{Isoc} ^{*} (\PP, T, X/K) } 
\ar[r] ^-{(\theta ^{+},\, |\U)}
\ar[d] ^-{(\hdag \widetilde{T})}
& { \mathrm{Isoc} ^{*} (\PP ^{(0)}, T ^{(0)}, X^{(0)}K) 
\times _{\mathrm{Isoc} ^{*} (\U ^{(0)}, Y ^{(0)}/K)}
\mathrm{Isoc} ^{*} (\U, Y/K)} 
\ar[d] ^-{((\hdag \widetilde{T}^{(0)}), |\widetilde{\U} )}
\\ 
{\mathrm{Isoc} ^{*} (\PP, \widetilde{T}, X/K) } 
\ar[r] ^-{(\widetilde{\theta} ^{+},\, |\widetilde{\U})}
& 
{\mathrm{Isoc} ^{*} (\PP ^{(0)}, \widetilde{T} ^{(0)}, X^{(0)}K) 
\times _{\mathrm{Isoc} ^{*} (\widetilde{\U} ^{(0)}, \widetilde{Y} ^{(0)}/K)}
\mathrm{Isoc} ^{*} (\widetilde{\U}, \widetilde{Y}/K) } .}
\end{equation}
D'après \cite{caro_surcoherent}, les foncteurs extensions commutent aux images inverses extraordinaires et foncteurs cohomologiques à support 
propre dans un sous-schéma fermé. D'après Virrion (voir \cite{virrion}), les foncteurs extensions commutent aux foncteurs duaux. 
On en déduit la commutativité à isomorphisme canonique près du diagramme \ref{a+|U-plfid-prop-diag1}.
Comme d'après \ref{dagT'pl-fid} (resp. \ref{a+|U-plfid-lemm}) le foncteur de gauche (resp. du bas) 
est pleinement fidèle, il en est alors de même de celui du haut. 
\end{proof}

\subsection{Une équivalence de catégories induite par le foncteur {\og extension-image inverse\fg}}

Dans cette section, nous reprenons les notations et hypothèses de la section \ref{formel-II}.

\begin{lemm}
\label{lemmTh2-cas lisse}
Le foncteur canonique 
\begin{equation}
(a ^{+},\, (\hdag \widetilde{T}) )\,:\,
\mathrm{Isoc} ^{*} (\PP, T, X/K) 
\to 
\mathrm{Isoc} ^{*} (\PP ^{(0)}, T ^{(0)}, X^{(0)}/K) 
\times _{\mathrm{Isoc} ^{*} (\PP ^{(0)}, \widetilde{T} ^{(0)}, X^{(0)}/K)}
\mathrm{Isoc} ^{*} (\PP, \widetilde{T}, X/K)
\end{equation}
est pleinement fidèle. 
\end{lemm}

\begin{proof}
Comme le foncteur $(\hdag \widetilde{T})$ est fidèle, $(a ^{+},\, (\hdag \widetilde{T}) )$ est fidèle. 
Comme $(\hdag \widetilde{T})$ est pleinement fidèle (voir le théorème \ref{dagT'pl-fid}) et comme $(\hdag \widetilde{T} ^{(0)}) $ est fidèle, on vérifie que la fidélité de
$(a ^{+},\, (\hdag \widetilde{T}) )$ est pleine.
\end{proof}

\begin{lemm}
\label{lemm2Th2-cas lisse}
On suppose $X$ lisse. 
Soit $\E \in \mathrm{Isoc} ^{*} (\PP, T, X/K) $. Alors 
$\E \riso \mathrm{Im} \left ( \DD _{T}\circ \DD _{\widetilde{T}} \circ (\hdag \widetilde{T})( \E) \to (\hdag \widetilde{T})( \E) \right )$. 
\end{lemm}

\begin{proof}
On note $\alpha $ les morphismes induits par fonctorialité par $Id \to (\hdag \widetilde{T}  )$. 
Comme $X$ est lisse, d'après \cite[6.1.4]{caro_devissge_surcoh}, $\E$ vérifie la propriété $P _{\PP, T}$ (voir la définition \cite[6.1.1]{caro_devissge_surcoh}). 
Via \ref{surcoh=>hol}, il en résulte que $\DD _{T} (\E (\hdag \widetilde{T})) $ est un $\D ^{\dag} _{\PP} (\hdag T) _{\Q}$-{\it module} surcohérent. 
Comme la catégorie $\mathrm{Isoc} ^{*} (\PP, T, X/K)$ est stable par le foncteur $\DD _{T}$ (voir \ref{Xlisse-Eqsp+}), 
on en déduit de même que
$\DD _{T} \circ \DD _{\widetilde{T} } (\E (\hdag \widetilde{T})) 
\riso \DD _{T} \circ  (\hdag \widetilde{T}) (\DD _{T} (\E) )$ est un $\D ^{\dag} _{\PP} (\hdag T) _{\Q}$-module surcohérent.
On définit $\FF $ comme égale à l'image du morphisme canonique composé de $\D ^{\dag} _{\PP} (\hdag T) _{\Q}$-modules surcohérents :
\begin{equation}
\label{contagi-def-beta}
\FF:= \mathrm{Im} \left (\DD _{T} \circ \DD _{\widetilde{T} } (\E (\hdag \widetilde{T})) \overset{\alpha}{\longrightarrow} (\hdag \widetilde{T} ) \circ \DD _{T} \circ \DD _{\widetilde{T} } (\E (\hdag \widetilde{T})) \riso 
\DD _{\widetilde{T} } \circ \DD _{\widetilde{T} } (\E (\hdag \widetilde{T})) \riso \E (\hdag \widetilde{T})\right ).
\end{equation}
Comme le morphisme \ref{contagi-def-beta} est un isomorphisme en dehors de $\widetilde{T} $, il en est de même du morphisme canonique
$\rho\,:\, (\hdag \widetilde{T} ) (\FF) \to \E (\hdag \widetilde{T})$ de $\D ^{\dag} _{\PP} (\hdag \widetilde{T} ) _{\Q}$-modules cohérents (qui factorise par définition l'inclusion canonique $\FF \subset \E (\hdag \widetilde{T})$). 
Par \cite[4.3.12]{Be1}, ce morphisme 
$\rho$
est donc un isomorphisme.

D'après l'étape $1)$ de la preuve de \ref{dagT'pl-fid}, 
comme $\E, \DD _{T} (\E) \in \mathrm{Isoc} ^{**} (\PP, T, X/K) $,
les morphismes canoniques $\alpha \,:\, \E \to (\hdag \widetilde{T}) ( \E)$ 
et
$\alpha \,:\, \DD _{T} (\E) \to (\hdag \widetilde{T}) ( \DD _{T} (\E))$ 
sont injectifs.
Grâce à \ref{dim-hom-m-dag}, 
en lui appliquant $\DD _{T} $, on obtient la surjection $\alpha ^{*}\,:\,   \DD _{T} (\hdag \widetilde{T}) ( \DD _{T} (\E)) \twoheadrightarrow \DD _{T} \DD _{T} (\E)$.
Comme le morphisme $(\hdag \widetilde{T}) (\alpha ^{*})$ est un isomorphisme en dehors de $\widetilde{T} $, celui-ci est un isomorphisme.
On obtient par fonctorialité la commutativité des carrés du bas du diagramme suivant: 
\begin{equation}
\label{contagiosité-2}
\xymatrix{
& 
{(\hdag \widetilde{T}) \circ \DD _{T} \circ \DD _{\widetilde{T}} ((\hdag \widetilde{T}) ( \E))  } 
\ar[r] ^-{\sim}
&
{\DD _{\widetilde{T}} \circ \DD _{\widetilde{T}} ((\hdag \widetilde{T}) ( \E))}
\ar[rd] ^-{\sim}
\\
{(\hdag \widetilde{T}) \circ \DD _{T} \circ \DD _{\widetilde{T}} \circ (\hdag \widetilde{T}) ( \E)} 
\ar[r] ^-{\sim}
\ar@{=}[ru] 
&
{(\hdag \widetilde{T})\circ \DD _{T} \circ (\hdag \widetilde{T}) (\DD _{T} ( \E))}  
\ar[r] ^-{\sim} _-{(\hdag \widetilde{T}) (\alpha ^{*})}
&
{(\hdag \widetilde{T})\circ \DD _{T}  (\DD _{T} ( \E))}  
\ar[r] ^-{\sim}
&
{(\hdag \widetilde{T}) ( \E)}
\\
{\DD _{T} \circ \DD _{\widetilde{T}} \circ (\hdag \widetilde{T}) ( \E)} 
\ar[u] ^-{\alpha}
\ar[r] ^-{\sim}
&
{\DD _{T} \circ (\hdag \widetilde{T}) (\DD _{T} ( \E))} 
\ar@{->>}[r] _-{\alpha ^{*}}
\ar[u] ^-{\alpha}
& 
{\DD _{T}  ( \DD _{T} ( \E))} 
\ar[u] ^-{\alpha}
\ar[r] ^-{\sim}
& 
{\E,} 
\ar@{_{(}->}[u] ^-{\alpha}
}
\end{equation}
où les isomorphismes horizontaux du carré du gauche (resp. de droite) sont induits par les isomorphismes 
$\DD _{\widetilde{T}} \circ (\hdag \widetilde{T}) ( \E) \riso  (\hdag \widetilde{T}) (\DD _{T} ( \E))$.
Pour valider la commutativité du trapèze (en haut) de \ref{contagiosité-2}, comme tous termes sont des 
$\D ^{\dag} _{\PP} (\hdag \widetilde{T}) _{\Q}$-modules cohérents, 
il suffit alors via \cite[4.3.12]{Be1} de l'établir en dehors de $\widetilde{T} $, ce qui est immédiat car le foncteur $(\hdag \widetilde{T})$ (resp. $\DD _{\widetilde{T}}$) est canoniquement isomorphe en dehors de $\widetilde{T} $ à l'identité (resp. à $\DD _{T}$) et car les morphismes de la forme $\alpha$ ou $\alpha ^{*}$ sont alors canoniquement égaux à l'identité.

On constate que le composé 
de \ref{contagi-def-beta} correspond au composé passant par la gauche puis le haut de \ref{contagiosité-2}.
On obtient alors 
le diagramme commutatif ci-dessous:
\begin{equation}
\label{contagiosité-3}
\xymatrix{
{\DD _{T} \circ \DD _{\widetilde{T}} (\E (\hdag \widetilde{T}))} 
\ar@{=}[d] 
\ar@{->>}[r]
& 
{\FF} 
\ar@{^{(}->}[r] ^-{}
&
{ \E (\hdag \widetilde{T})}
\ar@{=}[d]
\\ 
{\DD _{T} \circ \DD _{\widetilde{T}} ((\hdag \widetilde{T}) ( \E))} 
\ar@{->>}[r]
& 
{\E} 
\ar@{^{(}->}[r] ^-{\alpha}
&
{ (\hdag \widetilde{T}) ( \E),}
}
\end{equation}
dont la flèche surjective du bas à gauche est le composé du bas de \ref{contagiosité-2}.
Il en résulte l'isomorphisme canonique 
$\FF\riso \E$. 
\end{proof}

\begin{lemm}
\label{EDEsurcoh}
Soient $\E _{1},\E _{2}$ deux $\D ^{\dag} _{\PP}(\hdag T) _{\Q}$-modules surcohérents tels que 
$\DD _{T} (\E _{1}),\DD _{T} (\E _{2})$ soient $\D ^{\dag} _{\PP}(\hdag T) _{\Q}$-surcohérents.
Soit $\phi\,:\, \E _{1} \to \E _{2}$ un morphisme $\D ^{\dag} _{\PP}(\hdag T) _{\Q}$-linéaire. 
Notons $\E$ le noyau de $\phi$ et $\FF$ son image. 
Alors $\E$, $\mathcal{H} ^{0} \DD _{T}(\E)$, $\FF$, $\mathcal{H} ^{0} \DD _{T}(\FF)$ sont surcohérents et holonomes en tant que 
$\D ^{\dag} _{\PP}(\hdag T) _{\Q}$-modules.
\end{lemm}

\begin{proof}
Le Théorème \ref{surcoh=>hol} implique que $\E _{1}$ et $\E _{2}$ sont $\D ^{\dag} _{\PP}(\hdag T) _{\Q}$-holonomes.
Avec \ref{Stab-hol-dual} et \ref{hol-suiteexacte}, on valide alors la partie du lemme qui concerne la $\D ^{\dag} _{\PP}(\hdag T) _{\Q}$-holonomie.
De plus, comme la surcohérence est stable par noyau et image, $\E$ et $\FF$ sont
$\D ^{\dag} _{\PP}(\hdag T) _{\Q}$-surcohérents.
Comme $\mathcal{H} ^{1} \DD _{T}(\E _{2}/\FF)=0$ (voir \ref{dim-hom-m-dag}), en appliquant $ \mathcal{H} ^{0} \DD _{T} $ à l'injection
$\FF \hookrightarrow \E _{2}$, on obtient la surjection 
$ \mathcal{H} ^{0} \DD _{T} (\E _{2} ) \twoheadrightarrow  \mathcal{H} ^{0} \DD _{T} (\FF )$.
En appliquant le foncteur $\mathcal{H} ^{0} \DD _{T}$ à la suite exacte 
$0\to \E \to \E _{1} \to \FF \to 0$, 
comme 
$\mathcal{H} ^{-1} \DD _{T}(\E)=0$, 
on obtient la suite exacte 
$0\to \mathcal{H} ^{0} \DD _{T} (\FF)  
\to 
\mathcal{H} ^{0} \DD _{T} (\E _{1} ) 
\to 
\mathcal{H} ^{0} \DD _{T} (\E) \to 0$.
Il en résulte que $\mathcal{H} ^{0} \DD _{T} (\FF )$ est l'image du morphisme 
$ \mathcal{H} ^{0} \DD _{T} (\phi) \,:\,\mathcal{H} ^{0} \DD _{T} (\E _{2}) \to \mathcal{H} ^{0} \DD _{T} (\E _{1} )$.
Le module $\mathcal{H} ^{0} \DD _{T} (\FF )$ est donc $\D ^{\dag} _{\PP}(\hdag T) _{\Q}$-surcohérent.
Via la dernière suite exacte écrite ci-dessus, on en déduit que $\mathcal{H} ^{0} \DD _{T} (\E)$ est $\D ^{\dag} _{\PP}(\hdag T) _{\Q}$-surcohérent.
\end{proof}

\begin{lemm}
\label{a+a!=a+a+}
On suppose $X ^{(0)}$ lisse.
Soit $\widetilde{\E} ^{(0)} \in \mathrm{Isoc} ^{*} (\PP ^{(0)}, \widetilde{T} ^{(0)}, X^{(0)}/K) =\mathrm{Isoc} ^{**} (\PP ^{(0)}, \widetilde{T} ^{(0)}, X^{(0)}/K) $ (car $X ^{(0)}$ est lisse).
On dispose alors de l'isomorphisme canonique
$a _{1+} \circ a _{1} ^{+} ( \widetilde{\E} ^{(0)}) \riso a _{1+} \circ a _{1} ^{!} ( \widetilde{\E} ^{(0)})$.
\end{lemm}

\begin{proof}
D'après le théorème de désingularisation de de Jong (voir \cite{dejong} ou \cite{Ber-alterationdejong}), 
on dispose d'un morphisme surjectif, projectif, génériquement fini et étale
$\alpha\,:\, X'' \to X ^{(1)}$ 
tel que $X''$ soit lisse et $\alpha ^{-1} ( T ^{(1)} \cap X ^{(1)})$ soit un diviseur de $X ^{\prime \prime}$. Il existe ainsi un morphisme projectif et lisse de $\V$-schémas formels séparés et lisses de la forme $f ^{\prime \prime}\,:\,\PP'' \to \PP ^{(1)}$, une immersion fermée $X'' \hookrightarrow \PP''$ tels que leur composée donne le morphisme canonique $X'' \to \PP ^{(1)}$.
Posons $T '' := (f ^{\prime \prime })  ^{-1} (\widetilde{T} ^{(0)})$, 
$\alpha ^{!}:=\R \underline{\Gamma} ^\dag _{X''} \circ f ^{\prime \prime !}$,
$\alpha ^{+}:= \DD _{T''}\circ \R \underline{\Gamma} ^\dag _{X''} \circ f ^{\prime \prime !} \circ \DD _{\widetilde{T} ^{(0)}}$,
$\alpha _{+}:= f '' _{+}$
(notations analogues à \ref{stabIsoc*inv-rema}).
Comme $f ''$ est propre, on dispose du morphisme construit par adjonction:
$ \alpha _{+} \circ \alpha ^{!} \to Id$. En dualisant et via le théorème de dualité relative, on obtient le morphisme
$Id \to  \alpha _{+} \circ \alpha ^{+}$.
Or, comme $a _{1} \circ \alpha $ est un morphisme de variétés lisses, d'après les isomorphismes de commutation de \ref{Xlisse-Eqsp+} (ou voir \cite{caro-construction}), on dispose alors de l'isomorphisme canonique
$\alpha ^{+} \circ a _{1} ^{+} ( \widetilde{\E} ^{(0)}) 
\riso 
\alpha ^{!} \circ  a _{1} ^{!} ( \widetilde{\E} ^{(0)})$.
On en déduit la construction du morphisme composé: 
\begin{equation}
\notag
a _{1+} \circ a _{1} ^{+} ( \widetilde{\E} ^{(0)}) 
\to 
a _{1+} \circ \alpha _{+}\circ \alpha ^{+} \circ a _{1} ^{+} ( \widetilde{\E} ^{(0)}) 
\riso 
a _{1+} \circ \alpha _{+}\circ \alpha ^{!} \circ  a _{1} ^{!} ( \widetilde{\E} ^{(0)})
\to 
 a _{1+} \circ a _{1} ^{!} ( \widetilde{\E} ^{(0)}).
 \end{equation}
 Pour établir que ce morphisme est un isomorphisme, il suffit de le vérifier en dehors du diviseur 
 $\widetilde{T} ^{(0)}$, ce qui nous ramène au cas où le diviseur $\widetilde{T} ^{(0)}$ est vide
 et en particulier au cas où $X ^{(1)}$ est aussi lisse. Dans ce cas,
 on dispose de l'isomorphisme canonique
 $a _{1} ^{+} ( \widetilde{\E} ^{(0)}) \riso a _{1} ^{!} ( \widetilde{\E} ^{(0)})$ 
 et 
 le morphisme canonique 
$Id
\to 
\alpha _{+}\circ \alpha ^{+} 
\riso 
\alpha _{+} \alpha ^{!} 
\to 
Id$
de foncteurs définis sur $\mathrm{Isoc} ^{\dag \dag} (\PP ^{(1)}, X^{(1)}/K)$
 est un isomorphisme (car, comme il s'agit d'un morphisme d'isocristaux surconvergents, il suffit de le voir sur un ouvert dense de $X ^{(1)}$ où le morphisme induit par $\alpha$ est fini et étale).
\end{proof}

\begin{lemm}
\label{lemm-surcoh-ann-th2}
On suppose $X ^{(0)}$ lisse.
Soit $(\E ^{(0)}, \widetilde{\E}, \rho) \in \mathrm{Isoc} ^{*} (\PP ^{(0)}, T ^{(0)}, X^{(0)}/K) 
\times _{\mathrm{Isoc} ^{*} (\PP ^{(0)}, \widetilde{T} ^{(0)}, X^{(0)}/K)}
\mathrm{Isoc} ^{*} (\PP, \widetilde{T}, X/K)$.

Alors $\widetilde{\E}$, $\DD _{\widetilde{T}} (\widetilde{\E})$, $\DD _{T} (\widetilde{\E})$,
$\DD _{T} \circ \DD _{\widetilde{T}} (\widetilde{\E})$ sont $\D ^{\dag} _{\PP}(\hdag T) _{\Q}$-surcohérents.
De plus, 
pour tout $l\not = 0$, $\mathcal{H} ^{l}\DD _{T} (\widetilde{\E})=0$,
$\mathcal{H} ^{l} \DD _{T} \circ \DD _{\widetilde{T}} (\widetilde{\E})=0$.
\end{lemm}

\begin{proof}
Par définition,
$\mathcal{L}oc (\widetilde{\E})$ est de la forme $ (a ^{+} (\widetilde{\E}), \theta)$
(voir \ref{DefLoc}).
Par hypothèse, on dispose de l'isomorphisme $\rho\,:\, a ^{+} (\widetilde{\E}) \riso \E ^{(0)} (\hdag \widetilde{T} ^{(0)})$.
Par définition de \ref{DefLoc} du foncteur $ \mathcal{R}ecol$, via les Propositions \ref{desc-fini-ét} et \ref{a+a!=a+a+}, 
$\widetilde{\E}$ est alors isomorphe au noyau d'un morphisme de la forme
$\phi \,:\, a _{+} (\E ^{(0)} (\hdag \widetilde{T} ^{(0)}) ) \to \tilde{a} _{+} \circ a _{1} ^{!} (\E ^{(0)} (\hdag \widetilde{T} ^{(0)})) $.

Or, comme $X ^{(0)}$ est lisse, il résulte de \cite[6.1.4]{caro_devissge_surcoh} que $\E ^{(0)} $ vérifie la propriété $P _{\PP ^{(0)}, T ^{(0)}}$.
Par stabilité de cette propriété par image inverse extraordinaire, foncteur de localisation et image directe par un morphisme propre 
(voir \cite[6.1.3]{caro_devissge_surcoh}),
il en résulte que les deux termes du morphisme $\phi$ vérifient la propriété $P _{\PP, T }$.
Le Lemme \ref{EDEsurcoh} nous permet de conclure.
\end{proof}

\begin{nota}
\label{nota-Ttilde-T}
On note $\mathrm{Isoc} ^{*} (\PP, \widetilde{T} \supset T, X/K)$ la sous-catégorie pleine de
$\mathrm{Isoc} ^{*} (\PP, \widetilde{T} , X/K)$ des modules $\E$ tels qu'il existe $\G \in 
\mathrm{Isoc} ^{\dag \dag} (\U,Y/K)$ et un isomorphisme de la forme
$\E |\U \riso (\hdag \widetilde{T} \cap U) (\G)$.

On dispose de la factorisation 
$(\hdag \widetilde{T})\,:\,\mathrm{Isoc} ^{*} (\PP, T, X/K) \to  
\mathrm{Isoc} ^{*} (\PP, \widetilde{T} \supset T, X/K)$.
\end{nota}

\begin{theo}
\label{Th2}
On suppose $X ^{(0)}$ lisse. 
Avec les notations de \ref{nota-Ttilde-T}, le foncteur canonique 
\begin{equation}
(a ^{+},\, (\hdag \widetilde{T}) )\,:\,
\mathrm{Isoc} ^{*} (\PP, T, X/K) 
\to 
\mathrm{Isoc} ^{*} (\PP ^{(0)}, T ^{(0)}, X^{(0)}/K) 
\times _{\mathrm{Isoc} ^{*} (\PP ^{(0)}, \widetilde{T} ^{(0)}, X^{(0)}/K)}
\mathrm{Isoc} ^{*} (\PP, \widetilde{T} \supset T, X/K)
\end{equation}
est une équivalence de catégories. On dispose de plus d'un foncteur canonique
quasi-inverse noté $\mathcal{R}ecol$. 
\end{theo}

\begin{proof}
I) Le foncteur $(a ^{+},\, (\hdag \widetilde{T}) )$ est pleinement fidèle en vertu du Lemme \ref{lemmTh2-cas lisse}.

II) Construisons à présent le foncteur $\mathcal{R}ecol$ quasi-inverse de $(a ^{+},\, (\hdag \widetilde{T}) )$.

Soit $(\E ^{(0)}, \widetilde{\E}, \rho) \in \mathrm{Isoc} ^{*} (\PP ^{(0)}, T ^{(0)}, X^{(0)}/K) 
\times _{\mathrm{Isoc} ^{*} (\PP ^{(0)}, \widetilde{T} ^{(0)}, X^{(0)}/K)}
\mathrm{Isoc} ^{*} (\PP, \widetilde{T} \supset T, X/K)$.
De manière analogue à \ref{contagi-def-beta}, on construit grâce à \ref{lemm-surcoh-ann-th2}
le morphisme canonique de $\D ^{\dag} _{\PP}(\hdag T) _{\Q}$-modules surcohérents
$\phi \,:\, \DD _{T} \circ \DD _{\widetilde{T} } (\widetilde{\E}) \to  \widetilde{\E}$ tel que de plus $\DD _{T}(\phi)$ soit aussi un morphisme
de $\D ^{\dag} _{\PP}(\hdag T) _{\Q}$-modules surcohérents.
On définit alors un $\D ^{\dag} _{\PP}(\hdag T) _{\Q}$-module $\E$ en posant
\begin{equation}
\label{recol-th2}
\E:= \mathcal{R}ecol (\E ^{(0)}, \widetilde{\E}, \rho):=\mathrm{Im} (\phi).
\end{equation}
D'après le Lemme \ref{EDEsurcoh}, $\E$ et $\DD _{T}(\E)$ sont des $\D ^{\dag} _{\PP}(\hdag T) _{\Q}$-modules surcohérents.
De plus, en se rappelant de la définition \ref{nota-Ttilde-T}, 
il résulte de \ref{lemm2Th2-cas lisse} que $\E |\U \in \mathrm{Isoc} ^{*} (\U,Y/K) $, i.e. est dans l'image essentielle du foncteur $\sp _{Y\hookrightarrow \U, +}$.
Nous avons donc prouvé que $\E \in \mathrm{Isoc} ^{*} (\PP, T, X/K)$ (et même que $\E \in \mathrm{Isoc} ^{\dag \dag} (\PP, T, X/K)$).

Comme $(\hdag \widetilde{T}) (\phi)$ est un isomorphisme, 
la flèche canonique
$\E (\hdag \widetilde{T}) \to \widetilde{\E}$ est aussi un isomorphisme.
Par pleine fidélité de $(\hdag \widetilde{T} ^{(0)})$, on en déduit que 
$(a ^{+},\, (\hdag \widetilde{T}) ) \circ \mathcal{R}ecol (\E ^{(0)} , \widetilde{\E}, \rho)= (a ^{+},\, (\hdag \widetilde{T}) ) (\E) \riso  (\E ^{(0)},\widetilde{\E},\rho) $.
Comme le foncteur $(a ^{+},\, (\hdag \widetilde{T}) )$ est pleinement fidèle, il en résulte que 
les foncteurs $(a ^{+},\, (\hdag \widetilde{T}) )$ et $\mathcal{R}ecol$ sont quasi-inverses.
\end{proof}

\begin{rema}
\label{rema-recol-factdagdag}
D'après la preuve de \ref{Th2}, le foncteur $\mathcal{R}ecol$ de \ref{recol-th2} se factorise de la manière suivante :
$$\mathcal{R}ecol\,:\,
\mathrm{Isoc} ^{*} (\PP ^{(0)}, T ^{(0)}, X^{(0)}/K) 
\times _{\mathrm{Isoc} ^{*} (\PP ^{(0)}, \widetilde{T} ^{(0)}, X^{(0)}/K)}
\mathrm{Isoc} ^{*} (\PP, \widetilde{T} \supset T, X/K)
\to 
\mathrm{Isoc} ^{\dag \dag} (\PP, T, X/K) .$$
\end{rema}

\begin{coro}
\label{Isoc*=dagdag}
Soit $(\PP, T,X)$ un triplet lisse en dehors du diviseur (voir la définition \ref{defi-triplet}).
On bénéficie des égalités 
$\mathrm{Isoc} ^{*} (\PP, T, X/K)=\mathrm{Isoc} ^{**} (\PP, T, X/K)=\mathrm{Isoc} ^{\dag \dag} (\PP, T, X/K)$.
De plus, le foncteur $\DD _{T}$ préserve ces catégories.
\end{coro}

\begin{proof}
Grâce au théorème de désingularisation de de Jong, il existe un diagramme de la forme 
\ref{formel-II-diag} avec $X ^{(0)}$ lisse.
Les égalités résultent alors aussitôt du théorème \ref{Th2} et de la remarque \ref{rema-recol-factdagdag}.
La dernière assertion en découle grâce à la deuxième remarque de \ref{DDIsoc*=mod}.
\end{proof}

\section{Équivalence entre isocristaux partiellement surconvergent et surcohérents}

\subsection{Cas de la descente finie et étale en dehors des singularités surconvergentes}
En reprenant les notations et hypothèses de \ref{nota-diag}, supposons de plus $X^{(0)}$ lisse sur $k$.

\begin{prop}
\label{Isoc*=Isocdag-bullet}
Avec les notations de \ref{defi-Isocbullet*}, 
on dispose de l'équivalence canonique de catégories de la forme
$\sp _{X^{(\bullet)} \hookrightarrow \PP ^{(\bullet)},T ^{(\bullet)},+}\,:\,
\mathrm{Isoc} ^{\dag} (Y^{(\bullet)},X^{(\bullet)}/K)
\cong
\mathrm{Isoc} ^{*} (\PP^{(\bullet)}, T^{(\bullet)}, X^{(\bullet)}/K)$.
\end{prop}

\begin{proof}
Via le lemme \ref{rema-isocDagDagDense} et la remarque \ref{rema-isocDagDense}, en remarquant que l'adhérence d'une composante
irréductible de $Y ^{(0)}$ est encore lisse, 
on se ramène au cas où $Y$ (resp. $Y^{(0)}$) est intègre et dense dans $X$ (resp. $X^{(0)}$).

I) Construction du foncteur $\sp _{X^{(\bullet)} \hookrightarrow \PP ^{(\bullet)},T ^{(\bullet)},+}\,:\,
\mathrm{Isoc} ^{\dag} (Y^{(\bullet)},X^{(\bullet)}/K)
\to
\mathrm{Isoc} ^{*} (\PP^{(\bullet)}, T^{(\bullet)}, X^{(\bullet)}/K)$.

Soit $(E ^{(0)}, \epsilon) \in \mathrm{Isoc} ^{\dag} (Y^{(\bullet)},X^{(\bullet)}/K)$.
Comme $X^{(0)}$ est lisse, posons $\E ^{(0)} := \sp _{X^{(0)} \hookrightarrow \PP ^{(0)}, T ^{(0)},+} (E ^{(0)})$.
Posons $\FF _{1}:= a _{1} ^{+} (\E ^{(0)})$ et $\FF _{2}:= a _{2} ^{+} (\E ^{(0)})$.
D'après le théorème de désingularisation de de Jong (voir \cite{dejong} ou \cite{Ber-alterationdejong}), 
on dispose d'un morphisme surjectif, projectif, génériquement fini et étale
$\alpha\,:\, X'' \to X ^{(1)}$ 
tel que $X''$ soit lisse et $\alpha ^{-1} ( T ^{(1)} \cap X ^{(1)})$ soit un diviseur de $X ^{\prime \prime}$. Il existe ainsi un morphisme projectif et lisse de $\V$-schémas formels séparés et lisses de la forme $\PP'' \to \PP ^{(1)}$, une immersion fermée $X'' \hookrightarrow \PP''$ tels que leur composée donne le morphisme canonique $X'' \to \PP ^{(1)}$.

i) Construction de $\theta ''\,:\,\alpha ^{+} (\FF _{2}) \riso \alpha ^{+} (\FF _{1})$.

Pour $i=1,2$, comme $X''$ et $X ^{(0)}$ sont lisses, par commutativité des foncteurs de la forme 
$\sp _{X^{(0)} \hookrightarrow \PP ^{(0)}, T ^{(0)},+}$ aux images inverses et duaux (voir \cite{caro-construction}), 
on obtient l'isomorphisme canonique 
$( \alpha \circ a _{i} ) ^{+} \circ \E ^{(0)}  \riso 
\sp _{X '' \hookrightarrow \PP '', T '',+}  \circ  
\phantom{}{^{\vee}}\circ  ( \alpha \circ a _{i} ) ^{*}  \circ \phantom{}{^{\vee}} (E ^{(0)})$,
où $\vee$ désigne le foncteur dual dans les catégories respectives des isocristaux partiellement surconvergents.
Or, on dispose de l'isomorphisme 
de bidualité $( \alpha \circ a _{i} ) ^{*}  (E ^{(0)}) \riso \phantom{}{^{\vee}}\circ  ( \alpha \circ a _{i} ) ^{*}  \circ \phantom{}{^{\vee}} (E ^{(0)})$.
De plus, 
par transitivité de l'image inverse, 
$\alpha ^{+} (\FF _{i}) \riso 
( \alpha \circ a _{i} ) ^{+} \circ \E ^{(0)}$.
Ainsi, $\alpha ^{+} (\FF _{i}) \riso \sp _{X '' \hookrightarrow \PP '', T '',+}  \circ  ( \alpha \circ a _{i} ) ^{*}  (E ^{(0)})$.

L'isomorphisme structural de recollement 
$\epsilon \,:\,a _{2} ^{*}  (E ^{(0)}) \riso a _{1} ^{*}  (E ^{(0)})$
induit l'isomorphisme 
$\epsilon ''\,:\,( \alpha \circ a _{2} ) ^{*}  (E ^{(0)}) \riso ( \alpha \circ a _{1} ) ^{*}  (E ^{(0)})$.
Il en résulte l'isomorphisme canonique $\theta ''\,:\,\alpha ^{+} (\FF _{2}) \riso \alpha ^{+} (\FF _{1})$ s'inscrivant dans le diagramme commutatif
ci-dessous :
\begin{equation}
\label{Isoc*=Isocdag-bullet-diag1}
\xymatrix {
{\alpha ^{+} (\FF _{2}) } 
\ar@{=}[r] 
\ar@{.>}[d] ^-{\sim} _{\theta ''}
&
{\alpha ^{+} \circ a _{2} ^{+} (\E ^{(0)})}
\ar[r] ^-{\sim}
\ar@{.>}[d] ^-{\sim} _{\theta ''}
& 
{ \sp _{X '' \hookrightarrow \PP '', T '',+}  \circ  ( \alpha \circ a _{2} ) ^{*}  (E ^{(0)})} 
\ar[d] _-{\sim} ^{\sp _{X '' \hookrightarrow \PP '', T '',+}  (\epsilon'')}
\\ 
{\alpha ^{+} (\FF _{1})} 
\ar@{=}[r] 
&
{\alpha ^{+} \circ a _{1} ^{+} (\E ^{(0)})}
\ar[r] ^-{\sim}
& 
{ \sp _{X '' \hookrightarrow \PP '', T '',+}  \circ  ( \alpha \circ a _{1} ) ^{*}  (E ^{(0)})} 
}
\end{equation}
dont les isomorphismes horizontaux sont construits ci-dessus.

ii) Notons $\widehat{\epsilon} \,:\, b _{2} ^{*} (\widehat{E} ^{(0)} ) \riso b _{1} ^{*} (\widehat{E} ^{(0)} )$ l'isomorphisme de
$\mathrm{Isoc} ^{\dag} (Y^{(1)},Y^{(1)}/K)$ induit canoniquement par $\epsilon$.
On obtient l'isomorphisme canonique 
$\widehat{\theta}\,:\,\FF _{2} | \U ^{(1)} \riso \FF _{1} | \U ^{(1)}$ défini par commutativité du diagramme ci-dessous :
\begin{equation}
\label{Isoc*=Isocdag-bullet-diag2}
\xymatrix{
{\FF _{2} | \U ^{(1)} } 
\ar[r] ^-{\sim} 
\ar@{.>}[d] ^-{\sim} _{\widehat{\theta}}
&
{b _{2} ^{+} (\E ^{(0)} |\U ^{(0)})}
\ar[r] ^-{\sim} 
& 
{\sp _{Y ^{(1)} \hookrightarrow \U ^{(1)},+} \circ
b _{2} ^{*} (\widehat{E} ^{(0)} )} 
\ar[d] _-{\sim} ^{\sp _{Y ^{(1)} \hookrightarrow \U ^{(1)},+} (\widehat{\epsilon})}
\\ 
{\FF _{1} | \U ^{(1)} } 
\ar[r] ^-{\sim} 
&
{b _{1} ^{+} (\E ^{(0)} |\U ^{(0)})}
\ar[r] ^-{\sim} 
& 
{\sp _{Y ^{(1)} \hookrightarrow \U ^{(1)},+} \circ
b _{1} ^{*} (\widehat{E} ^{(0)} )} 
}
\end{equation}
dont les isomorphismes horizontaux sont construits de manière identique à ceux 
de \ref{Isoc*=Isocdag-bullet-diag1} (car les compactifications partielles $Y ^{(1)}$ et $Y ^{(0)}$ de $Y ^{(1)}$ et $Y ^{(0)}$ sont lisses).

iii)
Les morphismes $\theta''$, $\widehat{\theta}$ induisent canoniquement le même morphisme 
de la forme $\alpha ^{+} (\FF _{2})|\U ''\riso \alpha ^{+} (\FF _{1})|\U ''$. On obtient donc 
le morphisme 
$(\theta'',\widehat{\theta})$
de
$\mathrm{Isoc} ^{*} (\PP ^{\prime \prime}, T ^{\prime \prime}, X^{\prime \prime}K) 
\times _{\mathrm{Isoc} ^{*} (\U ^{\prime \prime}, Y ^{\prime \prime}/K)}
\mathrm{Isoc} ^{*} (\U ^{(1)}, Y ^{(1)}/K)$.
Or, d'après \ref{a+|U-plfid-prop}
le foncteur canonique  
\begin{equation}
\notag
(\alpha ^{+},\, |\U ^{(1)})\,:\,
\mathrm{Isoc} ^{*} (\PP ^{(1)}, T ^{(1)}, X^{(1)}K) 
\to 
\mathrm{Isoc} ^{*} (\PP ^{\prime \prime}, T ^{\prime \prime}, X^{\prime \prime}K) 
\times _{\mathrm{Isoc} ^{*} (\U ^{\prime \prime}, Y ^{\prime \prime}/K)}
\mathrm{Isoc} ^{*} (\U ^{(1)}, Y ^{(1)}/K)
\end{equation}
est pleinement fidèle. Il en résulte l'isomorphisme canonique 
$\theta\,:\,\FF _{2}\riso \FF _{1}$, i.e., $\theta\,:\,a _{2} ^{+} (\E ^{(0)}) \riso a _{1} ^{+} (\E ^{(0)})$.
Par fidélité du foncteur restriction $|\U ^{(2)}$, cet isomorphisme vérifie la condition de cocycle. 
On a donc construit le foncteur $\sp _{+}\,:\, \mathrm{Isoc} ^{\dag} (Y^{(\bullet)},X^{(\bullet)}/K)\to \mathrm{Isoc} ^{*} (\PP^{(\bullet)}, T^{(\bullet)}, X^{(\bullet)}/K)$
en posant 
$\sp _{X^{(\bullet)} \hookrightarrow \PP ^{(\bullet)},T ^{(\bullet)},+} (E ^{(0)}, \epsilon) 
:= (\sp _{X^{(0)} \hookrightarrow \PP ^{(0)}, T ^{(0)},+} (E ^{(0)}),\theta)$.

II) Le foncteur $\sp _{X^{(\bullet)} \hookrightarrow \PP ^{(\bullet)},T ^{(\bullet)},+}$ induit une équivalence de catégories.

 i) Comme $\sp _{X^{(0)} \hookrightarrow \PP ^{(0)}, T ^{(0)},+}$ est fidèle, il en est alors de même 
 de $\sp _{X^{(\bullet)} \hookrightarrow \PP ^{(\bullet)},T ^{(\bullet)},+} $.
Vérifions à présent que celle-ci est pleine. 

Soient 
$(E ^{(0)}, \epsilon),  (\widetilde{E} ^{(0)}, \widetilde{\epsilon})\in 
\mathrm{Isoc} ^{\dag} (Y^{(\bullet)},X^{(\bullet)}/K)$
et
$\psi \,:\, \sp _{X^{(\bullet)} \hookrightarrow \PP ^{(\bullet)},T ^{(\bullet)},+} (E ^{(0)}, \epsilon) \to
\sp _{X^{(\bullet)} \hookrightarrow \PP ^{(\bullet)},T ^{(\bullet)},+} (\widetilde{E} ^{(0)}, \widetilde{\epsilon})$.
Nous conservons les notations de l'étape I) de la preuve concernant les morphismes $\theta ''$, $\widehat{\theta}$, $\theta$, $\epsilon ''$, $\widehat{\epsilon}$ déduits à partir de $\epsilon$ ; de même en ajoutant des tildes.
Ainsi,
$(\E ^{(0)}, \theta):= \sp _{X^{(\bullet)} \hookrightarrow \PP ^{(\bullet)},T ^{(\bullet)},+} (E ^{(0)}, \epsilon)$,
$(\widetilde{\E} ^{(0)}, \widetilde{\theta}):= \sp _{X^{(\bullet)} \hookrightarrow \PP ^{(\bullet)},T ^{(\bullet)},+} (\widetilde{E} ^{(0)}, \widetilde{\epsilon})$.
Comme le foncteur $\sp _{X^{(0)} \hookrightarrow \PP ^{(0)}, T ^{(0)},+}$ est pleinement fidèle, 
il existe un morphisme $\phi\,:\, E ^{(0)}\to \widetilde{E} ^{(0)}$ qui induit le morphisme 
$\psi\,:\,\E ^{(0)} \to \widetilde{\E} ^{(0)}$.
Il reste à établir que $\phi$ commute aux données de recollement. 
Considérons le diagramme ci-après : 
\begin{equation}
\label{Isoc*=Isocdag-bullet-diag3}
\xymatrix @R=0,3cm{
&
{\alpha ^{+} \circ a _{2} ^{+} (\widetilde{\E} ^{(0)})}  
\ar'[d][dd] ^-{\widetilde{\theta} ''}
\ar[rr] 
&&
{\sp _{X '' \hookrightarrow \PP '', T '',+}  \circ  ( \alpha \circ a _{2} ) ^{*}  (\widetilde{E} ^{(0)})}
\ar[dd] ^-{\sp _{X '' \hookrightarrow \PP '', T '',+}  (\widetilde{\epsilon}'')}
\\ 
{\alpha ^{+} \circ a _{2} ^{+} (\E ^{(0)})} 
\ar[rr] 
\ar[ru] ^-{\psi}
\ar[dd] _(0.6){\theta ''} 
&&
{\sp _{X '' \hookrightarrow \PP '', T '',+}  \circ  ( \alpha \circ a _{2} ) ^{*}  (E ^{(0)})}
\ar[ru] ^-{\phi}
\ar[dd] ^(0.6){\sp _{X '' \hookrightarrow \PP '', T '',+}  (\epsilon'')}
\\
&
{\alpha ^{+} \circ a _{1} ^{+} (\widetilde{\E}^{(0)})} 
\ar'[r][rr] 
&&
{\sp _{X '' \hookrightarrow \PP '', T '',+}  \circ  ( \alpha \circ a _{1} ) ^{*}  (\widetilde{E} ^{(0)})}
\\
{\alpha ^{+} \circ a _{1} ^{+} (\E ^{(0)})} 
\ar[ru] _-{\psi}
\ar[rr] 
&&
{\sp _{X '' \hookrightarrow \PP '', T '',+}  \circ  ( \alpha \circ a _{1} ) ^{*}  (E ^{(0)}),}
\ar[ru] _-{\phi}
}
\end{equation}
où les flèches vers l'arrière sont celles canoniquement induites par $\phi$ ou $\psi$.
Comme $\psi$ commute aux données de recollement, le carré de gauche est commutatif. 
Les carrés horizontaux le sont par fonctorialité, tandis que ceux de devant et de derrière le sont 
par définition (voir \ref{Isoc*=Isocdag-bullet-diag1}).
Comme les cinq autres carrés du cube \ref{Isoc*=Isocdag-bullet-diag3} sont commutatifs, il en est de même du carré de droite. 
Comme le foncteur $\sp _{X '' \hookrightarrow \PP '', T '',+} $ est fidèle, ce carré reste commutatif sans 
$\sp _{X '' \hookrightarrow \PP '', T '',+} $.

De même, en construisant le cube déduit de \ref{Isoc*=Isocdag-bullet-diag2} par fonctorialité, 
on vérifie que le morphisme de $\mathrm{Isoc} ^{\dag} (Y^{(0)},Y^{(0)}/K)$ induit par $\phi$ commute aux données de
recollement. Par fidélité du foncteur 
$( \alpha ^{*}, j ^{(0)*})$ (voir le théorème \ref{plfid-dag}), il en résulte que $\phi$ commute bien aux données de recollement.

ii) Comme le foncteur $\sp _{X^{(0)} \hookrightarrow \PP ^{(0)}, T ^{(0)},+}$ est essentiellement surjectif, on vérifie de façon analogue
qu'il en est de même de $\sp _{X^{(\bullet)} \hookrightarrow \PP ^{(\bullet)},T ^{(\bullet)},+}$. 

\end{proof}

\begin{coro}
\label{Isoc*=dagdagFét}
Avec les notations et hypothèses de la section, on dispose d'une équivalence canonique de 
catégories 
$\sp _{X \hookrightarrow \PP, T,+}\,:\, 
\mathrm{Isoc} ^{\dag} (Y,X/K)
\cong 
\mathrm{Isoc} ^{*} (\PP, T, X/K)$.
\end{coro}

\begin{proof}
En utilisant le lemme \ref{rema-isocDagDagDense} et la remarque \ref{rema-isocDagDense}, 
on se ramène au cas où $Y$ (resp. $Y^{(0)}$) est intègre et dense dans $X$ (resp. $X^{(0)}$).

I) Construction de $\sp _{X \hookrightarrow \PP, T,+}$.

Notons $\mathcal{L}oc\,:\,\mathrm{Isoc} ^{\dag} (Y,X/K)\to \mathrm{Isoc} ^{\dag} (Y^{(\bullet)},X^{(\bullet)}/K)$ 
le foncteur canonique (construit de manière analogue à \ref{DefLoc}). 
D'après le théorème de descente de Shiho (voir \cite[7.3]{shiho-logRC-RCII}),
ce foncteur $\mathcal{L}oc$ est une équivalence de catégories. 
On obtient le foncteur $\sp _{X \hookrightarrow \PP, T,+}$ en posant
$\sp _{X \hookrightarrow \PP, T,+}:= \mathcal{R}ecol \circ \sp _{X^{(\bullet)} \hookrightarrow \PP ^{(\bullet)},T ^{(\bullet)},+}
\circ \mathcal{L}oc$. 
De plus, d'après \ref{Isoc*=Isocdag-bullet} (resp. \ref{desc-fini-ét})
le foncteur $\sp _{X^{(\bullet)} \hookrightarrow \PP ^{(\bullet)},T ^{(\bullet)},+}$ (resp. $\mathcal{R}ecol$)
est aussi une équivalence de catégories.
Le foncteur $\sp _{X \hookrightarrow \PP, T,+}$ est donc une équivalence de catégories.

II) Cette équivalence est canonique, i.e., est indépendant du morphisme $(\PP ^{(0)}, T ^{(0)}, X ^{(0)}, Y ^{(0)}) \to (\PP, T,X,Y)$  
de quadruplets lisses en dehors du diviseur
(voir la définition \ref{defi-triplet}) tel que $X ^{(0)}$ soit intègre, lisse et dont les conditions correspondantes du paragraphe 
\ref{nota-diag-pre} soient satisfaites.

En effet, soit $(f',a',b')\,:\,(\PP ^{\prime(0)}, T ^{\prime(0)}, X ^{\prime(0)}, Y ^{\prime(0)}) \to (\PP, T,X,Y)$  
un second morphisme de quadruplets lisses en dehors du diviseur
tel que $X ^{\prime(0)}$ soit intègre et lisse, $Y ^{\prime(0)}$ soit dense dans $X ^{\prime(0)}$ et dont les conditions du paragraphe 
\ref{nota-diag-pre} soient satisfaites.
Concernant ce second choix, on reprend les notations analogues à \ref{nota-diag} mais en ajoutant des primes.
On note $\sp _{X^{\prime (\bullet)} \hookrightarrow \PP ^{\prime(\bullet)},T ^{\prime (\bullet)},+}$ le foncteur 
défini via \ref{Isoc*=Isocdag-bullet} en remplaçant les données de descente sur $X ^{(0)}$ par celle sur $X ^{\prime (0)}$.
De plus, afin de faire la distinction entre les deux choix, notons à titre provisoire
$\mathcal{L}oc ^{\prime}\,:\,\mathrm{Isoc} ^{*} (\PP, T, X/K)
\to 
\mathrm{Isoc} ^{*} (\PP^{\prime (\bullet)}, T^{\prime (\bullet)}, X^{\prime(\bullet)}/K)$
et 
$\mathcal{L}oc ^{\prime}\,:\,\mathrm{Isoc} ^{\dag} (Y,X/K)\to \mathrm{Isoc} ^{\dag} (Y^{\prime (\bullet)},X^{\prime (\bullet)}/K)$
les foncteurs habituellement notés $\mathcal{L}oc$.
De même, on notera $\mathcal{R}ecol ^{\prime}\,:\, 
\mathrm{Isoc} ^{*} (\PP^{\prime (\bullet)}, T^{\prime (\bullet)}, X^{\prime(\bullet)}/K)
\to 
\mathrm{Isoc} ^{*} (\PP, T, X/K)$ (afin de le différencier du 
premier foncteur $\mathcal{R}ecol $).

En posant
$\sp ^{\prime}_{X \hookrightarrow \PP, T,+}:= \mathcal{R}ecol ^{\prime}\circ \sp _{X^{\prime (\bullet)} \hookrightarrow \PP ^{\prime(\bullet)},T ^{\prime (\bullet)},+}
\circ \mathcal{L}oc ^{\prime}$,
il s'agit ainsi d'établir que les deux foncteurs
$\sp _{X \hookrightarrow \PP, T,+}$
et $\sp ^{\prime}_{X \hookrightarrow \PP, T,+}$
sont canoniquement isomorphes. 

i) On se ramène au cas où $(f',a',b')\,:\,(\PP ^{\prime(0)}, T ^{\prime(0)}, X ^{\prime(0)}, Y ^{\prime(0)}) 
\to (\PP, T,X,Y)$  
se factorise par 
$(\PP ^{\prime(0)}, T ^{\prime(0)}, X ^{\prime(0)}, Y ^{\prime(0)}) 
\to
(\PP ^{(0)}, T ^{(0)}, X ^{(0)}, Y ^{(0)}) $.

Avec le lemme \ref{rema-isocDagDagDense} (et la remarque \ref{rema-isocDagDense}), quitte à le remplacer par une composante irréductible, on peut supposer $X ^{\prime (0)} \times _{X} X ^{(0)}$ intègre. 
D'après le théorème de désingularisation de de Jong, il existe un morphisme projectif, surjectif, génériquement fini et étale $X ^{\prime \prime (0)} \to X ^{\prime (0)} \times _{X} X ^{(0)}$ avec $X ^{\prime \prime (0)}$ lisse (comme ce morphisme est projectif, il se prolonge en un morphisme
de $\V$-schémas formels séparés et lisses de la forme $\PP ^{\prime \prime (0)} \to \PP ^{\prime (0)} \times _{\PP} \PP ^{(0)}$). 
Il existe ainsi un diviseur $\widetilde{T}$ contenant $T$ tel que 
$\widetilde{Y}:= X \setminus \widetilde{T}$ soit dense dans $X$ et tel que 
le morphisme $X ^{\prime \prime (0)} \to X$ soit fini et étale en dehors de $\widetilde{T}$. 

Or, d'après \ref{Loc-Recol-hdag}
le foncteur $\mathcal{R}ecol$ de \ref{def-recol}
commute aux extensions (i.e. aux foncteurs de la forme $(\hdag \widetilde{T})$).
De plus, il en est de même des foncteurs de la forme 
$\mathcal{L}oc\,:\,\mathrm{Isoc} ^{\dag} (Y,X/K)\to \mathrm{Isoc} ^{\dag} (Y^{(\bullet)},X^{(\bullet)}/K)$ 
ou de la forme 
$\sp _{X^{(\bullet)} \hookrightarrow \PP ^{(\bullet)},T ^{(\bullet)},+}$.
Or, d'après \ref{dagT'pl-fid}, les foncteurs extensions $(\hdag \widetilde{T})$ sont pleinement fidèles (lorsque 
$ X \setminus \widetilde{T} $ est dense dans $Y$). 
On se ramène donc à supposer le morphisme $X ^{\prime \prime (0)} \to X$ fini et étale en dehors
de $T$. Comme ce morphisme se factorise par $X ^{\prime \prime (0)} \to X ^{\prime(0)}$
et par $X ^{\prime \prime (0)} \to X ^{(0)}$ (et de même au niveau des schémas formels), on a donc vérifié l'étape $i)$.

ii) Fin de la preuve. 

a) Notons $\alpha ^{(0)}\,:\, X ^{\prime (0)}\to X ^{(0)}$.
On remarque que le morphisme $\alpha ^{(0)}$ induit la factorisation $\alpha ^{(1)}\,:\, X ^{\prime (1)}\to X ^{(1)}$
telle que, pour $i:=1,2$, $\alpha ^{(0)}\circ a ' _{i}= a _{i} \circ \alpha ^{(1)}$.
Soit 
$(E ^{(0)},\epsilon)\in \mathrm{Isoc} ^{\dag} (Y^{(\bullet)},X^{(\bullet)}/K)$.
L'isomorphisme de recollement $\epsilon$ induit alors canoniquement un isomorphisme de recollement sur $\alpha ^{(0)*}(E ^{(0)})$ que l'on note $\epsilon '$. 
On obtient le foncteur 
$\mathcal{L}oc ^{(0)}\,:\,\mathrm{Isoc} ^{\dag} (Y^{(\bullet)},X^{(\bullet)}/K)\to \mathrm{Isoc} ^{\dag} (Y^{\prime (\bullet)},X^{\prime (\bullet)}/K)$
défini en posant 
$\mathcal{L}oc ^{(0)} (E ^{(0)},\epsilon) := (\alpha ^{(0)*}(E ^{(0)}),\epsilon')$.
On dispose de plus de l'isomorphisme canonique 
$\mathcal{L}oc ^{(0)}  \circ \mathcal{L}oc  \riso \mathcal{L}oc ^{\prime}$.
De même, on construit le foncteur
$\mathcal{L}oc ^{(0)}\,:\,\mathrm{Isoc} ^{*} (\PP^{(\bullet)}, T^{(\bullet)}, X^{(\bullet)}/K)
\to 
\mathrm{Isoc} ^{*} (\PP^{\prime (\bullet)}, T^{\prime (\bullet)}, X^{\prime(\bullet)}/K)$
ainsi que de l'isomorphisme canonique
$\mathcal{L}oc ^{(0)}  \circ \mathcal{L}oc  \riso \mathcal{L}oc^{\prime}$.

b) Soient 
$E \in \mathrm{Isoc} ^{\dag} (Y,X/K)$ et 
$(E ^{(0)},\epsilon) := \mathcal{L}oc (E)\in \mathrm{Isoc} ^{\dag} (Y^{(\bullet)},X^{(\bullet)}/K)$.
Comme $\alpha ^{(0)}$ est un morphisme de $k$-variétés lisses,
on dispose de l'isomorphisme canonique
$$\alpha ^{(0)+} \circ \sp _{X^{(0)} \hookrightarrow \PP ^{(0)},T ^{(0)},+} (E ^{(0)})
\riso \sp _{X^{\prime (0)} \hookrightarrow \PP ^{\prime(0)},T ^{\prime (0)},+} \circ \alpha ^{(0)*} (E ^{(0)}).$$
Pour vérifier que cet isomorphisme commute aux isomorphismes de recollement induits par $\epsilon$, par fidélité du foncteur
$|\U ^{\prime(1)}$, on se ramène au cas où le diviseur $T$ est vide (i.e. $Y=X$), ce qui est immédiat.
Ainsi, on obtient l'isomorphisme  
\begin{equation}
\label{Isoc*=dagdagFét-iso1}
\mathcal{L}oc ^{(0)} \circ \sp _{X^{(\bullet)} \hookrightarrow \PP ^{(\bullet)},T ^{(\bullet)},+}( \mathcal{L}oc (E))
\riso \sp _{X^{\prime (\bullet)} \hookrightarrow \PP ^{\prime(\bullet)},T ^{\prime (\bullet)},+} \circ \mathcal{L}oc ^{(0)}  
(\mathcal{L}oc (E)).
\end{equation}
Pour vérifier que $\sp _{X \hookrightarrow \PP, T,+} (E)$
et $\sp ^{\prime}_{X \hookrightarrow \PP, T,+}(E)$
sont canoniquement isomorphes, il suffit d'établir l'isomorphisme canonique fonctoriel
$\mathcal{L}oc ^{\prime} \circ \sp _{X \hookrightarrow \PP, T,+}(E)
\riso
\mathcal{L}oc ^{\prime} \circ \sp ^{\prime}_{X \hookrightarrow \PP, T,+}(E)$.
Or, 
comme
$\mathcal{L}oc ^{(0)}  \circ \mathcal{L}oc \riso \mathcal{L}oc ^{\prime} $
et 
$\mathcal{L}oc\circ \mathcal{R}ecol \riso Id  $ (voir la preuve de \ref{desc-fini-ét}),
le terme
$\mathcal{L}oc ^{\prime} \circ \sp _{X \hookrightarrow \PP, T,+}(E )$ est canoniquement isomorphe au terme de gauche
de \ref{Isoc*=dagdagFét-iso1}.
De même, $\mathcal{L}oc ^{\prime} \circ \sp ^{\prime}_{X \hookrightarrow \PP, T,+}(E)$
est canoniquement isomorphe à celui de droite.

\end{proof}

\begin{rema}
\label{cas-lisse-sp+egaux}
Avec les notations de \ref{Isoc*=dagdagFét}, lorsque $X$ est lisse, le foncteur
$\sp _{X \hookrightarrow \PP, T,+}$ est par construction (e.g. voir l'étape II) de la preuve de \ref{Isoc*=dagdagFét})
canoniquement isomorphe à celui construit dans le cas lisse (voir \cite{caro-construction}).
\end{rema}

\begin{vide}
\label{sp-comm-a+-fini-etale}
Avec les notations de \ref{Isoc*=dagdagFét}, pour tout 
$E \in \mathrm{Isoc} ^{\dag} (Y,X/K)$, 
par construction de $\sp _{X \hookrightarrow \PP, T,+}$ (voir l'étape I) de la preuve de \ref{Isoc*=dagdagFét})
on bénéficie de l'isomorphisme canonique :
\begin{equation}
\label{sp-comm-a+-fini-etale-iso} 
a ^{+}\circ \sp _{X \hookrightarrow \PP, T,+}(E) 
\riso 
\sp _{X ^{(0)}\hookrightarrow \PP^{(0)}, T^{(0)},+}\circ a ^{*}(E) .
\end{equation}

\end{vide}

\begin{vide}
\label{sp+-comm-hdag}
Soit $T '$ un diviseur contenant $T$. Notons $Y' := X \setminus T'$, 
$Y ^{\prime (0)} := X ^{(0)} \setminus T^{\prime (0)}$, 
$j ' \,:\, Y' \hookrightarrow X$, $j ^{\prime (0)} \,:\, Y^{\prime (0)} \hookrightarrow X^{(0)}$ les immersions ouvertes induites etc (i.e. on rajoute des primes).
Considérons le diagramme 
\begin{equation}
\label{sp+-comm-hdag-diag}
\xymatrix{
{\mathrm{Isoc} ^{\dag} (Y,X/K)} 
\ar[r] ^-{\mathcal{L}oc}
\ar[d] ^-{j ^{\prime \dag}}
& 
{\mathrm{Isoc} ^{\dag} (Y^{(\bullet)},X^{(\bullet)}/K)} 
\ar[rr] ^-{\small \sp _{X^{(\bullet)} \hookrightarrow \PP ^{(\bullet)},T ^{(\bullet)},+}}
\ar[d] ^-{j ^{\prime (0)\dag}}
&&
{\mathrm{Isoc} ^{*} (\PP^{(\bullet)}, T^{(\bullet)}, X^{(\bullet)}/K)}
\ar[r] ^-{\mathcal{R}ecol}
\ar[d] ^-{(\hdag T ^{\prime(0)})}
&
{\mathrm{Isoc} ^{*} (\PP, T, X/K)}
\ar[d] ^-{(\hdag T ^{\prime})}
\\ 
{\mathrm{Isoc} ^{\dag} (Y',X/K)} 
\ar[r] ^-{\mathcal{L}oc}
& 
{\mathrm{Isoc} ^{\dag} (Y^{\prime (\bullet)},X^{(\bullet)}/K)} 
\ar[rr] ^-{\small \sp _{X^{(\bullet)} \hookrightarrow \PP ^{(\bullet)},T ^{\prime(\bullet)},+}}
&&
{\mathrm{Isoc} ^{*} (\PP^{(\bullet)}, T^{\prime (\bullet)}, X^{(\bullet)}/K)}
\ar[r] ^-{\mathcal{R}ecol}
&
{\mathrm{Isoc} ^{*} (\PP, T', X/K)}
}
\end{equation}
où $ (\hdag T ^{\prime(0)})$ (resp. $(\hdag T ^{\prime(0)})$) désigne le foncteur canonique induit par
$ (\hdag T ^{\prime(0)})$ (resp. $(\hdag T ^{\prime(0)})$).
Par commutation canonique des foncteurs de la forme $\sp _{X^{\prime(0)} \hookrightarrow \PP^{\prime(0)}, T^{\prime(0)},+} $ (i.e., les foncteurs construits dans \cite{caro-construction} lorsque la compactification partielle est lisse)
aux foncteurs d'extension du diviseur
(i.e. de la forme $j ^{\prime \dag}$ ou $(\hdag T ^{\prime})$), on vérifie que le carré central du diagramme \ref{sp+-comm-hdag-diag}
est commutatif à isomorphisme canonique près.
Comme il en est de même des deux autres carrés (immédiat pour le carré de gauche et voir \ref{Loc-Recol-hdag} pour celui de droite),
on bénéficie ainsi de l'isomorphisme canonique : 
\begin{equation}
\label{sp+-comm-hdag-iso}
(\hdag T ^{\prime})\circ \sp _{X \hookrightarrow \PP, T,+}
\riso 
\sp _{X \hookrightarrow \PP, T',+}\circ j ^{\prime \dag}.
\end{equation}

\end{vide}

\subsection{Cas général}

\begin{lemm}
\label{lemm-diag-eqcat-dagdag-ii}
Avec les notations de la section \ref{formel-II}, 
le foncteur $\sp _{X \hookrightarrow \PP, \widetilde{T},+}$
construit en \ref{Isoc*=dagdagFét}
induit l'équivalence de catégories :
$\sp _{X \hookrightarrow \PP, \widetilde{T},+}
\,:\,
\mathrm{Isoc} ^{\dag} (\widetilde{Y}\supset Y,X/K) 
\cong 
\mathrm{Isoc} ^{*} (\PP, \widetilde{T}\supset T, X/K)$.
\end{lemm}

\begin{proof}
Considérons le diagramme ci-dessous 
\begin{equation}
\label{eq1-lemm-diag-eqcat-dagdag-ii}
\xymatrix{
{\mathrm{Isoc} ^{\dag} (\widetilde{Y}\supset Y,X/K) }
\ar@{.>}[d] ^-{\sp _{X \hookrightarrow \PP, \widetilde{T},+}} 
\ar@{^{(}->}[r] ^-{}
& 
{\mathrm{Isoc} ^{\dag} (\widetilde{Y},X/K) }
\ar[d] ^-{\sp _{X \hookrightarrow \PP, \widetilde{T},+}} 
\ar[r] ^-{}
& 
{\mathrm{Isoc} ^{\dag} (\widetilde{Y},Y/K)}
\ar[d] ^-{\sp _{Y \hookrightarrow \U, \widetilde{T} \cap U,+}} 
\\
{\mathrm{Isoc} ^{*} (\PP, \widetilde{T}\supset T, X/K)} 
\ar@{^{(}->}[r] ^-{}
& 
{\mathrm{Isoc} ^{*} (\PP, \widetilde{T}, X/K)} 
\ar[r] ^-{}
& 
{\mathrm{Isoc} ^{*} (\U, \widetilde{T}\cap U, Y/K)}  
}
\end{equation}
dont la factorisation de gauche est immédiate par définition et dont le foncteur de droite est bien défini car $Y$ est lisse (voir la construction
de \cite{caro-construction}). 
La pleine fidélité de cette factorisation est immédiate (compte tenu de \ref{Isoc*=dagdagFét}).
Il reste à établir son essentielle surjectivité. 
Soit $\E \in \mathrm{Isoc} ^{*} (\PP, \widetilde{T}\supset T, X/K)$.
Comme le foncteur $\sp _{X \hookrightarrow \PP, \widetilde{T},+}$
au milieu du diagramme \ref{eq1-lemm-diag-eqcat-dagdag-ii} est une équivalence de catégories, 
il existe un isocristal $E $ de $\mathrm{Isoc} ^{\dag} (\widetilde{Y},X/K)$ tel que 
$\E \riso \sp _{X \hookrightarrow \PP, \widetilde{T},+} (E)$.
Notons $E _Y$ l'isocristal de $\mathrm{Isoc} ^{\dag} (\widetilde{Y},Y/K)$ induit par $E$.
Il découle de la remarque \ref{cas-lisse-sp+egaux} 
que le carré de droite du diagramme \ref{eq1-lemm-diag-eqcat-dagdag-ii}
est commutatif (à isomorphisme canonique près).
Ainsi $\E |\U \riso \sp _{Y \hookrightarrow \U, \widetilde{T} \cap U,+} (E _{Y})$.
Comme par hypothèse  $\E |\U$ provient par extension d'un élément de $\mathrm{Isoc} ^{*} (\U, Y/K)$, 
alors $E _Y$ provient par extension d'un élément de $\mathrm{Isoc} ^{\dag} (Y,Y/K)$
(on utilise pour cela la commutativité de type celle du diagramme commutatif de gauche de \ref{Xlisse-Eqsp+-diag}).
D'après le théorème de contagiosité de Kedlaya \cite[5.3.7]{kedlaya-semistableI}, 
il en résulte que $E$ provient d'un isocristal de $\mathrm{Isoc} ^{\dag} (Y,X/K)$, i.e.
$E \in \mathrm{Isoc} ^{\dag} (\widetilde{Y}\supset Y,X/K) $.
\end{proof}

\begin{theo}
\label{diag-eqcat-dagdag-ii}
Soit $(\PP, T,X,Y)$ un quadruplet lisse en dehors du diviseur (voir la définition \ref{defi-triplet}).
On dispose d'une équivalence de catégories canonique
notée $\sp _{X \hookrightarrow \PP, T,+}\,:\,\mathrm{Isoc} ^{\dag} (Y,X/K)\cong
\mathrm{Isoc} ^{*} (\PP, T, X/K) $.
\end{theo}

\begin{proof}
Avec le lemme \ref{rema-isocDagDagDense} (et la remarque \ref{rema-isocDagDense}), il ne coûte pas cher de supposer 
 $Y$ intègre et dense dans $X$.
Grâce au théorème de désingularisation de de Jong, il existe un diviseur $\widetilde{T}$ contenant $T$
et un diagramme de la forme \ref{formel-II-diag} satisfaisant aux conditions requises de \ref{formel-II} et
tel que de plus $X ^{(0)}$ soit lisse. 

Avec les notations de \ref{formel-II}, on définit le foncteur $\sp _{X \hookrightarrow \PP, T,+}$ comme celui rendant commutatif le diagramme ci-dessous :
\begin{equation}
\label{Th2-cas lisse-diag}
\xymatrix{
{\mathrm{Isoc} ^{\dag} (Y,X/K)} 
\ar[r] ^-{(a ^{*}, \widetilde{j} ^{\dag})}
\ar@{.>}[d] ^-{\sp _{X \hookrightarrow \PP, T,+}}
& 
{\mathrm{Isoc} ^{\dag} (Y^{(0)},X^{(0)}/K) 
\times _{\mathrm{Isoc} ^{\dag} (\widetilde{Y}^{(0)},X^{(0)}/K)} \mathrm{Isoc} ^{\dag} (\widetilde{Y}\supset Y,X/K) } 
\ar[d] ^-{\sp _{X ^{(0)} \hookrightarrow \PP^{(0)}, T^{(0)},+} \times \sp _{X \hookrightarrow \PP, \widetilde{T},+}}
\\ 
{\mathrm{Isoc} ^{*} (\PP, T, X/K)} 
& 
{\mathrm{Isoc} ^{*} (\PP ^{(0)}, T ^{(0)}, X^{(0)}/K) 
\times _{\mathrm{Isoc} ^{*} (\PP ^{(0)}, \widetilde{T} ^{(0)}, X^{(0)}/K)}
\mathrm{Isoc} ^{*} (\PP, \widetilde{T}\supset T, X/K),}
\ar[l] ^-{\mathcal{R}ecol} 
}
\end{equation}
dont le foncteur du bas a été construit en \ref{recol-th2}.
D'après \ref{EqCat-rig}
le foncteur canonique 
$(a ^{*}, \widetilde{j} ^{\dag})$
est une équivalence de catégories. 
De plus, 
d'après \ref{lemm-diag-eqcat-dagdag-ii},
il en est de même du foncteur $\sp _{X \hookrightarrow \PP, \widetilde{T},+}$.
Comme $X ^{(0)} $ est lisse, c'est encore le cas du foncteur
$\sp _{X ^{(0)} \hookrightarrow \PP^{(0)}, T^{(0)},+}$. 
Le foncteur de droite est donc une équivalence de catégories. 
Le théorème \ref{Th2} nous permet de conclure qu'il en est de même de $\sp _{X \hookrightarrow \PP, T,+}$.

Il reste à établir que ce foncteur ne dépend du choix d'un diagramme de la forme \ref{formel-II-diag}
et satisfaisant aux conditions requises. 
Traitons d'abord l'indépendance en $\widetilde{T}$. Soit $\widetilde{T}'$ un second diviseur tel que 
$a ^{-1} (X \setminus \widetilde{T}') \to 
(X \setminus \widetilde{T}')$ soit fini et étale. Quitte à considérer $\widetilde{T} \cup \widetilde{T}'$, on peut supposer
$\widetilde{T}' \supset \widetilde{T}$. Notons $\widetilde{j} ^\prime \,:\, X \setminus \widetilde{T}' \subset X$ l'inclusion canonique.
D'après \ref{sp+-comm-hdag-iso}, on dispose de l'isomorphisme canoniques 
$\sp _{X \hookrightarrow \PP, \widetilde{T}',+} \circ \widetilde{j} ^{\prime \dag} 
\riso (\hdag \widetilde{T}') \circ \sp _{X \hookrightarrow \PP, \widetilde{T},+}$.
Or, d'après \ref{Th2} (voir sa preuve), 
on bénéficie des isomorphismes canoniques $\mathcal{L}oc  \circ \mathcal{R}ecol\riso Id$, 
avec soit $\mathcal{L}oc = (a ^{+} ,(\hdag  \widetilde{T}))$ ou 
$\mathcal{L}oc = (a ^{ +} ,(\hdag  \widetilde{T} ^{\prime}))$ et le foncteur quasi-inverse $\mathcal{R}ecol$ correspondant.
On en déduit l'isomorphisme canonique $\mathcal{R}ecol \circ (Id \times (\hdag  \widetilde{T}') )\riso \mathcal{R}ecol $.
Il en résulte l'indépendance par rapport à $\widetilde{T}$.

Soient $(\widetilde{Y}^{\prime(0)}$, $Y^{\prime(0)}$, $X^{\prime(0)} $, $\PP ^{\prime(0)})$ 
quadruplet s'inscrivant dans un diagramme de la forme
\ref{formel-II-diag} et satisfaisant aux conditions qui y sont requises (pour les notations, on remplace {\og $(0)$\fg}
par {\og $\prime (0)$\fg}).
De manière analogue à l'étape II.i) de la preuve de \ref{Isoc*=dagdagFét}, 
quitte à désingulariser $X ^{\prime (0)} \times _{X} X ^{(0)}$ à la de Jong et quitte à agrandir $\widetilde{T}$, 
on se ramène à supposer que 
l'on dispose des morphismes de quadruplets lisses en dehors du diviseur
$(\PP ^{\prime(0)}, T ^{\prime(0)}, X ^{\prime(0)}, Y ^{\prime(0)}) 
\to
(\PP ^{(0)}, T ^{(0)}, X ^{(0)}, Y ^{(0)}) $
et
$(\PP ^{\prime(0)}, \widetilde{T} ^{\prime(0)}, X ^{\prime(0)}, \widetilde{Y} ^{\prime(0)}) 
\to
(\PP ^{(0)}, \widetilde{T} ^{(0)}, X ^{(0)}, \widetilde{Y} ^{(0)}) $
qui factorisent 
le diagramme \ref{formel-II-diag} où on remplace {\og $(0)$\fg} par {\og $\prime (0)$\fg}.
Notons $\alpha \,:\,X^{\prime(0)} \to X^{(0)}$ le morphisme de la factorisation
et $a'= a \circ \alpha \,:\, X^{\prime(0)} \to X$. 
On dispose de l'isomorphisme canonique 
$\alpha ^{+} \circ \sp _{X ^{(0)} \hookrightarrow \PP^{(0)}, T^{(0)},+}
\riso 
\sp _{X ^{\prime(0)} \hookrightarrow \PP^{\prime(0)}, T^{\prime (0)},+} \circ \alpha ^{*}$ 
(voir \cite{caro-construction}).
Comme $\mathcal{L}oc  \circ \mathcal{R}ecol\riso Id$ avec $\mathcal{L}oc = (a ^{+} ,(\hdag  \widetilde{T}))$ ou 
$\mathcal{L}oc = (a ^{\prime +} ,(\hdag  \widetilde{T}))$, on vérifie 
l'isomorphisme canonique $\mathcal{R}ecol \circ (\alpha ^{+} \times Id)\riso \mathcal{R}ecol $.
D'où le résultat.
\end{proof}

\begin{rema}
Afin d'établir \ref{diag-eqcat-dagdag-ii}, une idée aurait été de
construire un foncteur quasi-inverse $\mathcal{R}ecol$ {\og à la \ref{recol-th2}\fg} du foncteur du bas du diagramme :
\begin{equation}
\label{Rem-Th2-cas lisse-diag}
\xymatrix{
{\mathrm{Isoc} ^{\dag} (Y,X/K)} 
\ar[r] ^-{\widetilde{j} ^{\dag}} _-{\cong}
& 
{\mathrm{Isoc} ^{\dag} (\widetilde{Y}\supset Y,X/K) } 
\ar[d] ^-{ \sp _{X \hookrightarrow \PP, \widetilde{T},+}} _-{\cong}
\\ 
{\mathrm{Isoc} ^{*} (\PP, T, X/K)} 
\ar[r] ^-{(\hdag \widetilde{T})}
& 
{\mathrm{Isoc} ^{*} (\PP, \widetilde{T}\supset T, X/K).}
}
\end{equation}
Le problème est qu'il n'est pas immédiat que les objets de $\mathrm{Isoc} ^{*} (\PP, \widetilde{T}\supset T, X/K)$ sont 
$\D ^\dag _{\PP} (\hdag T) _\Q$-surcohérents (ainsi que leur duaux). 
Ce problème technique est seulement résolue {\og par descente\fg} via le lemme \ref{lemm-surcoh-ann-th2}, ce qui explique 
pourquoi nous avons recours au morphisme $a$.
\end{rema}

Avant d'établir via la proposition \ref{sp+-comm-dual} que cette équivalence de catégories commutent aux foncteurs duaux,
vérifions d'abord sa commutation aux images inverses via la proposition \ref{theta*com-sp} : 
\begin{prop}
\label{theta*com-sp}
Soit $\theta  \,:\, (\PP', T',X',Y')\to (\PP, T,X,Y)$ un morphisme de quadruplets lisses en dehors du diviseur.
On dispose de l'isomorphisme canonique : 
\begin{equation}
\label{iso-theta*com-sp}
\sp _{X '\hookrightarrow \PP', T',+} \circ \theta ^{*} 
\riso 
\theta ^{+}\circ  \sp _{X \hookrightarrow \PP, T,+} .
\end{equation}
\end{prop}

\begin{proof}
Posons $\theta = (\phi, \alpha, \beta)$.
Avec le lemme \ref{rema-isocDagDagDense} (et la remarque \ref{rema-isocDagDense}), on se ramène au cas où $Y$ (resp. $Y'$) est intègre et dense dans $X$ (resp. $X'$).
Lorsque $X$ et $X'$ sont lisses, ce théorème de commutation est déjà connu (voir \cite{caro-construction}).
Pour en déduire le cas général, l'idée est de se ramener au cas lisse grâce au théorème de pleine fidélité \ref{a+|U-plfid-prop} de la manière suivante : 
grâce au théorème de désingularisation de de Jong, 
il existe un morphisme de quadruplets lisses en dehors du diviseur
de la forme $\eta= (f , a, b )  \,:\, (\PP ^{(0)}, T^{(0)},X^{(0)},Y^{(0)})\to (\PP, T,X,Y)$
tel que $X^{(0)}$ soit lisse, $f$ soit projectif et lisse, $a$ soit projectif, surjectif, génériquement fini et étale, $ f ^{-1} ( T) = T ^{(0)}$
et 
$T ^{(0)} \cap X ^{(0)}$ soit un diviseur (à croisement normaux) de $X ^{(0)}$. 
De même en utilisant le théorème de désingularisation de de Jong (appliqué à $X ' \times _{X} X ^{(0)}$ élément de 
$(\PP', T',X',Y')\times _{(\PP, T,X,Y)} (\PP ^{(0)}, T^{(0)},X^{(0)},Y^{(0)})$ et à son diviseur complémentaire de $Y
 ' \times _{Y} Y ^{(0)}$), on construit 
un quadruplet $(\PP ^{\prime (0)}, T^{\prime (0)},X^{\prime (0)},Y^{\prime (0)})$ lisse en dehors du diviseur
avec $X^{\prime (0)}$ lisse et deux morphismes
$\eta ^{\prime}= (f ^{\prime}, a^{\prime}, b ^{\prime})  \,:\, (\PP ^{\prime (0)}, T^{\prime (0)},X^{\prime (0)},Y^{\prime (0)})\to (\PP', T',X',Y')$
et 
$\theta ^{(0)}= (\phi ^{(0)}, \alpha^{(0)}, \beta ^{(0)})  \,:\, (\PP ^{\prime (0)}, T^{\prime (0)},X^{\prime (0)},Y^{\prime (0)})\to
(\PP ^{(0)}, T^{(0)},X^{(0)},Y^{(0)})$
tels que 
$\theta  \circ \eta ^{\prime} = \eta  \circ \theta ^{(0)}$
et tels que 
$f^{\prime }$ soit projectif et lisse, $a^{\prime}$ soit projectif, surjectif, génériquement fini et étale,
$ (f ^{\prime}) ^{-1} ( T') = T ^{\prime (0)}$
et $T ^{\prime (0)} \cap X ^{\prime (0)}$ soit un diviseur (à croisement normaux) de $X ^{\prime(0)}$. 
Pour la commodité du lecteur, voici le diagramme illustrant nos notations : 
\begin{equation}
\label{carré-cartésien-X'0}
\xymatrix @R=0,3cm{
& 
{ Y ^{(0)} } 
\ar'[d][dd] ^-{b}
\ar[rr] 
&&
{ X ^{(0)} } 
\ar'[d][dd] ^-{a}
\ar[rr] 
&&
{\PP ^{(0)}}
\ar[dd] ^-{f}
\\ 
{Y ^{\prime (0)}} 
\ar[rr] 
\ar[ru] ^-{\beta ^{(0)}}
\ar[dd] ^-{b ^{\prime}}
&&
{X ^{\prime (0)}} 
\ar[rr] 
\ar[ru] ^-{\alpha ^{(0)}}
\ar[dd] ^(0.6){a'}
&&
{\PP ^{\prime (0)}} 
\ar[ru] ^-{\phi ^{(0)}}
\ar[dd] ^(0.6){f'}
\\
& 
{Y } 
\ar'[r][rr]
&&
{X } 
\ar'[r][rr]
&&
{\PP.} 
\\
{Y '} 
\ar[ru] _-{\beta}
\ar[rr]
&&
{X '} 
\ar[ru] _-{\alpha}
\ar[rr]
&&
{\PP'}
\ar[ru] _-{\phi}
}
\end{equation}

Par pleine fidélité du foncteur extension (voir le théorème \ref{dagT'pl-fid}), 
quitte à agrandir $T$ et $T'$, on se ramène au cas où $b $ et $b ^{\prime}$ sont finis et étales.
Notons $\U:= \PP \setminus T$ et $\U':= \PP' \setminus T'$.

D'après le cas de la compactification partielle lisse (voir \cite{caro-construction}), $\sp _{Y '\hookrightarrow \U',+} \circ \beta ^{*} \riso \beta ^{+} \circ \sp _{Y \hookrightarrow \U,+} $. D'où : 
$$ |\U' \circ \sp _{X '\hookrightarrow \PP', T',+} \circ \theta ^{*} 
\riso 
\sp _{Y '\hookrightarrow \U',+} \circ \beta ^{*} \circ |\U
\riso
\beta ^{+} \circ \sp _{Y \hookrightarrow \U,+} \circ  |\U
\riso 
 |\U' \circ \theta ^{+}\circ  \sp _{X \hookrightarrow \PP, T,+}.$$

Or, d'après \ref{a+|U-plfid-prop}
le foncteur
$(a ^{\prime+},\, |\U')$
(de \ref{a+|U-plfid-prop})
est pleinement fidèle.
Il reste donc à construire un isomorphisme compatible de la forme
$a ^{\prime +}\circ\sp _{X '\hookrightarrow \PP', T',+} \circ \theta ^{*} 
\riso 
a ^{\prime+}\circ\theta ^{+}\circ  \sp _{X \hookrightarrow \PP, T,+} $.
Comme $X ^{(0)}$ et $X ^{\prime (0)}$ sont lisses, 
$\sp _{X ^{\prime (0)}\hookrightarrow \PP^{\prime (0)}, T^{\prime (0)},+} \circ (\theta ^{(0)} )^{*} 
\riso 
(\theta ^{(0)}) ^{+}\circ  \sp _{X ^{(0)} \hookrightarrow \PP^{(0)}, T^{(0)},+} .$
Or, d'après \ref{sp-comm-a+-fini-etale}, comme $b $ et $b ^{\prime}$ sont finis et étales, 
on dispose des isomorphismes 
$a ^{+}\circ \sp _{X \hookrightarrow \PP, T,+}
\riso 
\sp _{X ^{(0)}\hookrightarrow \PP^{(0)}, T^{(0)},+}\circ a ^{*}$
et
$a ^{\prime+}\circ \sp _{X' \hookrightarrow \PP', T',+}
\riso 
\sp _{X ^{\prime (0)}\hookrightarrow \PP^{\prime (0)}, T^{\prime (0)},+}\circ a ^{\prime*}$.
Il en résulte les isomorphismes canoniques
\begin{gather}
\notag
a ^{\prime +}\circ\sp _{X '\hookrightarrow \PP', T',+} \circ \theta ^{*} 
\riso
\sp _{X ^{\prime (0)}\hookrightarrow \PP^{\prime (0)}, T^{\prime (0)},+}\circ a ^{\prime*}\circ \theta ^{*} 
\riso 
\sp _{X ^{(0)}\hookrightarrow \PP^{(0)}, T^{(0)},+}\circ \theta ^{(0)*} \circ a ^{*}
\\
\label{isoa0+-theta*com-sp}
\riso 
\theta  ^{(0) +}\circ  \sp _{X ^{(0)} \hookrightarrow \PP^{(0)}, T^{(0)},+} \circ a ^{*}
\riso 
\theta  ^{(0)+}\circ a ^{(0)+}\circ \sp _{X \hookrightarrow \PP, T,+}
\riso
a ^{(0)+}\circ\theta ^{+}\circ  \sp _{X \hookrightarrow \PP, T,+}.
\end{gather}
D'où le résultat.

\end{proof}

\section{Indépendance canonique de la catégorie des isocristaux partiellement surcohérents}

\subsection{Cas de la compactification partielle lisse}
Dans cette section, nous étudions lorsque $X$ est lisse l'indépendance en $\PP$ et $T$ de la catégorie
$(F\text{-})\mathrm{Isoc} ^{\dag \dag} (\PP, T, X/K)$. Commençons d'abord par quelques rappels et notations.

\begin{nota}
\label{nota-surhol-T}
Soient $\PP$ un $\V$-schéma formel séparé et lisse, $X$ un sous-schéma fermé de $P$, 
$T$ un diviseur de $P$ tel que $T \cap X$ soit un diviseur de $X$.
Soit $(\PP _{\alpha}) _{\alpha \in \Lambda}$ un recouvrement d'ouverts de $\PP$ tel que 
$X _\alpha:= P \cap P _{\alpha}$ soit affine.
Pour tout $\alpha \in \Lambda$, choisissons $\X _\alpha$ 
des $\V$-schémas formels lisses relevant $X _\alpha$
Rappelons que grâce à Elkik (\cite{elkik} de tels relèvements existent bien.

\begin{itemize}
\item Nous avons construit en \cite[2.5.2]{caro-construction}, 
la catégorie $(F\text{-})\mathrm{Coh} (X,\, (\X _\alpha) _{\alpha \in \Lambda},\, T\cap X/K)$ 
dont les objets sont les familles $(\E _\alpha) _{\alpha \in \Lambda}$
de $(F\text{-})\smash{\D} ^{\dag} _{\X _{\alpha} } (\hdag T  \cap X _{\alpha}) _{\Q}$-modules cohérents 
munie d'une donnée de recollement (en fait, pour définir ces données de recollements, il faut pour cela choisir des relèvements de $X _{\alpha} \cap X _{\beta}$ etc ; pour plus de détail, voir \cite[2.5.2]{caro-construction}).

\item De même, notons $(F\text{-})\mathrm{Isoc} ^{\dag \dag} (X,\, (\X _\alpha) _{\alpha \in \Lambda},\, T\cap X/K)$
la catégorie dont les objets sont les familles
$(\E _\alpha) _{\alpha \in \Lambda}$
de
$(F\text{-})\smash{\D} ^{\dag} _{\X _{\alpha} } (\hdag T  \cap X _{\alpha}) _{\Q}$-modules cohérents, $\O _{\X _{\alpha} } (\hdag T  \cap X _{\alpha}) _{\Q}$-cohérents munie d'une donnée de recollement.

\end{itemize}

\end{nota}

\begin{vide}
[Rappels]
\label{prop-donnederecol-dag}
Avec les notations de \ref{nota-6.2.1dev} et \ref{nota-surhol-T}, 
d'après la preuve de \cite[2.5.4]{caro-construction}, on dispose des deux foncteurs quasi-inverses canoniques
\begin{gather}
\label{def-Loc-dag}
\mathcal{L}oc\, :
\mathrm{Coh}  (X,\PP,T/K)
\rightarrow
\mathrm{Coh} (X,\, (\X _\alpha) _{\alpha \in \Lambda},\, T\cap X/K), 
\\
\label{def-Recol-dag}
\mathcal{R}ecol
\,:\,
\mathrm{Coh} (X,\, (\X _\alpha) _{\alpha \in \Lambda},\, T\cap X/K)
\rightarrow
\mathrm{Coh}  (X,\PP,T/K).
\end{gather}
D'après \cite[2.5.7 et 2.5.9]{caro-construction}, 
la catégorie $(F\text{-})\mathrm{Isoc} ^{\dag} (Y,X/K)$ des $(F\text{-})$isocristaux surconvergents sur $(Y,X)$ est canoniquement équivalente à $(F\text{-})\mathrm{Isoc} ^{\dag \dag} (X,\, (\X _\alpha) _{\alpha \in \Lambda},\, T\cap X/K)$. 
En composant cette équivalence avec le foncteur $\mathcal{R}ecol$ de \ref{def-Recol-dag}, on obtient alors le foncteur pleinement fidèle $\sp _{X\hookrightarrow \PP, T,+}\,:\, 
(F\text{-})\mathrm{Isoc} ^{\dag} (Y,X/K) \to \mathrm{Coh}  (X,\PP,T/K)$ (voir les premières lignes de la preuve de \cite[2.5.10]{caro-construction}). 
Grâce à \cite[6.1.4]{caro_devissge_surcoh}, on vérifie que son image essentielle est égale à 
$(F\text{-})\mathrm{Isoc} ^{\dag \dag} (\PP, T, X/K)$.
\end{vide}

\begin{lemm}
\label{lemme-coh-PXTindtP}
Soient $\PP$ un $\V$-schéma formel séparé et lisse, $X$ un sous-schéma fermé lisse de $P$, $T, T'$ deux diviseurs de $P$
tels que $T\cap X= T'\cap X$. 
On obtient alors les égalités :\\
$(F\text{-})\mathrm{Coh} (X,\, \PP,\, T)=(F\text{-})\mathrm{Coh} (X,\, \PP,\, T')$,
$(F\text{-})\mathrm{Isoc} ^{\dag \dag} (\PP, T, X/K)=(F\text{-})\mathrm{Isoc} ^{\dag \dag} (\PP, T', X/K)$.
\end{lemm}

\begin{proof}
La validation de chacune des égalités étant de nature locale en $\PP$, on peut donc supposer $\PP$ affine. 
Choisissons alors $u$ : $\X \hookrightarrow \PP$ une immersion fermée de $\V$-schémas formels lisses relevant $X \hookrightarrow P$. De plus, comme les modules ou les complexes que nous considérons sont à support dans $X$, le lemme est local en $X$. On se ramène alors au cas où $X$ est intègre. Dans ce cas, soit $T \cap X$ est un diviseur de $X$, soit
$T$ contient $X$. Lorsque $T$ contient $X$, toutes les catégories considérées sont réduites à l'élément nul (cela résulte de \cite[4.3.12]{Be1}). 
Supposons à présent que $T \cap X$ soit un diviseur de $X$.

Traitons d'abord l'égalité $(F\text{-})\mathrm{Coh} (X,\, \PP,\, T)=(F\text{-})\mathrm{Coh} (X,\, \PP,\, T')$.
Soit $\E\in(F\text{-})\mathrm{Coh} (X,\, \PP,\, T)$. 
Notons $u _{T+}$ le foncteur image directe par $u$ à singularités surconvergentes le long de $T$, i.e. le foncteur : 
$u _{T,+}\,:\, (F\text{-})D ^{\mathrm{b} }  _{\mathrm{coh}} (\D ^{\dag } _{\X} (\hdag T\cap X) _{\Q})\to 
(F\text{-})D ^{\mathrm{b} }  _{\mathrm{coh}} (\D ^{\dag } _{\PP} (\hdag T) _{\Q} )$. 
Dans ce cas, d'après le théorème de Berthelot-Kashiwara (voir \cite[5.3.3]{Beintro2}), il existe un $\D ^{\dag } _{\X} (\hdag T\cap X) _{\Q}$-module cohérent $\FF$ tel que $\E \riso u _{T+} (\FF)$.
Or, d'après \cite[1.1.9]{caro_courbe-nouveau}, $u _{T+} (\FF) \riso u _{T'+} (\FF)$. Il en résulte que $\E$ est un $\D ^{\dag } _{\PP} (\hdag T') _{\Q}$-module cohérent.
Ainsi, $\E \in (F\text{-})\mathrm{Coh} (X,\, \PP,\, T')$, i.e., 
$ (F\text{-})\mathrm{Coh} (X,\, \PP,\, T) \subset  (F\text{-})\mathrm{Coh} (X,\, \PP,\, T')$. 
Par symétrie, on obtient l'inclusion inverse. 

Vérifions à présent la deuxième égalité. Soit $\E$ un élément de $(F\text{-})\mathrm{Isoc} ^{\dag \dag} (\PP, T, X/K)$. 
Comme $X$ est lisse, d'après \ref{prop-donnederecol-dag}, 
$(F\text{-})\mathrm{Isoc} ^{\dag \dag} (\PP, T, X/K)$ 
est la catégorie des $(F\text{-})\D ^\dag _{\PP} (\hdag T) _\Q$-module cohérent à support dans $X$
et dans l'image essentielle du foncteur $\sp _{X \hookrightarrow \PP, T,+}$ (voir \cite{caro-construction}).
Plus précisément, comme on dispose du relèvement $u$, cela signifie 
qu'il existe un $\D ^{\dag } _{\X} (\hdag T\cap X) _{\Q}$-module cohérent $\FF$, 
$\O _{\X} (\hdag T\cap X) _{\Q}$-cohérent tel que $\E \riso u _{T+} (\FF)$.
Comme $T'\cap X =T\cap X$ et comme $u _{T+} (\FF) \riso u _{T'+} (\FF)$,
cela implique alors que $\E\in (F\text{-})\mathrm{Isoc} ^{\dag \dag} (\PP, T', X/K)$. On obtient alors par symétrie l'inclusion inverse.

\end{proof}

Le lemme suivant améliore \cite[3.2.6]{caro_surcoherent} ou \cite[4.10]{caro_surholonome}
car le morphisme n'est plus supposé forcément lisse. 
\begin{lemm}
\label{coh-PXTindtP}
  Soit $f$ : $\PP' \rightarrow \PP$ un morphisme de $\V$-schémas formels séparés et lisses,
  $X$ un sous-schéma fermé lisse de $P'$ tel que le morphisme induit $X \rightarrow P$ soit une immersion fermée,
  $Y$ un ouvert de $X$, 
  $T$ un diviseur de $P$ (resp. $T' $ un diviseur de $P'$) tels que $Y = X \setminus T$
  (resp. $Y = X \setminus T'$).

\begin{enumerate}
\item 
\label{coh-PXTindtP-i} 
Pour tout $\E \in (F\text{-})\mathrm{Coh} (X,\, \PP,\, T)$,
pour tout $\E '\in (F\text{-})\mathrm{Coh} (X,\, \PP',\, T')$,
pour tout $j \in \Z\setminus \{0\}$,
$$\mathcal{H} ^j (\R \underline{\Gamma} ^\dag _X f ^! (\E) ) =0,
\mathcal{H} ^j (f_+(\E')) =0.$$

\item \label{coh-PXTindtP-ii} Les foncteurs $\R \underline{\Gamma} ^\dag _X  f ^! $ et
$f _+$ induisent alors des équivalences quasi-inverses 
entre les catégories $(F\text{-})\mathrm{Coh} (X,\, \PP,\, T)$ et $(F\text{-})\mathrm{Coh} (X,\, \PP',\, T')$
(resp. entre les catégories $(F\text{-})\mathrm{Isoc} ^{\dag \dag} (\PP, T, X/K)$ et $(F\text{-})\mathrm{Isoc} ^{\dag \dag} (\PP', T', X/K)$).
\end{enumerate}
\end{lemm}

\begin{proof}
Afin d'alléger les notations, le cas avec structure de Frobenius étant analogue, nous omettrons d'indiquer {\og $(F\text{-})$\fg} dans toutes les catégories. 
Il ne coûte pas cher de supposer $P$ et $P'$ intègres. 
De même, comme $X$ est lisse et comme les modules sont à support dans $X$,
on se ramène au cas où $X$ intègre. 
Nous distinguons alors deux cas : soit $Y$ est vide soit $T \cap X = T'\cap X$ est un diviseur de $X$. 
Le premier cas est immédiat car les catégories qui interviennent sont dans ce cas nulles.
Traitons maintenant le deuxième cas. Dans ce cas $f ^{-1}( T)$ est un diviseur de $P'$. 
Par \ref{lemme-coh-PXTindtP}, comme $f ^{-1}( T) \cap X = T' \cap X$, on se ramène à supposer 
$T' = f ^{-1}( T)$.

Fixons $(\PP _{\alpha}) _{\alpha \in \Lambda}$ un recouvrement d'ouverts affines de $\PP$. 
On note 
$X _\alpha := X \cap P _\alpha$,
$\PP' _\alpha := f ^{-1} (\PP _\alpha)$, 
Pour tout $\alpha\in \Lambda $, choisissons
$\X _\alpha$ 
des $\V$-schémas formels lisses relevant $X _\alpha$
Soient $u ' _{\alpha }$ :
$\X  _{\alpha} \rightarrow \PP ' _{\alpha}$
des relèvements de
$X  _{\alpha} \rightarrow P ' _{\alpha}$
On note $f _\alpha$ : $\PP ' _\alpha \rightarrow \PP_\alpha$
le morphisme induit par $f$.
On pose $u  _\alpha:=f _{\alpha} \circ u ' _\alpha $ : $\X _\alpha \rightarrow \PP _\alpha$.
Notons que lorsqu'il faudra dans la suite de cette preuve vérifier les commutations aux données de recollement, il eût fallu
choisir d'autres relèvements (e.g. de $X _{\alpha} \cap X _{\beta}$, $X _{\alpha} \cap X _{\beta}\cap X _{\gamma}$), mais nous laissons au lecteur le soin de l'écrire.

$\bullet $ Soit $\E \in \mathrm{Coh} (X,\, \PP,\, T)$.
On dispose des isomorphismes canoniques (le deuxième résulte du dernier point de \cite[1.15]{caro_surholonome}) :
\begin{equation}
\label{coh-PXTindtP-iso}
\R \underline{\Gamma} ^\dag _X f ^! (\E) | \PP ' _\alpha 
\riso 
\R \underline{\Gamma} ^\dag _{X _\alpha} f _\alpha ^! (\E | \PP  _\alpha)
\riso 
u ' _{\alpha+} \circ u ^{\prime !}  _{\alpha}  \circ f _\alpha ^! (\E | \PP  _\alpha) 
\riso u ' _{\alpha+} \circ u ^{!}  _{\alpha}  (\E | \PP  _\alpha)  .
\end{equation}
Comme cela est local en $\PP'$, il en résulte alors que, pour tout entier $j \not = 0$, on ait
$\mathcal{H} ^j (\R \underline{\Gamma} ^\dag _X f ^! (\E) ) =0$ 
et 
$\R \underline{\Gamma} ^\dag _X f ^! (\E) \in \mathrm{Coh} (X,\, \PP',\, T')$.

L'équivalence de catégories 
$\mathcal{L}oc$ : $\mathrm{Coh} (X,\, \PP ,\, T)
\cong
\mathrm{Coh} (X,\, (\X _\alpha) _{\alpha \in \Lambda},\, T\cap X)$
de \ref{prop-donnederecol-dag} est définie en posant
$\mathcal{L}oc (\E )=  (u ^{!}  _{\alpha}  (\E | \PP  _\alpha))_\alpha$, ce dernier étant muni de la donnée de recollement canonique
(voir \cite[2.5.4]{caro-construction}).
De même avec des primes. 
Or, en appliquant $u ^{\prime !}  _{\alpha}  $ à \ref{coh-PXTindtP-iso},
via le théorème Berthelot-Kashiwara (voir \cite[5.3.3]{Beintro2}),
on obtient l'isomorphisme canonique
$u ^{\prime !}  _{\alpha} (\R \underline{\Gamma} ^\dag _X f ^! (\E) | \PP ' _\alpha )
\riso 
u ^{!}  _{\alpha}  (\E | \PP  _\alpha)  $, cet isomorphisme commutant aux données de recollement respectives.
On dispose ainsi du diagramme commutatif (à isomorphisme canonique près) :
\begin{equation}
\label{coh-PXTindtP-diagloc}
\xymatrix {
{\mathrm{Coh} (X,\, \PP ,\, T) }
\ar[r] ^-{\mathcal{L}oc} _-{\cong}
\ar[d] ^-{\R \underline{\Gamma} ^\dag _X f ^! } 
& 
{\mathrm{Coh} (X,\, (\X _\alpha) _{\alpha \in \Lambda},\, T\cap X) } 
\ar@{=}[d]
\\
{\mathrm{Coh} (X,\, \PP ',\, T') }
\ar[r] ^-{\mathcal{L}oc} _-{\cong}
& 
{\mathrm{Coh} (X,\, (\X _\alpha) _{\alpha \in \Lambda},\, T'\cap X) .} 
}
\end{equation}

$\bullet $ 
Soit $\E ' \in \mathrm{Coh} (X,\, \PP',\, T')$. 
D'après le théorème de Berthelot-Kashiwara (voir \cite[5.3.3]{Beintro2}), comme $\E'  | \PP' _\alpha $ est à support dans 
$X _{\alpha}$, $u_{\alpha } ^{\prime!} (\E'  | \PP' _\alpha )$ est un 
$\D ^{\dag} _{\X _{\alpha}} (\hdag T '\cap X _{\alpha}  )_{\Q}$-module cohérent et 
$\E'  | \PP' _\alpha  \riso u '_{\alpha +} \circ u_{\alpha } ^{\prime !} (\E'  | \PP' _\alpha )$.
Il en résulte l'isomorphisme $f _+ ( \E' ) | \PP _\alpha \riso f _{\alpha +} (\E'  | \PP' _\alpha ) \riso u_{\alpha +} \circ u_{\alpha } ^{\prime!} (\E'  | \PP' _\alpha )$.
Comme cela est local en $\PP$, on en déduit, pour $j \not =0$,
$ \mathcal{H} ^j (f_+(\E')) =0$ et $f_+(\E') \in \mathrm{Coh} (X,\, \PP,\, T)$.
 
L'équivalence de catégories 
$\mathcal{R}ecol$ :
$\mathrm{Coh} (X,\, (\X _\alpha) _{\alpha \in \Lambda},\, T \cap X)
\cong
\mathrm{Coh} (X,\, \PP ' ,\, T')$
de \ref{prop-donnederecol-dag} est
définie par $(\E _\alpha)_\alpha \mapsto (u ' _{\alpha+} (\E _\alpha))_\alpha$, où $(u ' _{\alpha+} (\E _\alpha))_\alpha$ est muni de la structure de recollement canonique, ce qui définit un objet de
$\mathrm{Coh} (X,\, \PP ' ,\, T')$ (voir \cite[2.5.4]{caro-construction}). 
On dispose aussi du foncteur $\mathcal{R}ecol$ :
$\mathrm{Coh} (X,\, (\X _\alpha) _{\alpha \in \Lambda},\, T' \cap X)
\cong
\mathrm{Coh} (X,\, \PP ,\, T)$, 
défini par 
$(\E _\alpha) _\alpha \mapsto (u  _{\alpha+} (\E _\alpha))_\alpha$.
Il en résulte le diagramme canonique commutatif (à isomorphisme canonique près) : 
\begin{equation}
\label{coh-PXTindtP-diagrecol}
\xymatrix {
{\mathrm{Coh} (X,\, (\X _\alpha) _{\alpha \in \Lambda},\, T\cap X) } 
\ar[r] ^-{\mathcal{R}ecol} _-{\cong}
\ar@{=}[d]
& 
{\mathrm{Coh} (X,\, \PP ,\, T) }
\\
{\mathrm{Coh} (X,\, (\X _\alpha) _{\alpha \in \Lambda},\, T' \cap X) } 
\ar[r] ^-{\mathcal{R}ecol} _-{\cong}
& 
{\mathrm{Coh} (X,\, \PP ',\, T') .} 
\ar[u] ^-{ f _+}
}
\end{equation}

$\bullet$ Comme les foncteurs $\mathcal{L}oc$ et $\mathcal{R}ecol$ sont quasi-inverses (voir \ref{prop-donnederecol-dag}),
via les diagrammes commutatifs (à isomorphisme canonique près)
\ref{coh-PXTindtP-diagloc} et \ref{coh-PXTindtP-diagrecol}, 
il en résulte que les foncteurs
$\R \underline{\Gamma} ^\dag _X f ^! $ et
$f _+$ induisent des équivalences quasi-inverses entre les catégories
$\mathrm{Coh} (X,\, \PP,\, T)$ et $\mathrm{Coh} (X,\, \PP',\, T')$.

$\bullet$ Concernant l'équivalence de catégories entre 
$\mathrm{Isoc} ^{\dag \dag} (\PP, T, X/K)$ et $\mathrm{Isoc} ^{\dag \dag} (\PP', T', X/K)$,
on procède de la même façon : il s'agit de remplacer dans la preuve respectivement
$\mathrm{Coh} (X,\, (\X _\alpha) _{\alpha \in \Lambda},\, T\cap X) $ par 
$\mathrm{Isoc} ^{\dag\dag } (X,\, (\X _\alpha) _{\alpha \in \Lambda},\, T\cap X/K)$ 
et 
$\mathrm{Coh} (X,\, \PP,\, T)$ par $\mathrm{Isoc} ^{\dag \dag} (\PP, T, X/K)$
et de même avec des primes.

\end{proof}

\subsection{Isomorphisme entre image inverse et image inverse extraordinaire d'un isocristal partiellement surcohérent}

\begin{lemm}
\label{adj-carr}
Considérons le diagramme commutatif suivant 
\begin{equation}
\label{adj-carr-diag}
\xymatrix {
{\PP ^{\prime (0)}} 
\ar[d] ^-{f'}  \ar[r] ^-{\phi ^{(0)}} 
& 
{\PP ^{(0)}} 
\ar[d] ^-{f} 
\\
{\PP ^{\prime}} 
\ar[r] ^-{\phi} 
& 
{\PP,} 
}
\end{equation}
où $\phi$ et $\phi ^{(0)}$ sont des morphismes propres de $\V$-schémas formels lisses,
$f$ et $f'$ sont des morphismes lisses $\V$-schémas formels lisses.
Soient $T$ un diviseur de $P$ tel que $T':= \phi ^{-1} (T)$ soit un diviseur de $P'$. 
On pose $ T ^{(0)}:= f ^{-1} (T)$ et 
$ T ^{\prime (0)}:= f ^{\prime -1} (T')$.
Soit $\E ' \in D ^\mathrm{b} _\mathrm{coh} (\D ^\dag _{\PP'} (\hdag T') _{\Q})$.
On dispose alors des morphismes canoniques :
$\phi ^{(0)} _{+} \circ f ^{\prime! } (\E ')\to f ^{! }\circ \phi  _{+} (\E ')$
et
$ f ^{+}\circ \phi ^{(0)} _{+} (\E ')
\to
\phi  _{+} \circ f ^{\prime+} (\E ')$.
De plus, lorsque le diagramme \ref{adj-carr-diag} est cartésien, ces deux morphismes sont des isomorphismes.
\end{lemm}

\begin{proof}
Le deuxième morphisme que l'on doit établir découle par dualité et via l'isomorphisme de dualité relative du premier.
Traitons donc celui-ci. Notons $\PP'' := \PP ' \times _{\PP} \PP ^{(0)} $, $\iota \,:\, \PP ^{\prime (0)} \to \PP'$,
$f''\,:\, \PP '' \to \PP'$ et $\phi '' \,:\, \PP'' \to \PP ^{(0)}$ les morphismes canoniques. 
Comme $\phi ''$ et $\phi ^{(0)}$ sont propres, $\iota$ l'est aussi. 
On dispose alors du morphisme d'adjonction $\iota _{+} \circ \iota ^{!} \to Id$. 
D'après la proposition \cite[3.1.8]{caro_surcoherent}, 
on dispose de l'isomorphisme canonique 
$\phi '' _{+} \circ f ^{\prime \prime !} (\E ') 
\riso 
f ^{ !} \circ \phi _{+} (\E ') $.
On en déduit les morphismes 
$$\phi ^{(0)} _{+} \circ f ^{\prime! } (\E ') 
\riso 
\phi '' _{+} \circ \iota _{+} \circ \iota ^{!} \circ  f ^{\prime\prime ! } (\E ')
\to 
\phi '' _{+} \circ   f ^{\prime\prime ! } (\E ')
\riso
 f ^{! }\circ \phi  _{+} (\E ').$$
Comme lorsque le diagramme \ref{adj-carr-diag} est cartésien
le morphisme d'adjonction $\iota _{+} \circ \iota ^{!} \to Id$ est un isomorphisme, on valide la dernière assertion.
\end{proof}

\begin{nota}
\label{nota-adj}
Soient $\theta = (\phi, \alpha, \beta) \,:\, (\PP', T',X',Y')\to (\PP, T,X,Y)$,
$\eta= (f , a, b )  \,:\, (\PP ^{(0)}, T^{(0)},X^{(0)},Y^{(0)})\to (\PP, T,X,Y)$,
$\theta ^{(0)}= (\phi ^{(0)}, \alpha^{(0)}, \beta ^{(0)})  \,:\, (\PP ^{\prime (0)}, T^{\prime (0)},X^{\prime (0)},Y^{\prime (0)})\to
(\PP ^{(0)}, T^{(0)},X^{(0)},Y^{(0)})$
et 
$\eta ^{\prime}= (f ^{\prime}, a^{\prime}, b ^{\prime})  \,:\, (\PP ^{\prime (0)}, T^{\prime (0)},X^{\prime (0)},Y^{\prime (0)})\to (\PP', T',X',Y')$
quatre morphismes de quadruplets lisses en dehors du diviseur
tels que 
$\theta  \circ \eta ^{\prime} = \eta  \circ \theta ^{(0)}$.
On obtient ainsi le parallélépipède rectangle de \ref{carré-cartésien-X'0}.
Afin de donner un sens à la notion d'image directe (du côté formel), on suppose que les quatre rectangles de \ref{carré-cartésien-X'0}
(i.e., les faces de devant, de derrière, du haut, du bas) vérifient les hypothèses analogues à celles que satisfait
\ref{nota-diag-deux}, notamment $T'= \phi ^{-1} (T)$, $T ^{\prime (0)}= \phi ^{-1} (T ^{(0)})$, 
$T ^{ (0)}= f ^{-1} (T)$ etc.

\end{nota}

\begin{lemm}
\label{lemm-iso-adj-cube}
On garde les notations et hypothèses de \ref{nota-adj}.
Pour tout $\E' \in \mathrm{Isoc} ^{*} (\PP', T', X'/K)$, on dispose des morphismes canoniques fonctoriels en $\E'$ : 
\begin{equation}
\label{iso-adj-cube}
\alpha ^{(0)} _{+} \circ  a ^{\prime !}(\E')\to a ^{!} \circ \alpha _{+} (\E ') ,\
a ^{+} \circ \alpha _{+} (\E ') \to \alpha ^{(0)} _{+} \circ  a ^{\prime +}(\E').
\end{equation}
Lorsque la face de gauche correspondante de \ref{nota-diag-deux} est cartésienne, ces morphismes sont des isomorphismes.
\end{lemm}

\begin{proof}
Via l'isomorphisme de dualité relative et de bidualité, le deuxième morphisme se construit en appliquant le foncteur dual au premier (utilisé, grâce à \ref{Isoc*=dagdag}, pour $\DD _{T'} (\E')$ à la place de $\E'$).
On définit le premier morphisme comme égal au composé :
\begin{gather}
\notag
\alpha ^{(0)} _{+} \circ  a ^{\prime !}(\E') = 
\phi ^{(0)} _{+} \circ  \R \underline{\Gamma}  ^{\dag } _{X ^{\prime(0)} } \circ  f ^{\prime !}(\E')
\to 
\phi ^{(0)} _{+} \circ  \R \underline{\Gamma}  ^{\dag } _{(\phi ^{(0)} ) ^{-1} (X ^{(0)} )} \circ  f ^{\prime !}(\E')
\riso
\\
\label{const-morp-adj-cube}
 \R \underline{\Gamma}  ^{\dag } _{X ^{(0)} }\circ \phi ^{(0)} _{+}  \circ  f ^{\prime !}(\E')
 \to 
  \R \underline{\Gamma}  ^{\dag } _{X ^{(0)} }\circ f ^{! }\circ \phi  _{+} (\E')
  = a ^{!} \circ \alpha _{+} (\E ') ,
\end{gather}
dont le dernier morphisme se déduit par fonctorialité de \ref{adj-carr}. 
  
  Supposons à présent que la face de gauche de \ref{nota-diag-deux} est cartésienne et vérifions alors que ce morphisme est un isomorphisme.
  Grâce à \cite[4.3.12]{Be1}, on se ramène au cas où $T$ est vide, i.e., la face du milieu de \ref{nota-diag-deux} 
  est égale à celle de gauche et est en particulier cartésienne.
  On en déduit alors un diagramme commutatif analogue à \ref{nota-diag-deux} où $\PP ^{\prime (0)}$ est remplacé par 
  $ \PP ^{\prime} \times _{\PP } \PP ^{(0)}$ (en effet, on dispose d'une immersion fermée canoniquement induite :
  $X ^{\prime} \times _{X } X ^{(0)} \hookrightarrow \PP ^{\prime} \times _{\PP } \PP ^{(0)}$. 
  Notons alors $\iota \,:\, 
  \PP ^{\prime (0)} \to \PP ^{\prime} \times _{\PP } \PP ^{(0)}$ le morphisme canonique.
  
  Comme $\iota$ est propre (n'oublions pas que par hypothèse les morphismes $f$, $f'$, $\phi$, $\phi ^{(0)}$ sont propres), d'après \ref{*dagdag-coh-PXTindtP}, $\iota _{+}$ et $\iota ^{!}$ induisent des équivalences quasi-inverses
  entre $\mathrm{Isoc} ^{*} ( \PP ^{\prime (0)} , X^{\prime (0)} /K)$ et 
  $\mathrm{Isoc} ^{*} (\PP ^{\prime} \times _{\PP } \PP ^{(0)}, X^{\prime (0)} /K)$.
  On se ramène ainsi au cas où $ \PP ^{\prime (0)} = \PP ^{\prime} \times _{\PP } \PP ^{(0)}$.
  Dans ce cas, d'après \ref{adj-carr}, le dernier morphisme de \ref{const-morp-adj-cube} est un isomorphisme.
  De plus, $(\phi ^{(0)} ) ^{-1} (X ^{(0)} )= P ' \times _{P} X ^{(0)}$. Or comme $\E'$ est à support dans $X'$, 
  $ f ^{\prime !}(\E')$ est à support dans $(f ') ^{-1}( X')=X ' \times _{P} P ^{(0)}$.
  Comme $(P ' \times _{P} X ^{(0)}) \cap (X ' \times _{P} P ^{(0)})=X ^{\prime} \times _{X } X ^{(0)}= X ^{\prime(0)}$, 
  il en résulte que 
  $\R \underline{\Gamma}  ^{\dag } _{(\phi ^{(0)} ) ^{-1} (X ^{(0)} )} \circ  f ^{\prime !}(\E')$
  est à support dans $X ^{\prime(0)}$.
  Le premier morphisme de \ref{const-morp-adj-cube} est donc un isomorphisme.
  D'où le résultat.
\end{proof}

\begin{lemm}
\label{theo-adj-fini-étale}
Soit $\theta = (\phi, \alpha, \beta) \,:\, (\PP', T',X',Y')\to (\PP, T,X,Y)$ un morphisme de quadruplets lisses en dehors du diviseur
tel que $\phi$ soit propre et lisse, $T'= \phi ^{-1} (T)$, $\alpha$ soit propre, surjectif et $\beta$ soit fini et étale.
Le foncteur $ \alpha _{+} \,:\,
\mathrm{Isoc} ^{*} (\PP', T', X'/K)
\to 
\mathrm{Isoc} ^{*} (\PP, T, X/K)$
est adjoint à gauche de $\alpha ^{+}$.
\end{lemm}

\begin{proof}
0) Avec le lemme \ref{rema-isocDagDagDense} (et la remarque \ref{rema-isocDagDense}), 
on se ramène au cas où $Y$ (resp. $Y'$) est intègre et dense dans $X$ (resp. $X'$).
Grâce au théorème de désingularisation de de Jong, 
il existe alors un quadruplet $(\PP ^{(0)}, T^{(0)},X^{ (0)},Y^{(0)})$ lisse en dehors du diviseur
avec $X^{ (0)}$ lisse et un morphisme
$\eta = (f , a, b )  \,:\, (\PP ^{ (0)}, T^{ (0)},X^{ (0)},Y^{ (0)})\to (\PP, T,X,Y)$
tels que  
$f$ soit projectif et lisse, $a$ soit projectif, surjectif, génériquement fini et étale, $ f ^{-1} ( T) = T ^{ (0)}$
et 
$T ^{(0)} \cap X ^{(0)}$ soit un diviseur (à croisement normaux) de $X ^{(0)}$. 

Or, d'après \ref{dagT'pl-fid}, pour tout diviseur $\widetilde{T} '$ contenant $T ^{\prime}$ et tel que $ X' \setminus \widetilde{T} '$ soit dense dans
$Y'$, le foncteur extension 
$(\hdag \widetilde{T} ')$ est pleinement fidèle. Quitte à agrandir $T$, on peut supposer que $b$ est fini et étale. 
Posons 
$Y^{\prime (0)} :=Y ' \times _{Y} Y ^{(0)}$, 
$X^{\prime \prime (0)} :=X ' \times _{X} X ^{(0)}$ 
$\PP^{\prime (0)} :=\PP ' \times _{\PP} \PP ^{(0)}$.
Notons $X^{\prime (0)}$ la normalisation de $X^{\prime \prime (0)}$,
$ T ^{ \prime (0)}:=f ^{\prime -1} ( T')$
$\theta ^{(0)}= (\phi ^{(0)}, \alpha^{(0)}, \beta ^{(0)})  \,:\, (\PP ^{\prime (0)}, T^{\prime (0)},X^{\prime (0)},Y^{\prime (0)})\to
(\PP ^{(0)}, T^{(0)},X^{(0)},Y^{(0)})$
et 
$\eta ^{\prime}= (f ^{\prime}, a^{\prime}, b ^{\prime})  \,:\, (\PP ^{\prime (0)}, T^{\prime (0)},X^{\prime (0)},Y^{\prime (0)})\to (\PP', T',X',Y')$
les morphismes canoniques induits par les projections.
On se retrouve ainsi dans la situation de \ref{nota-adj}
avec en outre 
$X ^{(0)}$ est lisse, $X ^{\prime (0)}$ est normal et la face de gauche du diagramme correspondant de \ref{nota-diag-deux} cartésienne.

{\it 1) Construction du morphisme $\alpha ^{(0)} _{+} \circ \alpha ^{(0)+} \to Id $.}

Comme $\alpha ^{(0)} $ est propre et surjectif avec $X ^{(0)}$ lisse et $X ^{\prime (0)}$ normal, 
comme $\beta ^{(0)} $ est fini et étale,
Tsuzuki a construit dans ce contexte (voir le chapitre \cite[5]{tsumono}), le foncteur image directe
$\alpha ^{(0)} _{*} \,:\,  \mathrm{Isoc} ^{\dag} (Y^{\prime (0)},X^{\prime (0)}/K) \to  \mathrm{Isoc} ^{\dag} (Y^{(0)},X^{(0)}/K)$.
De plus, il y a vérifié que ce foncteur est adjoint à droite et adjoint à gauche du foncteur image inverse
$(\alpha ^{(0)} ) ^{*} \,:\,    \mathrm{Isoc} ^{\dag} (Y^{(0)},X^{(0)}/K)
\to
\mathrm{Isoc} ^{\dag} (Y^{\prime (0)},X^{\prime (0)}/K) $.
Pour alléger les notations, on pose 
$\sp _{+}:= \sp _{X^{(0)} \hookrightarrow \PP ^{(0)}, T ^{(0)},+}$
et
$\sp ' _{+}:= \sp _{X^{\prime(0)} \hookrightarrow \PP ^{\prime(0)}, T ^{\prime(0)},+}$.
Soient $E ^{(0)}\in  \mathrm{Isoc} ^{\dag} (Y^{(0)},X^{(0)}/K)$,
$E ^{\prime (0)} \in \mathrm{Isoc} ^{\dag} (Y^{\prime (0)},X^{\prime (0)}/K)$.

Comme $\sp _{+}$ est pleinement fidèle et comme 
$\alpha ^{(0)} _{*} $ est adjoint à droite de $(\alpha ^{(0)} ) ^{*}$, on obtient 
l'isomorphisme canonique fonctoriel en $E ^{(0)}$ et $E ^{\prime (0)}$ :
\begin{equation}
\label{Hom-adj-finiét-1}
\mathrm{Hom} _{ \mathrm{Isoc} ^{\dag} (Y^{\prime (0)},X^{\prime (0)}/K)} 
((\alpha ^{(0)} ) ^{*} (E ^{(0)}), E ^{\prime (0)})
\riso
\mathrm{Hom} _{\mathrm{Isoc} ^{*} (\PP^{(0)}, T^{(0)}, X^{(0)}/K)} 
(\sp _{+} (E ^{ (0)}), \sp _{+} \circ \alpha ^{(0)} _{*} (E ^{\prime (0)}))
\end{equation}

D'après \ref{theta*com-sp}, on dispose de l'isomorphisme $\sp '_{+} \circ (\alpha ^{(0)} ) ^{*}  \riso (\alpha ^{(0)} ) ^{+} \circ \sp _{+}$.
De plus, d'après \ref{Isoc*a+}.\ref{Isoc*a+5}, 
le foncteur $\alpha ^{(0)} _{+}$ est adjoint à droite de $(\alpha ^{(0)} ) ^{+}$.
Comme le foncteur $\sp _{+} ^{\prime}$ est pleinement fidèle, 
on déduit de ces deux faits l'isomorphisme canonique fonctoriel en $E ^{(0)}$ et $E ^{\prime (0)}$ :
\begin{equation}
\label{Hom-adj-finiét-2}
\mathrm{Hom} _{ \mathrm{Isoc} ^{\dag} (Y^{\prime (0)},X^{\prime (0)}/K)} 
((\alpha ^{(0)} ) ^{*} (E ^{(0)}), E ^{\prime (0)})
\riso
\mathrm{Hom} _{\mathrm{Isoc} ^{*} (\PP^{(0)}, T^{(0)}, X^{(0)}/K)} 
(\sp _{+} (E ^{ (0)}), \alpha ^{(0)} _{+} \circ \sp ' _{+} (E ^{\prime (0)}))
\end{equation}
Il découle des isomorphismes \ref{Hom-adj-finiét-1} et \ref{Hom-adj-finiét-2}, 
l'isomorphisme canonique fonctoriel en $E ^{\prime (0)}$ :
\begin{equation}
\label{Hom-adj-finiét-3} 
\sp _{+} \circ \alpha ^{(0)} _{*} (E ^{\prime (0)}) \riso \alpha ^{(0)} _{+} \circ \sp ' _{+} (E ^{\prime (0)}).
\end{equation}

Soient $\E ^{(0)}\in \mathrm{Isoc} ^{*} (\PP^{(0)}, T^{(0)}, X^{(0)}/K)$,
$\E ^{\prime (0)} \in \mathrm{Isoc} ^{*} (\PP^{\prime(0)}, T^{\prime(0)}, X^{\prime(0)}/K)$.
Comme le foncteur $\sp _{+}$ est essentiellement surjectif, il existe $E ^{(0)} _{1}$, $E ^{(0)}_{2}$ (on en choisit deux afin
de valider l'indépendance par rapport au choix)
tel que 
$\E ^{(0)} \liso \sp _{+} (E ^{(0)} _{1})$ et 
 $\E ^{(0)} \riso \sp _{+} (E ^{(0)} _{2})$. 
 Comme $\alpha ^{(0)} _{*} $ est adjoint à gauche de $(\alpha ^{(0)} ) ^{*}$, 
 avec \ref{Hom-adj-finiét-3} et \ref{theta*com-sp}, on obtient le diagramme canonique :
\begin{equation}
\label{Hom-adj-finiét-4}
\xymatrix @R=0,3cm{
{\alpha ^{(0)} _{+} \circ \alpha ^{(0)+} \circ \sp _{+} (E ^{(0)} _{1}) } 
\ar[r] ^-{\sim}
\ar[d] ^-{\sim}
& 
{\alpha ^{(0)} _{+} \circ \alpha ^{(0)+}  (\E) } 
\ar[r] ^-{\sim}
& 
{\alpha ^{(0)} _{+} \circ \alpha ^{(0)+} \circ \sp _{+} (E ^{(0)} _{2}) } 
\ar[d] ^-{\sim}
\\ 
{ \sp _{+} \circ \alpha ^{(0)} _{*} \circ \alpha ^{(0)*}(E ^{(0)} _{1}) } 
\ar[d] ^-{\mathrm{adj}}
& 
{} 
& 
{ \sp _{+} \circ \alpha ^{(0)} _{*} \circ \alpha ^{(0)*}(E ^{(0)} _{2}) } 
\ar[d] ^-{\mathrm{adj}}
\\ 
{ \sp _{+} (E ^{(0)} _{1}) } 
\ar[r] ^-{\sim}
& 
{\E } 
\ar[r] ^-{\sim}
& 
{ \sp _{+} (E ^{(0)} _{2}) } 
}
\end{equation}
Comme $\sp _{+}$ est pleinement fidèle, le morphisme 
$ \sp _{+} (E ^{(0)} _{1}) \riso  \sp _{+} (E ^{(0)} _{2})$ provient 
d'un morphisme 
$E ^{(0)} _{1} \riso E ^{(0)} _{2}$.
Le diagramme \ref{Hom-adj-finiét-4} est donc commutatif. 
Le morphisme induit 
$\alpha ^{(0)} _{+} \circ \alpha ^{(0)+}  (\E ^{(0)})  \to \E ^{(0)} $
est donc canonique.

{\it 1 bis) De même, on construit le morphisme canonique 
$ \E ^{\prime(0)} \to   \alpha ^{(0)+} \circ \alpha ^{(0)} _{+}  (\E ^{\prime(0)}) $.}

{\it 2) Construction du morphisme canonique
$\alpha _{+} \circ \alpha ^{+} \to Id$.}

L'idée est d'utiliser la pleine fidélité de $(a ^{+}, |\U)$.
Soient $\E\in \mathrm{Isoc} ^{*} (\PP, T, X/K)$, $\E' \in \mathrm{Isoc} ^{*} (\PP', T', X'/K)$.
D'après \ref{lemm-iso-adj-cube},
on dispose de l'isomorphisme de changement de base $a ^{+} \circ \alpha _{+} (\E ') \to \alpha ^{(0)} _{+} \circ  a ^{\prime +}(\E')$. 
Avec le morphisme $\alpha ^{(0)} _{+} \circ \alpha ^{(0)+} \to Id $ construit à l'étape 1) de la preuve, 
on obtient le morphisme canonique : 
\begin{equation}
\label{Hom-adj-finiét-5}
a ^{+} \circ \alpha _{+} \circ \alpha ^{+} (\E)
\riso 
\alpha ^{(0)} _{+} \circ  a ^{\prime +} \circ \alpha ^{+} (\E)
\riso 
\alpha ^{(0)} _{+} \circ \alpha ^{(0)+} \circ a ^{+}  (\E) 
\overset{\mathrm{adj}}{\longrightarrow}
a ^{+}  (\E) .
\end{equation}
D'un autre côté, en notant $\U := \PP \setminus T$, 
d'après \ref{b+b+}, on bénéficie de l'isomorphisme d'adjonction
\begin{equation}
\label{Hom-adj-finiét-6}
|\U \circ \alpha _{+} \circ \alpha ^{+} (\E)
\riso 
\beta _{+} \circ \beta ^{+} (\E |\U )
\overset{\mathrm{adj}}{\longrightarrow}
\E |\U.
\end{equation}
D'après la proposition \ref{a+|U-plfid-prop}, le foncteur $(a ^{+}, |\U)$ est pleinement fidèle. 
Via les morphismes canoniques de  \ref{Hom-adj-finiét-5} et \ref{Hom-adj-finiét-6}, 
il en résulte le morphisme canonique 
$\alpha _{+} \circ \alpha ^{+} (\E) \to \E$.

2 bis) De même, on construit le morphisme
$ \E '\to \alpha ^{+}  \circ \alpha _{+}  (\E')$.

3) Pour en déduire que 
 $ \alpha _{+}$ est adjoint à gauche de $\alpha ^{+}$,
il suffit alors grâce à la proposition \cite[4.3.12]{Be1} de le 
vérifier en dehors des diviseurs. 
On se ramène ainsi à la situation déjà connue de \ref{b+b+}.
\end{proof}

\begin{theo}
\label{coro-theo-!=+}
Soit $\theta = (\phi, \alpha, \beta) \,:\, (\PP', T',X',Y')\to (\PP, T,X,Y)$ un morphisme de quadruplets lisses en dehors du diviseur.
Les foncteurs $ \alpha ^{+},  \alpha ^{!} \,:\,
\mathrm{Isoc} ^{*} (\PP, T, X/K) \to \mathrm{Isoc} ^{*} (\PP', T', X'/K) $
sont canoniquement isomorphisme.
\end{theo}

\begin{proof}
0) Avec le lemme \ref{rema-isocDagDagDense} (et la remarque \ref{rema-isocDagDense}), on se ramène au cas où $Y$ (resp. $Y'$) est intègre et dense dans $X$ (resp. $X'$).

1) Vérifions à présent le théorème lorsque 
$\phi$ soit propre et lisse, $T'= \phi ^{-1} (T)$, $\alpha$ soit propre, surjectif, génériquement fini et étale.

Grâce à \ref{dagT'pl-fid}, pour tout diviseur $\widetilde{T} '$ contenant $T ^{\prime}$
et tel que $ X' \setminus \widetilde{T} '$ soit dense dans
$Y'$, le foncteur d'extension des
scalaires $(\hdag \widetilde{T} ')$ est pleinement fidèle. Quitte à agrandir $T$, on se ramène ainsi 
au cas où $\beta$ est fini et étale. 
Dans ce cas, d'après \ref{theo-adj-fini-étale}, le foncteur $\alpha ^{+} $ est donc adjoint à droite de $\alpha _{+}$
Or, il en est de même de $\alpha ^{!}$.
D'où le résultat.

{\it 2) Cas Général.}

Grâce au théorème de désingularisation de de Jong, 
il existe $\eta= (f , a, b )  \,:\, (\PP ^{(0)}, T^{(0)},X^{(0)},Y^{(0)})\to (\PP, T,X,Y)$,
$\theta ^{(0)}= (\phi ^{(0)}, \alpha^{(0)}, \beta ^{(0)})  \,:\, (\PP ^{\prime (0)}, T^{\prime (0)},X^{\prime (0)},Y^{\prime (0)})\to
(\PP ^{(0)}, T^{(0)},X^{(0)},Y^{(0)})$
et 
$\eta ^{\prime}= (f ^{\prime}, a^{\prime}, b ^{\prime})  \,:\, (\PP ^{\prime (0)}, T^{\prime (0)},X^{\prime (0)},Y^{\prime (0)})\to (\PP', T',X',Y')$
trois morphismes de quadruplets lisses en dehors du diviseur
tels que 
$\theta  \circ \eta ^{\prime} = \eta  \circ \theta ^{(0)}$, 
$f$ et $f'$ soient projectifs et lisses,
$a$ et $a'$ soient projectifs, génériquement fini et étale, 
$X^{\prime (0)}$ et $X^{ (0)}$ soient lisses, 
$T ^{\prime (0)}= f ^{\prime -1} (T ^{(0)})$, 
$T ^{ (0)}= f ^{-1} (T)$,
$T ^{(0)} \cap X ^{(0)}$ soit un diviseur (à croisement normaux) de $X ^{(0)}$,
$T ^{\prime(0)} \cap X ^{\prime(0)}$ soit un diviseur (à croisement normaux) de $X ^{\prime(0)}$. 
Il résulte de l'étape 1), que les foncteurs $a ^{+}$ et $a ^{!}$ (resp. $a ^{\prime+}$ et $a ^{\prime !}$) sont canoniquement isomorphes.
Or, d'après la proposition \ref{a+|U-plfid-prop}, le foncteur $(a ^{\prime +}, |\U)$ est pleinement fidèle.
On se ramène donc au cas où $X$ et $X'$ sont lisses, ce qui est déjà connu (voir \cite{caro-construction}).
\end{proof}

\begin{rema}
\label{AvecFrob!=+facile}
Le but principal de cette section a été d'établir \ref{coro-theo-!=+}.
Lorsque l'on dispose d'une structure de Frobenius, cela se prouve plus facilement car on dispose alors du théorème
de pleine fidélité de Kedlaya de \cite{kedlaya_full_faithfull}.
En effet, grâce à ce théorème de pleine fidélité, il suffit de le voir en dehors des diviseurs, ce qui nous ramène à la situation de la compactification partiel lisse (i.e., celle de \cite{caro-construction}).

Cet isomorphisme de \ref{coro-theo-!=+} est un ingrédient fondamental dans la preuve 
du lemme \cite[6.3.1]{caro_devissge_surcoh} ou ici dans sa version étendue
\ref{casliss631dev}.
\end{rema}

Voici un corollaire du théorème \ref{coro-theo-!=+} :

\begin{prop}
\label{sp+-comm-dual}
Soit $(\PP, T,X,Y)$ un quadruplet lisse en dehors du diviseur (voir la définition \ref{defi-triplet}).
Soit $E \in \mathrm{Isoc} ^{\dag} (Y,X/K)$. On dispose de l'isomorphisme canonique dans $\mathrm{Isoc} ^{*} (\PP, T, X/K)$ : 
$$\sp _{X \hookrightarrow \PP, T,+} (E ^{\vee})
\riso 
\DD _{T}\circ \sp _{X \hookrightarrow \PP, T,+} (E ), $$
où $\vee$ désigne le dual dans $\mathrm{Isoc} ^{\dag} (Y,X/K)$
 et 
$\sp _{X \hookrightarrow \PP, T,+}$ a été défini au théorème
\ref{diag-eqcat-dagdag-ii}.
\end{prop}

\begin{proof}
Avec le lemme \ref{rema-isocDagDagDense} (et la remarque \ref{rema-isocDagDense}), on se ramène au cas où $X$ est intègre
et $Y$ est dense dans $X$.
Lorsque $X$ est lisse, cette proposition a été établie dans \cite{caro-construction}.
Pour se ramener au cas où $X$ est lisse, comme d'habitude, on utilise le théorème de pleine fidélité \ref{a+|U-plfid-prop} : 
d'après le théorème de désingularisation de de Jong, il existe un diviseur $\widetilde{T}$ contenant $T$
et un diagramme de la forme \ref{formel-II-diag} satisfaisant aux conditions requises de \ref{formel-II} et
tel que de plus $X ^{(0)}$ soit lisse. On reprend les notations correspondantes. 
Avec le théorème de bidualité, on déduit du théorème \ref{coro-theo-!=+}
l'isomorphisme de commutation 
$a ^{+}\circ \DD _{T}\circ \sp _{X \hookrightarrow \PP, T,+} (E ) \riso \DD _{T ^{(0)}}\circ a ^{+}\circ \sp _{X \hookrightarrow \PP, T,+} (E )$.
Or, d'après \ref{sp-comm-a+-fini-etale-iso},
$a ^{+}\circ \sp _{X \hookrightarrow \PP, T,+}(E) 
\riso 
\sp _{X ^{(0)}\hookrightarrow \PP^{(0)}, T^{(0)},+}\circ a ^{*}(E)$.
D'où: 
\begin{equation}
\label{sp+-comm-dual-iso1}
a ^{+}\circ \DD _{T}\circ \sp _{X \hookrightarrow \PP, T,+} (E ) 
\riso 
\DD _{T ^{(0)}}\circ \sp _{X ^{(0)}\hookrightarrow \PP^{(0)}, T^{(0)},+}\circ a ^{*}(E).
\end{equation}

De même, il découle de \ref{sp-comm-a+-fini-etale-iso} (et de la commutation du dual à $a ^{*}$) l'isomorphisme 
$a ^{+}\circ \sp _{X \hookrightarrow \PP, T,+}(E ^{\vee}) 
\riso 
\sp _{X ^{(0)}\hookrightarrow \PP^{(0)}, T^{(0)},+} ((a ^{*}(E)) ^{\vee})$.
Or, comme $X ^{(0)}$ est lisse, le foncteur $\sp _{X ^{(0)}\hookrightarrow \PP^{(0)}, T^{(0)},+} $ commute au dual.
D'où : 
\begin{equation}
\label{sp+-comm-dual-iso2}
a ^{+}\circ \sp _{X \hookrightarrow \PP, T,+}(E ^{\vee}) 
\riso 
\DD _{T ^{(0)}}\circ \sp _{X ^{(0)}\hookrightarrow \PP^{(0)}, T^{(0)},+}\circ a ^{*}(E).
\end{equation}
Par \ref{sp+-comm-dual-iso1} et \ref{sp+-comm-dual-iso2}, on obtient l'isomorphisme canonique
$a ^{+}\circ \DD _{T}\circ \sp _{X \hookrightarrow \PP, T,+} (E )  \riso 
a ^{+}\circ \sp _{X \hookrightarrow \PP, T,+}(E ^{\vee})$.
Celui-ci, en dehors de $T ^{(0)}$ est induit via $b ^{+}$ par l'isomorphisme canonique
$\DD \circ \sp _{Y \hookrightarrow \U, +} (j ^{*}(E) )  \riso 
\sp _{Y \hookrightarrow \U, +}(j ^{*}(E) ^{\vee})$.
D'où le résultat par pleine fidélité de $(a ^{+}, |\U)$ (voir \ref{a+|U-plfid-prop}).

\end{proof}

\subsection{Propriétés de finitude des isocristaux partiellement surcohérents}

\begin{lemm}
\label{casliss631dev}
  Soient $f$ : $\PP' \rightarrow \PP $ un morphisme de $\V$-schémas formels séparés et lisses,
  $u ' $ : $X'\hookrightarrow P'$ une immersion fermée avec $X'$ intègre, 
   $T'$ un diviseur de $P'$ ne contenant pas $X'$ et tel que $Y':= X'\setminus T'$ soit lisse.
 Soit $\E' \in (F\text{-})\mathrm{Isoc} ^{\dag \dag} (\PP', T', X'/K)$.

\begin{enumerate}
\item
\label{casliss631dev-i)}
  Il existe alors un diagramme commutatif de la forme
  \begin{equation}
  \label{diag-casliss631dev}
  \xymatrix {
  {X''} \ar[r] ^-{u''} \ar[d] _{a '} & {\P ^N _{P'}} \ar[r] \ar[d] &
  {\widehat{\P} ^N _{\PP'}} \ar[r] ^{\widehat{\P} ^N _f} \ar[d] ^-{q'} & {\widehat{\P} ^N _{\PP}} \ar[d]^-{q} \\
  {X'} \ar[r] ^-{u'} & {P'} \ar[r] & {\PP'} \ar[r] ^-f & {\PP,}
   }
\end{equation}
  où $X''$ est lisse sur $k$, $q$ et $q'$ sont les projections canoniques, $u''$ est une immersion fermée, $a  ^{\prime -1} (T'\cap X')$ est un diviseur à croisements normaux strict de $X''$,
  $a '$ est propre, surjectif, génériquement fini et étale. 
De plus,
$\E'$ est alors un facteur direct de $q ' _+  \circ a ^{\prime !} (\E ')$.

\item 
\label{casliss631dev-iv)} 
Supposons $f  \circ u' $ propre et qu'il existe un diviseur $T$ de $P$ tel que $T' = f ^{-1}(T)$.
Alors $f _+ (\E') ,\DD _{T} \circ f _{+} (\E') \in (F\text{-})D ^\mathrm{b} _\mathrm{surcoh} (\D ^\dag _{\PP} (\hdag T) _{\Q})$.

\item 
\label{casliss631dev-v)} 
Supposons que $f  \circ u' $ soit propre, 
qu'il existe un diviseur $T$ de $P$ tel que $T' = f ^{-1}(T)$
et que le morphisme induit $Y' \to P$ soit une immersion.
Alors $f _+ (\E') \in (F\text{-})\mathrm{Isoc} ^{\dag \dag} (\PP,T, X/K)$, 
où $X$ est l'adhérence de $Y'$ dans $P$.
\end{enumerate}
\end{lemm}

\begin{proof}
  Grâce au théorème de désingularisation de de Jong, il existe une variété quasi-projective lisse $X''$,
un morphisme projective génériquement fini et étale $a '$ : $X'' \rightarrow X'$
tel que $(a '  ) ^{-1} (T' \cap X')$
soit un diviseur à croisements normaux strict de $X''$.
Il existe donc une immersion de la forme $X'' \hookrightarrow \P ^N _{k}$.
Comme $a '$ est propre, l'immersion induite $X '' \hookrightarrow \P ^N _{X'}$
est fermée. D'où l'immersion fermée canonique :
$X '' \hookrightarrow \widehat{\P} ^N _{\PP'}$,
ce qui donne l'existence du diagramme commutatif \ref{diag-casliss631dev} satisfaisant aux propriétés requises.

Établissons à présent \ref{casliss631dev-i)}) (on  pourra comparer avec la preuve de \cite[6.3.1]{caro_devissge_surcoh}).
 Posons $T'':= q ^{\prime -1} (T')$.
Comme $q '$ est propre, le foncteur $q ' _{+} $ est adjoint à gauche de $q ^{\prime !}$ (dans la catégorie des complexes surcohérents)
et on bénéficie de l'isomorphisme de dualité relative
$q ' _{+} \circ \DD _{T''}  \riso \DD _{T'} \circ q ' _{+} $.
On en déduit les morphismes canoniques :
$q ' _{+} \circ \R \underline{\Gamma} _{X''} ^\dag \circ q ^{\prime !} (\E')
\to \E'$
et $\E' \to q ' _{+} \circ \DD _{T''}\circ \R \underline{\Gamma} _{X''} ^\dag \circ q ^{\prime !} \circ \DD _{T'}(\E')$.
D'après nos notations, $a ^{\prime !} (\E ') =  \R \underline{\Gamma} _{X''} ^\dag \circ q ^{\prime !} (\E')$ 
et $a ^{\prime +} (\E ')=\DD _{T''} \circ \R \underline{\Gamma} _{X''} ^\dag \circ q ^{\prime !} \circ \DD _{T'}(\E')$.
Or, d'après \ref{coro-theo-!=+}, $a ^{\prime !} (\E ') \riso a ^{\prime +} (\E ')$.
On obtient ainsi la suite 
$\E ' \to q ' _+  \circ a ^{\prime !} (\E ') \to \E'$.
Pour vérifier que la composée est un isomorphisme, par fidélité du foncteur restriction à un ouvert dense, on se ramène
au cas où $a '$ est fini et étale, ce qui est immédiat.

Vérifions maintenant \ref{casliss631dev-iv)}) : 
notons $\tilde{f} :=\widehat{\P} ^N _f$, $\widetilde{\PP}:=\widehat{\P} ^N _{\PP}$
$\widetilde{T}:=q ^{-1} (T)$, $\E'':=a ^{\prime !} (\E ') $.
Comme $f  \circ u' $ est propre, $\P ^N _{f }\circ u'' $
  est alors une immersion fermée. 
De plus, comme $T'' \cap X''=\widetilde{T}\cap X''$,
via 
\ref{coh-PXTindtP}.\ref{coh-PXTindtP-ii}
on vérifie $\tilde{f} _+  (\E'')\in (F\text{-})\mathrm{Isoc} ^{\dag \dag} (\widetilde{\PP}, \widetilde{T}, X''/K)$. 
Comme $q $ est propre, $q _+$ préserve la surcohérence et commute aux foncteurs duaux. 
Or, par transitivité de l'image directe 
$f _+  \circ q '_+ (\E'') \riso q _+ \circ \tilde{f} _+  (\E'')$.
Il en résulte que $f _+  \circ q '_+ (\E'')$ et $\DD _{T} \circ f _+  \circ q '_+ (\E'')$ sont dans $(F\text{-})D ^\mathrm{b} _\mathrm{surcoh} (\D ^\dag _{\PP} (\hdag T) _{\Q})$. Comme $\E'$ est un facteur direct de $q ' _+(\E'')$, alors 
$f _+ (\E')$ (resp. $\DD _{T} \circ f _+ (\E')$) est un facteur direct de
$f _+  \circ q '_+ (\E'')$ (resp. $\DD _{T} \circ f _+  \circ q '_+ (\E'')$).
D'où $f _+ (\E'), \DD _{T} \circ f _+ (\E') \in (F\text{-})D ^\mathrm{b} _\mathrm{surcoh} (\D ^\dag _{\PP} (\hdag T) _{\Q})$.

Traitons à présent \ref{casliss631dev-v)}. 
D'après \ref{casliss631dev-iv)}, on sait déjà que $f _+ (\E')\in (F\text{-})D ^\mathrm{b} _\mathrm{surcoh} (\D ^\dag _{\PP} (\hdag T) _{\Q})$.
Via \cite[4.3.12]{Be1} (et aussi \ref{Isoc*=dagdag} et \ref{rema-isocDagDagDense}), 
il ne reste plus qu'à vérifier que $f _+ (\E')|\U \in \mathrm{Isoc} ^{\dag \dag} (\U,Y'/K)$,
où $\U := \PP \setminus T$.
Cela résulte du cas de la compactification lisse déjà traité dans 
\ref{coh-PXTindtP}.
D'où le résultat.
\end{proof}

\begin{prop}
\label{isoctilstaiinvlis}
Soient $(\PP, T,X)$ un triplet lisse en dehors du diviseur 
et $\E \in (F\text{-})\mathrm{Isoc} ^{\dag \dag}( \PP , T , X /K)$. Alors $\E$ vérifie la propriété $P _{\PP, T}$ (voir la définition \cite[6.1.1]{caro_devissge_surcoh}).
\end{prop}

\begin{proof}
Via \ref{casliss631dev}.\ref{casliss631dev-i)}, on se ramène au cas où $X$ est lisse. Dans ce cas, cela découle de \cite[6.1.4]{caro_devissge_surcoh}.
\end{proof}

Donnons une application de \ref{isoctilstaiinvlis} :
\begin{prop}
\label{Recol-cont}
Soient $\PP$ un $\V$-schéma formel séparé et lisse, $T \subset T'$ deux diviseurs de $P$, $X$ un sous-schéma fermé de $P$. 
On pose $\U:= \PP \setminus T$, $\U':= \PP \setminus T'$, $Y:= X \setminus T$, $Y':= X \setminus T'$. 
On suppose de plus $Y$ lisse et $Y'$ dense dans $Y$.
Soit $\E ' \in \mathrm{Isoc} ^{*} (\PP, T ', X/K) $. 

S'il existe $\E \in \mathrm{Isoc} ^{*} (\PP, T, X/K) $ tel que 
$\E '  \riso (\hdag T ' ) (\E)$, alors 
$\E \riso \mathrm{Im} \left ( \DD _{T}\circ \DD _{T '} (\E ')\to \E ' \right )$. 
\end{prop}

\begin{proof}
D'après \ref{isoctilstaiinvlis}, $\E$ et $\DD _{T} (\E)$ vérifient la propriété $P _{\PP, T}$.
Il en résulte que $\E ' \riso (\hdag T') (\E)$ et
$ \DD _{T}\circ \DD _{T '} (\E ') \riso \DD _{T} \circ \DD _{T'} \circ (\hdag T') (\E)\riso\DD _{T} \circ  (\hdag T') \circ \DD _{T} (\E )$ sont $\D ^{\dag} _{\PP} (\hdag T) _{\Q}$-surcohérents.
On en déduit, en posant $\FF:= \mathrm{Im} \left ( \DD _{T}\circ \DD _{T '} (\E ')\to \E ' \right )$, 
que $\FF$ est $\D ^{\dag} _{\PP} (\hdag T) _{\Q}$-surcohérent.
Or, d'après le cas lisse \ref{lemm2Th2-cas lisse}, 
$\FF |\U \in \mathrm{Isoc} ^{*} (\U, Y/K) $. En d'autres termes, 
$\FF |\U$ est dans l'image essentielle du foncteur 
$\sp _{Y \hookrightarrow \U, +}$.
Ainsi $\FF \in \mathrm{Isoc} ^{*} (\PP, T, X/K) $.

L'inclusion canonique $\FF \subset \E'$ induit l'isomorphisme 
canonique $(\hdag T' )(\FF)  \riso \E'$ (en effet, via \cite[4.3.12]{Be1}, il suffit de le vérifier en dehors de $T'$). 
Par pleine fidélité du foncteur $(\hdag T' )$ (voir le théorème \ref{dagT'pl-fid}), il en résulte que $\FF \riso \E$.
\end{proof}

\begin{vide}
\label{contagiosité}
Soient $\PP$ un $\V$-schéma formel séparé et lisse, $T \subset T'$ deux diviseurs de $P$, $X$ un sous-schéma fermé de $P$. 
On pose $\U:= \PP \setminus T$, $\U':= \PP \setminus T'$, $Y:= X \setminus T$, $Y':= X \setminus T'$, 
$j\,:\, Y \subset X$ et $ j' \,:\, Y' \subset X$ les immersions ouvertes canoniques. 
On suppose de plus $Y$ lisse et $Y'$ dense dans $Y$.

$\bullet$ D'après le théorème de contagiosité de Kedlaya \cite[5.3.7]{kedlaya-semistableI}, 
le foncteur du bas du carré
\begin{equation}
\label{fonct-contagiosité}
\xymatrix @C=1,5 cm{
{\mathrm{Isoc} ^{*} (\PP, T , X/K)} 
\ar[r] ^-{( (\hdag T'), |\U)}
& 
{\mathrm{Isoc} ^{*} (\PP, T ', X/K) \times _{\mathrm{Isoc} ^{*} (\U, T ', Y/K)} \mathrm{Isoc} ^{*} (\U, Y/K)} 
\\ 
{\mathrm{Isoc} ^{\dag} (Y,X/K) } 
\ar[r] ^-{( j ^{\prime \dag}, j ^{*})}
\ar[u] ^-{\sp _{X \hookrightarrow \PP, T,+}}
& 
{\mathrm{Isoc} ^{\dag} (Y',X/K) \times _{\mathrm{Isoc} ^{\dag} (Y',Y/K)} \mathrm{Isoc} ^{\dag} (Y,Y/K) } 
\ar[u] _-{(\sp _{X \hookrightarrow \PP, T',+}, \sp _{Y \hookrightarrow \U,+})}
 }
\end{equation}
est une équivalence de catégories. 
Comme ce carré est commutatif à isomorphisme canonique près, comme 
d'après le théorème \ref{diag-eqcat-dagdag-ii} les foncteurs verticaux sont des équivalences de catégories, 
il en est de même du foncteur $( (\hdag T'), |\U)$ du haut.

$\bullet$ Grâce à \ref{Recol-cont}, on construit explicitement le foncteur $\mathcal{R} ecol$ quasi-inverse de $( (\hdag T'), |\U)$ en posant,
pour tout objet $( \E', \FF _{\U}, \rho)\in \mathrm{Isoc} ^{*} (\PP, T ', X/K) \times _{\mathrm{Isoc} ^{*} (\U, T ', Y/K)} \mathrm{Isoc} ^{*} (\U, Y/K)$,
$$\mathcal{R} ecol ( \E', \FF _{\U}, \rho) := 
\mathrm{Im} \left ( \DD _{T}\circ \DD _{T '} (\E ')\to \E ' \right ).$$

$\bullet$ En fait, pour tout $( \E', \FF _{\U}, \rho)$ tels que
$\E' $, 
$\DD _{T} \circ \DD _{T'} (\E') $ soient $\D ^{\dag} _{\PP} (\hdag T) _{\Q}$-surcohérents, on vérifie directement (i.e. sans utiliser le théorème de contagiosité de Kedlaya) que $\mathcal{R} ecol ( \E', \FF _{\U}, \rho)\in \mathrm{Isoc} ^{*} (\PP, T , X/K)$ et 
$( \E', \FF _{\U}, \rho) \riso ( (\hdag T'), |\U) \circ \mathcal{R} ecol ( \E', \FF _{\U}, \rho)$.
En effet, avec ces hypothèses, $\FF := \mathcal{R} ecol ( \E', \FF _{\U}, \rho)$ est $\D ^{\dag} _{\PP} (\hdag T) _{\Q}$-surcohérent.
Or, d'après le cas lisse \ref{lemm2Th2-cas lisse}, 
$\FF |\U \in \mathrm{Isoc} ^{*} (\U, Y/K) $. D'où $\FF \in \mathrm{Isoc} ^{*} (\PP, T , X/K)$.
Enfin, l'isomorphisme $( (\hdag T'), |\U) (\FF) \riso ( \E', \FF _{\U}, \rho)$ est immédiat par construction.
\end{vide}

\subsection{Indépendance}

\begin{lemm}
\label{gen-coh-PXTindtP}
Soit $\theta =(f, a,Id) \,:\, (\PP', T',X',Y)\to (\PP, T,X,Y)$ un morphisme de quadruplets lisses en dehors du diviseur
tel que $Y$ soit dense dans $X'$, $a$ soit propre et $f ^{-1} (T)$ soit un diviseur de $P'$ (e.g., si $f$ est lisse).

\begin{enumerate}

\item \label{gen-coh-PXTindtP-i} 
Soient $\E \in (F\text{-})\mathrm{Isoc} ^{\dag \dag} (\PP, T, X/K)$,
$\E '\in (F\text{-})\mathrm{Isoc} ^{\dag \dag} (\PP', T', X'/K)$). 
Pour tout $l \in \Z\setminus \{0\}$, 
$$\mathcal{H} ^l (\R \underline{\Gamma} ^\dag _{X'} f ^! (\E) ) =0,
\hspace{2cm}
\mathcal{H} ^l (f_+(\E')) =0.$$

\item \label{gen-coh-PXTindtP-ii} Les foncteurs $\R \underline{\Gamma} ^\dag _{X'}  f ^! $ et
$f _+$ induisent des équivalences quasi-inverses 
entre les catégories $(F\text{-})\mathrm{Isoc} ^{\dag \dag} (\PP, T, X/K)$ et $(F\text{-})\mathrm{Isoc} ^{\dag \dag} (\PP', T', X'/K)$.

\end{enumerate}
\end{lemm}

\begin{proof}
Soient  
$\E \in (F\text{-})\mathrm{Isoc} ^{\dag \dag} (\PP, T, X/K)$,
$\E' \in (F\text{-})\mathrm{Isoc} ^{\dag \dag} (\PP', T', X'/K)$.
Notons $\U:= \PP \setminus T$ et $\U':= \PP '\setminus T'$.

I) Vérifions d'abord le lemme lorsque $f$ et $a$ sont l'identité, i.e. $\PP'=\PP$ et $X'=X$. 
Par symétrie, il suffit d'établir l'inclusion 
$(F\text{-})\mathrm{Isoc} ^{\dag \dag} (\PP, T, X/K) \subset (F\text{-})\mathrm{Isoc} ^{\dag \dag} (\PP, T', X/K)$.
On procède de manière analogue à \cite[6.3.4]{caro_devissge_surcoh} :
choisissons un diagramme vérifiant les conditions de \ref{diag-casliss631dev} (avec $f =Id$) pour le diviseur $T'$. 
Remarquons que comme $T' \cap X = T \cap X$, alors $a ^{-1} (T \cap X)$ est aussi un diviseur de $X''$.
Posons $T'' _{1} := q ^{\prime-1} (T)$, 
$T'' _{2}:= q ^{\prime} (T')$,
$\E'' := \R \underline{\Gamma} _{X''} ^\dag \circ q ^{\prime !} (\E)$.
D'après \ref{casliss631dev}.\ref{casliss631dev-i)},  
$\E'' \in (F\text{-})\mathrm{Isoc} ^{\dag \dag} (\widehat{\P} ^N _{\PP}, T''_{1}, X''/K)$ 
et
$\E$ est un facteur direct de $q ' _+(\E'' )$.
Or, d'après le cas de la compactification partielle lisse (voir \ref{coh-PXTindtP}), 
comme $X'' \setminus T'' _{1}= X'' \setminus T'' _{2}$,
on obtient l'égalité
$(F\text{-})\mathrm{Isoc} ^{\dag \dag} (\widehat{\P} ^N _{\PP}, T''_{1}, X''/K)
=
(F\text{-})\mathrm{Isoc} ^{\dag \dag} (\widehat{\P} ^N _{\PP}, T''_{2}, X''/K)$.
En particulier, $\E''$ est $\D ^{\dag} _{\widehat{\P} ^N _{\PP}} (\hdag T''_{2}) _{\Q}$-surcohérent.
Par préservation de la surcohérence par image directe par un morphisme propre, il en résulte que 
$q ' _+ (\E'' )$ est $\D ^{\dag} _{\PP} (\hdag T') _{\Q}$-surcohérent.
Comme $\E$ est un facteur direct de $q ' _+(\E'' )$, 
celui-ci est aussi $\D ^{\dag} _{\PP} (\hdag T') _{\Q}$-surcohérent.
On en conclut la preuve grâce au corollaire \ref{Isoc*=dagdag}. 

II) Traitons à présent le cas général. 
Comme $a$ est propre, l'immersion ouverte $Y \hookrightarrow X '\setminus f ^{-1} (T) $ est aussi fermée.
Comme $Y$ est dense dans $X'$, on en déduit $X '\setminus f ^{-1} (T) =Y$. 
Ainsi, $f ^{-1} (T)$ est un diviseur de $P'$ tel que $X '\setminus f ^{-1} (T) =Y$.
d'après l'étape I), on se ramène à supposer $T'= f ^{-1} (T)$.

II.1) Grâce à \ref{casliss631dev}.\ref{casliss631dev-v)} (resp. grâce à \ref{stabIsoc*inv}), on obtient
l'assertion \ref{gen-coh-PXTindtP-i} pour l'image directe (resp. l'autre foncteur) et 
$f _+(\E')\in (F\text{-})\mathrm{Isoc} ^{\dag \dag} (\PP, T, X/K)$ (resp. $\R \underline{\Gamma} ^\dag _{X'}  f ^! (\E) 
\in  (F\text{-})\mathrm{Isoc} ^{\dag \dag} (\PP', T', X'/K)$).

II.2) Vérifions à présent que l'on dispose de l'isomorphisme canonique :
$\R \underline{\Gamma} _{X'}  ^\dag \circ f ^{!} \circ f _+ (\E ') \riso \E'$.
D'après le lemme \ref{casliss631dev}.\ref{casliss631dev-i)},
il existe un diagramme de la forme \ref{diag-casliss631dev} et satisfaisant aux conditions requises de ce lemme.
Avec ces notations, comme le foncteur $( a ^{\prime !}, |\U ')$ est pleinement fidèle (voir \ref{a+|U-plfid-prop}), 
il suffit de construire deux isomorphismes compatibles
$ a ^{\prime !} \circ \R \underline{\Gamma} _{X'}  ^\dag \circ f ^{!} \circ f _+ (\E ') \riso  a ^{\prime !} (\E')$
et
$|\U ' \circ \R \underline{\Gamma} _{X'}  ^\dag \circ f ^{!} \circ f _+ (\E ') \riso \E' |\U '$.
Le deuxième isomorphisme résulte aussitôt du cas de la compactification partielle lisse (voir \ref{coh-PXTindtP}).
Établissons à présent le premier. 
On rappelle que 
$a ^{\prime !} = \R \underline{\Gamma} _{X ^{\prime \prime}} ^\dag \circ q ^{\prime!} $.
Notons $\tilde{f} :=\widehat{\P} ^N _f\,:\, \widehat{\P} ^{N} _{\PP'} \to \widehat{\P} ^{N} _{\PP}$.
D'après  \ref{casliss631dev}.\ref{casliss631dev-iv)}, $f _{+} (\E ') $ est cohérent. 
D'après la proposition \cite[3.1.8]{caro_surcoherent} (il y est implicite que les foncteurs images inverses extraordinaires 
sont ceux définis en \cite[4.3.2.2]{Beintro2}),
on dispose de l'isomorphisme 
$\widetilde{f} _{+} \circ q ^{\prime!} (\E ') \riso q ^{!}\circ f _{+} (\E ') $.
On en déduit les isomorphismes :
\begin{equation}
\notag 
a ^{\prime !} \circ \R \underline{\Gamma} _{X'}  ^\dag \circ f ^{!} \circ f _+ (\E ')
\riso
\R \underline{\Gamma} _{X ^{\prime \prime}} ^\dag \circ \widetilde{f} ^{!}  \circ q ^{!} \circ f _+ (\E ')
\riso
\R \underline{\Gamma} _{X ^{\prime \prime}} ^\dag \circ \widetilde{f} ^{!}  \circ\widetilde{f} _{+} \circ q ^{\prime!} (\E ')
\riso
\R \underline{\Gamma} _{X ^{\prime \prime}} ^\dag \circ \widetilde{f} ^{!}  \circ\widetilde{f} _{+} \circ a ^{\prime !} (\E ')
\riso
a ^{\prime !} (\E '),
\end{equation}
le dernier isomorphisme, comme $X ^{\prime \prime}$ est lisse, résultant du cas de la compactification partielle lisse \ref{coh-PXTindtP}. 

II.3) De même, on établit l'isomorphisme canonique et $f _+  \circ \R \underline{\Gamma} _{X'}  ^\dag \circ f ^{!} (\E) \riso \E$.
D'où le résultat. 
\end{proof}

\begin{rema}
\label{rema-pashdag}
D'après \ref{gen-coh-PXTindtP} et avec ses notations, 
le morphisme canonique 
$\R \underline{\Gamma} ^\dag _{X'} f ^! (\E) \to 
(\hdag T') \circ \R \underline{\Gamma} ^\dag _{X'} f ^! (\E) $ est un isomorphisme.
Ainsi, il est inutile d'indiquer le foncteur $(\hdag T') $.
\end{rema}

\begin{lemm}
\label{p1!=p2!}
Soit $(\PP, T,X,Y)$ un quadruplet lisse en dehors du diviseur (voir la définition \ref{defi-triplet}).
Notons $p _{1}\,:\, \PP \times  \PP \to \PP$ et $p _{2}\,:\, \PP \times  \PP \to \PP$
les projections respectives à gauche et à droite.

Pour tout $\E \in (F\text{-})\mathrm{Isoc} ^{\dag \dag} (\PP, T, X/K)$, 
on dispose d'une équivalence canonique de catégories: 
\begin{equation}
\label{p1!=p2!-iso}
\R \underline{\Gamma} ^\dag _X p _{1} ^! (\E)
\riso
\R \underline{\Gamma} ^\dag _X p _{2} ^! (\E).
\end{equation}

\end{lemm}

\begin{proof}
Dans un premier temps, supposons $X$ lisse.
Notons $\U := \PP \setminus T$ et $H:= (P \times P ) \setminus (U \times U)$.
Comme $X$ est lisse, d'après \cite{caro-construction}, les foncteurs de la forme 
$\sp _{X \hookrightarrow \PP, T,+}$ commutent aux images inverses et extension. 
Avec la remarque \ref{rema-pashdag}, on obtient les isomorphismes canoniques
$$\R \underline{\Gamma} ^\dag _X p _{1} ^!  \circ \sp _{X \hookrightarrow \PP, T,+} 
\riso 
\sp _{X \hookrightarrow \PP \times \PP, H,+} 
\riso 
\R \underline{\Gamma} ^\dag _X p _{2} ^!  \circ \sp _{X \hookrightarrow \PP, T,+} 
.$$
Comme le foncteur $\sp _{X \hookrightarrow \PP, T,+} $ est essentiellement surjectif sur
$(F\text{-})\mathrm{Isoc} ^{\dag \dag} (\PP, T, X/K)$, on obtient 
l'isomorphisme canonique \ref{p1!=p2!-iso} lorsque $X$ est lisse.

Passons à présent au cas général. 
De manière analogue à l'étape II.2) de la preuve du lemme \ref{gen-coh-PXTindtP},
on se ramène au cas où $X$ est lisse grâce au théorème de désingularisation de de Jong 
et au théorème de pleine fidélité \ref{a+|U-plfid-prop}.
\end{proof}

\begin{defi}
\label{defi-(d)plong}
Un {\og couple $(Y,X)$ de $k$-variétés $d$-plongeable \fg} est la donnée 
d'une $k$-variété $X$, d'un ouvert $Y$ de $X$ tels qu'il existe
un $\V$-schéma formel $\PP $ séparé et lisse,
un diviseur $T$ de $P$ et une immersion fermée $X \hookrightarrow \PP$
vérifiant $Y = X \setminus T$.
Un morphisme de couples $(Y',X')\to (Y,X)$ de $k$-variétés $d$-plongeables est un morphisme de variétés
$a \,:\,X'\to X$ tel que $a (Y') \subset Y$.
\end{defi}

\begin{prop}
[et Définition]
  \label{defi-part-surcoh}
Soient $X$ une $k$-variété et $Y$ un ouvert lisse de $X$ tels que $(Y,X)$ soit $d$-plongeable (voir \ref{defi-(d)plong}).
Choisissons $\PP $ un $\V$-schéma formel séparé et lisse,
$T$ un diviseur de $P$ tels qu'il existe une immersion fermée $X \hookrightarrow \PP$
induisant l'égalité $Y = X \setminus T$. 

\begin{itemize}
\item La catégorie $(F\text{-})\mathrm{Isoc} ^{\dag \dag} (\PP, T, X/K)$ ne dépend, à équivalence canonique de catégories près, ni du choix de l'immersion fermée
$X \hookrightarrow \PP$ et ni de celui du diviseur $T$ tel que $Y = X \setminus T$ (mais seulement du couple $(Y,X)$ de $k$-variétés $d$-plongeable).

\item On note alors
$(F\text{-})\mathrm{Isoc}^{\dag \dag} (Y,X/K)$
sans ambiguïté à la place de $(F\text{-})\mathrm{Isoc} ^{\dag \dag} (\PP, T, X/K)$.
Ses objets sont les {\og $(F\text{-})$isocristaux partiellement surcohérents sur $(Y,X)$\fg} 
ou simplement {\og $(F\text{-})$isocristaux surcohérents sur $(Y,X)$\fg}.

\end{itemize}

\end{prop}

\begin{proof}
D'après lemme \ref{rema-isocDagDagDense}, on peut supposer que $Y$ est dense dans $X$.
Soit un deuxième choix $(X \hookrightarrow \PP',T')$.
Posons $\PP '':= \PP \times \PP'$. Soient $q\,:\, \PP'' \to \PP$, $q'\,:\, \PP'' \to \PP'$ les projections canoniques
et $T'':= q ^{-1} (T)$.
D'après \ref{gen-coh-PXTindtP},
les foncteurs $q _+$ et $\R \underline{\Gamma} _X ^\dag q ^! $
(resp. $q ' _+$ et $\R \underline{\Gamma} _X ^\dag q ^{\prime !} $)
induisent alors des équivalences quasi-inverses entre les catégories
$(F\text{-})\mathrm{Isoc} ^{\dag \dag} (\PP, T, X/K)$
et
$(F\text{-})\mathrm{Isoc} ^{\dag \dag} (\PP'', T'', X/K)$
(resp. $(F\text{-})\mathrm{Isoc} ^{\dag \dag} (\PP', T', X/K)$
et
$(F\text{-})\mathrm{Isoc} ^{\dag \dag} (\PP'', T'', X/K)$).
D'où le résultat.
\end{proof}

\begin{vide}
Soit $(Y,X)$ un couple de $k$-variétés $d$-plongeable.
D'après lemme \ref{rema-isocDagDagDense}, si $\overline{Y}$ désigne l'adhérence de $Y$ dans $X$ alors
$(F\text{-})\mathrm{Isoc}^{\dag \dag} (Y,\overline{Y}/K)=(F\text{-})\mathrm{Isoc}^{\dag \dag} (Y,X/K)$.

Si $\PP $ est un $\V$-schéma formel séparé et lisse,
$T$ est un diviseur de $P$ tels qu'il existe une immersion fermée $X \hookrightarrow \PP$
induisant l'égalité $Y = X \setminus T$, l'équivalence de catégories de \ref{prop-donnederecol-dag}
\begin{equation}
\label{cong-Xlisse}
  \sp _{X \hookrightarrow \PP, T,+}\, :\,
(F\text{-})\mathrm{Isoc}^{\dag} (Y,X/K) \cong (F\text{-})\mathrm{Isoc}^{\dag \dag} (Y,X/K),
\end{equation}
ne dépend canoniquement pas des choix faits et sera simplement notée $\sp _{(Y,X),+}$.
\end{vide}

\begin{prop}
[et Définition]
  \label{defi-part-surcoh-propre}
  Soit $Y$ une $k$-variété lisse.
On suppose qu'il existe une $k$-variété propre $X$ telle que $(Y,X)$ soit $d$-plongeable (voir \ref{defi-(d)plong}).
Choisissons $\PP $ un $\V$-schéma formel séparé et lisse,
$T$ un diviseur de $P$ tels qu'il existe une immersion fermée $X \hookrightarrow \PP$
induisant l'égalité $Y = X \setminus T$. 

\begin{itemize}
\item La catégorie $(F\text{-})\mathrm{Isoc} ^{\dag \dag} (\PP, T, X/K)$ ne dépend, à équivalence canonique de catégories près, ni du choix de la compactification propre $X$, ni du choix de l'immersion fermée
$X \hookrightarrow \PP$ et ni de celui du diviseur $T$ tel que $Y = X \setminus T$. 
\item On note alors
$(F\text{-})\mathrm{Isoc}^{\dag \dag} (Y/K)$
à la place de $(F\text{-})\mathrm{Isoc} ^{\dag \dag} (\PP, T, X/K)$.
Ses objets sont les {\og $(F\text{-})$isocristaux surcohérents sur $Y$\fg}.
\end{itemize}

\end{prop}

\begin{proof}
Soit un deuxième choix $(X '\hookrightarrow \PP',T')$.
Il suffit de considérer l'adhérence $X''$ de $Y$ dans $X \times X'$, l'immersion fermée induite
$X'' \hookrightarrow \PP \times \PP'$. En effet, comme les projection 
$X'' \to X$ et $X'' \to X'$ sont propres, on conclut la preuve grâce 
à \ref{gen-coh-PXTindtP}.
\end{proof}

\begin{nota}
\label{coro-isoctilstaiinvlis}

Soit $\theta = (f,a,b)\,:\, (\PP', T',X',Y')\to (\PP, T,X,Y)$ un morphisme de quadruplets lisses en dehors du diviseur.
De manière analogue à \cite[6.4.3]{caro_devissge_surcoh}, on déduit de \ref{defi-part-surcoh}  que les factorisations 
\begin{align}
\notag
\theta ^{!}\,:=\,(\hdag T') \circ \R \underline{\Gamma} ^\dag _{X'} \circ f ^! [-d _{X'/X}] \,:\,&(F\text{-})\mathrm{Isoc} ^{**} (\PP, T, X/K) \to (F\text{-})\mathrm{Isoc} ^{**} (\PP', T', X'/K),\\
\notag
\theta ^{+}\,:=\,\DD _{T'}\circ \R \underline{\Gamma} ^\dag _{X'}\circ (\hdag T') \circ f ^! [-d _{X'/X}] \circ \DD _{T} \,:\,&
(F\text{-})\mathrm{Isoc} ^{*} (\PP, T, X/K)\to (F\text{-})\mathrm{Isoc} ^{*} (\PP', T', X'/K)),
\end{align}
de \ref{stabIsoc*inv}
ne dépendent que du morphisme 
de couples $(X',Y')\to (Y,X)$ de $k$-variétés $d$-plongeables induit par $a$. 
D'après \ref{coro-theo-!=+}, ceux-ci sont isomorphes. 
On le notera alors sans ambiguïté $a ^{*}$.
\end{nota}

\section{Applications}

\subsection{Descriptions explicites des équivalence de catégories entre isocristaux surconvergents et surcohérents}

\begin{nota}
\label{notaSP+}
Soient $P ^\dag$ un $\V$-schéma formel faible lisse (voir la définition dans \cite{meredith-weakformalschemes}) de fibre spéciale $P$, $T $ un diviseur de $P $, $U ^\dag$ l'ouvert de
$P ^\dag $ complémentaire de $T $, $j$ : $U ^\dag \hookrightarrow P ^\dag$ l'immersion
  ouverte et $v$ : $Y ^\dag \hookrightarrow U ^\dag$ une immersion fermée
de $\V$-schémas formels faibles.
On suppose en outre $Y ^\dag$ affine et lisse. 
On note $U$, $Y$ les fibres spéciales respectives de $U ^\dag$, $Y ^\dag$ et $X $ l'adhérence schématique de
$Y $ dans $P $. 
\end{nota}

\begin{vide}
[Rappels]
Avec les notations de \ref{notaSP+}
Comme $Y ^{\dag}$ est affine et lisse, d'après \cite[5.1.1]{caro_devissge_surcoh}, 
on dispose d'un foncteur canonique pleinement fidèle noté $\sp _{*}$ de la catégorie 
$\mathrm{Isoc} ^{\dag} (Y/K)$ des isocristaux surconvergents sur $Y$ dans celle des 
$\smash{\D}  _{Y ^\dag ,\Q}$-modules globalement de présentation,
$\O _{Y ^\dag ,\Q}$-cohérents.

D'après \cite[5.1.3]{caro_devissge_surcoh}, on définit alors un $ \smash{\D} ^\dag _{\PP} (\hdag T ) _\Q$-module globalement de présentation finie à support dans $X $ à partir d'un isocristal surconvergent $E$ sur $Y$ en posant:
\begin{equation}
\label{sp+ff}
\sp _{Y ^{\dag}  \hookrightarrow U ^{\dag},T  , + }(E) :=\smash{\D} ^\dag _{\PP} (\hdag T ) _\Q \otimes _{j _* \smash{\D} _{U^\dag, \Q} } j_* v _+ (\sp _{*} (E)).
\end{equation}

\end{vide}

Nous aurons besoin du lemme \ref{lemm-pass-cvasurcv-isopart1} ci-dessous lors de la preuve du théorème \ref{sp+essent}.

\begin{lemm}
\label{lemm-pass-cvasurcv-isopart1}
Avec les notations de \ref{notaSP+}, 
soit $\E $ un $\smash{\D}  _{Y ^\dag ,\Q}$-module globalement de présentation, $\O _{Y ^\dag ,\Q}$-cohérent
 qui soit associé via le foncteur $\sp _{*}$ à un élément de $\mathrm{Isoc} ^{\dag} (Y/K)$ (voir \cite[5.1.1]{caro_devissge_surcoh}).
Posons $E := \Gamma (Y ^{\dag} , \E)$, $D _{Y ^{\dag}  ,\Q}: = \Gamma (Y ^{\dag} , \D _{Y ^{\dag}  ,\Q})$,
$D  _{U ^{\dag}\leftarrow Y ^{\dag}  ,\Q} =\Gamma (Y ^{\dag} , \D _{U ^{\dag}\leftarrow Y ^{\dag}  ,\Q})$,
$v _{+} (E) := D _{U ^{\dag}\leftarrow Y ^{\dag}  ,\Q} \otimes _{D _{Y ^{\dag}  ,\Q}} E$.

On bénéficie de l'isomorphisme canonique 
\begin{equation}
\label{pass-cvasurcv-isopart1}
  \sp _{Y ^\dag \hookrightarrow U ^\dag, T ,+} (E)
\liso
\D ^\dag _{\PP } (\hdag T ) _\Q
\otimes _{D _{U^\dag ,K}} v_+(E).
\end{equation}
\end{lemm}

\begin{proof}
On dispose, pour tout $D _{U^\dag }$-module $M$, d'un morphisme
$j_* \D _{U^\dag,\Q} \otimes _{D _{U^\dag ,K}} M \rightarrow
j_* (\D _{U^\dag,\Q} \otimes _{D _{U^\dag ,K}} M )$
fonctoriel en $M$.
Lorsque $M= D _{U^\dag , K}$, celui-ci est un isomorphisme.
En appliquant ces deux foncteurs à une présentation finie de $v_+(E)$,
on obtient un morphisme entre deux présentations finies. 
Par le lemme des cinq, il en résulte
l'isomorphisme
$j_* \D _{U^\dag,\Q} \otimes _{D _{U^\dag ,K}} v_+(E) \riso
j_* (\D _{U^\dag,\Q} \otimes _{D _{U^\dag ,K}}v_+(E) )$.
Or, par (passage de droite à gauche de) \cite[2.4.1]{caro_devissge_surcoh},
$v_+ (\E) \riso \D _{U^\dag} \otimes _{D _{U^\dag }} v_+(E)$.
Donc, $j_* \D _{U^\dag,\Q} \otimes _{D _{U^\dag ,K}} v_+(E)  \riso j_* v_+ (\E) $.
D'où le résultat par extension via $j_* \D _{U^\dag,\Q} \to \D ^\dag _{\PP } (\hdag T ) _\Q$.
\end{proof}

\begin{rema}
Rappelons que l'idée de travailler sur des schémas faiblement formels résultait de deux remarques : 
d'une part une théorie de $\D$-modules sur des schémas faiblement formels a été formulée par Mebkhout et Narvaez-Macarro dans 
\cite{Mebkhout-Narvaez-Macarro_D} et d'autre part Noot-Huyghe a établi dans \cite{huyghe-comparaison} des théorèmes de comparaison entre les deux approches (celle de Berthelot et celle de Mebkhout et Narvaez-Macarro).
Néanmoins, le défaut de la construction \ref{sp+ff} est (jusqu'à preuve du contraire) de ne pas avoir de généralisation
pour $\mathrm{Isoc} ^{\dag} (Y,X/K)$ à la place de $\mathrm{Isoc} ^{\dag} (Y/K)$
lorsque la compactification partielle $X$ de $Y$ n'est pas propre.
\end{rema}

\begin{theo}
\label{spYdag+to}
On garde les notations de \ref{notaSP+}.
Pour tout $E \in \mathrm{Isoc} ^{\dag }(Y/K)  $, on a 
$\sp _{Y ^{\dag}  \hookrightarrow U ^{\dag},T  , + }(E) \in  \mathrm{Isoc} ^{\dag \dag }(\PP, T,X/K)$.
\end{theo}

\begin{proof}
Soit $E \in \mathrm{Isoc} ^{\dag }(Y/K)$. Modulo quelques différences décrites ci-après, il s'agit de reprendre la première partie de la preuve de \cite[1.3.5]{caro-2006-surcoh-surcv}
  de la façon suivante : les autres points utilisés de la preuve (e.g. \cite[5.2.3, 5.2.6, 7.2.4]{caro_devissge_surcoh}) ne nécessitant ni de structure de Frobenius ni d'hypothèse de propreté de $P$, il reste à effectuer les deux changements suivants. 
On remplace la référence de \cite[6.3.1]{caro_devissge_surcoh} par \ref{casliss631dev}.
Enfin, avec les notations de la preuve de \cite[1.3.5]{caro-2006-surcoh-surcv}, il reste à vérifier sans utiliser le théorème de pleine fidélité de Kedlaya que $\widetilde{E} '$ provient d'un isocristal surconvergent sur $V'$. 
Pour cela, grâce au théorème de contagiosité de Kedlaya (voir \cite[5.3.7]{kedlaya-semistableI} ou \ref{contagiosité}),
il suffit de le vérifier en dehors de $T'$, ce qui nous ramène au cas immédiat de la compactification lisse.

Terminons par une remarque : 
 en fait, grâce à \ref{Isoc*=dagdag}, on peut simplifier la preuve de \cite[1.3.5]{caro-2006-surcoh-surcv}
 en supprimant dans \cite[1.3.5.1]{caro-2006-surcoh-surcv} la condition concernant le dual, i.e. concernant 
$\DD _{T}\R \underline{\Gamma} ^{\dag} _{Z} \E$.
\end{proof}

On déduit du théorème \ref{spYdag+to} le corollaire suivant :
\begin{coro}
\label{sp+dag-indt}
Soit $E \in \mathrm{Isoc} ^{\dag }(Y/K)  $. L'objet $\sp _{Y ^{\dag}  \hookrightarrow U ^{\dag},T  , + }(E) $ est canoniquement indépendant
du choix de l'immersion fermée $Y ^{\dag} \hookrightarrow U ^{\dag}$ relevant $Y \hookrightarrow U$.
\end{coro}

\begin{proof}
Soient $v _{1}, v _{2}\,:\, Y ^{\dag} \hookrightarrow U ^{\dag}$ deux relèvements de  
$Y \hookrightarrow U$. Notons $w=(v _{1}, v _{2}) \,:\, Y ^{\dag} \hookrightarrow U ^{\dag} \times  U ^{\dag}$, 
$H:= (P \times P ) \setminus (U\times U)$ le diviseur de $P\times P$,
$Z$ l'adhérence de $Y$ dans $P \times P$,
$p _{1}\,:\, \PP \times  \PP \to \PP$ et $p _{2}\,:\, \PP \times  \PP \to \PP$
les projections respectives à gauche et à droite.
D'après les propositions \cite[5.2.3, 5.2.6]{caro_devissge_surcoh} (pour utiliser la deuxième proposition, on a aussi besoin de la propriété de finitude de \ref{spYdag+to}), 
les foncteurs de la forme $\sp _{+}$ commute aux images inverses, i.e., 
on dispose des isomorphismes canoniques 
\begin{gather}
\notag
\sp _{Y ^{\dag} \underset{w}{\hookrightarrow} U ^{\dag} \times  U ^{\dag}, H, +} 
\riso 
(\hdag H) \circ \sp _{Y ^{\dag} \underset{w}{\hookrightarrow} U ^{\dag} \times  U ^{\dag}, p ^{-1} _{1} (T), +} 
\riso
(\hdag H) \circ \R \underline{\Gamma} ^{\dag} _{Z} p _{1} ^{!}
\circ 
\sp _{Y ^{\dag} \underset{v _{1}}{\hookrightarrow}   U ^{\dag}, T, +} 
\riso 
\R \underline{\Gamma} ^{\dag} _{Z} p _{1} ^{!}
\circ 
\sp _{Y ^{\dag} \underset{v _{1}}{\hookrightarrow}   U ^{\dag}, T, +} ,
\end{gather}
le dernier isomorphisme résultant de la remarque \ref{rema-pashdag}.
De même :
$\sp _{Y ^{\dag} \underset{w}{\hookrightarrow} U ^{\dag} \times  U ^{\dag}, H, +}
\riso \R \underline{\Gamma} ^{\dag} _{Z} p _{2} ^{!}
\circ 
\sp _{Y ^{\dag} \underset{v _{2}}{\hookrightarrow}   U ^{\dag}, T, +}.$
On en déduit alors, via \ref{p1!=p2!} et \ref{gen-coh-PXTindtP}, l'isomorphisme canonique
\begin{equation}
\label{sp+dag-indt-iso}
\tau _{v _{1}, v _{2}}\,:\,
\sp _{Y ^{\dag} \underset{v _{2}}{\hookrightarrow}   U ^{\dag}, T, +} (E)\riso \sp _{Y ^{\dag} \underset{v _{1}}{\hookrightarrow}   U ^{\dag}, T, +} (E).
\end{equation}
\end{proof}

\begin{rema}
\label{sp+dag-indt-rema}
Soient $v _{1}, v _{2}, v _{3}\,:\, Y ^{\dag} \hookrightarrow U ^{\dag}$ trois relèvements de  
$Y \hookrightarrow U$.
On vérifie, par construction que l'isomorphisme canonique de 
\ref{sp+dag-indt-iso}, la condition de cocycle 
$\tau _{v _{1}, v _{2}} \circ \tau _{v _{2}, v _{3}}=\tau _{v _{1}, v _{3}}$
(pour cela, on fait intervenir $U ^{\dag} \times  U ^{\dag} \times U ^{\dag}$).
De plus, par construction $\tau _{v _{1}, v _{1}}=Id$.
\end{rema}

\begin{vide}
\label{nota-spY->UdagT}
Soient $P ^\dag$ un $\V$-schéma formel faible séparé et lisse de fibre spéciale $P$, $T $ un diviseur de $P $, $U ^\dag$ l'ouvert de
$P ^\dag $ complémentaire de $T $, $U$ la fibre spéciale de $U ^\dag$, 
 $j$ : $U ^\dag \hookrightarrow P ^\dag$ l'immersion ouverte, 
  $v$ : $Y \hookrightarrow U $ une immersion fermée de $k$-schémas lisses.
On note $X $ l'adhérence schématique de $Y $ dans $P $.

On construit alors par recollement le foncteur 
\begin{equation}
\label{spY->UdagT}
\sp _{Y \hookrightarrow U ^{\dag}, T,+}\,:\, 
\mathrm{Isoc} ^\dag (Y/K) \to \mathrm{Isoc} ^{\dag\dag} (\PP, T, X /K)
\end{equation}
de la manière suivante :  
choisissons $(P ^{\dag} _{\alpha}) _{\alpha}$ un recouvrement d'ouverts affines de $P ^{\dag}$ tel que 
$T \cap P _{\alpha}$ soit défini par une équation. En particulier, 
$U ^{\dag} _{\alpha}:= P ^{\dag} _{\alpha} \cap U ^{\dag}$ 
est affine. On note 
$Y  _{\alpha}:= P _{\alpha} \cap Y $, $T _{\alpha} : = P _{\alpha} \cap T $, $X _{\alpha} : = P _{\alpha} \cap X$.
D'après un théorème d'Elkik (voir \cite{elkik}), il existe un $\V$-schéma formel faible affine et lisse 
$Y ^{\dag} _{\alpha}$ relevant $Y _{\alpha}$. 
On choisit $v _{\alpha}\,:\, Y ^{\dag} _{\alpha} \hookrightarrow U ^{\dag} _{\alpha}$
une immersion fermée relevant
$Y _{\alpha} \hookrightarrow U  _{\alpha}$
Soit $E \in \mathrm{Isoc} ^\dag (Y/K) $. Il lui correspond
une famille $(E _{\alpha}) _{\alpha}$  d'isocristaux surconvergents sur $Y _{\alpha}$ 
munie d'une donnée de recollement (e.g., en remarquant que la structure de Frobenius est inutile, voir la section \cite[1.4.4]{caro-2006-surcoh-surcv}).
On note $\E _{\alpha} := \sp _{Y ^{\dag}   _{\alpha}\underset{v _{\alpha}}{\hookrightarrow} U ^{\dag}  _{\alpha}, T _{\alpha},+} (E _{\alpha})$.
Grâce à \ref{sp+dag-indt} et à sa remarque \ref{sp+dag-indt-rema}, on vérifie (de manière analogue à \cite[2.5]{caro-construction})
que la donnée de recollement 
de $E _{\alpha}$ induit une donnée de recollement sur la famille $(\E _{\alpha}) _{\alpha}$ de 
$\D ^{\dag} _{\PP _{\alpha}} (\hdag T _{\alpha}) _{\Q}$-modules cohérents à support dans $X _{\alpha}$. 
Cette famille $(\E _{\alpha}) _{\alpha}$ se recolle donc en un $\D ^{\dag} _{\PP} (\hdag T ) _{\Q}$-module cohérent à support dans $X $, que l'on note $\E$. 
D'après le théorème \ref{spYdag+to}, 
$\E _{\alpha} \in   \mathrm{Isoc} ^{\dag \dag }(\PP_{\alpha}, T_{\alpha},X_{\alpha}/K)$.
Il en résulte que $\E  \in   \mathrm{Isoc} ^{\dag \dag }(\PP, T,X/K)$ (en effet, le fait que $\E$ soit un $\D ^{\dag} _{\PP} (\hdag T ) _{\Q}$-module surcohérent
à support dans $X$ est local en $\PP$, de même que celui que $\E |\U$ soit dans l'image essentielle de $\sp _{Y\hookrightarrow \U,+}$).

\end{vide}

\begin{prop}
\label{5.2.6-recol}
Considérons le diagramme commutatif :
\begin{equation}
\label{formel-II-diag-copiebis}
\xymatrix{
{Y^{\prime}} \ar@{^{(}->}[r]  ^-{v ^{\prime}} \ar[d] ^-{b}
& 
{U ^{\dag \prime} } \ar@{^{(}->}[r] ^-{j^{\prime}} \ar[d] ^-{g}
& 
{P ^{\dag \prime} }  \ar[d] ^-{f}
\\ 
{Y} \ar@{^{(}->}[r] ^-{v}
&
 {U  ^{\dag} } \ar@{^{(}->}[r] ^-{j} & {P  ^{\dag} ,}
}
\end{equation}
où $f$ et $g$ sont des morphismes lisses de $\V$-schémas formels faibles séparés et lisses,
$b$ est un morphisme de $k$-variétés lisses, 
$j$ et $j^{\prime}$ sont des immersions ouvertes,
$v$ et $v ^{\prime}$ sont des immersions fermées.
On suppose qu'il existe un diviseur $T$ de $P$
(resp. $T ^{\prime}$ de $P ^{\prime}$) tel que $U ^{\dag}= P ^{\dag} \setminus T$ 
(resp. $U ^{\prime\dag }= P ^{\prime\dag }\setminus T^{\prime}$).
Notons $X '$ l'adhérence de $Y'$ dans $P'$.
On dispose alors, pour tout $E \in \mathrm{Isoc} ^\dag (Y/K)$, de l'isomorphisme canonique 
\begin{equation}
\sp _{Y ^{\prime}  \hookrightarrow U ^{\prime\dag},T ^{\prime} , + }  \circ  b ^{*} (E) 
\riso
(\hdag T') \circ \R \underline{\Gamma} ^{\dag} _{X '} \circ  f ^{!} \circ  \sp _{Y  \hookrightarrow U ^{\dag},T  , + } (E) [-d _{X' /X}].
\end{equation}
\end{prop}

\begin{proof}
Cela résulte par recollement de \cite[5.2.3 et 5.2.6]{caro_devissge_surcoh} 
(en effet, la proposition \cite[5.2.6]{caro_devissge_surcoh} peut être utilisée grâce à la propriété
de finitude \ref{spYdag+to}).
\end{proof}

Le théorème suivant répond positivement à la conjecture \cite[5.3.3.(a)]{caro_devissge_surcoh} : 
\begin{theo}
 \label{sp+essent}
Avec les notations de \ref{nota-spY->UdagT}, on suppose $X $ lisse.
Le foncteur $\sp _{X  \hookrightarrow \PP,T  , + }\,:\,  \mathrm{Isoc} ^{\dag} (Y,X /K)  \to \mathrm{Isoc} ^{\dag\dag} (Y,X /K) $
est une équivalence de catégories
s'inscrivant dans le diagramme canonique 
\begin{equation}
\xymatrix{
{\mathrm{Isoc} ^\dag (Y/K)}
\ar[rr] ^-{}
\ar[rd] _-{\sp _{Y  \hookrightarrow U ^{\dag},T  , + }} 
&&
{\mathrm{Isoc} ^{\dag} (Y,X /K) } 
\ar[ld] ^-{ \sp _{X  \hookrightarrow \PP,T  , + }} _-{\cong}
\\ 
& 
{\mathrm{Isoc} ^{\dag \dag} (Y,X /K)} 
& 
}
\end{equation}
commutatif à isomorphisme canonique près. 
\end{theo}

\begin{proof}
Dans un premier temps, supposons $P ^{\dag}$ affine et $T$ défini par une équation locale. 
De plus, supposons qu'il existe une immersion fermée
$u\,:\, X ^{\dag} \hookrightarrow P ^{\dag}$
de $\V$-schémas formels faibles lisses 
relevant $X \hookrightarrow P$,
une immersion fermée 
$v\,:\, Y ^{\dag} \hookrightarrow U ^{\dag}$
relevant $Y \hookrightarrow U ^{\dag}$, 
une immersion ouverte
$\lambda \,:\,
Y ^{\dag} \hookrightarrow X ^{\dag}$
tels que 
$j \circ v = u \circ \lambda$.

Soit $\E $ un $\smash{\D}  _{Y ^\dag ,\Q}$-module globalement de présentation, $\O _{Y ^\dag ,\Q}$-cohérent
qui soit associé via le foncteur $\sp _{*}$ à un élément de $\mathrm{Isoc} ^{\dag} (Y/K)$.
Or, on bénéficie d'après \cite[2.7.3.(ii)]{huyghe-comparaison}
des injections
$j_* \D ^\dag _{Y ^{\dag}  ,\Q} \rightarrow \D ^\dag _{\X  } (\hdag T ) _\Q$
et
$j_* \O _{Y ^{\dag}  ,\Q} \rightarrow \O _{\X  } (\hdag T ) _\Q$.
On obtient alors la flèche : 
$\O _{\X  } (\hdag T ) _\Q \otimes _{j_* \O _{Y ^{\dag},\Q}} j _{*}(\E) \to 
\D ^\dag _{\X  } (\hdag T ) _\Q \otimes _{j_* \D _{Y ^{\dag}  ,\Q}} j _{*}(\E)$.
Comme 
$\D ^\dag _{\X  } (\hdag T ) _\Q \otimes _{j_* \D _{Y ^{\dag}  ,\Q}} j _{*}(\E)$
est un 
$\D ^\dag _{\X  } (\hdag T ) _\Q$-module cohérent dont la restriction sur $\Y$ est 
un isocristal convergent sur $Y$ (voir le corollaire \cite[5.1.3]{caro_devissge_surcoh}), 
$\D ^\dag _{\X  } (\hdag T ) _\Q \otimes _{j_* \D _{Y ^{\dag}  ,\Q}} j _{*}(\E)$
est donc un isocristal sur $Y$ surconvergent le long de $T$ (voir le théorème de Berthelot énoncé dans \cite[2.2.12]{caro_courbe-nouveau}).
En particulier, ce dernier est $\O _{\X  } (\hdag T ) _\Q$-cohérent.
Comme la flèche 
$\O _{\X  } (\hdag T ) _\Q \otimes _{j_* \O _{Y ^{\dag},\Q}} j _{*}(\E) \to 
\D ^\dag _{\X  } (\hdag T ) _\Q \otimes _{j_* \D _{Y ^{\dag}  ,\Q}} j _{*}(\E)$
est un morphisme de 
$\O _{\X  } (\hdag T ) _\Q$-modules cohérents qui est un isomorphisme en dehors de $T$,
celui-ci est un isomorphisme. 
On obtient ainsi un foncteur 
$\mathrm{Isoc} ^{\dag} (Y/K) \to
\mathrm{Isoc} ^{\dag} (Y,X/K)\cong \mathrm{Isoc} ^{\dag \dag} (\X, T, X/K)$
défini par 
$\E \to  \D ^\dag _{\X  } (\hdag T ) _\Q \otimes _{j_* \D _{Y ^{\dag}  ,\Q}} j _{*}(\E).$

D'après \ref{lemm-pass-cvasurcv-isopart1} est avec ses notations, on dispose de l'isomorphisme
canonique 
$ \sp _{Y ^{\dag}  \hookrightarrow U ^{\dag},T  , + }(E) \riso
\smash{\D} ^\dag _{\PP} (\hdag T ) _\Q \otimes _{D _{U^\dag, \Q} }  v _+ (E)$.
Comme $j \circ v = u \circ \lambda$, le morphisme canonique 
$D _{Y ^{\dag}  ,\Q} \to \D ^\dag _{\X  } (\hdag T ) _\Q$
induit 
$ v ^{*}D _{Y ^{\dag}  ,\Q} 
\to 
u ^{*} \D ^\dag _{\X  } (\hdag T ) _\Q$.
D'où :
$D _{U ^{\dag}\leftarrow Y ^{\dag}  ,\Q}  \to  \D ^\dag _{\PP \leftarrow \X  } (\hdag T ) _\Q$.
Il en résulte la flèche 
$$v _{+} (E) \to 
u _{+} (\D ^\dag _{\X  } (\hdag T ) _\Q \otimes _{j_* \D _{Y ^{\dag}  ,\Q}} j _{*}(\E))
= \sp _{X \hookrightarrow \PP, T,+} 
(\D ^\dag _{\X  } (\hdag T ) _\Q \otimes _{j_* \D _{Y ^{\dag}  ,\Q}} j _{*}(\E)).$$ 
On en déduit par extension le morphisme canonique 
\begin{equation}
\label{fleche-sp->sp}
 \sp _{Y ^{\dag}  \hookrightarrow U ^{\dag},T  , + }(E)  
 \to 
 \sp _{X \hookrightarrow \PP, T,+} 
(\D ^\dag _{\X  } (\hdag T ) _\Q \otimes _{j_* \D _{Y ^{\dag}  ,\Q}} j _{*}(\E)).
\end{equation}
Cette flèche \ref{fleche-sp->sp} est un morphisme de
$\D ^\dag _{\PP  } (\hdag T ) _\Q$-modules cohérents dont la restriction en dehors de $T$ est un isomorphisme (cela résulte de \cite[5.1.3 et 5.2.4]{caro_devissge_surcoh}).
 D'après \cite[4.3.12]{Be1}, cela implique que 
le morphisme \ref{fleche-sp->sp} est en fait un isomorphisme.
D'où le résultat dans le premier cas. 

Passons à présent au cas général. 
Soit $(P ^{\dag} _{\alpha}) _{\alpha}$ un recouvrement d'ouverts affines de $P ^{\dag}$ tel que 
$T \cap P _{\alpha}$ soit défini par une équation. 
On note 
$U ^{\dag} _{\alpha}:= P ^{\dag} _{\alpha} \cap U ^{\dag}$ , 
$Y  _{\alpha}:= P _{\alpha} \cap Y $, $T _{\alpha} : = P _{\alpha} \cap T $, $X _{\alpha} : = P _{\alpha} \cap X$.
D'après un théorème d'Elkik (voir \cite{elkik}), comme $X _{\alpha}$ est affine et lisse, 
il existe un $\V$-schéma formel faible affine et lisse 
$X ^{\dag} _{\alpha}$ relevant $X _{\alpha}$. 
On note $j _{\alpha}\,:\, P ^{\dag} _{\alpha} \subset P ^{\dag}$ l'immersion ouverte canonique.
Choisissons 
$u _{\alpha}\,:\, X ^{\dag} _{\alpha} \hookrightarrow P ^{\dag} _{\alpha}$ une immersion fermée relevant
$X  _{\alpha} \hookrightarrow P _{\alpha}$.
En choisissant un isomorphisme 
$Y ^{\dag} _{\alpha} \riso X ^{\dag} _{\alpha}  \cap U ^{\dag} _{\alpha} $, 
on obtient 
une immersion ouverte $\lambda _{\alpha}\,:\,Y ^{\dag} _{\alpha} \hookrightarrow X ^{\dag} _{\alpha} $
et une immersion fermée de la forme
$v _{\alpha}\,:\, Y ^{\dag} _{\alpha} \hookrightarrow U ^{\dag} _{\alpha}$
telles que
$u _{\alpha } \circ \lambda _{\alpha} = j _{\alpha}\circ v _{\alpha}$.
On se ramène ainsi par recollement au premier cas traité. 
\end{proof}

\begin{theo}
\label{spYdag+}
Avec les notations de \ref{nota-spY->UdagT}, on suppose $P$ propre.
Le foncteur $ \sp _{Y  \hookrightarrow U ^{\dag},T  , + }$ induit alors l'équivalence de catégories : 
\begin{equation}
\label{spYdag+-Eq1}
\sp _{Y   \hookrightarrow U ^{\dag},T  , + }\,:\,\mathrm{Isoc} ^{\dag }(Y/K)  \cong \mathrm{Isoc} ^{\dag \dag }(Y/K).
\end{equation}
Pour tout $E \in \mathrm{Isoc} ^{\dag }(Y/K)$, 
on dispose de l'isomorphisme canonique :
\begin{equation}
\label{spYdag+-Eq2}
\sp _{Y   \hookrightarrow U ^{\dag},T  , + }(F ^{*} (E))
\riso 
F ^{*} \sp _{Y   \hookrightarrow U ^{\dag},T  , + }(E).
\end{equation}
Plus généralement, l'équivalence de catégories \ref{spYdag+-Eq1} s'étend par recollement au cas d'une variété lisse sur $k$ quelconque $Y$. 
On la note alors  $\sp _{Y+}\,:\,\mathrm{Isoc} ^{\dag }(Y/K)  \cong \mathrm{Isoc} ^{\dag \dag }(Y/K)$.
\end{theo}

\begin{proof}
I) Vérifions d'abord l'équivalence \ref{spYdag+-Eq1}.
Commençons par quelques notations. Grâce au théorème de désingularisation de de Jong (voir \cite{dejong}), 
il existe un diviseur $\widetilde{T}$ de $P$ tel que, 
en posant $\widetilde{Y} = X \setminus  \widetilde{T}$,
on obtienne le diagramme commutatif 
\begin{equation}
\label{formel-II-diag-copie}
\xymatrix{
{\widetilde{Y}^{(0)}} \ar@{^{(}->}[r]  ^-{l ^{(0)}} \ar[d] ^-{c}
& 
{Y^{(0)}} \ar@{^{(}->}[r]  ^-{j ^{(0)}} \ar[d] ^-{b}
& 
{X^{(0)} } \ar@{^{(}->}[r] ^-{u^{(0)}} \ar[d] ^-{a}
& 
{P ^{\dag (0)} }  \ar[d] ^-{f}
\\ 
{\widetilde{Y}} \ar@{^{(}->}[r] ^-{l}
&
{Y} \ar@{^{(}->}[r] ^-{j }
&
 {X } \ar@{^{(}->}[r] ^-{u} & {P  ^{\dag} ,}
}
\end{equation}
où les deux carrés de gauche sont cartésiens, 
$f$ est un morphisme propre et lisse de $\V$-schémas formels faibles séparés et lisses,
$a$ est un morphisme propre, surjectif de $k$-variétés avec $X ^{(0)}$ lisse, 
$b$ est un morphisme de $k$-variétés lisses, 
$c$ est un morphisme fini et étale, 
$j$ et $j^{\prime}$ sont des immersions ouvertes, $u$ et $u^{(0)}$ sont des immersions fermées.
On note $\widetilde{j}\,:\, \widetilde{Y}\hookrightarrow X$ et $\widetilde{j}\,:\, \widetilde{Y} ^{(0)}\hookrightarrow X ^{(0)}$ les immersions ouvertes. 
On note $U ^{\dag}:= P ^{\dag} \setminus T$, 
$T ^{(0)}:=f ^{-1}(T)$, $U ^{\dag (0)}:= P ^{\dag(0)}\setminus T^{(0)}$ et $g \,:\, U ^{\dag(0)} \to U ^{\dag}$ le morphisme induit par $f$.
De même en rajoutant des tildes (e.g. on pose $\widetilde{U} ^{\dag} :=  P ^{\dag} \setminus \widetilde{T}$ etc.).
On reprend de plus les notations analogues à celles de \ref{formel-II} concernant $a _{1}$, $a _{2}$ etc.

1) Vérifions que le foncteur $ \sp _{Y  \hookrightarrow U ^{\dag},T  , + }$ est pleinement fidèle. 
Pour cela, considérons le diagramme ci-dessous :
\begin{equation}
\label{spYdag+-carre2bis}
\xymatrix{
{\mathrm{Isoc} ^{*} (\PP, T, X/K)}
\ar[r] ^-{(a ^{!}, |\U)}
& 
{\mathrm{Isoc} ^{*} (\PP ^{(0)},T^{(0)}, X^{(0)}/K)
\times _{\mathrm{Isoc} ^{*} (\U ^{(0)}, Y^{(0)}/K)} 
\mathrm{Isoc} ^{*} (\U, Y/K)} 
 \\ 
{\mathrm{Isoc} ^{\dag }(Y/K)} 
 \ar[r] ^-{(a ^{*},  j ^{*})}
 \ar[u] ^-{\sp _{Y   \hookrightarrow U ^{\dag},T  , + }}
 & 
 {\mathrm{Isoc} ^{\dag }(Y^{(0)}/K) 
 \times _{\mathrm{Isoc} ^{\dag }(Y  ^{(0)},Y  ^{(0)} /K) } 
 \mathrm{Isoc} ^{\dag }(Y,Y/K).} 
 \ar[u] _-{(\sp _{Y^{(0)} \hookrightarrow U ^{\dag (0)},T^{(0)},+}, 
 \sp _{Y   \hookrightarrow \U, + })}
 }
\end{equation}
Il résulte de \ref{5.2.6-recol} et de \cite[5.2.4]{caro_devissge_surcoh} 
que ce diagramme est commutatif à isomorphisme canonique près. 
Or, d'après le théorème \ref{plfid-dag} (resp. \ref{a+|U-plfid-prop} et \ref{coro-theo-!=+}), 
le foncteur du bas (resp. du haut) est pleinement fidèle. 
D'après  \cite{caro-construction}, les foncteurs $ \sp _{Y   \hookrightarrow \U, + }$ 
et $\sp _{X^{(0)} \hookrightarrow \PP ^{(0)},T^{(0)},+}$
sont pleinement fidèles. 
De plus, comme $X ^{(0)}$ est propre et lisse,  d'après \ref{sp+essent}, 
$\sp _{Y^{(0)} \hookrightarrow U ^{\dag (0)},T^{(0)},+}
\riso \sp _{X^{(0)} \hookrightarrow \PP ^{(0)},T^{(0)},+}$. 
Le foncteur de droite de \ref{spYdag+-carre2bis} est donc pleinement fidèle. 
D'où le résultat.

2) Construisons à présent le foncteur
$\sp _{Y^{(\bullet)} \hookrightarrow \widetilde{U} ^{\dag (\bullet)},\widetilde{T} ^{(\bullet)},+}\,:\,
\mathrm{Isoc} ^{\dag} (\widetilde{Y}^{(\bullet)},X^{(\bullet)}/K)
\to
\mathrm{Isoc} ^{*} (\PP^{(\bullet)}, \widetilde{T}^{(\bullet)}, X^{(\bullet)}/K)$.
Soit $(E ^{(0)}, \epsilon)\in \mathrm{Isoc} ^{\dag} (\widetilde{Y}^{(\bullet)},X^{(\bullet)}/K)$.
On définit l'isomorphisme $\theta$ en demandant au diagramme suivant d'être commutatif : 
\begin{equation}
\label{epsilon-theta-sp+}
\xymatrix{
{a _{2} ^{!} \sp _{\widetilde{Y}^{(0)} \hookrightarrow \widetilde{U} ^{\dag (0)},\widetilde{T} ^{(0)},+} (E ^{(0)})} 
\ar@{.>}[r] ^-{\theta} _-{\sim}
& 
{a _{1} ^{!} \sp _{\widetilde{Y}^{(0)} \hookrightarrow \widetilde{U} ^{\dag (0)},\widetilde{T} ^{(0)},+} (E ^{(0)})} 
 \\ 
{\sp _{\widetilde{Y}^{(1)} \hookrightarrow \widetilde{U} ^{\dag (1)},\widetilde{T} ^{(1)},+} (a _{2} ^{*}  E ^{(0)})} 
\ar[u] ^-{\sim}
\ar[r] _-{\epsilon} ^-{\sim}
& 
{\sp _{\widetilde{Y}^{(1)} \hookrightarrow \widetilde{U} ^{\dag (1)},\widetilde{T} ^{(1)},+} (a _{1} ^{*}  E ^{(0)}),} 
\ar[u] ^-{\sim}
}
\end{equation}
les isomorphismes verticaux résultant de \ref{5.2.6-recol}. L'isomorphisme $\theta$ vérifie la condition de cocycle (par \cite[4.3.12]{Be1}, il suffit de le voir en dehors des diviseurs, ce qui est immédiat).
On construit ainsi le foncteur $\sp _{\widetilde{Y}^{(\bullet)} \hookrightarrow \widetilde{U} ^{\dag (\bullet)},\widetilde{T} ^{(\bullet)},+}$ 
en posant 
$\sp _{\widetilde{Y}^{(\bullet)} \hookrightarrow \widetilde{U} ^{\dag (\bullet)},\widetilde{T} ^{(\bullet)},+} (E ^{(0)}, \epsilon)
:= (\sp _{\widetilde{Y}^{(0)} \hookrightarrow \widetilde{U} ^{\dag (0)},\widetilde{T} ^{(0)},+} (E ^{(0)}), \theta)$.

3) Ce foncteur 
$\sp _{\widetilde{Y}^{(\bullet)} \hookrightarrow \widetilde{U} ^{\dag (\bullet)},\widetilde{T} ^{(\bullet)},+}$ 
est une équivalence de catégories.
Preuve : comme $X ^{(0)}$ est lisse, le foncteur $\sp _{\widetilde{Y}^{(0)} \hookrightarrow \widetilde{U} ^{\dag (0)},\widetilde{T} ^{(0)},+} $ 
induit une équivalence de catégories (cela résulte du théorème \ref{sp+essent}).
Or, d'après l'étape 1) (établit pour une situation générale), le foncteur 
$\sp _{\widetilde{Y}^{(1)} \hookrightarrow \widetilde{U} ^{\dag (1)},\widetilde{T} ^{(1)},+}$ est pleinement fidèle. 
En considérant le diagramme \ref{epsilon-theta-sp+}, cela implique que la donnée de recollement $\epsilon$ sur $E ^{(0)}$
équivaut à une donnée de recollement $\theta$ sur $\sp _{\widetilde{Y}^{(0)} \hookrightarrow \widetilde{U} ^{\dag (0)},\widetilde{T} ^{(0)},+} (E ^{(0)})$, 
ce qui termine la preuve de 3). 
 
4) Le foncteur 
$\sp _{\widetilde{Y}   \hookrightarrow \widetilde{U} ^{\dag},\widetilde{T}  , + }$ est une équivalence de catégories.
À cette fin, considérons le diagramme
\begin{equation}
\label{spYdag+-carre1}
\xymatrix{
{\mathrm{Isoc} ^{\dag \dag }(\widetilde{Y}/K)} 
\ar@{=}[r]
&
{\mathrm{Isoc} ^{*} (\PP, \widetilde{T}, X/K)}
\ar[r] ^-{\mathcal{L}oc}
& 
{\mathrm{Isoc} ^{*} (\PP^{(\bullet)}, \widetilde{T}^{(\bullet)}, X^{(\bullet)}/K)} 
\\
& 
{\mathrm{Isoc} ^{\dag }(\widetilde{Y}/K)} 
\ar[u] ^-{\sp _{\widetilde{Y}   \hookrightarrow \widetilde{U} ^{\dag},\widetilde{T}  , + }}
\ar[r] ^-{\mathcal{L}oc}
& 
{\mathrm{Isoc} ^{\dag} (\widetilde{Y}^{(\bullet)},X^{(\bullet)}/K) .}
\ar[u] _-{\sp _{\widetilde{Y} ^{(\bullet)} \hookrightarrow \widetilde{U} ^{\dag (\bullet)},\widetilde{T} ^{(\bullet)},+}} 
}
\end{equation}
Via \ref{5.2.6-recol} et \ref{coro-theo-!=+}, on vérifie que ce diagramme est commutatif à isomorphisme canonique près. 
Or, on vient d'établir que le foncteur de droite est une équivalence de catégories. 
Avec le théorème de descente de Shiho de \cite[7.3]{shiho-logRC-RCII} (resp. via \ref{desc-fini-ét} et \ref{coro-theo-!=+}), 
il en est de même du foncteur du bas (resp. du haut), 
ce qui conclut la preuve de 4).

5) 
On dispose de l'équivalence de
catégories 
$\sp _{\widetilde{Y}   \hookrightarrow \widetilde{U} ^{\dag},\widetilde{T}  , + }\,:\,
\mathrm{Isoc} ^{\dag }(\widetilde{Y}\supset Y/K)
\riso
\mathrm{Isoc} ^{*} (\PP, \widetilde{T} \supset T, X/K)$, 
où $\mathrm{Isoc} ^{\dag }(\widetilde{Y}\supset Y/K)$ désigne plus simplement 
$\mathrm{Isoc} ^{\dag }(\widetilde{Y}\supset Y,X/K)$.
En effet, il s'agit de procéder de manière analogue à la preuve de \ref{lemm-diag-eqcat-dagdag-ii}
en utilisant l'étape 4), \ref{5.2.6-recol} et \ref{sp+essent}.

6) Fin de la preuve de l'étape I).
Considérons à présent le diagramme ci-dessous :
\begin{equation}
\label{spYdag+-carre2}
\xymatrix{
{\mathrm{Isoc} ^{*} (\PP, T, X/K)}
\ar[r] ^-{(a ^{!}, (\hdag \widetilde{T}))}
& 
{\mathrm{Isoc} ^{*} (\PP ^{(0)},T^{(0)}, X^{(0)}/K)
\times _{\mathrm{Isoc} ^{*} (\PP ^{(0)}, \widetilde{T}^{(0)}, X^{(0)}/K)} 
\mathrm{Isoc} ^{*} (\PP, \widetilde{T}\supset T, X/K)} 
 \\ 
{\mathrm{Isoc} ^{\dag }(Y/K)} 
 \ar[r] ^-{(a ^{*},  \widetilde{j} ^{\dag})}
 \ar[u] ^-{\sp _{Y   \hookrightarrow U ^{\dag},T  , + }}
 & 
 {\mathrm{Isoc} ^{\dag }(Y^{(0)}/K) 
 \times _{\mathrm{Isoc} ^{\dag }(\widetilde{Y } ^{(0)} /K) } 
 \mathrm{Isoc} ^{\dag }(\widetilde{Y}\supset Y/K).} 
 \ar[u] _-{(\sp _{Y^{(0)} \hookrightarrow U ^{\dag (0)},T^{(0)},+}, 
 \sp _{\widetilde{Y}   \hookrightarrow \widetilde{U} ^{\dag},\widetilde{T}  , + })}
 }
\end{equation}
Il résulte de \ref{5.2.6-recol} que ce diagramme est commutatif à isomorphisme canonique près. 
Or, d'après le théorème \ref{EqCat-rig} (resp. \ref{Th2} et \ref{coro-theo-!=+}), 
le foncteur du bas (resp. du haut) est une équivalence de catégories. 
D'après l'étape 5) (resp. grâce au théorème \ref{sp+essent}) le foncteur
$\sp _{\widetilde{Y}   \hookrightarrow \widetilde{U} ^{\dag},\widetilde{T}  , + }$
(resp. $\sp _{Y^{(0)} \hookrightarrow U ^{\dag (0)},T^{(0)},+}$)
est une équivalence de catégories. 
Cela implique qu'il en est de même du foncteur de droite de \ref{spYdag+-carre2}.
Comme cela est le cas pour les trois autres foncteurs de \ref{spYdag+-carre2},
le foncteur 
$\sp _{Y   \hookrightarrow U ^{\dag},T  , + }$ est donc une équivalence de catégories.

II) Traitons à présent l'isomorphisme \ref{spYdag+-Eq2}. Posons $\E := \sp _{Y   \hookrightarrow U ^{\dag},T  , + }(E)$.
Notons $p _{1},p _{2}\,:\, \PP \times \PP \to \PP$ respectivement les projections à gauche et à droite, 
$\gamma _{F}:= (Id, F) \,:\, P \hookrightarrow P \times P$. 
Ainsi, $F ^{*} (\E) = p _{1+} \circ \R \underline{\Gamma} ^\dag _{P} p ^{!} _{2} (\E)$, 
où l'on voit $P$ comme sous-schéma fermé de $P\times P$ via $\gamma _{F}$.
Les propositions 
\ref{5.2.6-recol} et \ref{gen-coh-PXTindtP} nous permettent de conclure.

III) La dernière assertion résulte de ce que l'on vient d'établir à l'étape I) et
de la procédure de recollement de la section \cite[1.4]{caro-2006-surcoh-surcv}
et via \cite[2.2.3--4]{caro-2006-surcoh-surcv} dont les démonstrations restent valables sans structure de Frobenius.
\end{proof}

\subsection{Dévissage en isocristaux surconvergents}

Soit $\PP$ un $\V$-schéma formel séparé et lisse.
Si $S$ est un sous-schéma de $P$, on notera $\overline{S}$ l'adhérence schématique de $S$ dans $P$.

\begin{lemm}
\label{lemm-theodev2}
  Soient $T$ un diviseur de $P$, $Y $ un sous-schéma fermé de $P \setminus T $,
  $Y _1$ une composante irréductible de $Y$ et
  $\E $ un $(F\text{-})\D ^\dag _{\PP} (\hdag T) _\Q$-module surcohérent à support dans $\overline{Y}$.

  Il existe un diviseur $T '\supset T$ de $P$ tel que,
  en notant $Y' := \overline{Y} \setminus T'$,
  $Y'$ soit intègre, lisse, inclus et dense dans $Y _1$ et
  $\E (\hdag T') \in (F\text{-})\mathrm{Isoc} ^{\dag \dag}( \PP, T', \overline{Y}/K) = (F\text{-})\mathrm{Isoc} ^{\dag \dag}( Y', \overline{Y'}/K)$.
\end{lemm}

\begin{proof}
Il s'agit de reprendre la preuve de \cite[3.1.1]{caro-2006-surcoh-surcv} en remplaçant dans cette preuve respectivement la référence \cite[2.3.2]{caro-2006-surcoh-surcv} par \ref{Isoc*=dagdag} 
et la référence \cite[2.2.17]{caro_courbe-nouveau} par \ref{borel9.3} (qui donne une version sans structure de Frobenius du théorème \cite[2.2.17]{caro_courbe-nouveau}). 
\end{proof}

Grâce à \ref{diag-eqcat-dagdag-ii}, on étend la notion de dévissage en $F$-isocristaux surconvergents au cas sans structure de Frobenius et au cas partiellement surconvergent de la manière suivante : 
\begin{defi}
\label{defi-dev}
Soient $T$ un diviseur de $P$, $Y $ un sous-schéma fermé de $P \setminus T $,
  $\E \in (F\text{-})D ^\mathrm{b} _\mathrm{coh} (\smash{\D} ^\dag _{\PP} (\hdag T) _\Q)$
  à support dans $\overline{Y}$.
Le complexe $\E$ {\og se dévisse en isocristaux surconvergents\fg} 
s'il existe des diviseurs $T _1, \dots,T _{r}$ contenant $T$ avec $T _r =T$ tels que,
en notant $T _0: =\overline{Y}$ et, 
pour $i=0,\dots , r-1$, $X _{i}=T _0 \cap \dots \cap T _i $,
\begin{enumerate}
\item  le schéma $Y _i:= X _{i}\setminus T _{i+1}$
est lisse ;
\item 
$  \mathcal{H} ^j (\R \underline{\Gamma} ^\dag _{X _i}  (\hdag T _{i+1}) (\E) )
  \in (F\text{-})\mathrm{Isoc} ^{\dag \dag}( \PP, T _{i+1}, X _i/K)
  =(F\text{-})\mathrm{Isoc} ^{\dag \dag}( Y_i,X _{i}/K) $, pour tout entier $j$.
 \end{enumerate}
On note $(F\text{-})D ^\mathrm{b} _\mathrm{\text{dév}} (\smash{\D} ^\dag _{\PP} (\hdag T) _\Q)$ la sous-catégorie pleine de $(F\text{-})D ^\mathrm{b} _\mathrm{coh} (\smash{\D} ^\dag _{\PP} (\hdag T) _\Q)$ des complexes dévissables en isocristaux surconvergents. 
\end{defi}

Le théorème suivant améliore \cite{caro_devissge_surcoh,caro-2006-surcoh-surcv} car 
nous nous sommes affranchi de l'utilisation d'une structure de Frobenius :

\begin{theo}
\label{caradev}
Soit $T$ un diviseur de $P$. On dispose de l'inclusion 
$(F\text{-})D ^\mathrm{b} _\mathrm{surcoh} (\smash{\D} ^\dag _{\PP} (\hdag T) _\Q) \subset (F\text{-})D ^\mathrm{b} _\mathrm{\text{\rm dév}} (\smash{\D} ^\dag _{\PP} (\hdag T) _\Q)$.
\end{theo}

\begin{proof}
On procède par récurrence sur $Y$ en utilisant le lemme précédent \ref{lemm-theodev2}.
\end{proof}

Lorsque l'on dispose d'une structure de Frobenius, de manière identique à \cite[2.3.18]{caro-Tsuzuki-2008}, on vérifie : 
\begin{theo}
\label{theo-overcoh=dev}
Soient $T$ un diviseur de $P$, 
$\E \in F \text{-}D ^\mathrm{b} _\mathrm{coh} (\D ^\dag _{\PP} (\hdag T) _{\Q})$.
Les assertions suivantes sont équivalentes :
\begin{enumerate}
  \item Le $F$-complexe $\E$ est $\D ^\dag _{\PP} (\hdag T) _{\Q}$-surcohérent ;
  \item Le $F$-complexe $\E$ est $\D ^\dag _{\PP,\Q}$-surcohérent ;
   \item Le $F$-complexe $\E$ est surholonome ;
  \item Le $F$-complexe $\E$ est dévissable en $F$-isocristaux surconvergents.
\end{enumerate}
\end{theo}

\begin{rema}
\label{rema-sans-Frob-faux}
Le théorème \ref{theo-overcoh=dev} est faux sans structure de Frobenius. 
En effet, en reprenant le contre-exemple de Berthelot donné à la fin de \cite{Be1}, il existe des $\D ^\dag _{\PP} (\hdag T) _{\Q}$-modules 
$\O _{\PP} (\hdag T) _{\Q}$-cohérents qui ne sont pas $\D ^\dag _{\PP,\Q}$-cohérents. 
Or, d'après \cite[6.1.4]{caro_devissge_surcoh}, les $\D ^\dag _{\PP} (\hdag T) _{\Q}$-modules cohérents $\O _{\PP} (\hdag T) _{\Q}$-cohérents sont $\D ^\dag _{\PP} (\hdag T) _{\Q}$-surcohérents.
On en déduit qu'un $\D ^\dag _{\PP} (\hdag T) _{\Q}$-module surcohérent n'est pas forcément $\D ^\dag _{\PP,\Q}$-cohérent ni a fortiori $\D ^\dag _{\PP,\Q}$-surcohérent.
\end{rema}

\bibliographystyle{smfalpha}

\newcommand{\etalchar}[1]{$^{#1}$}
\providecommand{\bysame}{\leavevmode ---\ }
\providecommand{\og}{``}
\providecommand{\fg}{''}
\providecommand{\smfandname}{et}
\providecommand{\smfedsname}{\'eds.}
\providecommand{\smfedname}{\'ed.}
\providecommand{\smfmastersthesisname}{M\'emoire}
\providecommand{\smfphdthesisname}{Th\`ese}


\bigskip
\noindent Daniel Caro\\
Laboratoire de Mathématiques Nicolas Oresme\\
Université de Caen
Campus 2\\
14032 Caen Cedex\\
France.\\
email: daniel.caro@math.unicaen.fr

\end{document}